\documentclass{amsart}
\usepackage{graphicx}
\usepackage[all]{xy}
\usepackage{latexsym}
\usepackage{amssymb}
\usepackage{amsmath}
\usepackage{color}

\xyoption{arc}

 \newfam\cyrfam

  \font\tencyr=wncyr10

  \font\sevencyr=wncyr7

  \font\fivecyr=wncyr5

  \textfont\cyrfam=\tencyr \scriptfont\cyrfam=\sevencyr

    \scriptscriptfont\cyrfam=\fivecyr

  \newfam\cyifam

  \font\tencyi=wncyi10

  \font\sevencyi=wncyi7

  \font\fivecyi=wncyi5

  \textfont\cyifam=\tencyi \scriptfont\cyifam=\sevencyi

    \scriptscriptfont\cyifam=\fivecyi

\def\id{{\mbox{1 \hskip -7pt 1}}}
\newcommand{\sgn}{{\mathit s  \mathit g\mathit  n}}
 \newcommand{\lon}{\longrightarrow}
 \newcommand{\bu}{\bullet}
 
 \newcommand{\rar}{\rightarrow}

\newcommand{\p}{{\partial}}
\newcommand{\pc}{{\partial}_\centerdot}
\newcommand{\Id}{{\mathrm{Id}}}

\newcommand{\Der}{\mathrm{Der}}
\newcommand{\Tw}{\mathsf{Tw}}
\newcommand{\TW}{\mathsf{Tw}}
\newcommand{\wTW}{\widetilde{\mathsf{Tw}}}
\newcommand{\tw}{\mathsf{tw}}

\newcommand{\RGra}{\mathcal{R} \mathcal{G} ra}
\newcommand{\GRav}{\mathcal{G} \mathcal{R} av}
\newcommand{\Grav}{\mathcal{G} rav}
\newcommand{\Gra}{\mathcal{G} ra}

 \newcommand{\Z}{{\mathbb Z}}
 \newcommand{\bS}{{\mathbb S}}

 \newcommand{\R}{{\mathbb R}}
 \newcommand{\N}{{\mathbb N}}
 \newcommand{\K}{{\mathbb K}}

 \newcommand{\ot}{\otimes}

\newcommand{\sd}{{\p}}

\newcommand{\Def}{\mathsf{Def}}

\newcommand{\HGC}{\mathsf{HGC}}
\newcommand{\wt}{\widetilde}
\newcommand{\GC}{\mathsf{GC}}

\newcommand{\grt}{\fg\fr\ft}

\newcommand{\PROP}{\mathsf{PROP}}

 \newcommand{\Beq}{\begin{equation}}
 \newcommand{\Eeq}{\end{equation}}
 \newcommand{\Beqr}{\begin{eqnarray}}
 \newcommand{\Eeqr}{\end{eqnarray}}
 \newcommand{\Beqrn}{\begin{eqnarray*}}
 \newcommand{\Eeqrn}{\end{eqnarray*}}
 \newcommand{\Ba}{\begin{array}}
 \newcommand{\Ea}{\end{array}}
 \newcommand{\Bi}{\begin{itemize}}
 \newcommand{\Ei}{\end{itemize}}
 \newcommand{\Bc}{\begin{center}}
 \newcommand{\Ec}{\end{center}}
 \newcommand{\fg}{{\mathfrak g}}

\newcommand{\fr}{{\mathfrak r}}
\newcommand{\fs}{{\mathfrak s}}

\newcommand{\ft}{{\mathfrak t}}

 \newcommand{\f}{{\mathcal O}}
 \newcommand{\cA}{{\mathcal A}}
 \newcommand{\cB}{{\mathcal B}}
 \newcommand{\cC}{{\mathcal C}}
 \newcommand{\caD}{{\mathcal D}}
 \newcommand{\cE}{{\mathcal E}}
 \newcommand{\cF}{{\mathcal F}}
 \newcommand{\cG}{{\mathcal G}}
 \newcommand{\caH}{{\mathcal H}}

 \newcommand{\caL}{{\mathcal L}}
 \newcommand{\cM}{{\mathcal M}}
 
 \newcommand{\cP}{{\mathcal P}}
 
 \newcommand{\cR}{{\mathcal R}}
 \newcommand{\cS}{{\mathcal S}}
 \newcommand{\cT}{{\mathcal T}}

 \newcommand{\cV}{{\mathcal V}}

 \newcommand{\al}{\alpha}
 \newcommand{\be}{\beta}
 \newcommand{\ga}{\gamma}
 
 \newcommand{\Ga}{\Gamma}

 \newcommand{\la}{\lambda}

 \newcommand{\Img}{{\mathsf I\mathsf m}\, }
 \newcommand{\Hom}{{\mathrm H\mathrm o\mathrm m}}
 
 \newcommand{\sip}{\smallskip}
 \newcommand{\bip}{\bigskip}
 \newcommand{\mip}{\vspace{2.5mm}}
\newcommand{\wh}{\widehat}

\newcommand{\LB}{\mathcal{L}\mathit{ieb}}
\newcommand{\LoB}{\mathcal{L}\mathit{ieb}^\diamond}

\newcommand{\LBd}{\mathcal{L}\mathit{ieb}_{d}}
\newcommand{\LBcd}{\mathcal{L}\mathit{ieb}_{c,d}}

\newcommand{\HoLBd}{\mathcal{H}\mathit{olieb}_{d}}
\newcommand{\HoLBcd}{\mathcal{H}\mathit{olieb}_{c,d}}

\newcommand{\wHoLBcd}{\widehat{\mathcal{H}\mathit{olieb}}_{c,d}}

\newcommand{\LoBcd}{\mathcal{L}\mathit{ieb}_{c,d}^\diamond}

\newcommand{\HoLoBd}{\mathcal{H}\mathit{olieb}_{d}^\diamond}
\newcommand{\HoLoBcd}{\mathcal{H}\mathit{olieb}_{c,d}^\diamond}
\newcommand{\HoLoBpcd}{\mathcal{H}\mathit{olieb}_{c,d}^{\diamond +}}
\newcommand{\HoLoB}{\mathcal{H}\mathit{olieb}^\diamond}

\newcommand{\HoLB}{\mathcal{H}\mathit{olieb}}

\newcommand{\Holie}{\mathcal{H} \mathit{olie}}
\newcommand{\Lie}{\mathcal{L} \mathit{ie}}

\newcommand{\HoqLBcd}{\mathcal{H}\mathit{oqlieb}_{c,d}}

\theoremstyle{plain}
\swapnumbers

\newtheorem{prop-def}[theorem]{Proposition-definition}

\newtheorem{f-theorem}{Formality Theorem}[section]
\newtheorem{main-theorem}{Main~Theorem}[section]
\newtheorem{section-theorem}{Theorem}[section]

\theoremstyle{definition}

\begin{document}

 \sloppy

 \newenvironment{proo}{\begin{trivlist} \item{\sc {Proof.}}}
  {\hfill $\square$ \end{trivlist}}

\long\def\symbolfootnote[#1]#2{\begingroup%
\def\thefootnote{\fnsymbol{footnote}}\footnote[#1]{#2}\endgroup}

\title{Twisting of properads}

\author{Sergei~A. Merkulov}
\address{Sergei~Merkulov:  Department of Mathematics, University of Luxembourg,  Grand Duchy of Luxembourg }
\email{sergei.merkulov@uni.lu}

\begin{abstract} We study T.\ Willwacher's twisting endofunctor $\tw$ in the category of dg prop(erad)s $\cP$ under the operad of (strongly homotopy) Lie algebras, $i:\Lie\rar \cP$. It is proven that if $\cP$ is a properad under properad of Lie bialgebras $\LB$, then the associated twisted properad  $\tw\cP$ becomes in general a properad under quasi-Lie bialgebras (rather than under $\LB$). This result implies that the cyclic cohomology of any cyclic homotopy associative algebra has in general an induced structure of a quasi-Lie bialgebra. We show that the cohomology of the twisted properad $\tw\LB$ is highly non-trivial --- it contains the cohomology of the so called haired graph complex introduced and studied recently in the context of the theory of long knots and the theory of moduli spaces $\cM_{g,n}$ of algebraic curves of arbitrary genus $g$ with $n$ punctures.

\sip

Using a polydifferential functor from the category of props to the category of operads, we introduce  the notion of a Maurer-Cartan element of a strongly homotopy Lie bialgebra, and use it to construct
 a new  twisting endofunctor $\TW$ in the category dg prop(erad)s $\cP$ under $\HoLB$, the minimal resolution of $\LB$. We prove that $\TW \HoLB$ is quasi-isomorphic to $\LB$, and establish its relation to the homotopy theory of triangular Lie bialgebras.
 It is proven that the dg Lie algebra $\Def(\HoLB\stackrel{}{\rar} \cP)$ controlling deformations of $i$ acts on $\TW\cP$ by derivations. In some important examples this dg Lie algebra has a rich  and interesting cohomology
(containing, for example, the Grothendieck-Teichm\"uller Lie algebra).

\sip

Finally, we introduce a diamond version $\TW^\lozenge$ of $\TW$ which works in the category of dg properads under {\em involutive}\, (strongly homotopy) Lie bialgebras, and discuss its applications in string topology.

\sip

%\sip
%\noindent {\sc Mathematics Subject Classifications} (2000). 17B37, 16W30, 51M20.

%\noindent {\sc Key words}. Properads, Lie bialgebras,  string topology
\end{abstract}

 \maketitle
 \markboth{}{}

{\small
{\small
\tableofcontents
}
}

{\Large
\section{\bf Introduction}
}

\bip

In this paper we study some new aspects of a well-known twisting endofunctor $\tw$ \cite{W}  in a certain subcategory of the category of properads, and introduce a new one. Both of them have applications in several areas of modern research --- the string topology, the theory of moduli spaces of algebraic curves, the theory of cyclic strongly homotopy associative algebras etc --- which we discuss below.

\sip

Let $\cP$ be a properad under the operad $\Lie_d$ of (degree $d\in \Z$ shifted) Lie algebras, that is, the one equipped with a morphism
$$
i: \Lie_d \lon \cP.
$$
Thomas Willwacher introduced in \cite{W} a {\it twisting endofunctor}\,
$$
\tw: \cP \lon \tw\cP
$$
 in the category of such properads; the twisted properad $\tw\cP$ is obtained from $\cP$ by adding to it a new generator degree $d$ generator $\Ba{c}\resizebox{1.5mm}{!}{\begin{xy}
 <0mm,0.5mm>*{};<0mm,4mm>*{}**@{-},
 <0mm,0mm>*{\bullet};<0mm,0mm>*{}**@{},
 \end{xy}}\Ea$ with no inputs and precisely one output encoding the defining property
 of a Maurer-Cartan element of a generic $\caL ie_d$-algebra. This twisting construction originated in the formality theory of the chain operad of chains of little disks operad  \cite{Ko2, LV, W} and found many other important applications including the Deligne conjecture \cite{DW}, the homotopy theory of configuration spaces \cite{CW} and the theory of moduli spaces of algebraic curves \cite{MW1, Me2}. The twisted properad comes equipped with a canonical morphism
 $$
 \Lie_d \lon \tw\cP
 $$
so the twisting construction is indeed an endofunctor in the category, $\PROP_{\Lie_d}$, of operads under $\Lie_d$. It can be naturally extended \cite{W} to the category $\PROP_{\Holie_d}$ of properads under $\Holie_d$, the minimal resolution of $\Lie_d$. It is proven in \cite{DW} that $\tw \Holie_d$ is quasi-isomorphic to $\Lie_d$; the cohomologies of twisted versions of some other classical operads under $\Lie_d$ have been computed in \cite{DW, DSV}.

\sip

The main purpose of this paper is to study the restriction of  T.\ Willwacher's  twisting endofunctor $\tw$,
 $$
 \tw: \PROP_{\LBcd} \lon \PROP_{\Lie_d}
 $$
 to the subcategory $\PROP_{\LBcd} \subset \PROP_{\Lie_d}$  of properads under the properad $\LBcd$ of (degree shifted) Lie bialgebras, i.e.\ the ones which come equipped with a non-trivial morphism of properads,
\Beq\label{1: i from LBcd to P}
i: \LBcd\lon \cP.
\Eeq
and then to appropriately modify $\tw \rightsquigarrow \Tw$ such that the new functor $\Tw$ which becomes an {\em endofuntor}\, of $ \PROP_{\LBcd}$. Everything will work, of course, in the category $\PROP_{\HoLBcd}$ of properads under
$\HoLBcd$, the minimal resolution of $\LBcd$.

\sip

The notion of Lie bialgebra was introduced by Vladimir Drinfeld in \cite{Dr1}  in the context of the theory of Yang-Baxter equations  and the deformation theory  of universal enveloping algebras. This notion and its involutive version have found since many applications in algebra, string topology, contact homology, theory of associators and the theory of Riemann surfaces with punctures. The properad $\LBcd$ controls Lie bialgebras with Lie bracket of degree $1-d$ and Lie-cobracket of degree $1-c$; the properads $\LBcd$ with the same parity of $c+d\in \Z$ are isomorphic to each other up to degree shift so that there are essentially two different types of such properads, even and odd ones.
 %\cite{Tu} [Turaev, CS, Cielback etc, Schedler, ET...]

\sip

In general, given $\cP\in \PROP_{\LBcd}$,  the associated twisted properad $\tw\cP$ is no more a properad under $\LBcd$ ---  a generic Maurer-Cartan element of the Lie bracket need not to respect the Lie cobracket. Surprisingly enough, the cohomology of any twisted properad $\tw\cP$ always comes equipped with an induced from $i$ morphism of properads,
\Beq\label{1: i^q map}
i^q: q\LB_{c-1,d}\lon H^\bu(\tw\cP),
\Eeq
 where $q\LB_{c-1,d}$ is the properad of (degree shifted) {\it quasi}-Lie bialgebras which have been also introduced by Vladimir Drinfeld in \cite{Dr1} in the context of the theory of quantum groups. The map  $i^q$ is described explicitly in   Theorem  {\ref{2: qLien to twP}} below. In the special case when $\cP$ is the properad of ribbon graphs $\RGra_d$ introduced in \cite{MW1}, the map $i^q$ has been found (in a slightly different but equivalent form) in \cite{Me2}. As $\tw\RGra_d$ acts canonically (almost by its very construction), on the reduced cyclic cohomology $H^\bu(Cyc(A))$ of an arbitrary cyclic strongly homotopy associative algebra $A$ (equipped with the degree $-d$ scalar product), we deduce a new observation that {\em $H^\bu(Cyc(A))$ is always a quasi-Lie bialgebra}, see \S {\ref{3: subsec on qLieb on Cyc(A)}} for full details. It is worth noting that  the twisted properad of ribbon graphs $\tw \RGra_d$  controls \cite{Me2}  the totality of compactly supported cohomology groups $\prod_{n\geq 1, 2g+n\geq 3} H_c^{\bu- d(2g-2+n)}(\cM_{g,n})$ of moduli spaces $\cM_{g,n}$ of genus $g$ algebraic curves with $m$ boundaries and $n$ punctures, and the associated map
$$
i^q: q\LB_{d-1,d}\lon H(\tw\RGra_d)\simeq \prod_{n\geq 1, 2g+n\geq 3} H_c(\cM_{g,n})
$$
is non-trivial on infinitely many elements of $q\LB_{-1,0}$ (see \S 3.9 in \cite{Me2}).

\sip

 The deformation complex
$$
\Def\left(q\LB_{c-1,d}\stackrel{i^q}{\lon} H^\bu(\tw\cP)\right)
$$
of the morphism $i^q$ has, in general, a much richer cohomology than the complex $\Def(\Lie_{d}\stackrel{i}{\lon} \tw\cP)$; moreover,  that cohomology comes always equipped with a morphism of cohomology groups,
$$
H^\bu(\GC^{\geq 2}_{c+d-1})\lon H^\bu\left(
\Def\left(q\LB_{c-1,d}\stackrel{i^q}{\lon} H^\bu(\tw\cP)\right)\right)
$$
where $\GC^{\geq 2}_n$ is the famous Maxim Kontsevich graph complex \cite{Ko1} (more precisely, its extension allowing graphs with bivalent vertices). The case $c+d=3$ is of special interest as the dg Lie algebra
$H^\bu(\GC^{\geq 2}_{2})$ contains the Grothendieck-Teichm\"uller Lie algebra \cite{W}. The case $c+d=2$ is also of interest as it corresponds to the odd Kontsevich graph complex $H^\bu(\GC_1^{\geq 2})$ which contains a rich subspace generated by trivalent graphs.
 
\sip

Another important example of a dg properad in the subcategory $\PROP_{\LBcd}$ is the properad
 $\LBcd$ itself; the cohomology $H^\bu(\tw\LBcd)$ of the associated twisted properad  is highly non-trivial:
we show in \S {\ref{3: subsec on reduced twisting of (prop)erads under HoLBcd}}  that, for any natural number $N\geq 1$, there is an injection of cohomology groups,
 $$
 H^\bu(\HGC_{c+d}^N) \lon H^{\bu + dN}(\tw \LB_{c,d})
 $$
where $\HGC_d^N$ is a version of the Kontsevich graph complex $\GC_d$ with $N$ labelled hairs which has been  introduced and studied recently in the context of the theory of moduli spaces $\cM_{g,n}$ of algebraic curves of arbitrary genus $g$ with $n$ punctures \cite{CGP} and the theory of long knots \cite{FTW}.

\sip

The functor $\tw$ is {\em not}\, an endofunctor of the category $\PROP_{\HoLBcd}$
and, contrary to the canonical projection $\tw \Holie_d \rar \Holie_d$, the analogous ``forgetful" map
$$
\tw \HoLBcd \lon \HoLBcd
$$
is {\em not}\, a quasi-isomorphism (as the above result with the haired graph complex demonstrates). 
In \S 4 we introduce a new twisting endofunctor
 $$
 \Tw: \PROP_{\HoLBcd} \lon \PROP_{\HoLBcd}
 $$  which is an enlargement of $\tw$ in the sense of the number of new generators (hence the notation) and which fixes both these ``not"s .
The key point is to introduce the correct notion of a Maurer-Cartan element of a generic $\HoLBcd$-algebra. The idea is to use a polydifferential functor \cite{MW1}
$$
\Ba{rccc}
\f: & \text{\sf Category of dg props} & \lon & \text{\sf Category of dg operads}\\
&     \cP   &\lon & \f\cP
\Ea
$$
whose main property is that,
given any dg prop $\cP$ and its any representation
 in a dg vector space $V$, the dg operad $\f\cP$ comes equipped canonically with an induced representation
 in the  graded commutative tensor algebra ${\odot^\bu} V$ given in terms of polydifferential --- with respect to the standard multiplication in ${\odot^\bu} V$ --- operators. The point is that the dg operad $\f\HoLBcd$ comes equipped with a highly non-trivial morphism morphism of dg properads which was discovered in \cite{MW1},
 $$
 \Holie_{c+d}^+ \lon \f\HoLBcd,
 $$
and whose image brings into play {\em all}\, generators of the properad $\HoLBcd$, not just the
ones spanning the sub-properad $\Holie_d$ of $\HoLBcd$. Here the symbol $+$ means a slight extension of $\Holie_d$ which take cares about deformations of the differential in representation spaces \cite{Me1}. It makes sense to talk about Maurer-Cartan elements of representations of $\Holie_{c+d}^+$ as usual, and hence it makes sense to talk about {\em Maurer-Cartan elements $\ga\in {\odot^{\bu\geq 1}}( V[c])$ of an arbitrary $\HoLBcd$-algebra}\, $V$ via the above morphism. Rather surprisingly, these MC elements $\ga$ can be used to twist not only the dg operad $\f\HoLBcd$ but  the dg properad $\HoLBcd$ itself giving thereby rise to a new twisting endofunctor $\Tw$ on the category $\PROP_{\HoLBcd}$! As explained in \S 4, the twisting endofunctor $\Tw$ adds to a generic properad $\cP$ under $\HoLBcd$ infinitely many new (skew)symmetric generators,
$$
\Ba{c}\resizebox{14mm}{!}{\begin{xy}
 <0mm,0mm>*{\bu};
 <-0.6mm,0.44mm>*{};<-8mm,5mm>*{^1}**@{-},
 <-0.4mm,0.7mm>*{};<-4.5mm,5mm>*{^2}**@{-},
 <0mm,5mm>*{\ldots},
 <0.4mm,0.7mm>*{};<4.5mm,5mm>*{^{}}**@{-},
 <0.6mm,0.44mm>*{};<8mm,5mm>*{^m}**@{-},
 \end{xy}}\Ea =
 (-1)^{c|\sigma|}
 \Ba{c}\resizebox{17mm}{!}{\begin{xy}
 <0mm,0mm>*{\bu};
 <-0.6mm,0.44mm>*{};<-11mm,7mm>*{^{\sigma(1)}}**@{-},
 <-0.4mm,0.7mm>*{};<-4.5mm,7mm>*{^{\sigma(2)}}**@{-},
 <0mm,5mm>*{\ldots},
 <0.4mm,0.7mm>*{};<4.5mm,6mm>*{}**@{-},
 <0.6mm,0.44mm>*{};<10mm,7mm>*{^{\sigma(m)}}**@{-},
 \end{xy}}\Ea \ \ \ \forall \sigma\in \bS_m, \ m\geq 1.
$$
with the $m=1$ corolla $\Ba{c}\resizebox{1.5mm}{!}{\begin{xy}
 <0mm,0.5mm>*{};<0mm,4mm>*{}**@{-},
 <0mm,0mm>*{\bullet};<0mm,0mm>*{}**@{},
 \end{xy}}\Ea$  corresponding to the original functor $\tw$. We show in \S {\ref{4: subsec on trinagular}} that the quotient of $\Tw \HoLBcd$ by the ideal generated by  $\Ba{c}\resizebox{1.5mm}{!}{\begin{xy}
 <0mm,0.5mm>*{};<0mm,4mm>*{}**@{-},
 <0mm,0mm>*{\bullet};<0mm,0mm>*{}**@{},
 \end{xy}}\Ea$  gives us a  dg free properad closely related to the properad of strongly homotopy {\em triangular}\, Lie bialgebras.

 \sip

It is proven in \S 4 that the natural projection
$$
\Tw\HoLBcd \lon \HoLBcd
$$
is a quasi-isomorphism. We also prove for any $\cP\in \PROP_{\HoLBcd}$ t the dg Lie algebra
 $\Def(\HoLBcd\stackrel{i}{\rar} \cP)$ controlling deformations of the given morphism
 $
 i: \HoLBcd \rar \cP$
 acts on $\Tw\cP$ by derivations. In the case $\cP=\HoLBcd$ the associated complex
 $$
 \Def(\HoLBcd\stackrel{\Id}{\rar} \HoLBcd)
 $$
has cohomology equal (up to one rescaling class) to $H^\bu(\GC_{c+d+1}^{\geq 2})$ \cite{MW2} so that, for any dg properad $\cP$ under $\HoLBcd$ there is always a morphism of cohomology groups
$$
H^\bu(\GC_{c+d+1}^{\geq 2}) \lon H^\bu(\Def(\HoLBcd\stackrel{i}{\rar} \cP))
$$
and hence an action of $H^\bu(\GC_{c+d+1}^{\geq 2})$  on the cohomology of the twisted dg properad  $\Tw\cP$ by derivations (which can be in concrete cases homotopy trivial).

\sip

A similar trick via the polydifferential functor $\f$ works fine in the case of strongly homotopy {\em involutive}\,  Lie bialgebras; we denote the associate properad by $\LoBcd$ and its minimal resolution by $\HoLoBcd$. There is again a highly non-trivial morphism of dg properads \cite{MW1},
$$
\Holie^{\diamond +}_{c+d} \to \f\HoLoBcd,
$$
where $\Holie^{\diamond +}$ is a {\em diamond}\, extension of  $\Holie_d$ which was introduced and studied in \cite{CMW}. Maurer-Cartan elements of $\Holie^{\diamond +}$-algebras are defined in the standard way so that the above morphism of properads can be immediately translated into the notion of {\em Maurer-Cartan element of a $\HoLoBcd$-algebra}. However this time we obtain essentially nothing new: this approach just re-discovers the notion which has been introduced earlier in general in \cite{CFL}, and in the special case of $\LoBcd$-algebras of cyclic words in \cite{Ba1,Ba2}. Therefore the diamond extension $\Tw^\diamond$ of the twisting endofunctor $\Tw$ from \S 3 gives us essentially nothing new as well: we obtain just a properadic incarnation of the twisting constructions in \cite{Ba1,Ba2,CFL}.  This incarnation fits nicely the beautiful approach to string topology developed in \cite{NW} via the so called partition functions of closed manifolds. Full details are given in \S 5.

\sip

%%%%%%%%%%%%%%%%%%%%%%%%%%%%%%%%%%%%%%%%%%%%%%%

{\bf Notation}. We work over a field $\K$ of characteristic zero.
 The set $\{1,2, \ldots, n\}$ is abbreviated to $[n]$;  its group of automorphisms is
denoted by $\bS_n$;
%for a non-negative integer $d$, the symbol  $\sgn^d$ stands for the
%one-dimensional representation of $\bS_n$ on which $\bS_n$ acts trivially if $d$ is even
%and by sign if $d$ is odd.
the trivial (resp., sign) one-dimensional representation of
 $\bS_n$ is denoted by $\id_n$ (resp., by $\sgn_n$). We often abbreviate $\sgn_n^d:= \sgn_n^{\ot |d|}$, $d\in \Z$.
 %; we also set, for $d\in \Z$,
 %$$
 %\sgn^d_n:=\left\{ \Ba{ll} \id_n & \mathrm{if}\ d\ \mathrm{is\ even}\\
 %\sgn_n & \mathrm{if}\ d\ \mathrm{is\ odd}.
 %\Ea\right.
 %$$
The cardinality of a finite set $A$ is denoted by $\# A$.

\sip

We work throughout in the category of $\Z$-graded vector spaces over a field $\K$
of characteristic zero.
If $V=\oplus_{i\in \Z} V^i$ is a graded vector space, then
$V[d]$ stands for the graded vector space with $V[d]^i:=V^{i+d}$; the canonical isomorphism $V\rar V[d]$ is denoted by $\fs^d$.
% and and $s^d$ for the associated isomorphism $V\rar V[d]$;
for $v\in V^i$ we set $|v|:=i$.
%For a pair of graded vector spaces $V_1$ and $V_2$, the symbol $\Hom_i(V_1,V_2)$ stands
%for the space of homogeneous linear maps of degree $i$, and
%$\Hom(V_1,V_2):=\bigoplus_{i\in \Z}\Hom_i(V_1,V_2)$; for example, $s^k\in \Hom_{-k}(V,V[k])$.

\sip

For a
prop(erad) $\cP$ we denote by $\cP\{d\}$ a prop(erad) which is uniquely defined by
 the following property:
for any graded vector space $V$ a representation
of $\cP\{d\}$ in $V$ is identical to a representation of  $\cP$ in $V[d]$; in particular, one has for an
endomorphism properad $\cE nd_V\{-d\}=\cE nd_{V[d]}$.
Thus a map $\cP\{d\} \rar \cE nd_V$ is the same as $
\cP \rar  \cE nd_{V[d]}\equiv \cE nd_V\{-d\}$. The operad controlling  Lie algebras with Lie bracket of degree $-d$ is denoted by $\caL ie_{d+1}$ while its minimal resolution by $\caH \mathit{olie}_{d+1}$; thus $\caL ie_{d+1}$ is equal to $\caL \mathit{ie}\{d\}$ if ones uses the standard notation $\caL ie:=\Lie_1$ for the ordinary operad of Lie algebras.

%\sip

%{\tiny
%For a right (resp., left) module $V$ over a group $G$ we denote by $V_G$ (resp.\
%$_G\hspace{-0.5mm}V$)
% the $\K$-vector space of coinvariants:
%$V/\{g(v) - v\ |\ v\in V, g\in G\}$ and by $V^G$ (resp.\ $^GV$) the subspace
%of invariants: $\{\forall g\in G\ :\ g(v)=v,\ v\in V\}$. If $G$ is finite, then these
%spaces are canonically isomorphic as $char(\K)=0$.}

\sip

 We often used the following elements
 $$
\oint_{123}\hspace{-1mm}:=\sum_{k=1}^3 (123)^k\in \K[\bS_3], \ \ \ \ \ \
\text{Alt}^d_{\bS_n}:=\sum_{\sigma\in \bS_n} (-1)^{d|\sigma|} \sigma\in \K[\bS_n]
 $$
as linear operators on $\bS_3$- and, respectively, $\bS_n$-modules.

%If $V$ is a left $G_1$-module and a
%right $G_2$-module and actions of $G_1$ and $G_2$ commute, then  the associated
%spaces of coinvariants and invariants are denoted by $_{G_1}\hspace{-0.5mm}V_{G_2}$ and,
%respectively,
%by $^{G_1}\hspace{-0.5mm}V^{G_2}$.

%{\bf Acknowledgement}.
%It is a pleasure to thank ...  for valuable communications.

\bip

\bip

{\Large
\section{\bf Twisting of operads under $\Lie_d$ --- an overview in pictures}
}

\bip

\subsection{Introduction}  This section is a more or less self-contained exposition of Thomas Willwacher's construction \cite{W} of the twisting endofunctor in the category of operads under the operad of Lie algebras. For purely pedagogical purposes, we consider here a new intermediate step based on the ``plus" endofunctor from \cite{Me1}. We try to show all the details (including elementary ones) with emphasis on the action of the deformation complexes on twisted operads. Many calculations in the later sections, where we discuss some new material, are more tedious but analogous to the ones reviewed here.

\subsection{Reminder about $\Holie_d$}
Recall that the operad of degree shifted Lie algebras is defined, for any integer  $d\in \Z$, as a quotient,

$$
\Lie_{d}:=\cF ree\langle e\rangle/\langle\cR\rangle,
$$
of the free prop generated by an  $\bS$-module $e=\{e(n)\}_{n\geq 2}$ with
 all $e(n)=0$ except\footnote{When representing elements of all operads and props
  below as (decorated) graphs we tacitly assume that all edges and legs are {\em directed}\, along the flow going from the bottom of the graph to the top.}
$$
e(2):= \sgn_2^{d}\ot \id_1[d-1]=\mbox{span}\left\langle
\Ba{c}\begin{xy}
 <0mm,0.66mm>*{};<0mm,3mm>*{}**@{-},
 <0.39mm,-0.39mm>*{};<2.2mm,-2.2mm>*{}**@{-},
 <-0.35mm,-0.35mm>*{};<-2.2mm,-2.2mm>*{}**@{-},
 <0mm,0mm>*{\bu};<0mm,0mm>*{}**@{},
   %<0mm,0.66mm>*{};<0mm,3.4mm>*{^1}**@{},
   <0.39mm,-0.39mm>*{};<2.9mm,-4mm>*{^{_2}}**@{},
   <-0.35mm,-0.35mm>*{};<-2.8mm,-4mm>*{^{_1}}**@{},
\end{xy}\Ea
=(-1)^{d}
\Ba{c}\begin{xy}
 <0mm,0.66mm>*{};<0mm,3mm>*{}**@{-},
 <0.39mm,-0.39mm>*{};<2.2mm,-2.2mm>*{}**@{-},
 <-0.35mm,-0.35mm>*{};<-2.2mm,-2.2mm>*{}**@{-},
 <0mm,0mm>*{\bu};<0mm,0mm>*{}**@{},
   %<0mm,0.66mm>*{};<0mm,3.4mm>*{^1}**@{},
   <0.39mm,-0.39mm>*{};<2.9mm,-4mm>*{^{_1}}**@{},
   <-0.35mm,-0.35mm>*{};<-2.8mm,-4mm>*{^{_2}}**@{},
\end{xy}\Ea
\right\rangle
$$
modulo the ideal generated by the following relation
\Beq\label{2: Lie operad Jacobi relation}
%%%%%%%%%%%%%% Lie %%%%%%%%%%%%%%%%%%%%%%%%
\oint_{123}\hspace{-1mm}\Ba{c}\resizebox{9mm}{!}{ \begin{xy}
 <0mm,0mm>*{\bu};<0mm,0mm>*{}**@{},
 <0mm,0.69mm>*{};<0mm,3.0mm>*{}**@{-},
 <0.39mm,-0.39mm>*{};<2.4mm,-2.4mm>*{}**@{-},
 <-0.35mm,-0.35mm>*{};<-1.9mm,-1.9mm>*{}**@{-},
 <-2.4mm,-2.4mm>*{\bu};<-2.4mm,-2.4mm>*{}**@{},
 <-2.0mm,-2.8mm>*{};<0mm,-4.9mm>*{}**@{-},
 <-2.8mm,-2.9mm>*{};<-4.7mm,-4.9mm>*{}**@{-},
    <0.39mm,-0.39mm>*{};<3.3mm,-4.0mm>*{^3}**@{},
    <-2.0mm,-2.8mm>*{};<0.5mm,-6.7mm>*{^2}**@{},
    <-2.8mm,-2.9mm>*{};<-5.2mm,-6.7mm>*{^1}**@{},
 \end{xy}}\Ea
 \equiv
 \Ba{c}\resizebox{9mm}{!}{ \begin{xy}
 <0mm,0mm>*{\bu};<0mm,0mm>*{}**@{},
 <0mm,0.69mm>*{};<0mm,3.0mm>*{}**@{-},
 <0.39mm,-0.39mm>*{};<2.4mm,-2.4mm>*{}**@{-},
 <-0.35mm,-0.35mm>*{};<-1.9mm,-1.9mm>*{}**@{-},
 <-2.4mm,-2.4mm>*{\bu};<-2.4mm,-2.4mm>*{}**@{},
 <-2.0mm,-2.8mm>*{};<0mm,-4.9mm>*{}**@{-},
 <-2.8mm,-2.9mm>*{};<-4.7mm,-4.9mm>*{}**@{-},
    <0.39mm,-0.39mm>*{};<3.3mm,-4.0mm>*{^3}**@{},
    <-2.0mm,-2.8mm>*{};<0.5mm,-6.7mm>*{^2}**@{},
    <-2.8mm,-2.9mm>*{};<-5.2mm,-6.7mm>*{^1}**@{},
 \end{xy}}\Ea
 +
\Ba{c}\resizebox{9mm}{!}{ \begin{xy}
 <0mm,0mm>*{\bu};<0mm,0mm>*{}**@{},
 <0mm,0.69mm>*{};<0mm,3.0mm>*{}**@{-},
 <0.39mm,-0.39mm>*{};<2.4mm,-2.4mm>*{}**@{-},
 <-0.35mm,-0.35mm>*{};<-1.9mm,-1.9mm>*{}**@{-},
 <-2.4mm,-2.4mm>*{\bu};<-2.4mm,-2.4mm>*{}**@{},
 <-2.0mm,-2.8mm>*{};<0mm,-4.9mm>*{}**@{-},
 <-2.8mm,-2.9mm>*{};<-4.7mm,-4.9mm>*{}**@{-},
    <0.39mm,-0.39mm>*{};<3.3mm,-4.0mm>*{^2}**@{},
    <-2.0mm,-2.8mm>*{};<0.5mm,-6.7mm>*{^1}**@{},
    <-2.8mm,-2.9mm>*{};<-5.2mm,-6.7mm>*{^3}**@{},
 \end{xy}}\Ea
 +
\Ba{c}\resizebox{9mm}{!}{ \begin{xy}
 <0mm,0mm>*{\bu};<0mm,0mm>*{}**@{},
 <0mm,0.69mm>*{};<0mm,3.0mm>*{}**@{-},
 <0.39mm,-0.39mm>*{};<2.4mm,-2.4mm>*{}**@{-},
 <-0.35mm,-0.35mm>*{};<-1.9mm,-1.9mm>*{}**@{-},
 <-2.4mm,-2.4mm>*{\bu};<-2.4mm,-2.4mm>*{}**@{},
 <-2.0mm,-2.8mm>*{};<0mm,-4.9mm>*{}**@{-},
 <-2.8mm,-2.9mm>*{};<-4.7mm,-4.9mm>*{}**@{-},
    <0.39mm,-0.39mm>*{};<3.3mm,-4.0mm>*{^1}**@{},
    <-2.0mm,-2.8mm>*{};<0.5mm,-6.7mm>*{^3}**@{},
    <-2.8mm,-2.9mm>*{};<-5.2mm,-6.7mm>*{^2}**@{},
 \end{xy}}\Ea
 =0.
\Eeq
Its minimal resolution $\Holie_d$ is a dg free operad whose (skew)symmetric generators,
\Beq\label{2: Lie_inf corolla}
\Ba{c}\resizebox{22mm}{!}{ \xy
(1,-5)*{\ldots},
(-13,-7)*{_1},
(-8,-7)*{_2},
(-3,-7)*{_3},
(7,-7)*{_{n-1}},
(13,-7)*{_n},
 (0,0)*{\bu}="a",
(0,5)*{}="0",
(-12,-5)*{}="b_1",
(-8,-5)*{}="b_2",
(-3,-5)*{}="b_3",
(8,-5)*{}="b_4",
(12,-5)*{}="b_5",
\ar @{-} "a";"0" <0pt>
\ar @{-} "a";"b_2" <0pt>
\ar @{-} "a";"b_3" <0pt>
\ar @{-} "a";"b_1" <0pt>
\ar @{-} "a";"b_4" <0pt>
\ar @{-} "a";"b_5" <0pt>
\endxy}\Ea
=(-1)^d
\Ba{c}\resizebox{23mm}{!}{\xy
(1,-6)*{\ldots},
(-13,-7)*{_{\sigma(1)}},
(-6.7,-7)*{_{\sigma(2)}},
%(-3,-7)*{_{\sigma(3)}},
%(7,-8)*{_{n-1}},
(13,-7)*{_{\sigma(n)}},
 (0,0)*{\bu}="a",
(0,5)*{}="0",
(-12,-5)*{}="b_1",
(-8,-5)*{}="b_2",
(-3,-5)*{}="b_3",
(8,-5)*{}="b_4",
(12,-5)*{}="b_5",
\ar @{-} "a";"0" <0pt>
\ar @{-} "a";"b_2" <0pt>
\ar @{-} "a";"b_3" <0pt>
\ar @{-} "a";"b_1" <0pt>
\ar @{-} "a";"b_4" <0pt>
\ar @{-} "a";"b_5" <0pt>
\endxy}\Ea,
\ \ \ \forall \sigma\in \bS_n,\ n\geq2,
\Eeq
have degrees $1+d-nd$.
 The differential in $\Holie_d$ is given by
\Beq\label{2: d in Lie_infty}
\delta\hspace{-3mm}
\Ba{c}\resizebox{21mm}{!}{\xy
(1,-5)*{\ldots},
(-13,-7)*{_1},
(-8,-7)*{_2},
(-3,-7)*{_3},
(7,-7)*{_{n-1}},
(13,-7)*{_n},
 (0,0)*{\bu}="a",
(0,5)*{}="0",
(-12,-5)*{}="b_1",
(-8,-5)*{}="b_2",
(-3,-5)*{}="b_3",
(8,-5)*{}="b_4",
(12,-5)*{}="b_5",
\ar @{-} "a";"0" <0pt>
\ar @{-} "a";"b_2" <0pt>
\ar @{-} "a";"b_3" <0pt>
\ar @{-} "a";"b_1" <0pt>
\ar @{-} "a";"b_4" <0pt>
\ar @{-} "a";"b_5" <0pt>
\endxy}\Ea
=
\sum_{A\varsubsetneq [n]\atop
\# A\geq 2}\pm
%
% \sum_{[m]=I_1\sqcup I_2\atop
% {|I_1|\geq 0, |I_2|\geq 1}}
% \sum_{[n]=J_1\sqcup J_2\atop
% {|J_1|\geq 1, |J_2|\geq 0}
%}\hspace{0mm}
%(-1)^{d(\# A  + \sigma(A))}
%
\Ba{c}\resizebox{19mm}{!}{\begin{xy}
<10mm,0mm>*{\bu},
<10mm,0.8mm>*{};<10mm,5mm>*{}**@{-},
<0mm,-10mm>*{...},
<14mm,-5mm>*{\ldots},
<13mm,-7mm>*{\underbrace{\ \ \ \ \ \ \ \ \ \ \ \ \  }},
<14mm,-10mm>*{_{[n]\setminus A}};
<10.3mm,0.1mm>*{};<20mm,-5mm>*{}**@{-},
<9.7mm,-0.5mm>*{};<6mm,-5mm>*{}**@{-},
<9.9mm,-0.5mm>*{};<10mm,-5mm>*{}**@{-},
<9.6mm,0.1mm>*{};<0mm,-4.4mm>*{}**@{-},
<0mm,-5mm>*{\bu};
<-5mm,-10mm>*{}**@{-},
<-2.7mm,-10mm>*{}**@{-},
<2.7mm,-10mm>*{}**@{-},
<5mm,-10mm>*{}**@{-},
<0mm,-12mm>*{\underbrace{\ \ \ \ \ \ \ \ \ \ }},
<0mm,-15mm>*{_{A}}.
\end{xy}}
\Ea
\Eeq
%where the summation is taken over the ordered subsets of $[n]$ and $\sigma(A)$ is the sign of the %permutation $[n]\rar A\sqcup {n]\setminus A$; the vertices are ordered in such a way that the upper %vertex comes first. %(1+d-Ad)(1+d+nd-Ad)=1+d-Ad+ ndA
If $d$ is even, all the signs above are equal to $-1$.

\subsection{``Plus" extension}\label{2: subsection on plus functor}
We also consider a dg operad $\HoLB_d^+$ which is an extension of $\HoLBd$ by an extra degree 1 generator $\Ba{c}\resizebox{1.7mm}{!}{\begin{xy}
 <0mm,-0.55mm>*{};<0mm,-3mm>*{}**@{-},
 <0mm,0.5mm>*{};<0mm,3mm>*{}**@{-},
 <0mm,0mm>*{\bullet};<0mm,0mm>*{}**@{},
 \end{xy}}\Ea$ and the differential given by an above formula with the summation running over all possible non-empty subsets $A\subset [n]$.

 \sip

 More generally, there is an endofunctor on the category of dg props (or dg operads)
 introduced in \cite{Me1}

 $$
 \Ba{rccc}
 ^+: & \text{category of dg props} & \lon & \text{category of dg props} \\
     &   (\cP,\sd)  & \lon & (\cP^+, \sd^+)
 \Ea
 $$
 defined as follows. For any dg prop $\cP$,
 let $\cP^+$ be the free prop generated by $\cP$ and one other operation
$\Ba{c}\resizebox{2.2mm}{!}{
\begin{xy}
 <0mm,-0.55mm>*{};<0mm,-4mm>*{}**@{-},
 <0mm,0.5mm>*{};<0mm,4mm>*{}**@{-},
 <0mm,0mm>*{\bullet};<0mm,0mm>*{}**@{},
 \end{xy}}\Ea$ of arity $(1,1)$ and of cohomological degree $+1$. On $\cP^+$ one defines a differential $\sd^+$ by setting its value on the new generator by
 $$
 \sd^+
 \Ba{c}\resizebox{2.0mm}{!}{
\begin{xy}
 <0mm,-0.55mm>*{};<0mm,-4mm>*{}**@{-},
 <0mm,0.5mm>*{};<0mm,4mm>*{}**@{-},
 <0mm,0mm>*{\bullet};<0mm,0mm>*{}**@{},
 \end{xy}}\Ea
  :=-
  \Ba{c}\resizebox{2.8mm}{!}{ \begin{xy}
 <0mm,0mm>*{};<0mm,-4mm>*{}**@{-},
 <0mm,0mm>*{};<0mm,8mm>*{}**@{-},
 <0mm,0mm>*{\bullet};
 <0mm,4mm>*{\bullet};
 \end{xy}
 }\Ea
 $$
 % we order vertices such that the bottom one comes first;
 % differential on vertices of graphs acts from "above"
 and on any other element $a\in \cP(m,n)$ (which we identify pictorially with the $(m,n)$-corolla
 whose vertex is decorated with $a$) by the formula
 \Beq\label{2: delta^+}
 \sd^+
 \Ba{c}\resizebox{15mm}{!}{
 \begin{xy}
 <0mm,0mm>*{\circ};<0mm,0mm>*{}**@{},
 <0mm,0mm>*{};<-8mm,5mm>*{}**@{-},
 <0mm,0mm>*{};<-4.5mm,5mm>*{}**@{-},
 <0mm,0mm>*{};<-1mm,5mm>*{\ldots}**@{},
 <0mm,0mm>*{};<4.5mm,5mm>*{}**@{-},
 <0mm,0mm>*{};<8mm,5mm>*{}**@{-},
   <0mm,0mm>*{};<-8.5mm,5.5mm>*{^1}**@{},
   <0mm,0mm>*{};<-5mm,5.5mm>*{^2}**@{},
   %<0mm,0mm>*{};<4.5mm,5.5mm>*{^{m\hspace{-0.5mm}-\hspace{-0.5mm}1}}**@{},
   <0mm,0mm>*{};<9.0mm,5.5mm>*{^m}**@{},
 <0mm,0mm>*{};<-8mm,-5mm>*{}**@{-},
 <0mm,0mm>*{};<-4.5mm,-5mm>*{}**@{-},
 <0mm,0mm>*{};<-1mm,-5mm>*{\ldots}**@{},
 <0mm,0mm>*{};<4.5mm,-5mm>*{}**@{-},
 <0mm,0mm>*{};<8mm,-5mm>*{}**@{-},
   <0mm,0mm>*{};<-8.5mm,-6.9mm>*{^1}**@{},
   <0mm,0mm>*{};<-5mm,-6.9mm>*{^2}**@{},
   %<0mm,0mm>*{};<4.5mm,-6.9mm>*{^{n\hspace{-0.5mm}-\hspace{-0.5mm}1}}**@{},
   <0mm,0mm>*{};<9.0mm,-6.9mm>*{^n}**@{},
 \end{xy}}\Ea:= \sd
\Ba{c}\resizebox{15mm}{!}{ \begin{xy}
 <0mm,0mm>*{\circ};<0mm,0mm>*{}**@{},
 <0mm,0mm>*{};<-8mm,5mm>*{}**@{-},
 <0mm,0mm>*{};<-4.5mm,5mm>*{}**@{-},
 <0mm,0mm>*{};<-1mm,5mm>*{\ldots}**@{},
 <0mm,0mm>*{};<4.5mm,5mm>*{}**@{-},
 <0mm,0mm>*{};<8mm,5mm>*{}**@{-},
   <0mm,0mm>*{};<-8.5mm,5.5mm>*{^1}**@{},
   <0mm,0mm>*{};<-5mm,5.5mm>*{^2}**@{},
   %<0mm,0mm>*{};<4.5mm,5.5mm>*{^{m\hspace{-0.5mm}-\hspace{-0.5mm}1}}**@{},
   <0mm,0mm>*{};<9.0mm,5.5mm>*{^m}**@{},
 <0mm,0mm>*{};<-8mm,-5mm>*{}**@{-},
 <0mm,0mm>*{};<-4.5mm,-5mm>*{}**@{-},
 <0mm,0mm>*{};<-1mm,-5mm>*{\ldots}**@{},
 <0mm,0mm>*{};<4.5mm,-5mm>*{}**@{-},
 <0mm,0mm>*{};<8mm,-5mm>*{}**@{-},
   <0mm,0mm>*{};<-8.5mm,-6.9mm>*{^1}**@{},
   <0mm,0mm>*{};<-5mm,-6.9mm>*{^2}**@{},
   %<0mm,0mm>*{};<4.5mm,-6.9mm>*{^{n\hspace{-0.5mm}-\hspace{-0.5mm}1}}**@{},
   <0mm,0mm>*{};<9.0mm,-6.9mm>*{^n}**@{},
 \end{xy}}
 \Ea
-
\overset{m-1}{\underset{i=0}{\sum}}
\Ba{c}\resizebox{16mm}{!}{
\begin{xy}
 <0mm,0mm>*{\circ};<0mm,0mm>*{}**@{},
 <0mm,0mm>*{};<-8mm,5mm>*{}**@{-},
 <0mm,0mm>*{};<-3.5mm,5mm>*{}**@{-},
 <0mm,0mm>*{};<-6mm,5mm>*{..}**@{},
 <0mm,0mm>*{};<0mm,5mm>*{}**@{-},
  <0mm,5mm>*{\bullet};
  <0mm,5mm>*{};<0mm,8mm>*{}**@{-},
  <0mm,5mm>*{};<0mm,9mm>*{^{i\hspace{-0.2mm}+\hspace{-0.5mm}1}}**@{},
<0mm,0mm>*{};<8mm,5mm>*{}**@{-},
<0mm,0mm>*{};<3.5mm,5mm>*{}**@{-},
 <0mm,0mm>*{};<6mm,5mm>*{..}**@{},
   <0mm,0mm>*{};<-8.5mm,5.5mm>*{^1}**@{},
   <0mm,0mm>*{};<-4mm,5.5mm>*{^i}**@{},
   <0mm,0mm>*{};<9.0mm,5.5mm>*{^m}**@{},
 <0mm,0mm>*{};<-8mm,-5mm>*{}**@{-},
 <0mm,0mm>*{};<-4.5mm,-5mm>*{}**@{-},
 <0mm,0mm>*{};<-1mm,-5mm>*{\ldots}**@{},
 <0mm,0mm>*{};<4.5mm,-5mm>*{}**@{-},
 <0mm,0mm>*{};<8mm,-5mm>*{}**@{-},
   <0mm,0mm>*{};<-8.5mm,-6.9mm>*{^1}**@{},
   <0mm,0mm>*{};<-5mm,-6.9mm>*{^2}**@{},
   %<0mm,0mm>*{};<4.5mm,-6.9mm>*{^{n\hspace{-0.5mm}-\hspace{-0.5mm}1}}**@{},
   <0mm,0mm>*{};<9.0mm,-6.9mm>*{^n}**@{},
 \end{xy}}\Ea
 + (-1)^{|a|}
\overset{n-1}{\underset{i=0}{\sum}}
 \Ba{c}\resizebox{16mm}{!}{\begin{xy}
 <0mm,0mm>*{\circ};<0mm,0mm>*{}**@{},
 <0mm,0mm>*{};<-8mm,-5mm>*{}**@{-},
 <0mm,0mm>*{};<-3.5mm,-5mm>*{}**@{-},
 <0mm,0mm>*{};<-6mm,-5mm>*{..}**@{},
 <0mm,0mm>*{};<0mm,-5mm>*{}**@{-},
  <0mm,-5mm>*{\bullet};
  <0mm,-5mm>*{};<0mm,-8mm>*{}**@{-},
  <0mm,-5mm>*{};<0mm,-10mm>*{^{i\hspace{-0.2mm}+\hspace{-0.5mm}1}}**@{},
<0mm,0mm>*{};<8mm,-5mm>*{}**@{-},
<0mm,0mm>*{};<3.5mm,-5mm>*{}**@{-},
 <0mm,0mm>*{};<6mm,-5mm>*{..}**@{},
   <0mm,0mm>*{};<-8.5mm,-6.9mm>*{^1}**@{},
   <0mm,0mm>*{};<-4mm,-6.9mm>*{^i}**@{},
   <0mm,0mm>*{};<9.0mm,-6.9mm>*{^n}**@{},
 <0mm,0mm>*{};<-8mm,5mm>*{}**@{-},
 <0mm,0mm>*{};<-4.5mm,5mm>*{}**@{-},
 <0mm,0mm>*{};<-1mm,5mm>*{\ldots}**@{},
 <0mm,0mm>*{};<4.5mm,5mm>*{}**@{-},
 <0mm,0mm>*{};<8mm,5mm>*{}**@{-},
   <0mm,0mm>*{};<-8.5mm,5.5mm>*{^1}**@{},
   <0mm,0mm>*{};<-5mm,5.5mm>*{^2}**@{},
   <0mm,0mm>*{};<4.5mm,5.5mm>*{^{m\hspace{-0.5mm}-\hspace{-0.5mm}1}}**@{},
   <0mm,0mm>*{};<9.0mm,5.5mm>*{^m}**@{},
 \end{xy}}\Ea.
\Eeq
where $\p$ is the original differential in $\cP$.
The dg prop $(\cP^+, \sd^+)$ is uniquely characterized by the property: there is a one-to-one correspondence between representations
 $$
 \rho: \cP^+ \lon \cE nd_V
 $$
of $(\cP^+, \sd^+)$ in a dg vector space $(V,d)$, and representations of $\cP$ in the same space $V$
but equipped with a deformed differential $d+d'$, where $d':=\rho(\begin{xy}
 <0mm,-0.55mm>*{};<0mm,-3mm>*{}**@{-},
 <0mm,0.5mm>*{};<0mm,3mm>*{}**@{-},
 <0mm,0mm>*{\bullet};<0mm,0mm>*{}**@{},
 \end{xy})$.

\subsection{From morphisms $\Holie_d^+$ to twisted morphisms from $\Holie_d$}\label{2: subsec on twisting d in operads} Given any dg operad $(\cA=\{\cA(n)\}_{n\geq 0},\sd)$,
it is well-known that any element $h\in \cA(1)$ defines a derivation $D_h$ of the (non-differential) operad $\cA$ by the formula analogous to (\ref{2: delta^+}),
$$
D_ha= h\circ_1 a - (-1)^{|h||a|} \sum_{i=1}^n a\circ_i h, \ \ \forall \ a\in \cA(n).
$$
Moreover, if $|h|=1$ and
$\sd h=-h\circ_1 h$, then the operator
$$
\pc:= \p+ D_h
$$
is also a differential in $\cA$ (which acts on $h$ by
$
\pc h= h\circ h
$).
%Indeed, for any $a\in \cA(n)$ one has $D_h h=+ 2 h\circ_1 h!$
%\Beqrn
%\p_h^2(a) &=& -\p(D_h a) - D_h(\p a) + D_h^2a\\
%&=& (\p h)\circ_1 a -  \sum_{i=1}^n a\circ_i (\p h) + (h\circ_1 h)\circ_1 a -  %\sum_{i=1}^n a\circ_i (h\circ_1 h)\\
%&=& 0.
%\Eeqrn
%\sip
Assume we have a morphism of dg operads
\Beq\label{2: g^+ map from Holie}
g^+: (\Holie_d^+, \delta^+) \lon (\cA,\sd)
\Eeq
Then the element $h:=g^+(\begin{xy}
 <0mm,-0.55mm>*{};<0mm,-3mm>*{}**@{-},
 <0mm,0.5mm>*{};<0mm,3mm>*{}**@{-},
 <0mm,0mm>*{\bullet};<0mm,0mm>*{}**@{},
 \end{xy})$
 satisfies all the conditions specified above so that the sum
 $$
 \pc:= \p + D_{g^+(\begin{xy}
 <0mm,-0.55mm>*{};<0mm,-3mm>*{}**@{-},
 <0mm,0.5mm>*{};<0mm,3mm>*{}**@{-},
 <0mm,0mm>*{\bullet};<0mm,0mm>*{}**@{},
 \end{xy})}
 $$
is a differential in $\cA$. Hence we have the following
\subsubsection{\bf Proposition}\label{2: Prop on Holie -> O with d+} {\it For any morphism of dg operads (\ref{2: g^+ map from Holie}) there is an associated morphism of dg operads,
$$
g: (\Holie_d, \delta) \lon (\cA, \pc)
$$
given by the restriction of $g^+$ to the generators of $\Holie_d$.}

\begin{proof} Abbreviating $C_n:=\left(\hspace{-2.5mm} \Ba{c}\resizebox{12mm}{!}{  \xy
(-0,6)*{}="1";
(0,+1)*{\bu}="C";
(-7,-7)*+{_1}="L1";
(-3,-7)*+{_2}="L2";
(2,-5)*{...};
(7,-7)*+{_n}="L3";
\ar @{-} "C";"L1" <0pt>
\ar @{-} "C";"L2" <0pt>
\ar @{-} "C";"L3" <0pt>
\ar @{-} "C";"1" <0pt>
 \endxy}
 \Ea
 \hspace{-2.5mm}\right)$, we have for any $n\geq 2$,

\ \ \ \
$
\pc g(C_n)\equiv \pc g^+(C_n)=\p g^+(C_n) + D_{g^+(\begin{xy}
 <0mm,-0.55mm>*{};<0mm,-2mm>*{}**@{-},
 <0mm,0.5mm>*{};<0mm,2mm>*{}**@{-},
 <0mm,0mm>*{\bullet};<0mm,0mm>*{}**@{},
 \end{xy})}(C_n)
 =g^+(\delta^+ C_n)  + D_{g^+(\begin{xy}
 <0mm,-0.55mm>*{};<0mm,-2mm>*{}**@{-},
 <0mm,0.5mm>*{};<0mm,2mm>*{}**@{-},
 <0mm,0mm>*{\bullet};<0mm,0mm>*{}**@{},
 \end{xy})}(C_n)= g(\delta C_n).
$
\end{proof}

\subsection{Twisting of $\Holie_d$ by a Maurer-Cartan element}\label{2: subsec on the def of twA} Let $\wt{\tw}\Holie_d$ be a dg free operad generated by degree $1+d-nd$ corollas  (\ref{2: Lie_inf corolla}) of type $(1,n)$ with $n\geq 2$, and also by an additional corolla  $\Ba{c}\resizebox{1.5mm}{!}{\begin{xy}
 <0mm,0.5mm>*{};<0mm,3.8mm>*{}**@{-},
 <0mm,0mm>*{\bu};<0mm,0mm>*{}**@{},
 \end{xy}}\Ea$ of type $(1,0)$ and of degree $d$. The differential is defined on the $(1,n\geq 2)$ generators  by the standard formula (\ref{2: d in Lie_infty}), while on the new generator it is defined as follows
 \Beq\label{2: delta on (1,0) MC corolla}
 \delta \Ba{c}\resizebox{1.5mm}{!}{\begin{xy}
 <0mm,0.5mm>*{};<0mm,4mm>*{}**@{-},
 <0mm,0mm>*{\bu};<0mm,0mm>*{}**@{},
 \end{xy}}\Ea= -
\sum_{{ k\geq 2}}
\frac{1}{k!}
\Ba{c}\resizebox{13mm}{!}{  \xy
(-5,8)*{}="1";
%
    %(-25,-5)*+{_1}="l1";
(-3,-5)*{...};
 <-5mm,-9mm>*{\underbrace{ \ \ \ \ \ \ \ \ \ \ \ \ \ \ \ \ \ }_{k}},
    (-5,+3)*{\bu}="L";
 (-14,-5)*{\bu}="B";
  (-8,-5)*{\bu}="C";
   (3,-5)*{\bu}="D";
     %<0mm,12mm>*{_m},
\ar @{-} "D";"L" <0pt>
\ar @{-} "C";"L" <0pt>
\ar @{-} "B";"L" <0pt>
\ar @{-} "1";"L" <0pt>
%
 %
 %\ar @{-} "l1";"L" <0pt>
 \endxy}
 \Ea
\Eeq

 \subsubsection{\bf Lemma}\label{2: Lemma on d^2 for MC Lie}  $\delta^2=0$, {\em i.e.\ it is indeed a differential in  $\wt{\tw}\Holie_d$}.
\begin{proof} We have (assuming that $d$ is even to simplify signs)
%but
% $d (a_\circ b)= (da)\circ b + (-1)^a a\circ (db)$... for $d$ even $a$ is odd... so at %second term we should get $-a\circ db$...$db=-...$... we have to be careful about the %ordering of vertices
\Beqrn
 \delta^2\hspace{-1mm} \Ba{c}\resizebox{1.5mm}{!}{\begin{xy}
 <0mm,0.5mm>*{};<0mm,4.2mm>*{}**@{-},
 <0mm,0mm>*{\bu};<0mm,0mm>*{}**@{},
 \end{xy}}\Ea
 &=&
  -\sum_{{ k\geq 2}}
\frac{1}{k!}
\delta\left(
\Ba{c}\resizebox{13mm}{!}{  \xy
(-5,8)*{}="1";
%
    %(-25,-5)*+{_1}="l1";
(-3,-5)*{...};
 <-5mm,-9mm>*{\underbrace{ \ \ \ \ \ \ \ \ \ \ \ \ \ \ \ \ \ }_{k}},
    (-5,+3)*{\bu}="L";
 (-14,-5)*{\bu}="B";
  (-8,-5)*{\bu}="C";
   (3,-5)*{\bu}="D";
     %<0mm,12mm>*{_m},
\ar @{-} "D";"L" <0pt>
\ar @{-} "C";"L" <0pt>
\ar @{-} "B";"L" <0pt>
\ar @{-} "1";"L" <0pt>
%
 %
 %\ar @{-} "l1";"L" <0pt>
 \endxy}
 \Ea\right)\\
 &=& +\sum_{{ k\geq 2}}
\frac{1}{k!}
\left(\sum_{k=k'+k''\atop k'\geq 2,k''\geq 1}
\frac{k!}{k'!k''!}
\Ba{c}\resizebox{16mm}{!}{  \xy
(-5,8)*{}="1";
(-3,-5)*{...};
(-15,-13)*{...};
 <-15mm,-17mm>*{\underbrace{ \ \ \ \ \  \ \ \ \ \ \ }_{k'}},
 <-3mm,-9mm>*{\underbrace{ \ \ \ \ \  \ \ \ \ \ \ }_{k''}},
    (-5,+3)*{\bu}="L";
 (-14,-5)*{\bu}="B";
 (-20,-13)*{\bu}="b1";
 (-11,-13)*{\bu}="b2";
  (-8,-5)*{\bu}="C";
   (3,-5)*{\bu}="D";
     %<0mm,12mm>*{_m},
\ar @{-} "D";"L" <0pt>
\ar @{-} "C";"L" <0pt>
\ar @{-} "B";"L" <0pt>
\ar @{-} "B";"b1" <0pt>
\ar @{-} "B";"b2" <0pt>
\ar @{-} "1";"L" <0pt>
%
 %
 %\ar @{-} "l1";"L" <0pt>
 \endxy}
 \Ea\right)
- \sum_{{ k\geq 2}}
\frac{1}{k!}
\left(\sum_{l\geq 2}
\frac{k}{l!}
\Ba{c}\resizebox{16mm}{!}{  \xy
(-5,8)*{}="1";
(-3,-5)*{...};
(-15,-13)*{...};
 <-15mm,-17mm>*{\underbrace{ \ \ \ \ \  \ \ \ \ \ \ }_{l}},
 <-3mm,-9mm>*{\underbrace{ \ \ \ \ \  \ \ \ \ \ \ }_{k-1}},
    (-5,+3)*{\bu}="L";
 (-14,-5)*{\bu}="B";
 (-20,-13)*{\bu}="b1";
 (-11,-13)*{\bu}="b2";
  (-8,-5)*{\bu}="C";
   (3,-5)*{\bu}="D";
     %<0mm,12mm>*{_m},
\ar @{-} "D";"L" <0pt>
\ar @{-} "C";"L" <0pt>
\ar @{-} "B";"L" <0pt>
\ar @{-} "B";"b1" <0pt>
\ar @{-} "B";"b2" <0pt>
\ar @{-} "1";"L" <0pt>
%
 %
 %\ar @{-} "l1";"L" <0pt>
 \endxy}
 \Ea\right)=0,
 \Eeqrn
 where the first summand comes from (\ref{2: d in Lie_infty}) and the second one from (\ref{2: delta on (1,0) MC corolla}).
 \end{proof}
A representation
$$
\rho: \wt{\tw}\Holie_d \lon \cE nd_V
$$
of $\wt{\tw}\Holie_d$ in a dg (appropriately filtered) vector space $(V,d)$ is given by
 a $\Holie_d$-algebra structure $\{\mu_n\}_{n\geq 1}$ on $V$,
$$
mu_1:=d,\ \  \mu_n:=\rho\left(\hspace{-2mm} \Ba{c}\resizebox{13mm}{!}{  \xy
(-0,6)*{}="1";
(0,+1)*{\bu}="C";
(-7,-7)*+{_1}="L1";
(-3,-7)*+{_2}="L2";
(2,-5)*{...};
(7,-7)*+{_n}="L3";
\ar @{-} "C";"L1" <0pt>
\ar @{-} "C";"L2" <0pt>
\ar @{-} "C";"L3" <0pt>
\ar @{-} "C";"1" <0pt>
 \endxy}
 \Ea
 \hspace{-2mm}\right): \odot^n (V[d]) \rar V[d+1], \ \ n\geq
$$
together with a special element
$
m:= \rho(\hspace{-1mm}\Ba{c}\resizebox{1.4mm}{!}{\begin{xy}
 <0mm,0.5mm>*{};<0mm,5mm>*{}**@{-},
 <0mm,0mm>*{\bu};<0mm,0mm>*{}**@{},
 \end{xy}}\Ea\hspace{-1mm})
$
satisfying the equation (a filtration on $V$ is assumed to be such that this infinite in general sum makes sense)
$$
dm + \sum_{k\geq 2} \frac{1}{k!} \mu_k(m,\ldots, m)=0.
$$
Such an element is called the Maurer-Cartan element of the given $\Holie_d$-algebra structure  on $V$.

\subsubsection{\bf Proposition}\label{2: map from HoLB+ to tilda twA} {\em There is a morphism of dg operads
$$
c^+: (\Holie_d^+,\delta^+)  \lon (\wt{\tw}\Holie_d,\delta)
$$
given on the generators as follows:}
\Beq\label{2: c^+ on Holie-corollas}
\Ba{c}\resizebox{1.5mm}{!}{\begin{xy}
 <0mm,-0.55mm>*{};<0mm,-3mm>*{}**@{-},
 <0mm,0.5mm>*{};<0mm,3mm>*{}**@{-},
 <0mm,0mm>*{\bu};<0mm,0mm>*{}**@{},
 \end{xy}}\Ea
 \stackrel{c^+}{\lon}
 \sum_{{ k\geq 1}}
\frac{1}{k!}
% \sum_{m=p+m_1+...+m_n, n,p\geq 1, m_{\bu}\geq 0 \atop n+p\geq 3}
\Ba{c}\resizebox{22mm}{!}{  \xy
(-25,8)*{}="1";
    (-25,-5)*+{_1}="l1";
(-3,-5)*{...};
 <-5mm,-9mm>*{\underbrace{ \ \ \ \ \ \ \ \ \ \ \ \ \ \ \ \ \ }_{k}},
    (-25,+3)*{\bu}="L";
 (-14,-5)*{\bu}="B";
  (-8,-5)*{\bu}="C";
   (3,-5)*{\bu}="D";
     %<0mm,12mm>*{_m},
\ar @{-} "D";"L" <0pt>
\ar @{-} "C";"L" <0pt>
\ar @{-} "B";"L" <0pt>
\ar @{-} "1";"L" <0pt>
 \ar @{-} "l1";"L" <0pt>
 \endxy}
 \Ea
 , \ \
\Ba{c}\resizebox{14mm}{!}{  \xy
(-0,6)*{}="1";
(0,+1)*{\bu}="C";
(-7,-7)*+{_1}="L1";
(-3,-7)*+{_2}="L2";
(2,-5)*{...};
(7,-7)*+{_n}="L3";
\ar @{-} "C";"L1" <0pt>
\ar @{-} "C";"L2" <0pt>
\ar @{-} "C";"L3" <0pt>
\ar @{-} "C";"1" <0pt>
 \endxy}
 \Ea
\stackrel{c^+}{\lon}
 \sum_{{ k\geq 0}}
\frac{1}{k!}
% \sum_{m=p+m_1+...+m_n, n,p\geq 1, m_{\bu}\geq 0 \atop n+p\geq 3}
\Ba{c}\resizebox{26mm}{!}{  \xy
(-25,8)*{}="1";
(-30,-7)*+{_1}="l2";
    (-27,-7)*+{_2}="l1";
 (-24,-6)*{...};
   (-20,-7)*+{_n}="ln";
(-3,-5)*{...};
 <-5mm,-9mm>*{\underbrace{ \ \ \ \ \ \ \ \ \ \ \ \ \ \ \ \ \ }_{k}},
    (-25,+3)*{\bu}="L";
 (-14,-5)*{\bu}="B";
  (-8,-5)*{\bu}="C";
   (3,-5)*{\bu}="D";
     %<0mm,12mm>*{_m},
\ar @{-} "D";"L" <0pt>
\ar @{-} "C";"L" <0pt>
\ar @{-} "B";"L" <0pt>
\ar @{-} "1";"L" <0pt>
 \ar @{-} "l1";"L" <0pt>
  \ar @{-} "l2";"L" <0pt>
   \ar @{-} "ln";"L" <0pt>
 \endxy}
 \Ea
\ \forall n\geq 2.
\Eeq

\begin{proof}
One has to check that $\delta \circ c^+ = c^+\circ \delta^+$. One has (assuming for simplicity of signs again that $d$ is even)
\Beqrn
\delta\circ c^+(\hspace{-1mm} \Ba{c}\resizebox{1.5mm}{!}{\begin{xy}
 <0mm,-0.55mm>*{};<0mm,-3mm>*{}**@{-},
 <0mm,0.5mm>*{};<0mm,3mm>*{}**@{-},
 <0mm,0mm>*{\bu};<0mm,0mm>*{}**@{},
 \end{xy}}\Ea\hspace{-1mm})
 \hspace{-2mm}
 &=&  %%%%%%%%%%%%%%%% level 1
 \hspace{-2mm}
 \sum_{{ k\geq 1}}
\frac{1}{k!}
% \sum_{m=p+m_1+...+m_n, n,p\geq 1, m_{\bu}\geq 0 \atop n+p\geq 3}
\ \delta \left(\hspace{-1mm} \Ba{c}\resizebox{22mm}{!}{  \xy
(-25,8)*{}="1";
    (-25,-5)*+{_1}="l1";
(-3,-5)*{...};
 <-5mm,-9mm>*{\underbrace{ \ \ \ \ \ \ \ \ \ \ \ \ \ \ \ \ \ }_{k}},
    (-25,+3)*{\bu}="L";
 (-14,-5)*{\bu}="B";
  (-8,-5)*{\bu}="C";
   (3,-5)*{\bu}="D";
     %<0mm,12mm>*{_m},
\ar @{-} "D";"L" <0pt>
\ar @{-} "C";"L" <0pt>
\ar @{-} "B";"L" <0pt>
\ar @{-} "1";"L" <0pt>
 \ar @{-} "l1";"L" <0pt>
 \endxy}
 \Ea\hspace{-1mm} \right)
=
-\sum_{{ k\geq 1}}
\frac{1}{k!}
% \sum_{m=p+m_1+...+m_n, n,p\geq 1, m_{\bu}\geq 0 \atop n+p\geq 3}
\ \left(\sum_{k+1=k'+k''\atop k'\geq 2,k''\geq 0}\frac{k!}{k'!k''!}
\Ba{c}\resizebox{21mm}{!}{  \xy
(-27,9)*{}="1";
(-27,-5)*+{_1}="d";
(-3,-5)*{...};
(-15,-13)*{...};
 <-15mm,-17mm>*{\underbrace{ \ \ \ \ \  \ \ \ \ \ \ }_{k'}},
 <-3mm,-9mm>*{\underbrace{ \ \ \ \ \  \ \ \ \ \ \ }_{k''}},
    (-27,+3)*{\bu}="L";
 (-22,-5)*{\bu}="B";
 (-20,-13)*{\bu}="b1";
 (-11,-13)*{\bu}="b2";
  (-8,-5)*{\bu}="C";
   (3,-5)*{\bu}="D";
     %<0mm,12mm>*{_m},
\ar @{-} "D";"L" <0pt>
\ar @{-} "C";"L" <0pt>
\ar @{-} "B";"L" <0pt>
\ar @{-} "B";"b1" <0pt>
\ar @{-} "B";"b2" <0pt>
\ar @{-} "1";"L" <0pt>
\ar @{-} "d";"L" <0pt>
%
 %
 %\ar @{-} "l1";"L" <0pt>
 \endxy}
 \Ea
 \right.\\
%%%%%%%%%%%%%%%%%%%%%%%%%%%% level 2 %%%%%%%%%%%%%%%%%%%%%%
&& \hspace{-3mm}+ \left.
 \sum_{k=k'+k''\atop k',k''\geq 1}\frac{k!}{k'!k''!}
\Ba{c}\resizebox{21mm}{!}{  \xy
(-27,9)*{}="1";
(-24,-13)*+{_1}="d";
(-3,-5)*{...};
(-15,-13)*{...};
 <-15mm,-17mm>*{\underbrace{ \ \ \ \ \  \ \ \ \ \ \ }_{k'}},
 <-3mm,-9mm>*{\underbrace{ \ \ \ \ \  \ \ \ \ \ \ }_{k''}},
    (-27,+3)*{\bu}="L";
 (-24,-5)*{\bullet}="B";
 (-20,-13)*{\bu}="b1";
 (-11,-13)*{\bu}="b2";
  (-8,-5)*{\bu}="C";
   (3,-5)*{\bu}="D";
     %<0mm,12mm>*{_m},
\ar @{-} "D";"L" <0pt>
\ar @{-} "C";"L" <0pt>
\ar @{-} "B";"L" <0pt>
\ar @{-} "B";"b1" <0pt>
\ar @{-} "B";"b2" <0pt>
\ar @{-} "1";"L" <0pt>
\ar @{-} "d";"B" <0pt>
%
 %
 %\ar @{-} "l1";"L" <0pt>
 \endxy}
 \Ea
\right)
+
\sum_{{ k\geq 1}}
\frac{1}{k!}\hspace{-1mm}
% \sum_{m=p+m_1+...+m_n, n,p\geq 1, m_{\bu}\geq 0 \atop n+p\geq 3}
\ \left(\sum_{l\geq 2}\frac{k}{l!}\hspace{-1mm}
\Ba{c}\resizebox{21mm}{!}{  \xy
(-27,9)*{}="1";
(-27,-5)*+{_1}="d";
(-3,-5)*{...};
(-15,-13)*{...};
 <-15mm,-17mm>*{\underbrace{ \ \ \ \ \  \ \ \ \ \ \ }_{l}},
 <-3mm,-9mm>*{\underbrace{ \ \ \ \ \  \ \ \ \ \ \ }_{k-1}},
    (-27,+3)*{\bu}="L";
 (-22,-5)*{\bullet}="B";
 (-20,-13)*{\bu}="b1";
 (-11,-13)*{\bu}="b2";
  (-8,-5)*{\bu}="C";
   (3,-5)*{\bu}="D";
     %<0mm,12mm>*{_m},
\ar @{-} "D";"L" <0pt>
\ar @{-} "C";"L" <0pt>
\ar @{-} "B";"L" <0pt>
\ar @{-} "B";"b1" <0pt>
\ar @{-} "B";"b2" <0pt>
\ar @{-} "1";"L" <0pt>
\ar @{-} "d";"L" <0pt>
%
 %
 %\ar @{-} "l1";"L" <0pt>
 \endxy}
 \Ea\hspace{-2mm}\right)
 \\
 &=&
-\hspace{-1mm} \sum_{k',k''\geq 1}\frac{1}{k'!k''!}\hspace{-2mm}
\Ba{c}\resizebox{21mm}{!}{  \xy
(-27,9)*{}="1";
(-24,-13)*+{_1}="d";
(-3,-5)*{...};
(-15,-13)*{...};
 <-15mm,-17mm>*{\underbrace{ \ \ \ \ \  \ \ \ \ \ \ }_{k'}},
 <-3mm,-9mm>*{\underbrace{ \ \ \ \ \  \ \ \ \ \ \ }_{k''}},
    (-27,+3)*{\bu}="L";
 (-24,-5)*{\bullet}="B";
 (-20,-13)*{\bu}="b1";
 (-11,-13)*{\bu}="b2";
  (-8,-5)*{\bu}="C";
   (3,-5)*{\bu}="D";
     %<0mm,12mm>*{_m},
\ar @{-} "D";"L" <0pt>
\ar @{-} "C";"L" <0pt>
\ar @{-} "B";"L" <0pt>
\ar @{-} "B";"b1" <0pt>
\ar @{-} "B";"b2" <0pt>
\ar @{-} "1";"L" <0pt>
\ar @{-} "d";"B" <0pt>
%
 %
 %\ar @{-} "l1";"L" <0pt>
 \endxy}
 \Ea
= - c^+\left(\hspace{-2mm}
  \Ba{c}\resizebox{2.8mm}{!}{ \begin{xy}
 <0mm,0mm>*{};<0mm,-4mm>*{}**@{-},
 <0mm,0mm>*{};<0mm,8mm>*{}**@{-},
 <0mm,0mm>*{\bullet};
 <0mm,4mm>*{\bullet};
 \end{xy}
 }\Ea\hspace{-2mm}\right)
=
 c^+\circ \delta^+\left(\hspace{-1mm}\Ba{c}\resizebox{1.5mm}{!}{\begin{xy}
 <0mm,-0.55mm>*{};<0mm,-3mm>*{}**@{-},
 <0mm,0.5mm>*{};<0mm,3mm>*{}**@{-},
 <0mm,0mm>*{\bullet};<0mm,0mm>*{}**@{},
 \end{xy}}\Ea\hspace{-1mm}\right).
\Eeqrn
%%%%%%%%%%%%%%%%%%%%%%%%%%%%%%%%%%%%%%%%%%%%%%%%%%%%%%%%%%%%%%%%
Similarly one checks the required equality for any $n\geq 2$
\Beqrn
\delta\circ c^+\hspace{-1mm}\left(\hspace{-3mm}
\Ba{c}\resizebox{12mm}{!}{  \xy
(-0,6)*{}="1";
(0,+1)*{\bu}="C";
(-7,-7)*+{_1}="L1";
(-3,-7)*+{_2}="L2";
(2,-5)*{...};
(7,-7)*+{_n}="L3";
\ar @{-} "C";"L1" <0pt>
\ar @{-} "C";"L2" <0pt>
\ar @{-} "C";"L3" <0pt>
\ar @{-} "C";"1" <0pt>
 \endxy}
 \Ea
\hspace{-3mm}\right)\hspace{-3mm}
&=& % 1st level
\hspace{-2mm}
 \sum_{{ k\geq 0}}
\frac{1}{k!}\, \delta\left(\hspace{-2mm}
\Ba{c}\resizebox{23mm}{!}{  \xy
(-25,8)*{}="1";
(-30,-7)*+{_1}="l2";
    (-27,-7)*+{_2}="l1";
 (-24,-6)*{...};
   (-20,-7)*+{_n}="ln";
(-3,-5)*{...};
 <-5mm,-9mm>*{\underbrace{ \ \ \ \ \ \ \ \ \ \ \ \ \ \ \ \ \ }_{k}},
    (-25,+3)*{\bu}="L";
 (-14,-5)*{\bu}="B";
  (-8,-5)*{\bu}="C";
   (3,-5)*{\bu}="D";
\ar @{-} "D";"L" <0pt>
\ar @{-} "C";"L" <0pt>
\ar @{-} "B";"L" <0pt>
\ar @{-} "1";"L" <0pt>
 \ar @{-} "l1";"L" <0pt>
  \ar @{-} "l2";"L" <0pt>
   \ar @{-} "ln";"L" <0pt>
 \endxy}
 \Ea\right)
 =
\hspace{-0.5mm} -\sum_{k\geq 1}\frac{1}{k!}\left(
 \sum_{k=k'+k''\atop k'\geq 0,k''\geq 1}\frac{k!}{k'!k''!}
\Ba{c}\resizebox{21mm}{!}{  \xy
(-27,9)*{}="1";
(-31,-13)*{_1}="l1";
 (-27,-12)*{...};
 (-24,-13)*+{_n}="ln";
(-3,-5)*{...};
(-15,-13)*{...};
 <-15mm,-17mm>*{\underbrace{ \ \ \ \ \  \ \ \ \ \ \ }_{k'}},
 <-3mm,-9mm>*{\underbrace{ \ \ \ \ \  \ \ \ \ \ \ }_{k''}},
    (-27,+3)*{\bu}="L";
 (-24,-5)*{\bullet}="B";
 (-20,-13)*{\bu}="b1";
 (-11,-13)*{\bu}="b2";
  (-8,-5)*{\bu}="C";
   (3,-5)*{\bu}="D";
     %<0mm,12mm>*{_m},
\ar @{-} "D";"L" <0pt>
\ar @{-} "C";"L" <0pt>
\ar @{-} "B";"L" <0pt>
\ar @{-} "B";"b1" <0pt>
\ar @{-} "B";"b2" <0pt>
\ar @{-} "1";"L" <0pt>
\ar @{-} "l1";"B" <0pt>
\ar @{-} "ln";"B" <0pt>
%
 %
 %\ar @{-} "l1";"L" <0pt>
 \endxy}
 \Ea
 \right.\\
 %%%%%%%%%%%%%%%%%%%%%%%%%%%%%%%% 2 level %%%%%%%%%%%%%%%%%%%
&& \hspace{-21mm}+\left.
\sum_{k=k'+k''\atop k'\geq 1, k''\geq 0}\frac{k!}{k'!k''!}\sum_{i=1}^n
\Ba{c}\resizebox{23mm}{!}{  \xy
(-27,9)*{}="1";
(-24,-13)*+{_i}="d";
(-24,-5)*+{}="n1";
(-31,-5)*+{}="n2";
(-27.5,-3)*{...};
(-3,-5)*{...};
(-15,-13)*{...};
 <-15mm,-17mm>*{\underbrace{ \ \ \ \ \  \ \ \ \ \ \ }_{k'}},
 <-3mm,-9mm>*{\underbrace{ \ \ \ \ \  \ \ \ \ \ \ }_{k''}},
 <-28mm,-8mm>*{\underbrace{ \ \ \ }_{[n]\setminus i}},
    (-27,+3)*{\bu}="L";
 (-18,-5)*{\bullet}="B";
 (-20,-13)*{\bu}="b1";
 (-11,-13)*{\bu}="b2";
  (-8,-5)*{\bu}="C";
   (3,-5)*{\bu}="D";
\ar @{-} "D";"L" <0pt>
\ar @{-} "C";"L" <0pt>
\ar @{-} "B";"L" <0pt>
\ar @{-} "B";"b1" <0pt>
\ar @{-} "B";"b2" <0pt>
\ar @{-} "1";"L" <0pt>
\ar @{-} "d";"B" <0pt>
\ar @{-} "n1";"L" <0pt>
\ar @{-} "n2";"L" <0pt>
 \endxy}
 \Ea\hspace{-2mm} \right)
 -\sum_{k\geq 2}\frac{1}{k!}
  \left(
  \sum_{k=k'+k''\atop k'\geq 2, k''\geq 0}\frac{k!}{k'!k''!}\sum_{i=1}^n
\Ba{c}\resizebox{23mm}{!}{  \xy
(-27,9)*{}="1";
%(-24,-13)*+{_i}="d";
(-24,-5)*+{_n}="n1";
(-31,-5)*+{_1}="n2";
(-27,-4)*{...};
(-3,-5)*{...};
(-15,-13)*{...};
 <-15mm,-17mm>*{\underbrace{ \ \ \ \ \  \ \ \ \ \ \ }_{k'}},
 <-3mm,-9mm>*{\underbrace{ \ \ \ \ \  \ \ \ \ \ \ }_{k''}},
 %<-28mm,-8mm>*{\underbrace{ \ \ \ }_{[n]\setminus i}},
%
    (-27,+3)*{\bu}="L";
 (-18,-5)*{\bullet}="B";
 (-20,-13)*{\bu}="b1";
 (-11,-13)*{\bu}="b2";
  (-8,-5)*{\bu}="C";
   (3,-5)*{\bu}="D";
\ar @{-} "D";"L" <0pt>
\ar @{-} "C";"L" <0pt>
\ar @{-} "B";"L" <0pt>
\ar @{-} "B";"b1" <0pt>
\ar @{-} "B";"b2" <0pt>
\ar @{-} "1";"L" <0pt>
%\ar @{-} "d";"B" <0pt>
\ar @{-} "n1";"L" <0pt>
\ar @{-} "n2";"L" <0pt>
 \endxy}
 \Ea\hspace{-2mm} \right)
 \\
%%%%%%%%%%%%%%%%%%%%%%%%%%%%%%%%%%%%%%%%%%% 3 level %%%%%%%%%%%%%%%%%%%%
&&\hspace{-21mm}
+\sum_{{ k\geq 1}}
\frac{1}{k!}
% \sum_{m=p+m_1+...+m_n, n,p\geq 1, m_{\bu}\geq 0 \atop n+p\geq 3}
\ \left(\sum_{l\geq 2}\frac{k}{l!}
\Ba{c}\resizebox{23mm}{!}{  \xy
(-27,9)*{}="1";
%(-24,-13)*+{_i}="d";
(-24,-5)*+{_n}="n1";
(-31,-5)*+{_1}="n2";
(-27,-4)*{...};
(-3,-5)*{...};
(-15,-13)*{...};
 <-15mm,-17mm>*{\underbrace{ \ \ \ \ \  \ \ \ \ \ \ }_{l}},
 <-3mm,-9mm>*{\underbrace{ \ \ \ \ \  \ \ \ \ \ \ }_{k-1}},
    (-27,+3)*{\bu}="L";
 (-18,-5)*{\bullet}="B";
 (-20,-13)*{\bu}="b1";
 (-11,-13)*{\bu}="b2";
  (-8,-5)*{\bu}="C";
   (3,-5)*{\bu}="D";
\ar @{-} "D";"L" <0pt>
\ar @{-} "C";"L" <0pt>
\ar @{-} "B";"L" <0pt>
\ar @{-} "B";"b1" <0pt>
\ar @{-} "B";"b2" <0pt>
\ar @{-} "1";"L" <0pt>
%\ar @{-} "d";"B" <0pt>
\ar @{-} "n1";"L" <0pt>
\ar @{-} "n2";"L" <0pt>
 \endxy}
 \Ea\hspace{-2mm}\right)
  - \sum_{k\geq 0}\frac{1}{k!}
\sum_{k=k'+k''\atop k',k''\geq 0}\frac{k!}{k'!k''!}\sum_{[n]=I'\sqcup I'' \atop \# I', \#I'' \geq 2}
\Ba{c}\resizebox{23mm}{!}{  \xy
(-27,9)*{}="1";
(-26,-13)*+{}="d1";
(-20,-13)*+{}="d2";
(-22,-11.2)*{...};
(-24,-5)*+{}="n1";
(-31,-5)*+{}="n2";
(-27.5,-3)*{...};
(-3,-5)*{...};
(-15,-13)*{...};
 <-11mm,-17mm>*{\underbrace{  \ \ \  \ \ \ \ \ \ }_{k'}},
 <-3mm,-9mm>*{\underbrace{ \ \ \ \ \  \ \ \ \ \ \ }_{k''}},
 <-28mm,-7.5mm>*{\underbrace{ \ \ \ }_{I''}},
  <-22mm,-14mm>*{\underbrace{ \ \ \ }_{I'}},
    (-27,+3)*{\bu}="L";
 (-18,-5)*{\bullet}="B";
 (-15,-13)*{\bu}="b1";
 (-8,-13)*{\bu}="b2";
  (-8,-5)*{\bu}="C";
   (3,-5)*{\bu}="D";
\ar @{-} "D";"L" <0pt>
\ar @{-} "C";"L" <0pt>
\ar @{-} "B";"L" <0pt>
\ar @{-} "B";"b1" <0pt>
\ar @{-} "B";"b2" <0pt>
\ar @{-} "1";"L" <0pt>
\ar @{-} "d1";"B" <0pt>
\ar @{-} "d2";"B" <0pt>
\ar @{-} "n1";"L" <0pt>
\ar @{-} "n2";"L" <0pt>
 \endxy}
 \Ea
 %%%%%%%%%%%%%%%%%%%%%%%%%%%%%%%%%%%%%%%%%%% 4 level %%%%%%%%%%%%%%%%%%
 \\
 &=&\hspace{-2mm} -c^+\left(\hspace{-2mm}
 \Ba{c}\resizebox{12mm}{!}{  \xy
 (-0,11)*{}="1";
(-0,6)*{\bu}="0";
(0,+1)*{\bu}="C";
(-7,-7)*+{_1}="L1";
(-3,-7)*+{_2}="L2";
(2,-5)*{...};
(7,-7)*+{_n}="L3";
\ar @{-} "C";"L1" <0pt>
\ar @{-} "C";"L2" <0pt>
\ar @{-} "C";"L3" <0pt>
\ar @{-} "C";"0" <0pt>
\ar @{-} "1";"0" <0pt>
 \endxy}
 \Ea
  \hspace{-3mm}
 +
\sum_{i=1}^n\hspace{-2mm}
 \Ba{c}\resizebox{15mm}{!}{  \xy
 (-0,5)*{}="1";
(3,-13)*+{_i}="d1";
(0,+1)*{\bu}="C";
(-7,-7)*+{_1}="L1";
(3.3,-7)*{\bu}="Li";
(-3,-7)*+{_2}="L2";
(0,-6)*{...};
(6,-6)*{...};
(10,-7)*+{_n}="L3";
\ar @{-} "C";"L1" <0pt>
\ar @{-} "C";"L2" <0pt>
\ar @{-} "C";"L3" <0pt>
\ar @{-} "C";"Li" <0pt>
\ar @{-} "C";"1" <0pt>
\ar @{-} "d1";"Li" <0pt>
 \endxy}
 \Ea   \hspace{-2mm}   \right) - 0 +  c^+\hspace{-1mm}\left(\delta\hspace{-3mm}
\Ba{c}\resizebox{12mm}{!}{  \xy
(-0,6)*{}="1";
(0,+1)*{\bu}="C";
(-7,-7)*+{_1}="L1";
(-3,-7)*+{_2}="L2";
(2,-5)*{...};
(7,-7)*+{_n}="L3";
\ar @{-} "C";"L1" <0pt>
\ar @{-} "C";"L2" <0pt>
\ar @{-} "C";"L3" <0pt>
\ar @{-} "C";"1" <0pt>
 \endxy}
 \Ea
\hspace{-2mm}\right)
=
  c^+\hspace{-1mm}\left(\delta^+\hspace{-3mm}
\Ba{c}\resizebox{12mm}{!}{  \xy
(-0,6)*{}="1";
(0,+1)*{\bu}="C";
(-7,-7)*+{_1}="L1";
(-3,-7)*+{_2}="L2";
(2,-5)*{...};
(7,-7)*+{_n}="L3";
\ar @{-} "C";"L1" <0pt>
\ar @{-} "C";"L2" <0pt>
\ar @{-} "C";"L3" <0pt>
\ar @{-} "C";"1" <0pt>
 \endxy}
 \Ea
\hspace{-2mm}\right).
\Eeqrn
\end{proof}
Hence by Proposition {\ref{2: Prop on Holie -> O with d+}}, the differential in the operad $\wt{\tw}\Holie_d$ can be twisted,
$$
\delta\rar
 \delta_{\centerdot}=\delta + D_{c^+(\hspace{-1.4mm}\Ba{c}\resizebox{1mm}{!}{\begin{xy}
 <0mm,-0.55mm>*{};<0mm,-3mm>*{}**@{-},
 <0mm,0.5mm>*{};<0mm,3mm>*{}**@{-},
 <0mm,0mm>*{\bullet};<0mm,0mm>*{}**@{},
 \end{xy}}\Ea\hspace{-1.4mm})}.
  $$
  The operad $\wt{\tw}\Holie_d$ equipped  with the twisted differential $\delta_{\centerdot}$ is denoted from now on by  $\tw\Holie_d$.

 \subsubsection{\bf Definition-proposition} The data $\tw\Holie_d:= \{\tw\Holie_d(n)\}_{n\geq 0}, \delta_{\centerdot})$
 is called the {\it twisted operad of strongly homotopy Lie algebras}.
 %\footnote{Note that  elements $\tw\Holie_d$ have at least one input, i.e.\ elements with zero inputs %are excluded (as a matter of convention, cf.\ \cite{W}).}.
 It comes equipped with a
 monomorphism
 $$
 c: (\Holie_d,\delta) \lon (\tw \Holie_d, \delta_{\centerdot})
 $$
 given on the generators of $\Holie_d$ by the second expression in formula (\ref{2: c^+ on Holie-corollas}).

\subsection{Twisting of operads under $\Holie_d$ \cite{W}}
Let $(\cA, \sd)$ be a dg operad equipped with a non-trivial morphism of operads
$$
i: \HoLB_d \lon \cA
$$
Such an operad is called an {\it operad under $\Holie_d$}. {\em Generic}\, elements of $\cA=\{\cA(n)\}_{n\geq 0}$ are denoted in this paper as decorated corollas with, say, white vertices (to distinguish them from generators of $\Holie_d$),
$$
\Ba{c}\resizebox{21mm}{!}{ \xy
(1,-5)*{\ldots},
(-13,-7)*{_1},
(-8,-7)*{_2},
(-3,-7)*{_3},
(7,-7)*{_{n-1}},
(13,-7)*{_n},
 (0,0)*{\circ}="a",
(0,5)*{}="0",
(-12,-5)*{}="b_1",
(-8,-5)*{}="b_2",
(-3,-5)*{}="b_3",
(8,-5)*{}="b_4",
(12,-5)*{}="b_5",
\ar @{-} "a";"0" <0pt>
\ar @{-} "a";"b_2" <0pt>
\ar @{-} "a";"b_3" <0pt>
\ar @{-} "a";"b_1" <0pt>
\ar @{-} "a";"b_4" <0pt>
\ar @{-} "a";"b_5" <0pt>
\endxy}\Ea \in \cA(n),\ \ \ n\geq 1.
$$
The images of the generators (\ref{2: Lie_inf corolla}) of $\Holie_d$ under the map $i$ are denoted by decorated corollas with vertices shown as $\circledcirc$ (to emphasize the special status of these elements of $\cA$),
$$
\Ba{c}\resizebox{22mm}{!}{ \xy
(1,-5)*{\ldots},
(-13,-7)*{_1},
(-8,-7)*{_2},
(-3,-7)*{_3},
(7,-7)*{_{n-1}},
(13,-7)*{_n},
 (0,0)*{\circledcirc}="a",
 (0,0)*{\bu},
(0,5)*{}="0",
(-12,-5)*{}="b_1",
(-8,-5)*{}="b_2",
(-3,-5)*{}="b_3",
(8,-5)*{}="b_4",
(12,-5)*{}="b_5",
\ar @{-} "a";"0" <0pt>
\ar @{-} "a";"b_2" <0pt>
\ar @{-} "a";"b_3" <0pt>
\ar @{-} "a";"b_1" <0pt>
\ar @{-} "a";"b_4" <0pt>
\ar @{-} "a";"b_5" <0pt>
\endxy}\Ea
:=i\left(
\Ba{c}\resizebox{21mm}{!}{ \xy
(1,-5)*{\ldots},
(-13,-7)*{_1},
(-8,-7)*{_2},
(-3,-7)*{_3},
(7,-7)*{_{n-1}},
(13,-7)*{_n},
 (0,0)*{\bu}="a",
(0,5)*{}="0",
(-12,-5)*{}="b_1",
(-8,-5)*{}="b_2",
(-3,-5)*{}="b_3",
(8,-5)*{}="b_4",
(12,-5)*{}="b_5",
\ar @{-} "a";"0" <0pt>
\ar @{-} "a";"b_2" <0pt>
\ar @{-} "a";"b_3" <0pt>
\ar @{-} "a";"b_1" <0pt>
\ar @{-} "a";"b_4" <0pt>
\ar @{-} "a";"b_5" <0pt>
\endxy}\Ea\right)\in \cA(n),\ \ \ n\geq 2.
$$
It is worth noting that some of these elements can stand for the zero vector in $\cA(n)$ as we do not assume in general that the map $i$ is an injection on every generator.

\sip

We define a dg operad $\widetilde{\tw}\cA=\{\widetilde{\tw}\cA(n)\}_{n\geq 0}$ as an operad generated freely by $\cA$ and one new generator $\Ba{c}\resizebox{2.5mm}{!}{ \xy
 (0,0)*{\bu}="a",
(0,4)*{}="0",
\ar @{-} "a";"0" <0pt>
\endxy}\Ea$ of type $(1,0)$ and of cohomological degree $d$. The differential $\sd$ in $\widetilde{\tw}\cA$ is equal to $\sd$ when acting on elements of $\cA$, and its action on the new generator is  defined by
\Beq\label{2: sd on (1,0) generator in A}
 \sd \Ba{c}\resizebox{1.4mm}{!}{\begin{xy}
 <0mm,0.5mm>*{};<0mm,5mm>*{}**@{-},
 <0mm,0mm>*{\bu};<0mm,0mm>*{}**@{},
 \end{xy}}\Ea= -
\sum_{{ k\geq 2}}
\frac{1}{k!}
\Ba{c}\resizebox{15mm}{!}{  \xy
(-5,8)*{}="1";
%
    %(-25,-5)*+{_1}="l1";
(-3,-5)*{...};
 <-5mm,-8mm>*{\underbrace{ \ \ \ \ \ \ \ \ \ \ \ \ \ \ \ \ \ }_{k}},
    (-5,+3)*{\circledcirc}="L";(-5,+3)*{\bu};
 (-14,-4)*{\bu}="B";
  (-8,-4)*{\bu}="C";
   (3,-4)*{\bu}="D";
     %<0mm,12mm>*{_m},
\ar @{-} "D";"L" <0pt>
\ar @{-} "C";"L" <0pt>
\ar @{-} "B";"L" <0pt>
\ar @{-} "1";"L" <0pt>
%
 %
 %\ar @{-} "l1";"L" <0pt>
 \endxy}
 \Ea
 \Eeq
There is a chain of operadic morphisms,
$$
i^+: (\Holie_d^+, \delta^+) \stackrel{c^+}{\lon} (\wt{\tw}\Holie_d, \delta) \stackrel{i}{\lon} (\wt{\tw}\cA, \sd)
$$
where the map $i$ is extended to the extra generator as the identity map.
\sip

Using Proposition {\ref{2: Prop on Holie -> O with d+}}, one concludes that the differential $\sd$  in $\widetilde{\tw}\cA$  can be twisted,
$$
\sd \rar  \sd_\centerdot:= \sd + D_{i^+(\hspace{-1.4mm}\Ba{c}\resizebox{1mm}{!}{\begin{xy}
 <0mm,-0.55mm>*{};<0mm,-3mm>*{}**@{-},
 <0mm,0.5mm>*{};<0mm,3mm>*{}**@{-},
 <0mm,0mm>*{\bullet};<0mm,0mm>*{}**@{},
 \end{xy}}\Ea\hspace{-1.4mm})}
$$
This makes the $\bS$-module $\widetilde{\tw}\cA$ into a new {\em d}g operad denoted from now on by $\tw\cA$.

\subsubsection{\bf Definition-proposition} For any dg operad $(\cA,\delta)$ under $\Holie_d$, the associated dg operad
$$
\tw\cA:= \{\tw\cA(n), \sd_{\centerdot})\}_{n\geq 0}
$$
 is called  the {\it twisted extension of $\cA$} or  the {\it twisted operad of $\cA$}. There is
 \Bi
\item[(i)] a morphism of dg operads
 $$
\iota:  (\Holie_d,\delta) \lon (\tw\cA, \sd_{\centerdot})
 $$
 which factors through the composition
 $$
 (\Holie_d,\delta) \stackrel{c}{\lon} (\tw\Holie_d), \delta_{\centerdot}) \stackrel{\tw(i)}\lon (\tw\cA, \sd_{\centerdot})
 $$
 and hence is given explicitly by
\Beq\label{2: i on Holie-corollas to Tw(O)}
\Ba{c}\resizebox{22mm}{!}{ \xy
(1,-5)*{\ldots},
(-13,-7)*{_1},
(-8,-7)*{_2},
(-3,-7)*{_3},
(7,-7)*{_{n-1}},
(13,-7)*{_n},
 (0,0)*{\bu}="a",
(0,5)*{}="0",
(-12,-5)*{}="b_1",
(-8,-5)*{}="b_2",
(-3,-5)*{}="b_3",
(8,-5)*{}="b_4",
(12,-5)*{}="b_5",
\ar @{-} "a";"0" <0pt>
\ar @{-} "a";"b_2" <0pt>
\ar @{-} "a";"b_3" <0pt>
\ar @{-} "a";"b_1" <0pt>
\ar @{-} "a";"b_4" <0pt>
\ar @{-} "a";"b_5" <0pt>
\endxy}\Ea
\ \
\stackrel{\iota}{\lon}
\ \
 \sum_{{ k\geq 0}}
\frac{1}{k!}
% \sum_{m=p+m_1+...+m_n, n,p\geq 1, m_{\bu}\geq 0 \atop n+p\geq 3}
\Ba{c}\resizebox{27mm}{!}{  \xy
(-25,8)*{}="1";
(-30,-7)*+{_1}="l2";
    (-27,-7)*+{_2}="l1";
 (-24,-6)*{...};
   (-20,-7)*+{_n}="ln";
(-3,-5)*{...};
 <-5mm,-9mm>*{\underbrace{ \ \ \ \ \ \ \ \ \ \ \ \ \ \ \ \ \ }_{k}},
    (-25,+3)*{\circledcirc}="L";(-25,+3)*{\bu};
 (-14,-5)*{\bu}="B";
  (-8,-5)*{\bu}="C";
   (3,-5)*{\bu}="D";
     %<0mm,12mm>*{_m},
\ar @{-} "D";"L" <0pt>
\ar @{-} "C";"L" <0pt>
\ar @{-} "B";"L" <0pt>
\ar @{-} "1";"L" <0pt>
 \ar @{-} "l1";"L" <0pt>
  \ar @{-} "l2";"L" <0pt>
   \ar @{-} "ln";"L" <0pt>
 \endxy}
 \Ea \ \ \ \forall\ n\geq 2.
\Eeq
\item[(ii)] a natural epimorphism of dg operads
$$
p: (\tw\cA, \sd_\centerdot) \lon (\cA, \sd)
$$
which sends the extra generator $\Ba{c}\resizebox{2.5mm}{!}{ \xy
 (0,0)*{\bu}="a",
(0,4)*{}="0",
\ar @{-} "a";"0" <0pt>
\endxy}\Ea$ to zero.
\Ei

%In applications it is often useful to study the part $\tw\cA(0)$  of $\tw\cA$ generated by elements  %with all inputs filled with the MC generator and the part $\{\tw\cA(n)\}_{n\geq 1}$ with at least one
%``free" input separately, and these parts get often in application separate names (see Examples in \S %?.? below). When we want to emphasize that we work only with a certain part of $\tw\cA$, we denote the %former part $\tw\cA(0)$ by $\tw_\bu(\cA)$, and its complementary part $\{\tw\cA(n)\}_{n\geq 1}$  by %$\tw_\circ\cA$.

\subsubsection{\bf Proposition \cite{DW}}\label{2: Prop of DW on exactness of Tw} {\em The endofunctor $\tw$ in the category of operads under
$\Holie_d$ is exact, i.e.\  any diagram
$$
\Holie_d \lon  \cA \stackrel{g}{\lon} \cA'
$$
with $g$ being a quasi-isomorphism, the map $\tw g$ in the associated  diagram
$$
\Holie_d \lon  \tw\cA \stackrel{\tw(g)}{\lon} \tw\cA'
$$
is also a quasi-isomorphism.}

\subsection{An action of the deformation complex of $\Holie_d\stackrel{i}{\rar} \cA$ on $\tw\cA$}
Given a dg operad $(\cA,\sd)$ under $\Holie_d$,
$$
\Ba{rccc}
i: & \HoLB_d  & \lon & \cA \\
&
\Ba{c}\resizebox{15mm}{!}{  \xy
(-0,6)*{}="1";
(0,+1)*{\bu}="C";
(-7,-7)*+{_1}="L1";
(-3,-7)*+{_2}="L2";
(2,-5)*{...};
(7,-7)*+{_n}="L3";
\ar @{-} "C";"L1" <0pt>
\ar @{-} "C";"L2" <0pt>
\ar @{-} "C";"L3" <0pt>
\ar @{-} "C";"1" <0pt>
 \endxy}
 \Ea
 &\lon &
 \Ba{c}\resizebox{15mm}{!}{  \xy
(-0,6)*{}="1";
(0,+1)*{\circledcirc}="C";(0,+1)*{\bu};
(-7,-7)*+{_1}="L1";
(-3,-7)*+{_2}="L2";
(2,-5)*{...};
(7,-7)*+{_n}="L3";
\ar @{-} "C";"L1" <0pt>
\ar @{-} "C";"L2" <0pt>
\ar @{-} "C";"L3" <0pt>
\ar @{-} "C";"1" <0pt>
 \endxy}
 \Ea
\Ea
$$
Consider a dg Lie algebra controlling deformations of the morphism $i$ (see \cite{MV} for several equivalent constructions of such a dg Lie algebra),
$$
\Def\left(\Holie_d \stackrel{i}{\rar} \cA\right)=\prod_{n\geq 2} \cA(n)\ot_{\bS_n} \sgn_n^{|d|}[d(1-n)]
$$
Its Maurer-Cartan elements are in 1-1 correspondence with morphisms $\Holie_d\rar \cA$ which are deformations of $i$; in particular the zero MC element corresponds to $i$ itself.
%The differential and the Lie bracket  can be described very explicitly (cf.\ \cite{MV}).
An element $\ga$ of the above complex can be represented pictorially as a collection of $(1,n)$-corollas,
$$
\ga=\{\hspace{-2mm} \Ba{c}\resizebox{12mm}{!}{  \xy
(-0,6)*{}="1";
(0,-9.4)*{\underbrace{\ \ \ \ \ \ \ \ \ \ \ \ }_n};
(0,+1)*{\circledast}="C";
(-7,-7)*+{}="L1";
(-3,-7)*+{}="L2";
(2,-5)*{...};
(7,-7)*+{}="L3";
\ar @{-} "C";"L1" <0pt>
\ar @{-} "C";"L2" <0pt>
\ar @{-} "C";"L3" <0pt>
\ar @{-} "C";"1" <0pt>
 \endxy}
 \Ea\hspace{-2mm} \}_{n\geq 2}, %\ \ \text{or, better, as a formal sum}\ \ \ \ga=\sum_{n\geq 2}\frac{1}{n!}\Ba{c}\resizebox{14mm}{!}{  \xy
%(-0,6)*{}="1";
%%
%(0,+1)*{\circledast}="C";
%(-7,-7)*+{_1}="L1";
%(-3,-7)*+{_2}="L2";
%(2,-5)*{...};
%(7,-7)*+{_n}="L3";
%
%\ar @{-} "C";"L1" <0pt>
%\ar @{-} "C";"L2" <0pt>
%\ar @{-} "C";"L3" <0pt>
%\ar @{-} "C";"1" <0pt>
% \endxy}
% \Ea
$$
of corollas whose vertices are decorated with  elements of $\cA(n)\ot_{\bS_n} \sgn_n^{|d|}$ and whose input legs are (skew)symmetrized (so that we can omit their labels); the degrees of decorations of vertices are shifted by $d(1-n)$. To distinguish these elements  from the generic elements of $\cA$
as well as from  the images of $\Holie_d$-generators under $i$, we denote the vertices of such corollas from now on  by $\circledast$. A formal sum of such corollas is homogeneous of degree $p$ if and only if the degree of each contributing $(1,n)$-corolla is equal to $p+d-dn$. The differential $\delta$ in the deformation complex $\Def\left(\Holie_d \stackrel{i}{\rar} \cA\right)$ can then be given explicitly by
\Beq\label{2: MC eqn in Def(Lie to A)}
%\sum_{n\geq 2}
\delta \hspace{-3mm}
\Ba{c}\resizebox{12mm}{!}{  \xy
(-0,6)*{}="1";
(0,-9.5)*{\underbrace{\ \ \ \ \ \ \ \ \ \ \ \ }_n};
(0,+1)*{\circledast}="C";
(-7,-7)*+{}="L1";
(-3,-7)*+{}="L2";
(2,-5)*{...};
(7,-7)*+{}="L3";
\ar @{-} "C";"L1" <0pt>
\ar @{-} "C";"L2" <0pt>
\ar @{-} "C";"L3" <0pt>
\ar @{-} "C";"1" <0pt>
 \endxy}
 \Ea
=
\sd \hspace{-3mm}
\Ba{c}\resizebox{12mm}{!}{  \xy
(-0,6)*{}="1";
(0,-9.5)*{\underbrace{\ \ \ \ \ \ \ \ \ \ \ \ }_n};
(0,+1)*{\circledast}="C";
(-7,-7)*+{}="L1";
(-3,-7)*+{}="L2";
(2,-5)*{...};
(7,-7)*+{}="L3";
\ar @{-} "C";"L1" <0pt>
\ar @{-} "C";"L2" <0pt>
\ar @{-} "C";"L3" <0pt>
\ar @{-} "C";"1" <0pt>
 \endxy}
 \Ea
 +
\sum_{[n]=[n']\sqcup [n'']\atop n'\geq 2,n''\geq 1}\left(
\pm
\Ba{c}\resizebox{19mm}{!}{  \xy
(-5,8)*{}="1";
(-3,-5)*{...};
(-12.5,-13)*{...};
 <-15mm,-17mm>*{\underbrace{ \ \ \ \ \  \ \ \ \ \ \ }_{n'}},
 <-3mm,-9mm>*{\underbrace{ \ \ \ \ \  \ \ \ \ \ \ }_{n''}},
    (-5,+3)*{\circledcirc}="L";(-5,+3)*{\bu};
 (-14,-5)*{\circledast}="B";
 (-20,-13)*{}="b1";
 (-17,-13)*{}="b2";
 (-10,-13)*{}="b3";
  (-8,-5)*{}="C";
   (3,-5)*{}="D";
     %<0mm,12mm>*{_m},
\ar @{-} "D";"L" <0pt>
\ar @{-} "C";"L" <0pt>
\ar @{-} "B";"L" <0pt>
\ar @{-} "B";"b1" <0pt>
\ar @{-} "B";"b2" <0pt>
\ar @{-} "B";"b3" <0pt>
\ar @{-} "1";"L" <0pt>
%
 %
 %\ar @{-} "l1";"L" <0pt>
 \endxy}
 \Ea
 %%%%%%%%%%%%%%%%%%%%%%%%%%%%%%%%%%%%%%%%%%%%%
\mp (-1)^{|\circledast|}
\Ba{c}\resizebox{19mm}{!}{  \xy
(-5,8)*{}="1";
(-3,-5)*{...};
(-12.5,-13)*{...};
 <-15mm,-17mm>*{\underbrace{ \ \ \ \ \  \ \ \ \ \ \ }_{n'}},
 <-3mm,-9mm>*{\underbrace{ \ \ \ \ \  \ \ \ \ \ \ }_{n''}},
    (-5,+3)*{\circledast}="L";
 (-14,-5)*{\circledcirc}="B";(-14,-5)*{\bu};
 (-20,-13)*{}="b1";
 (-17,-13)*{}="b2";
 (-10,-13)*{}="b3";
  (-8,-5)*{}="C";
   (3,-5)*{}="D";
     %<0mm,12mm>*{_m},
\ar @{-} "D";"L" <0pt>
\ar @{-} "C";"L" <0pt>
\ar @{-} "B";"L" <0pt>
\ar @{-} "B";"b1" <0pt>
\ar @{-} "B";"b2" <0pt>
\ar @{-} "B";"b3" <0pt>
\ar @{-} "1";"L" <0pt>
%
 %
 %\ar @{-} "l1";"L" <0pt>
 \endxy}
 \Ea
 \right)
\Eeq
where the rule of signs depends on $d$ and is read from (\ref{2: d in Lie_infty}); for $d$ even the first $\pm$-symbol is $+1$, while the second one is $-1$.
%the differential $\sd$ in the first term above acts on the decorations of $\circledast$-corollas.
%In fact the above equation should be understood as an infinite collection of equations, one for for all %summands with precisely  $N\geq 2$  incoming legs.

\sip

Let $(\Der(\tw\cA),\ [\ ,\ ])$ be the  Lie algebra of derivations of the non-differential operad $\tw\cA$. The differential $\sd_\centerdot$ in $\tw\cA$ is, of course, its MC element and hence makes $\Der(\tw\cA)$ into a dg Lie algebra with the differential given by the commutator $[\sd_\centerdot,\ ]$.

\subsubsection{\bf Proposition}\label{2: Theroem on Def action on TwA}{\it
There is a canonical morphism of dg Lie algebras
\Beq\label{2: Def to Der(A)}
\Ba{rccc}
\Phi: & \Def\left(\Holie_d \stackrel{i}{\rar} \cA\right) &\lon & \Der(\tw\cA)\\
      &    \ga & \lon &    \Phi_\ga
\Ea
\Eeq
where the derivation $\Phi_\ga$ is given on the generators by
$$
\Ba{rccc}
\Phi_\ga: &\tw\cA & \lon &  \tw\cA \ \ \ \ \ \ \ \ \ \ \ \\
& \Ba{c}\resizebox{1.4mm}{!}{\begin{xy}
 <0mm,0.5mm>*{};<0mm,5mm>*{}**@{-},
 <0mm,0mm>*{\bu};<0mm,0mm>*{}**@{},
 \end{xy}}\Ea
& \lon & \displaystyle
\sum_{{ k\geq 2}}
\frac{1-k}{k!}
\Ba{c}\resizebox{14mm}{!}{  \xy
(-5,8)*{}="1";
%
    %(-25,-5)*+{_1}="l1";
(-3,-5)*{...};
 <-5mm,-9mm>*{\underbrace{ \ \ \ \ \ \ \ \ \ \ \ \ \ \ \ \ \ }_{k}},
    (-5,+3)*{\circledast}="L";
 (-14,-5)*{\bu}="B";
  (-8,-5)*{\bu}="C";
   (3,-5)*{\bu}="D";
     %<0mm,12mm>*{_m},
\ar @{-} "D";"L" <0pt>
\ar @{-} "C";"L" <0pt>
\ar @{-} "B";"L" <0pt>
\ar @{-} "1";"L" <0pt>
%
 %
 %\ar @{-} "l1";"L" <0pt>
 \endxy}
 \Ea   \ \ \ \ \ \ \ \ \
  \\
%%%%%%%%%%%%%%%%%%%%%%%%%%%%%%%%%%%%%
& \Ba{c}\resizebox{14mm}{!}{  \xy
(-0,6)*{}="1";
(0,+1)*{\circ}="C";
(-7,-7)*+{_1}="L1";
(-3,-7)*+{_2}="L2";
(2,-5)*{...};
(7,-7)*+{_n}="L3";
\ar @{-} "C";"L1" <0pt>
\ar @{-} "C";"L2" <0pt>
\ar @{-} "C";"L3" <0pt>
\ar @{-} "C";"1" <0pt>
 \endxy}
 \Ea
&\lon &\displaystyle
\sum_{k\geq 1} \frac{1}{k!}
\left( -
 \Ba{c}\resizebox{17mm}{!}{  \xy
 (-0,11)*{}="1";
(-0,6)*{\circledast}="0";
(5,1)*{\bu}="r1";
(13,1)*{\bu}="r2";
(9,1)*{_{...}};
 <9mm,-3mm>*{\underbrace{ \ \ \ \ \ \ \ \ }_{k}},
(0,+1)*{\circ}="C";
(-7,-7)*+{_1}="L1";
(-3,-7)*+{_2}="L2";
(2,-5)*{...};
(7,-7)*+{_n}="L3";
\ar @{-} "C";"L1" <0pt>
\ar @{-} "C";"L2" <0pt>
\ar @{-} "C";"L3" <0pt>
\ar @{-} "C";"0" <0pt>
\ar @{-} "1";"0" <0pt>
\ar @{-} "r1";"0" <0pt>
\ar @{-} "r2";"0" <0pt>
 \endxy}
 \Ea
{ +}
  (-1)^{|\circledast||\circ|}\sum_{i=1}^n\hspace{-2mm}
 \Ba{c}\resizebox{17mm}{!}{  \xy
 (-0,5)*{}="1";
(3,-15)*+{_i}="d1";
(7,-14.5)*{\bu}="d2";
(15,-14.5)*{\bu}="d3";
(11,-14.5)*{_{...}};
 <11mm,-19mm>*{\underbrace{ \ \ \ \ \ \ \ \ }_{k}},
(0,+1)*{\circ}="C";
(-7,-7)*+{_1}="L1";
(3,-8)*{\circledast}="Li";
(-3,-7)*+{_2}="L2";
(0,-6)*{...};
(6,-6)*{...};
(10,-7)*+{_n}="L3";
\ar @{-} "C";"L1" <0pt>
\ar @{-} "C";"L2" <0pt>
\ar @{-} "C";"L3" <0pt>
\ar @{-} "C";"Li" <0pt>
\ar @{-} "C";"1" <0pt>
\ar @{-} "d1";"Li" <0pt>
\ar @{-} "d2";"Li" <0pt>
\ar @{-} "d3";"Li" <0pt>
 \endxy}
 \Ea
 \right)
 \Ea
$$
}
\begin{proof} Any derivation of $\tw \cA$ is uniquely determined by its values on the generators, i.e.\ on  $\Ba{c}\resizebox{1.5mm}{!}{\begin{xy}
 <0mm,0.5mm>*{};<0mm,4mm>*{}**@{-},
 <0mm,0mm>*{\bu};<0mm,0mm>*{}**@{},
 \end{xy}}\Ea$ and on every element  $\Ba{c}\resizebox{11mm}{!}{  \xy
(-0,6)*{}="1";
(0,+1)*{\circ}="C";
(-7,-7)*+{_1}="L1";
(-3,-7)*+{_2}="L2";
(2,-5)*{...};
(7,-7)*+{_n}="L3";
\ar @{-} "C";"L1" <0pt>
\ar @{-} "C";"L2" <0pt>
\ar @{-} "C";"L3" <0pt>
\ar @{-} "C";"1" <0pt>
 \endxy}
 \Ea$ of $\cA$. The first value can be chosen arbitrary, while the second ones are subject to the condition that they are derivations of the operad structure in $\cA$; as $\Phi_\ga$ applied  to $\Ba{c}\resizebox{11mm}{!}{  \xy
(-0,6)*{}="1";
(0,+1)*{\circ}="C";
(-7,-7)*+{_1}="L1";
(-3,-7)*+{_2}="L2";
(2,-5)*{...};
(7,-7)*+{_n}="L3";
\ar @{-} "C";"L1" <0pt>
\ar @{-} "C";"L2" <0pt>
\ar @{-} "C";"L3" <0pt>
\ar @{-} "C";"1" <0pt>
 \endxy}
 \Ea$ is precisely of the form $D_h$ discussed \S{\ref{2: subsec on twisting d in operads}}, we conclude that the above formulae do define a derivation of $\tw\cA$ as a non-differential operad. Hence it remains to show that $\Phi_\ga$ respects differentials in both dg Lie algberas, that is, satisfies the equation
 \Beq\label{2: compatibility of Der with d}
 \Phi_{\delta\ga}=[\p_\centerdot, \Phi_\ga].
\Eeq
It is straightforward to check that the operator equality (\ref{2: compatibility of Der with d}) holds true when applied to generators of  $\tw\cA$ if and only if one has the equality,
$$
\sum_{{ k\geq 1}}
\frac{1}{k!}
\Ba{c}\resizebox{22mm}{!}{  \xy
(-25,8)*{}="1";
(-25,-7)*+{}="l1";
(-8,-5)*{...};
 <-9mm,-10mm>*{\underbrace{ \ \ \ \ \ \ \ \ \ \ \ \ \ \ \ \ \ }_{k}},
    (-25,+1)*{(\delta\circledast)}="L";
 (-18,-5)*{\bu}="B";
  (-13,-5)*{\bu}="C";
   (-2,-5)*{\bu}="D";
     %<0mm,12mm>*{_m},
\ar @{-} "D";"L" <0pt>
\ar @{-} "C";"L" <0pt>
\ar @{-} "B";"L" <0pt>
\ar @{-} "1";"L" <0pt>
 \ar @{-} "l1";"L" <0pt>
 \endxy}
 \Ea
 =
\sum_{{ k\geq 1}}
 \Ba{c}\resizebox{22mm}{!}{  \xy
(-25,8)*{}="1";
(-25,-7)*+{}="l1";
(-8,-5)*{...};
 <-9mm,-10mm>*{\underbrace{ \ \ \ \ \ \ \ \ \ \ \ \ \ \ \ \ \ }_{k}},
    (-25,+1)*{(\sd_\centerdot\circledast)}="L";
 (-18,-5)*{\bu}="B";
  (-13,-5)*{\bu}="C";
   (-2,-5)*{\bu}="D";
     %<0mm,12mm>*{_m},
\ar @{-} "D";"L" <0pt>
\ar @{-} "C";"L" <0pt>
\ar @{-} "B";"L" <0pt>
\ar @{-} "1";"L" <0pt>
 \ar @{-} "l1";"L" <0pt>
 \endxy}
 \Ea
  +
 \sum_{k\geq 1\atop l\geq 2} \frac{(-1)^{|\circledast|}}{k!l!}
 \Ba{c}\resizebox{20mm}{!}{  \xy
 (0,11)*{}="1";
  (-5,1)*{}="l";
(-0,6)*{\circledcirc}="0";(-0,6)*{\bu};
(5,1)*{\bu}="r1";
(13,1)*{\bu}="r2";
(9,1)*{_{...}};
 <9mm,-3mm>*{\underbrace{ \ \ \ \ \ \ \ \ }_{k}},
(0,+1)*{\circledast}="C";
(-7,-7)*{_\bu}="L1";
(-3,-7)*{_\bu}="L2";
(2,-5)*{...};
(7,-7)*{_\bu}="L3";
<0mm,-11mm>*{\underbrace{ \ \ \ \ \ \ \ \ \ \ \ \ \ \ }_{l}},
\ar @{-} "C";"L1" <0pt>
\ar @{-} "C";"L2" <0pt>
\ar @{-} "C";"L3" <0pt>
\ar @{-} "C";"0" <0pt>
\ar @{-} "1";"0" <0pt>
\ar @{-} "l";"0" <0pt>
\ar @{-} "r1";"0" <0pt>
\ar @{-} "r2";"0" <0pt>
 \endxy}
 \Ea
$$
which is indeed the case due to (\ref{2: MC eqn in Def(Lie to A)}). The compatibility of the map $\Phi$ with Lie brackets is almost obvious.
\end{proof}
Thus the Lie algebra $H^\bu(\Def(\Holie_d \rar \cA))$ acts on the cohomology of the twisted operad $\tw\cA$ by derivations. For some operads (see Example in \S {\ref{2: Example Graphs_d operad}}) below) this cohomology Lie Lie algebra can be extremely rich and interesting.

\subsection{Twisting of $\Lie_d$} Assume, in the above notation, that $\cA$ is $\Lie_d$ and the morphism
$$
i: \Holie_d \lon \Lie_d
$$
is the canonical quasi-isomorphism. Then $\tw\Lie_d$ is a dg operad generated by
the degree $1-d$ corolla \hspace{-2mm}
$
\Ba{c}\resizebox{7mm}{!}{  \xy
(-5,5)*{}="1";
    (-5,+1)*{\bu}="L";
  (-8,-5)*+{_1}="C";
   (-2,-5)*+{_2}="D";
\ar @{-} "D";"L" <0pt>
\ar @{-} "C";"L" <0pt>
\ar @{-} "1";"L" <0pt>
%
 %
 %\ar @{-} "l1";"L" <0pt>
 \endxy}
 \Ea\hspace{-2mm} = (-1)^d\hspace{-2mm}
 \Ba{c}\resizebox{7mm}{!}{  \xy
(-5,5)*{}="1";
    (-5,+1)*{\bu}="L";
  (-8,-5)*+{_2}="C";
   (-2,-5)*+{_1}="D";
\ar @{-} "D";"L" <0pt>
\ar @{-} "C";"L" <0pt>
\ar @{-} "1";"L" <0pt>
%
 %
 %\ar @{-} "l1";"L" <0pt>
 \endxy}
 \Ea\hspace{-2mm}
$
(modulo the Jacobi relations) and the degree $d$ $(1,0)$-corolla  $\Ba{c}\resizebox{1.5mm}{!}{\begin{xy}
 <0mm,0.5mm>*{};<0mm,4mm>*{}**@{-},
 <0mm,0mm>*{\bu};<0mm,0mm>*{}**@{},
 \end{xy}}\Ea$.
The twisted differential $\delta_\centerdot$ is trivial on the first generator (due to the Jacobi identities)
$$
\delta_\centerdot\hspace{-3mm}
\Ba{c}\resizebox{8mm}{!}{  \xy
(-5,6)*{}="1";
    (-5,+1)*{\bu}="L";
  (-8,-5)*+{_1}="C";
   (-2,-5)*+{_2}="D";
\ar @{-} "D";"L" <0pt>
\ar @{-} "C";"L" <0pt>
\ar @{-} "1";"L" <0pt>
%
 %
 %\ar @{-} "l1";"L" <0pt>
 \endxy}
 \Ea
 \equiv
 \Ba{c}\resizebox{11.0mm}{!}{  \xy
(-9,8)*{}="1";
    (-9,+3)*{\bu}="L";
 (-14,-3.5)*{\bu}="B";
 (-18,-10)*+{_1}="b1";
 (-10,-10)*+{_2}="b2";
  (-3,-5)*{\bu}="C";
     %<0mm,12mm>*{_m},
\ar @{-} "C";"L" <0pt>
\ar @{-} "B";"L" <0pt>
\ar @{-} "B";"b1" <0pt>
\ar @{-} "B";"b2" <0pt>
\ar @{-} "1";"L" <0pt>
%
 %
 %\ar @{-} "l1";"L" <0pt>
 \endxy}
 \Ea
 +(-1)^d
  \Ba{c}\resizebox{11.5mm}{!}{  \xy
(-9,8)*{}="1";
    (-9,+3)*{\bu}="L";
 (-14,-3.5)*{\bu}="B";
 (-18,-10)*+{_1}="b1";
 (-10,-10)*{\bu}="b2";
  (-3,-5)*+{_2}="C";
     %<0mm,12mm>*{_m},
\ar @{-} "C";"L" <0pt>
\ar @{-} "B";"L" <0pt>
\ar @{-} "B";"b1" <0pt>
\ar @{-} "B";"b2" <0pt>
\ar @{-} "1";"L" <0pt>
%
 %
 %\ar @{-} "l1";"L" <0pt>
 \endxy}
 \Ea
 +
   \Ba{c}\resizebox{11.5mm}{!}{  \xy
(-9,8)*{}="1";
    (-9,+3)*{\bu}="L";
 (-14,-3.5)*{\bu}="B";
 (-18,-10)*+{_2}="b1";
 (-10,-10)*{\bu}="b2";
  (-3,-5)*+{_1}="C";
     %<0mm,12mm>*{_m},
\ar @{-} "C";"L" <0pt>
\ar @{-} "B";"L" <0pt>
\ar @{-} "B";"b1" <0pt>
\ar @{-} "B";"b2" <0pt>
\ar @{-} "1";"L" <0pt>
%
 %
 %\ar @{-} "l1";"L" <0pt>
 \endxy}
 \Ea=0,
$$
 while it acts  on the second generator  by
$$
 \delta_\centerdot \Ba{c}\resizebox{1.4mm}{!}{\begin{xy}
 <0mm,0.5mm>*{};<0mm,5mm>*{}**@{-},
 <0mm,0mm>*{\bu};<0mm,0mm>*{}**@{},
 \end{xy}}\Ea=
\frac{1}{2}
\Ba{c}\resizebox{6.0mm}{!}{  \xy
(-5,6)*{}="1";
    (-5,+1)*{\bu}="L";
  (-8,-5)*{\bu}="C";
   (-2,-5)*{\bu}="D";
\ar @{-} "D";"L" <0pt>
\ar @{-} "C";"L" <0pt>
\ar @{-} "1";"L" <0pt>
%
 %
 %\ar @{-} "l1";"L" <0pt>
 \endxy}
 \Ea.
 $$
 The twisted morphism $\tw(i): \Lie_d \lon \tw\Lie_d$ becomes in this case the obvious inclusion.
This example is important for understanding those twisted operads $\tw\cA$ for which the map $\Holie_d\rar \cA$ factors through the above projection,
$$
\Holie_d \lon \Lie_d \lon \cA.
$$
We call such dg operads {\em operads under $\Lie_d$}. The deformation complex of the epimorphism $i$ has almost trivial cohomology, $H^\bu(\Def(\HoLBd\rar \Lie_d))= \R[-1]$; the only non-trivial cohomology class acts on $\tw\LBd$ by rescaling its generators.

\subsection{Twisting of (prop)operads under $\Lie_d$}\label{2: subsec on twisting of (prop)erads under Lie} Let $(\cA, \sd)$ be a dg operad equipped with an operadic morphism
$$
\Ba{rccc}
\imath: & (\Lie_d, 0) & \lon & (\cA,\sd)\\
& \Ba{c}\begin{xy}
 <0mm,0.66mm>*{};<0mm,4mm>*{}**@{-},
 <0.39mm,-0.39mm>*{};<2.2mm,-2.2mm>*{}**@{-},
 <-0.35mm,-0.35mm>*{};<-2.2mm,-2.2mm>*{}**@{-},
 <0mm,0mm>*{\bu};<0mm,0mm>*{}**@{},
   %<0mm,0.66mm>*{};<0mm,3.4mm>*{^1}**@{},
   <0.39mm,-0.39mm>*{};<2.9mm,-4mm>*{^{_2}}**@{},
   <-0.35mm,-0.35mm>*{};<-2.8mm,-4mm>*{^{_1}}**@{},
\end{xy}\Ea & \lon & \Ba{c}\resizebox{8mm}{!}{  \xy
(-5,6)*{}="1";
    (-5,+1)*{\circledcirc}="L";(-5,+1)*{\bu};
  (-8,-5)*+{_1}="C";
   (-2,-5)*+{_2}="D";
\ar @{-} "D";"L" <0pt>
\ar @{-} "C";"L" <0pt>
\ar @{-} "1";"L" <0pt>
%
 %
 %\ar @{-} "l1";"L" <0pt>
 \endxy}
 \Ea
\Ea
$$
where $\Lie_d$ is understood as a differential operad with the trivial differential. Then $(\tw\cA,\sd_\centerdot)$ is an operad freely generated by  $\cA$ and one new generator  $\Ba{c}\resizebox{1.5mm}{!}{\begin{xy}
 <0mm,0.5mm>*{};<0mm,4mm>*{}**@{-},
 <0mm,0mm>*{\bu};<0mm,0mm>*{}**@{},
 \end{xy}}\Ea$ of degree $d$. The differential $\sd_\centerdot$ acts, by definition, on an element $a$ of $\cA(n)$ (identified with the $a$-decorated $(1,n)$-corolla) by a formula similar to (\ref{2: delta^+})
 \Beq\label{2: sd_c on generators of TwO}
\sd_\centerdot\hspace{-3mm}
\Ba{c}\resizebox{14mm}{!}{  \xy
(-0,6)*{}="1";
(0,+1)*{\circ}="C";
(-7,-7)*+{_1}="L1";
(-3,-7)*+{_2}="L2";
(2,-5)*{...};
(7,-7)*+{_n}="L3";
\ar @{-} "C";"L1" <0pt>
\ar @{-} "C";"L2" <0pt>
\ar @{-} "C";"L3" <0pt>
\ar @{-} "C";"1" <0pt>
 \endxy}
 \Ea
 :=
 \sd\hspace{-3mm}
\Ba{c}\resizebox{14mm}{!}{  \xy
(-0,6)*{}="1";
(0,+1)*{\circ}="C";
(-7,-7)*+{_1}="L1";
(-3,-7)*+{_2}="L2";
(2,-5)*{...};
(7,-7)*+{_n}="L3";
\ar @{-} "C";"L1" <0pt>
\ar @{-} "C";"L2" <0pt>
\ar @{-} "C";"L3" <0pt>
\ar @{-} "C";"1" <0pt>
 \endxy}
 \Ea
 +
 \Ba{c}\resizebox{15mm}{!}{  \xy
 (-0,11)*{}="1";
(-0,6)*{\circledcirc}="0";(-0,6)*{\bu};
(5,1)*{\bu}="r";
(0,+1)*{\circ}="C";
(-7,-7)*+{_1}="L1";
(-3,-7)*+{_2}="L2";
(2,-5)*{...};
(7,-7)*+{_n}="L3";
\ar @{-} "C";"L1" <0pt>
\ar @{-} "C";"L2" <0pt>
\ar @{-} "C";"L3" <0pt>
\ar @{-} "C";"0" <0pt>
\ar @{-} "1";"0" <0pt>
\ar @{-} "r";"0" <0pt>
 \endxy}
 \Ea
 -
  (-1)^{|a|}\sum_{i=1}^n\hspace{-2mm}
 \Ba{c}\resizebox{16mm}{!}{  \xy
 (-0,5)*{}="1";
(3,-15)*+{_i}="d1";
(9,-14.5)*{\bu}="d2";
(0,+1)*{\circ}="C";
(-7,-7)*+{_1}="L1";
(3,-8)*{\circledcirc}="Li";(3,-8)*{\bu};
(-3,-7)*+{_2}="L2";
(0,-6)*{...};
(6,-6)*{...};
(10,-7)*+{_n}="L3";
\ar @{-} "C";"L1" <0pt>
\ar @{-} "C";"L2" <0pt>
\ar @{-} "C";"L3" <0pt>
\ar @{-} "C";"Li" <0pt>
\ar @{-} "C";"1" <0pt>
\ar @{-} "d1";"Li" <0pt>
\ar @{-} "d2";"Li" <0pt>
 \endxy}
 \Ea
\Eeq
and on the extra generator as follows
\Beq\label{2: d_c on boxdot}
 \sd_\centerdot \Ba{c}\resizebox{1.4mm}{!}{\begin{xy}
 <0mm,0.5mm>*{};<0mm,5mm>*{}**@{-},
 <0mm,0mm>*{\bu};<0mm,0mm>*{}**@{},
 \end{xy}}\Ea=
\frac{1}{2}
\Ba{c}\resizebox{7mm}{!}{  \xy
(-5,6)*{}="1";
    (-5,+1)*{\circledcirc}="L";(-5,+1)*{\bu};
  (-8,-5)*{\bu}="C";
   (-2,-5)*{\bu}="D";
\ar @{-} "D";"L" <0pt>
\ar @{-} "C";"L" <0pt>
\ar @{-} "1";"L" <0pt>
%
 %
 %\ar @{-} "l1";"L" <0pt>
 \endxy}
 \Ea.
 \Eeq

By Proposition~{\ref{2: Prop of DW on exactness of Tw}}, the canonical epimorphism
$$
\tw \Holie_d \lon \tw\Lie_d
$$
is a quasi-isomorphism. Moreover, it is proven in \cite{DW} that the natural projections,
\Beq\label{2: quasi-iso twHoLB to HoLB, twLB to LB}
\tw\Holie_d \lon \Holie_d, \ \ \  \tw\Lie_d\lon \Lie_d
\Eeq
are quasi-isomorphisms as well.

\subsubsection{\bf Example: Twisting of $\cA ss$} The operad of associative algebras $\cA ss$ is obviously an operad under $\Lie_1$ and hence can be twisted. It is proven in \cite{CL} that the natural projection
$$
\tw \cA ss \lon \cA ss
$$
is a quasi-isomorphism.

\subsection{Example: M.\ Kontsevich's operad of graphs}\label{2: Example Graphs_d operad} Here is an example of the twisting procedure used in \cite{W} to reproduce an important dg operad of graphs $\cG raphs_d$ which has been invented by M.\ Kontsevich in \cite{Ko2} in the context of a new proof of the formality of the little disks operad, and which was further studied in \cite{LV,W}.
 By a {\em graph}\, $\Ga$ we understand  a 1-dimensional $CW$ complex whose 0-cells are called vertices and 1-cells are called edges; the set of vertices of $\Ga$ is denoted by $V(\Ga)$ and the set of edges by $E(\Ga)$.
 %A graph $\Ga$ is called {\em directed}\, if its edge $e\in E(\Ga)$ comes equipped with an orientation or, plainly speaking,
%with a choice of a direction.
Let $\Gra_d(n)$, $d\in \Z$, stand for the graded vector space generated by graphs $\Ga$ such that
\Bi
\item[(i)] $\Ga$ has precisely $n$ vertices which are labelled, that is an isomorphism $V(\Ga)\rar [n]$ is fixed;
\item[(ii)] $\Ga$ is equipped with an orientation which is for $d$ even is defined as an ordering of edges (up the sign action of $\bS_{\# E(\Ga)}$), while for $d$ odd it is defined as a choice of the direction on each edge (up to the sign action of $\bS_2$ whose generator flips the direction).
\item[(iii)] $\Ga$ is assigned the cohomological degree $(1-d)\# E(\Ga)$.
\Ei
For example,
 $$
\Ba{c}\resizebox{9mm}{!}{
\xy
(0,2)*{^1},
(8,2)*{^2},
 (0,0)*{\circ}="a",
(8,0)*{\circ}="b",
\ar @{-} "a";"b" <0pt>
\endxy}\Ea \in \cG ra_{d}(2), \ \ \
\Ba{c}\resizebox{10mm}{!}{
\xy
(-7,4)*{^1},
(-7,-4)*{_2},
(4,2)*{^3},
 {\ar@{-}(-5,4)*{\circ};(-5,-4)*{\circ}};
   {\ar@{-}(-5,4)*{\circ};(4,0)*{\circ}};
 {\ar@{-}(4,0)*{\circ};(-5,-4)*{\circ}};
\endxy}\Ea\in \cG ra_{d}(3)
$$
where for $d$ odd one should assume a choice of directions on edges (defined up to a flip
and multiplication by $-1$).
The $\Z$-graded vector space $\Gra_d(n)$ is an $\bS_n$-module with the permutation group acting on graphs by relabeling their vertices. The $\bS$-module
$$
\Gra_d:= \{
\Gra_d(n)\}
$$
is an operad \cite{W} with the
 operadic compositions
\Beq\label{3: operad comp in Gra}
\Ba{rccc}
\circ_i: & \cG ra (n)\ot \cG ra (m) & \lon &  \cG ra (m+n-1)\\
&  \Ga_1 \ot \Ga_2   &\lon & \Ga_1\circ_i \Ga_2
\Ea
\Eeq
defined as follows:  $\Ga_1\circ_i \Ga_2$ is the linear combination of graphs obtained by substituting the graph $\Ga_2$ into the $i$-labeled vertex of $\Ga_1$ and taking a sum over all possible re-attachments of dangling edges (attached earlier to that vertex) to the vertices of $\Ga_2$. Here is an example (for $d$ odd),
$$
\Ba{c}\resizebox{10mm}{!}{
\xy
(-5,2)*{^1},
(5,2)*{^2},
   {\ar@/^0.6pc/(-5,0)*{\bu};(5,0)*{\bu}};
 {\ar@/^0.6pc/(5,0)*{\bu};(-5,0)*{\bu}};
\endxy}
\Ea
 \ \ \circ_1\
 \Ba{c}\resizebox{4mm}{!}{
\xy
(-2,7)*{^1},
(-2,0)*{_2},
 {\ar@{->}(0,7)*{\bu};(0,0)*{\bu}};
\endxy}\Ea
=
 \Ba{c}\resizebox{12mm}{!}{
\xy
(-7,7)*{^1},
(-7,0)*{_2},
(5,2)*{^3},
 {\ar@{->}(-5,7)*{\bu};(-5,0)*{\bu}};
   {\ar@/^0.6pc/(-5,0)*{\bu};(5,0)*{\bu}};
 {\ar@/^0.6pc/(5,0)*{\bu};(-5,0)*{\bu}};
\endxy}\Ea
+
 \Ba{c}\resizebox{12mm}{!}{
\xy
(-7,0)*{^1},
(-7,-7)*{_2},
(5,2)*{^3},
 {\ar@{->}(-5,0)*{\bu};(-5,-7)*{\bu}};
   {\ar@/^0.6pc/(-5,0)*{\bu};(5,0)*{\bu}};
 {\ar@/^0.6pc/(5,0)*{\bu};(-5,0)*{\bu}};
\endxy}\Ea
+
 \Ba{c}\resizebox{11mm}{!}{
\xy
(-7,4)*{^1},
(-7,-4)*{_2},
(5,2)*{^3},
 {\ar@{->}(-5,4)*{\bu};(-5,-4)*{\bu}};
   {\ar@{->}(-5,4)*{\bu};(5,0)*{\bu}};
 {\ar@{->}(5,0)*{\bu};(-5,-4)*{\bu}};
\endxy}\Ea
+
 \Ba{c}\resizebox{11mm}{!}{
\xy
(-7,4)*{^1},
(-7,-4)*{_2},
(5,2)*{^3},
 {\ar@{->}(-5,4)*{\bu};(-5,-4)*{\bu}};
   {\ar@{<-}(-5,4)*{\bu};(5,0)*{\bu}};
 {\ar@{<-}(5,0)*{\bu};(-5,-4)*{\bu}};
\endxy}
\Ea
$$
There is a morphism of operads \cite{W}
$$
\Ba{ccc}
\Lie_d & \lon & \Gra_d\\
\Ba{c}
\xy
 <0mm,0.55mm>*{};<0mm,3.5mm>*{}**@{-},
 <0.5mm,-0.5mm>*{};<2.2mm,-2.2mm>*{}**@{-},
 <-0.48mm,-0.48mm>*{};<-2.2mm,-2.2mm>*{}**@{-},
 <0mm,0mm>*{\bu};<0mm,0mm>*{}**@{},
 <0.5mm,-0.5mm>*{};<2.7mm,-3.2mm>*{_2}**@{},
 <-0.48mm,-0.48mm>*{};<-2.7mm,-3.2mm>*{_1}**@{},
 \endxy\Ea
&\lon &
 \xy
(0,2)*{_{1}},
(5,2)*{_{2}},
 (0,0)*{\circ}="a",
(5,0)*{\circ}="b",
\ar @{-} "a";"b" <0pt>
\endxy
\Ea
$$
so that one can apply the twisting  endofunctor to  $\Gra_d$. The resulting dg operad $\tw \Gra_d$ is generated by graphs with two types of vertices, white ones which are labelled and black ones which are unlabelled and assigned the cohomological degree $d$, e.g.
$$
\Ba{c}\resizebox{10mm}{!}{
\xy
(-7,4)*{^2},
(-7,-4)*{_1},
%(4,2)*{^3},
 {\ar@{-}(-5,4)*{\circ};(-5,-4)*{\circ}};
   {\ar@{-}(-5,4)*{\circ};(4,0)*{\circ}};
 {\ar@{-}(4,0)*{\bu};(-5,-4)*{\circ}};
\endxy}\Ea\in \cG raphs_{d}(2)
$$
The differential acts on white vertices and black vertices by splitting them,
\Beq\label{2: d in Graphs}
 \xy
(0,2)*{_{i}},
 (0,0)*{\circ}="a",
\endxy
\rightsquigarrow
 \xy
(0,2)*{_{i}},
 (0,0)*{\circ}="a",
(5,0)*{\bu}="b",
\ar @{-} "a";"b" <0pt>
\endxy
,
\ \ \ \
 \xy
 (0,0)*{\bu}="a",
\endxy
\rightsquigarrow
 \xy
 (0,0)*{\bu}="a",
(5,0)*{\bu}="b",
\ar @{-} "a";"b" <0pt>
\endxy
\Eeq
and re-attaching edges.
The dg sub-operad of $\tw\Gra_d$ generated by graphs with at least one white vertex is denoted by $\cG raphs_d$. It is proven in \cite{Ko2,LV} that its cohomology $H^\bu(\cG raphs_d)$ is the operad of $d$-algebras. The case $d=2$ is of special interest as $2$-algebras are precisely the Gerstenhaber algebras which have many applications in algebra, geometry and mathematical physics.

\bip

%%%%%%%%%%%%%%%%%%%%%%%%%%%%%%%%%%%%%%%%%%%%%%%%%%%%%%%%%%%%%%%%%%%%%%%%%%%%%
%%%%%%%%%%%%%%%%%%%%%%%%%%%%%%%%%%%%%%%%%%%%%%%%%%%%%%%%%%%%%%%%%%%%%%%%%%%%%

{\Large
\section{\bf Partial twisting of properads under $\LB_d$ and quasi-Lie bialgebras}
}

\mip

\subsection{Reminder on the properads of (degree shifted) Lie bialgebras and quasi-Lie bialgebras} The properad of degree shifted Lie bialgebras is defined, for any pair of integer  $c,d\in \Z$,  as the quotient
$$
\LB_{c,d}:=\cF ree\langle E_0\rangle/\langle\cR\rangle,
$$
of the free prop generated by an  $\bS$-bimodule $E_0=\{E_0(m,n)\}_{m,n\geq 0}$ with
 all $E_0(m,n)=0$ except

  $$
E_0(2,1):=\id_1\ot \sgn_2^{c}[c-1]=\mbox{span}\left\langle
\Ba{c}\begin{xy}
 <0mm,-0.55mm>*{};<0mm,-2.5mm>*{}**@{-},
 <0.5mm,0.5mm>*{};<2.2mm,2.2mm>*{}**@{-},
 <-0.48mm,0.48mm>*{};<-2.2mm,2.2mm>*{}**@{-},
 <0mm,0mm>*{\bu};<0mm,0mm>*{}**@{},
 %<0mm,-0.55mm>*{};<0mm,-3.8mm>*{_1}**@{},
 <0.5mm,0.5mm>*{};<2.7mm,2.8mm>*{^{_2}}**@{},
 <-0.48mm,0.48mm>*{};<-2.7mm,2.8mm>*{^{_1}}**@{},
 \end{xy}\Ea
=(-1)^{c}
\Ba{c}\begin{xy}
 <0mm,-0.55mm>*{};<0mm,-2.5mm>*{}**@{-},
 <0.5mm,0.5mm>*{};<2.2mm,2.2mm>*{}**@{-},
 <-0.48mm,0.48mm>*{};<-2.2mm,2.2mm>*{}**@{-},
 <0mm,0mm>*{\bu};<0mm,0mm>*{}**@{},
 %<0mm,-0.55mm>*{};<0mm,-3.8mm>*{_1}**@{},
 <0.5mm,0.5mm>*{};<2.7mm,2.8mm>*{^{_1}}**@{},
 <-0.48mm,0.48mm>*{};<-2.7mm,2.8mm>*{^{_2}}**@{},
 \end{xy}\Ea
   \right\rangle
$$
$$
E_0(1,2):= \sgn_2^{d}\ot \id_1[d-1]=\mbox{span}\left\langle
\Ba{c}\begin{xy}
 <0mm,0.66mm>*{};<0mm,3mm>*{}**@{-},
 <0.39mm,-0.39mm>*{};<2.2mm,-2.2mm>*{}**@{-},
 <-0.35mm,-0.35mm>*{};<-2.2mm,-2.2mm>*{}**@{-},
 <0mm,0mm>*{\bu};<0mm,0mm>*{}**@{},
   <0.39mm,-0.39mm>*{};<2.9mm,-4mm>*{^{_2}}**@{},
   <-0.35mm,-0.35mm>*{};<-2.8mm,-4mm>*{^{_1}}**@{},
\end{xy}\Ea
=(-1)^{d}
\Ba{c}\begin{xy}
 <0mm,0.66mm>*{};<0mm,3mm>*{}**@{-},
 <0.39mm,-0.39mm>*{};<2.2mm,-2.2mm>*{}**@{-},
 <-0.35mm,-0.35mm>*{};<-2.2mm,-2.2mm>*{}**@{-},
 <0mm,0mm>*{\bu};<0mm,0mm>*{}**@{},
   <0.39mm,-0.39mm>*{};<2.9mm,-4mm>*{^{_1}}**@{},
   <-0.35mm,-0.35mm>*{};<-2.8mm,-4mm>*{^{_2}}**@{},
\end{xy}\Ea
\right\rangle
$$
by the ideal generated by the following relations
\Beq\label{3: R for LieB}
\cR:\left\{
\Ba{c}
\oint_{123}\hspace{-1mm} \Ba{c}\resizebox{8.1mm}{!}{
\begin{xy}
 <0mm,0mm>*{\bu};<0mm,0mm>*{}**@{},
 <0mm,-0.49mm>*{};<0mm,-3.0mm>*{}**@{-},
 <0.49mm,0.49mm>*{};<1.9mm,1.9mm>*{}**@{-},
 <-0.5mm,0.5mm>*{};<-1.9mm,1.9mm>*{}**@{-},
 <-2.3mm,2.3mm>*{\bu};<-2.3mm,2.3mm>*{}**@{},
 <-1.8mm,2.8mm>*{};<0mm,4.9mm>*{}**@{-},
 <-2.8mm,2.9mm>*{};<-4.6mm,4.9mm>*{}**@{-},
   <0.49mm,0.49mm>*{};<2.7mm,2.3mm>*{^3}**@{},
   <-1.8mm,2.8mm>*{};<0.4mm,5.3mm>*{^2}**@{},
   <-2.8mm,2.9mm>*{};<-5.1mm,5.3mm>*{^1}**@{},
 \end{xy}}\Ea
% +
%\Ba{c}\resizebox{7mm}{!}{\begin{xy}
% <0mm,0mm>*{\bu};<0mm,0mm>*{}**@{},
% <0mm,-0.49mm>*{};<0mm,-3.0mm>*{}**@{-},
% <0.49mm,0.49mm>*{};<1.9mm,1.9mm>*{}**@{-},
% <-0.5mm,0.5mm>*{};<-1.9mm,1.9mm>*{}**@{-},
% <-2.3mm,2.3mm>*{\bu};<-2.3mm,2.3mm>*{}**@{},
% <-1.8mm,2.8mm>*{};<0mm,4.9mm>*{}**@{-},
% <-2.8mm,2.9mm>*{};<-4.6mm,4.9mm>*{}**@{-},
%   <0.49mm,0.49mm>*{};<2.7mm,2.3mm>*{^2}**@{},
%   <-1.8mm,2.8mm>*{};<0.4mm,5.3mm>*{^1}**@{},
%   <-2.8mm,2.9mm>*{};<-5.1mm,5.3mm>*{^3}**@{},
% \end{xy}}\Ea
% +
%\Ba{c}\resizebox{7mm}{!}{\begin{xy}
% <0mm,0mm>*{\bu};<0mm,0mm>*{}**@{},
% <0mm,-0.49mm>*{};<0mm,-3.0mm>*{}**@{-},
% <0.49mm,0.49mm>*{};<1.9mm,1.9mm>*{}**@{-},
% <-0.5mm,0.5mm>*{};<-1.9mm,1.9mm>*{}**@{-},
% <-2.3mm,2.3mm>*{\bu};<-2.3mm,2.3mm>*{}**@{},
% <-1.8mm,2.8mm>*{};<0mm,4.9mm>*{}**@{-},
% <-2.8mm,2.9mm>*{};<-4.6mm,4.9mm>*{}**@{-},
%   <0.49mm,0.49mm>*{};<2.7mm,2.3mm>*{^1}**@{},
%   <-1.8mm,2.8mm>*{};<0.4mm,5.3mm>*{^3}**@{},
%   <-2.8mm,2.9mm>*{};<-5.1mm,5.3mm>*{^2}**@{},
% \end{xy}}\Ea
=0
 \ \ , \ \
%%%%%%%%%%%%%% Lie %%%%%%%%%%%%%%%%%%%%%%%%
\oint_{123}\hspace{-1mm} \Ba{c}\resizebox{8.4mm}{!}{ \begin{xy}
 <0mm,0mm>*{\bu};<0mm,0mm>*{}**@{},
 <0mm,0.69mm>*{};<0mm,3.0mm>*{}**@{-},
 <0.39mm,-0.39mm>*{};<2.4mm,-2.4mm>*{}**@{-},
 <-0.35mm,-0.35mm>*{};<-1.9mm,-1.9mm>*{}**@{-},
 <-2.4mm,-2.4mm>*{\bu};<-2.4mm,-2.4mm>*{}**@{},
 <-2.0mm,-2.8mm>*{};<0mm,-4.9mm>*{}**@{-},
 <-2.8mm,-2.9mm>*{};<-4.7mm,-4.9mm>*{}**@{-},
    <0.39mm,-0.39mm>*{};<3.3mm,-4.0mm>*{^3}**@{},
    <-2.0mm,-2.8mm>*{};<0.5mm,-6.7mm>*{^2}**@{},
    <-2.8mm,-2.9mm>*{};<-5.2mm,-6.7mm>*{^1}**@{},
 \end{xy}}\Ea
% +
%\Ba{c}\resizebox{8.4mm}{!}{ \begin{xy}
% <0mm,0mm>*{\bu};<0mm,0mm>*{}**@{},
% <0mm,0.69mm>*{};<0mm,3.0mm>*{}**@{-},
% <0.39mm,-0.39mm>*{};<2.4mm,-2.4mm>*{}**@{-},
% <-0.35mm,-0.35mm>*{};<-1.9mm,-1.9mm>*{}**@{-},
% <-2.4mm,-2.4mm>*{\bu};<-2.4mm,-2.4mm>*{}**@{},
% <-2.0mm,-2.8mm>*{};<0mm,-4.9mm>*{}**@{-},
% <-2.8mm,-2.9mm>*{};<-4.7mm,-4.9mm>*{}**@{-},
%    <0.39mm,-0.39mm>*{};<3.3mm,-4.0mm>*{^2}**@{},
%    <-2.0mm,-2.8mm>*{};<0.5mm,-6.7mm>*{^1}**@{},
%    <-2.8mm,-2.9mm>*{};<-5.2mm,-6.7mm>*{^3}**@{},
% \end{xy}}\Ea
% +
%\Ba{c}\resizebox{8.4mm}{!}{ \begin{xy}
% <0mm,0mm>*{\bu};<0mm,0mm>*{}**@{},
% <0mm,0.69mm>*{};<0mm,3.0mm>*{}**@{-},
% <0.39mm,-0.39mm>*{};<2.4mm,-2.4mm>*{}**@{-},
% <-0.35mm,-0.35mm>*{};<-1.9mm,-1.9mm>*{}**@{-},
% <-2.4mm,-2.4mm>*{\bu};<-2.4mm,-2.4mm>*{}**@{},
% <-2.0mm,-2.8mm>*{};<0mm,-4.9mm>*{}**@{-},
% <-2.8mm,-2.9mm>*{};<-4.7mm,-4.9mm>*{}**@{-},
%    <0.39mm,-0.39mm>*{};<3.3mm,-4.0mm>*{^1}**@{},
%    <-2.0mm,-2.8mm>*{};<0.5mm,-6.7mm>*{^3}**@{},
%    <-2.8mm,-2.9mm>*{};<-5.2mm,-6.7mm>*{^2}**@{},
% \end{xy}}\Ea
=0
 \\
%%%%%%%%%%%%%%%%%%%%%%% Lie[1]Bi %%%%%%%%%%%%%%%
 (-1)^{c+d}\Ba{c}\resizebox{5mm}{!}{\begin{xy}
 <0mm,2.47mm>*{};<0mm,0.12mm>*{}**@{-},
 <0.5mm,3.5mm>*{};<2.2mm,5.2mm>*{}**@{-},
 <-0.48mm,3.48mm>*{};<-2.2mm,5.2mm>*{}**@{-},
 <0mm,3mm>*{\bu};<0mm,3mm>*{}**@{},
  <0mm,-0.8mm>*{\bu};<0mm,-0.8mm>*{}**@{},
<-0.39mm,-1.2mm>*{};<-2.2mm,-3.5mm>*{}**@{-},
 <0.39mm,-1.2mm>*{};<2.2mm,-3.5mm>*{}**@{-},
     <0.5mm,3.5mm>*{};<2.8mm,5.7mm>*{^2}**@{},
     <-0.48mm,3.48mm>*{};<-2.8mm,5.7mm>*{^1}**@{},
   <0mm,-0.8mm>*{};<-2.7mm,-5.2mm>*{^1}**@{},
   <0mm,-0.8mm>*{};<2.7mm,-5.2mm>*{^2}**@{},
\end{xy}}\Ea
+(-1)^{cd}\left(
\Ba{c}\resizebox{7mm}{!}{\begin{xy}
 <0mm,-1.3mm>*{};<0mm,-3.5mm>*{}**@{-},
 <0.38mm,-0.2mm>*{};<2.0mm,2.0mm>*{}**@{-},
 <-0.38mm,-0.2mm>*{};<-2.2mm,2.2mm>*{}**@{-},
<0mm,-0.8mm>*{\bu};<0mm,0.8mm>*{}**@{},
 <2.4mm,2.4mm>*{\bu};<2.4mm,2.4mm>*{}**@{},
 <2.77mm,2.0mm>*{};<4.4mm,-0.8mm>*{}**@{-},
 <2.4mm,3mm>*{};<2.4mm,5.2mm>*{}**@{-},
     <0mm,-1.3mm>*{};<0mm,-5.3mm>*{^1}**@{},
     <2.5mm,2.3mm>*{};<5.1mm,-2.6mm>*{^2}**@{},
    <2.4mm,2.5mm>*{};<2.4mm,5.7mm>*{^2}**@{},
    <-0.38mm,-0.2mm>*{};<-2.8mm,2.5mm>*{^1}**@{},
    \end{xy}}\Ea
  + (-1)^{d}
\Ba{c}\resizebox{7mm}{!}{\begin{xy}
 <0mm,-1.3mm>*{};<0mm,-3.5mm>*{}**@{-},
 <0.38mm,-0.2mm>*{};<2.0mm,2.0mm>*{}**@{-},
 <-0.38mm,-0.2mm>*{};<-2.2mm,2.2mm>*{}**@{-},
<0mm,-0.8mm>*{\bu};<0mm,0.8mm>*{}**@{},
 <2.4mm,2.4mm>*{\bu};<2.4mm,2.4mm>*{}**@{},
 <2.77mm,2.0mm>*{};<4.4mm,-0.8mm>*{}**@{-},
 <2.4mm,3mm>*{};<2.4mm,5.2mm>*{}**@{-},
     <0mm,-1.3mm>*{};<0mm,-5.3mm>*{^2}**@{},
     <2.5mm,2.3mm>*{};<5.1mm,-2.6mm>*{^1}**@{},
    <2.4mm,2.5mm>*{};<2.4mm,5.7mm>*{^2}**@{},
    <-0.38mm,-0.2mm>*{};<-2.8mm,2.5mm>*{^1}**@{},
    \end{xy}}\Ea
  + (-1)^{d+c}
\Ba{c}\resizebox{7mm}{!}{\begin{xy}
 <0mm,-1.3mm>*{};<0mm,-3.5mm>*{}**@{-},
 <0.38mm,-0.2mm>*{};<2.0mm,2.0mm>*{}**@{-},
 <-0.38mm,-0.2mm>*{};<-2.2mm,2.2mm>*{}**@{-},
<0mm,-0.8mm>*{\bu};<0mm,0.8mm>*{}**@{},
 <2.4mm,2.4mm>*{\bu};<2.4mm,2.4mm>*{}**@{},
 <2.77mm,2.0mm>*{};<4.4mm,-0.8mm>*{}**@{-},
 <2.4mm,3mm>*{};<2.4mm,5.2mm>*{}**@{-},
     <0mm,-1.3mm>*{};<0mm,-5.3mm>*{^2}**@{},
     <2.5mm,2.3mm>*{};<5.1mm,-2.6mm>*{^1}**@{},
    <2.4mm,2.5mm>*{};<2.4mm,5.7mm>*{^1}**@{},
    <-0.38mm,-0.2mm>*{};<-2.8mm,2.5mm>*{^2}**@{},
    \end{xy}}\Ea
 + (-1)^{c}
\Ba{c}\resizebox{7mm}{!}{\begin{xy}
 <0mm,-1.3mm>*{};<0mm,-3.5mm>*{}**@{-},
 <0.38mm,-0.2mm>*{};<2.0mm,2.0mm>*{}**@{-},
 <-0.38mm,-0.2mm>*{};<-2.2mm,2.2mm>*{}**@{-},
<0mm,-0.8mm>*{\bu};<0mm,0.8mm>*{}**@{},
 <2.4mm,2.4mm>*{\bu};<2.4mm,2.4mm>*{}**@{},
 <2.77mm,2.0mm>*{};<4.4mm,-0.8mm>*{}**@{-},
 <2.4mm,3mm>*{};<2.4mm,5.2mm>*{}**@{-},
     <0mm,-1.3mm>*{};<0mm,-5.3mm>*{^1}**@{},
     <2.5mm,2.3mm>*{};<5.1mm,-2.6mm>*{^2}**@{},
    <2.4mm,2.5mm>*{};<2.4mm,5.7mm>*{^1}**@{},
    <-0.38mm,-0.2mm>*{};<-2.8mm,2.5mm>*{^2}**@{},
    \end{xy}}\Ea=0\right).
    \Ea
\right.
\Eeq
where the vertices are ordered implicitly  in such a way that the ones on the top come first.

\sip

V.\ Drinfeld introduced  \cite{Dr2} the notion of {\it quasi-Lie bialgebra}\ or {\it Lie quasi-bialgebra}.
The prop(erad) $q\LBcd$ controlling degree shifted quasi-Lie bialgebras can be defined, for any pair of integer  $c,d\in \Z$, as the quotient
$$
q\LB_{c,d}:=\cF ree\langle E_q\rangle/\langle\cR_q\rangle,
$$
of the free prop(erad) generated by an  $\bS$-bimodule $Q=\{Q(m,n)\}_{m,n\geq 0}$ with
 all $Q(m,n)=0$ except
\Beqrn
Q(2,1)&:=&\id_1\ot \sgn_2^{c}[c-1]=\mbox{span}\left\langle
\Ba{c}\begin{xy}
 <0mm,-0.55mm>*{};<0mm,-3mm>*{}**@{-},
 <0.5mm,0.5mm>*{};<2.2mm,2.2mm>*{}**@{-},
 <-0.48mm,0.48mm>*{};<-2.2mm,2.2mm>*{}**@{-},
 <0mm,0mm>*{\bu};<0mm,0mm>*{}**@{},
 %<0mm,-0.55mm>*{};<0mm,-3.8mm>*{_1}**@{},
 <0.5mm,0.5mm>*{};<2.7mm,2.8mm>*{^{_2}}**@{},
 <-0.48mm,0.48mm>*{};<-2.7mm,2.8mm>*{^{_1}}**@{},
 \end{xy}\Ea
=(-1)^{c}
\Ba{c}\begin{xy}
 <0mm,-0.55mm>*{};<0mm,-3mm>*{}**@{-},
 <0.5mm,0.5mm>*{};<2.2mm,2.2mm>*{}**@{-},
 <-0.48mm,0.48mm>*{};<-2.2mm,2.2mm>*{}**@{-},
 <0mm,0mm>*{\bu};<0mm,0mm>*{}**@{},
 %<0mm,-0.55mm>*{};<0mm,-3.8mm>*{_1}**@{},
 <0.5mm,0.5mm>*{};<2.7mm,2.8mm>*{^{_1}}**@{},
 <-0.48mm,0.48mm>*{};<-2.7mm,2.8mm>*{^{_2}}**@{},
 \end{xy}\Ea
   \right\rangle,\\
%%%%%%
Q(1,2)&:=& \sgn_2^{d}\ot \id_1[d-1]=\mbox{span}\left\langle
\Ba{c}\begin{xy}
 <0mm,0.66mm>*{};<0mm,3mm>*{}**@{-},
 <0.39mm,-0.39mm>*{};<2.2mm,-2.2mm>*{}**@{-},
 <-0.35mm,-0.35mm>*{};<-2.2mm,-2.2mm>*{}**@{-},
 <0mm,0mm>*{\bu};<0mm,0mm>*{}**@{},
   <0.39mm,-0.39mm>*{};<2.9mm,-4mm>*{^{_2}}**@{},
   <-0.35mm,-0.35mm>*{};<-2.8mm,-4mm>*{^{_1}}**@{},
\end{xy}\Ea
=(-1)^{d}
\Ba{c}\begin{xy}
 <0mm,0.66mm>*{};<0mm,3mm>*{}**@{-},
 <0.39mm,-0.39mm>*{};<2.2mm,-2.2mm>*{}**@{-},
 <-0.35mm,-0.35mm>*{};<-2.2mm,-2.2mm>*{}**@{-},
 <0mm,0mm>*{\bu};<0mm,0mm>*{}**@{},
   <0.39mm,-0.39mm>*{};<2.9mm,-4mm>*{^{_1}}**@{},
   <-0.35mm,-0.35mm>*{};<-2.8mm,-4mm>*{^{_2}}**@{},
\end{xy}\Ea
\right\rangle,\\
%%%%%%%%%
Q(3,0)&:=& (\sgn_3)^{\ot|c|}[2c-d-1]=\mbox{span}\left\langle
\Ba{c}\begin{xy}
 <0mm,-1mm>*{\bu};<-4mm,3mm>*{^{_1}}**@{-},
 <0mm,-1mm>*{\bu};<0mm,3mm>*{^{_2}}**@{-},
 <0mm,-1mm>*{\bu};<4mm,3mm>*{^{_3}}**@{-},
 \end{xy}\Ea= (-1)^{c|\sigma|}
 \Ba{c}\begin{xy}
 <0mm,-1mm>*{\bu};<-6mm,3mm>*{^{_{\sigma(1)}}}**@{-},
 <0mm,-1mm>*{\bu};<0mm,3mm>*{^{_{\sigma(2)}}}**@{-},
 <0mm,-1mm>*{\bu};<6mm,3mm>*{^{_{\sigma(3)}}}**@{-},
 \end{xy}\Ea\ \  \forall \sigma\in \bS_3
\right\rangle,
\Eeqrn
modulo the ideal generated by the following relations
\Beq\label{3: R for qLieB}
\cR_q:\left\{\hspace{-2mm}
\Ba{c}
\displaystyle
\oint_{123}
\left(\hspace{-2mm} \Ba{c}\resizebox{9.4mm}{!}{
\begin{xy}
 <0mm,0mm>*{\bu};<0mm,0mm>*{}**@{},
 <0mm,-0.49mm>*{};<0mm,-3.0mm>*{}**@{-},
 <0.49mm,0.49mm>*{};<1.9mm,1.9mm>*{}**@{-},
 <-0.5mm,0.5mm>*{};<-1.9mm,1.9mm>*{}**@{-},
 <-2.3mm,2.3mm>*{\bu};<-2.3mm,2.3mm>*{}**@{},
 <-1.8mm,2.8mm>*{};<0mm,4.9mm>*{}**@{-},
 <-2.8mm,2.9mm>*{};<-4.6mm,4.9mm>*{}**@{-},
   <0.49mm,0.49mm>*{};<2.7mm,2.3mm>*{^3}**@{},
   <-1.8mm,2.8mm>*{};<0.4mm,5.3mm>*{^2}**@{},
   <-2.8mm,2.9mm>*{};<-5.1mm,5.3mm>*{^1}**@{},
 \end{xy}}\Ea
+
\Ba{c}\resizebox{12.5mm}{!}{
\begin{xy}
(0,0)*{\bu};(-4,5)*{^{1}}**@{-},
(0,0)*{\bu};(0,5)*{^{2}}**@{-},
(0,0)*{\bu};(4,5)*{\bu}**@{-},
(4,5)*{\bu};(4,10)*{^{3}}**@{-},
(4,5)*{\bu};(8,0)*{_{\, 1}}**@{-},
\end{xy}
}
\Ea
\hspace{-2mm}
\right)
=0
 \ , \ \ \
%%%%%%%%%%%%%% Lie %%%%%%%%%%%%%%%%%%%%%%%%
\oint_{123}\hspace{-1mm} \Ba{c}\resizebox{10.0mm}{!}{ \begin{xy}
 <0mm,0mm>*{\bu};<0mm,0mm>*{}**@{},
 <0mm,0.69mm>*{};<0mm,3.0mm>*{}**@{-},
 <0.39mm,-0.39mm>*{};<2.4mm,-2.4mm>*{}**@{-},
 <-0.35mm,-0.35mm>*{};<-1.9mm,-1.9mm>*{}**@{-},
 <-2.4mm,-2.4mm>*{\bu};<-2.4mm,-2.4mm>*{}**@{},
 <-2.0mm,-2.8mm>*{};<0mm,-4.9mm>*{}**@{-},
 <-2.8mm,-2.9mm>*{};<-4.7mm,-4.9mm>*{}**@{-},
    <0.39mm,-0.39mm>*{};<3.3mm,-4.0mm>*{^3}**@{},
    <-2.0mm,-2.8mm>*{};<0.5mm,-6.7mm>*{^2}**@{},
    <-2.8mm,-2.9mm>*{};<-5.2mm,-6.7mm>*{^1}**@{},
 \end{xy}}\Ea =0
\ , \ \ \
\text{Alt}^c_{\bS_4}\hspace{-1.5mm}
\Ba{c}\resizebox{11.4mm}{!}{
\begin{xy}
(0,0)*{\bu};(-4,5)*{^{1}}**@{-},
(0,0)*{\bu};(0,5)*{^{2}}**@{-},
(0,0)*{\bu};(4,5)*{\bu}**@{-},
(4,5)*{\bu};(1.5,10)*{^{3}}**@{-},
(4,5)*{\bu};(6.5,10)*{^{4}}**@{-},
\end{xy}}
\Ea
%+
%(-1)^c
%
%\Ba{c}\resizebox{11.4mm}{!}{
%\begin{xy}
%(0,0)*{\bu};(-4,5)*{^{4}}**@{-},
%(0,0)*{\bu};(0,5)*{^{1}}**@{-},
%(0,0)*{\bu};(4,5)*{\bu}**@{-},
%(4,5)*{\bu};(1.5,10)*{^{2}}**@{-},
%(4,5)*{\bu};(6.5,10)*{^{3}}**@{-},
%\end{xy}}
%\Ea\hspace{-1mm}\right)
=0 \vspace{2mm}\\
%%%%%%%%%%%%%%%%%%%%%%% Lie[1]Bi %%%%%%%%%%%%%%%
\hspace{-1mm} \Ba{c}\resizebox{6mm}{!}{\begin{xy}
 <0mm,2.47mm>*{};<0mm,0.12mm>*{}**@{-},
 <0.5mm,3.5mm>*{};<2.2mm,5.2mm>*{}**@{-},
 <-0.48mm,3.48mm>*{};<-2.2mm,5.2mm>*{}**@{-},
 <0mm,3mm>*{\bu};<0mm,3mm>*{}**@{},
  <0mm,-0.8mm>*{\bu};<0mm,-0.8mm>*{}**@{},
<-0.39mm,-1.2mm>*{};<-2.2mm,-3.5mm>*{}**@{-},
 <0.39mm,-1.2mm>*{};<2.2mm,-3.5mm>*{}**@{-},
     <0.5mm,3.5mm>*{};<2.8mm,5.7mm>*{^2}**@{},
     <-0.48mm,3.48mm>*{};<-2.8mm,5.7mm>*{^1}**@{},
   <0mm,-0.8mm>*{};<-2.7mm,-5.2mm>*{^1}**@{},
   <0mm,-0.8mm>*{};<2.7mm,-5.2mm>*{^2}**@{},
\end{xy}}\Ea\hspace{-2mm}
+(-1)^{cd+c+d}\left(\hspace{-2mm}
\Ba{c}\resizebox{7mm}{!}{\begin{xy}
 <0mm,-1.3mm>*{};<0mm,-3.5mm>*{}**@{-},
 <0.38mm,-0.2mm>*{};<2.0mm,2.0mm>*{}**@{-},
 <-0.38mm,-0.2mm>*{};<-2.2mm,2.2mm>*{}**@{-},
<0mm,-0.8mm>*{\bu};<0mm,0.8mm>*{}**@{},
 <2.4mm,2.4mm>*{\bu};<2.4mm,2.4mm>*{}**@{},
 <2.77mm,2.0mm>*{};<4.4mm,-0.8mm>*{}**@{-},
 <2.4mm,3mm>*{};<2.4mm,5.2mm>*{}**@{-},
     <0mm,-1.3mm>*{};<0mm,-5.3mm>*{^1}**@{},
     <2.5mm,2.3mm>*{};<5.1mm,-2.6mm>*{^2}**@{},
    <2.4mm,2.5mm>*{};<2.4mm,5.7mm>*{^2}**@{},
    <-0.38mm,-0.2mm>*{};<-2.8mm,2.5mm>*{^1}**@{},
    \end{xy}}\Ea\hspace{-1mm}
  + (-1)^{d}
\Ba{c}\resizebox{7mm}{!}{\begin{xy}
 <0mm,-1.3mm>*{};<0mm,-3.5mm>*{}**@{-},
 <0.38mm,-0.2mm>*{};<2.0mm,2.0mm>*{}**@{-},
 <-0.38mm,-0.2mm>*{};<-2.2mm,2.2mm>*{}**@{-},
<0mm,-0.8mm>*{\bu};<0mm,0.8mm>*{}**@{},
 <2.4mm,2.4mm>*{\bu};<2.4mm,2.4mm>*{}**@{},
 <2.77mm,2.0mm>*{};<4.4mm,-0.8mm>*{}**@{-},
 <2.4mm,3mm>*{};<2.4mm,5.2mm>*{}**@{-},
     <0mm,-1.3mm>*{};<0mm,-5.3mm>*{^2}**@{},
     <2.5mm,2.3mm>*{};<5.1mm,-2.6mm>*{^1}**@{},
    <2.4mm,2.5mm>*{};<2.4mm,5.7mm>*{^2}**@{},
    <-0.38mm,-0.2mm>*{};<-2.8mm,2.5mm>*{^1}**@{},
    \end{xy}}\Ea\hspace{-1mm}
  + (-1)^{d+c}
\Ba{c}\resizebox{7mm}{!}{\begin{xy}
 <0mm,-1.3mm>*{};<0mm,-3.5mm>*{}**@{-},
 <0.38mm,-0.2mm>*{};<2.0mm,2.0mm>*{}**@{-},
 <-0.38mm,-0.2mm>*{};<-2.2mm,2.2mm>*{}**@{-},
<0mm,-0.8mm>*{\bu};<0mm,0.8mm>*{}**@{},
 <2.4mm,2.4mm>*{\bu};<2.4mm,2.4mm>*{}**@{},
 <2.77mm,2.0mm>*{};<4.4mm,-0.8mm>*{}**@{-},
 <2.4mm,3mm>*{};<2.4mm,5.2mm>*{}**@{-},
     <0mm,-1.3mm>*{};<0mm,-5.3mm>*{^2}**@{},
     <2.5mm,2.3mm>*{};<5.1mm,-2.6mm>*{^1}**@{},
    <2.4mm,2.5mm>*{};<2.4mm,5.7mm>*{^1}**@{},
    <-0.38mm,-0.2mm>*{};<-2.8mm,2.5mm>*{^2}**@{},
    \end{xy}}\Ea\hspace{-1mm}
 + (-1)^{c}
\Ba{c}\resizebox{7mm}{!}{\begin{xy}
 <0mm,-1.3mm>*{};<0mm,-3.5mm>*{}**@{-},
 <0.38mm,-0.2mm>*{};<2.0mm,2.0mm>*{}**@{-},
 <-0.38mm,-0.2mm>*{};<-2.2mm,2.2mm>*{}**@{-},
<0mm,-0.8mm>*{\bu};<0mm,0.8mm>*{}**@{},
 <2.4mm,2.4mm>*{\bu};<2.4mm,2.4mm>*{}**@{},
 <2.77mm,2.0mm>*{};<4.4mm,-0.8mm>*{}**@{-},
 <2.4mm,3mm>*{};<2.4mm,5.2mm>*{}**@{-},
     <0mm,-1.3mm>*{};<0mm,-5.3mm>*{^1}**@{},
     <2.5mm,2.3mm>*{};<5.1mm,-2.6mm>*{^2}**@{},
    <2.4mm,2.5mm>*{};<2.4mm,5.7mm>*{^1}**@{},
    <-0.38mm,-0.2mm>*{};<-2.8mm,2.5mm>*{^2}**@{},
    \end{xy}}\Ea\hspace{-2mm}\right)=0.
    \Ea
\right.
\Eeq
Its minimal resolution $\HoqLBcd$ is a free operad
$
\HoqLBcd:=\cF ree \left\langle E\right\rangle
$
 generated by an $\bS$-bimodule
 $
 E_q=\{E_q(m,n)\}_{m\geq 1, n\geq 0, m+n\geq 3}
 $
 with
 $$
{E}_q(m,n):=\sgn_m^{\ot |c|}\ot \sgn_n^{|d|}[cm+dn-1-c-d]\equiv\text{span}\left\langle\hspace{-1mm}
\Ba{c}\resizebox{15mm}{!}{\begin{xy}
 <0mm,0mm>*{\bu};<0mm,0mm>*{}**@{},
 <-0.6mm,0.44mm>*{};<-8mm,5mm>*{}**@{-},
 <-0.4mm,0.7mm>*{};<-4.5mm,5mm>*{}**@{-},
 <0mm,0mm>*{};<1mm,5mm>*{\ldots}**@{},
 <0.4mm,0.7mm>*{};<4.5mm,5mm>*{}**@{-},
 <0.6mm,0.44mm>*{};<8mm,5mm>*{}**@{-},
   <0mm,0mm>*{};<-10.5mm,5.9mm>*{^{\sigma(1)}}**@{},
   <0mm,0mm>*{};<-4mm,5.9mm>*{^{\sigma(2)}}**@{},
   <0mm,0mm>*{};<10.0mm,5.9mm>*{^{\sigma(m)}}**@{},
 <-0.6mm,-0.44mm>*{};<-8mm,-5mm>*{}**@{-},
 <-0.4mm,-0.7mm>*{};<-4.5mm,-5mm>*{}**@{-},
 <0mm,0mm>*{};<1mm,-5mm>*{\ldots}**@{},
 <0.4mm,-0.7mm>*{};<4.5mm,-5mm>*{}**@{-},
 <0.6mm,-0.44mm>*{};<8mm,-5mm>*{}**@{-},
   <0mm,0mm>*{};<-10.5mm,-6.9mm>*{^{\tau(1)}}**@{},
   <0mm,0mm>*{};<-4mm,-6.9mm>*{^{\tau(2)}}**@{},
   %<0mm,0mm>*{};<4.5mm,-6.9mm>*{^{n\hspace{-0.5mm}-\hspace{-0.5mm}1}}**@{},
   <0mm,0mm>*{};<10.0mm,-6.9mm>*{^{\tau(n)}}**@{},
 \end{xy}}\Ea\hspace{-3mm}
=(-1)^{c|\sigma|+d|\tau|}\hspace{-1mm}
\Ba{c}\resizebox{12mm}{!}{\begin{xy}
 <0mm,0mm>*{\bu};<0mm,0mm>*{}**@{},
 <-0.6mm,0.44mm>*{};<-8mm,5mm>*{}**@{-},
 <-0.4mm,0.7mm>*{};<-4.5mm,5mm>*{}**@{-},
 <0mm,0mm>*{};<-1mm,5mm>*{\ldots}**@{},
 <0.4mm,0.7mm>*{};<4.5mm,5mm>*{}**@{-},
 <0.6mm,0.44mm>*{};<8mm,5mm>*{}**@{-},
   <0mm,0mm>*{};<-8.5mm,5.5mm>*{^1}**@{},
   <0mm,0mm>*{};<-5mm,5.5mm>*{^2}**@{},
   %<0mm,0mm>*{};<4.5mm,5.5mm>*{^{m\hspace{-0.5mm}-\hspace{-0.5mm}1}}**@{},
   <0mm,0mm>*{};<9.0mm,5.5mm>*{^m}**@{},
 <-0.6mm,-0.44mm>*{};<-8mm,-5mm>*{}**@{-},
 <-0.4mm,-0.7mm>*{};<-4.5mm,-5mm>*{}**@{-},
 <0mm,0mm>*{};<-1mm,-5mm>*{\ldots}**@{},
 <0.4mm,-0.7mm>*{};<4.5mm,-5mm>*{}**@{-},
 <0.6mm,-0.44mm>*{};<8mm,-5mm>*{}**@{-},
   <0mm,0mm>*{};<-8.5mm,-6.9mm>*{^1}**@{},
   <0mm,0mm>*{};<-5mm,-6.9mm>*{^2}**@{},
   %<0mm,0mm>*{};<4.5mm,-6.9mm>*{^{n\hspace{-0.5mm}-\hspace{-0.5mm}1}}**@{},
   <0mm,0mm>*{};<9.0mm,-6.9mm>*{^n}**@{},
 \end{xy}}\Ea\hspace{-1mm}
 \right\rangle_{ \forall \sigma\in \bS_m \atop \forall\tau\in \bS_n}
$$
The differential in $\HoqLBcd$ is given on the generators by
  %a formula similar to (\ref{3: d in LBcd_infty}),
 \Beq\label{3: d in qLBcd_infty}
\delta
\Ba{c}\resizebox{14mm}{!}{\begin{xy}
 <0mm,0mm>*{\bu};<0mm,0mm>*{}**@{},
 <-0.6mm,0.44mm>*{};<-8mm,5mm>*{}**@{-},
 <-0.4mm,0.7mm>*{};<-4.5mm,5mm>*{}**@{-},
 <0mm,0mm>*{};<-1mm,5mm>*{\ldots}**@{},
 <0.4mm,0.7mm>*{};<4.5mm,5mm>*{}**@{-},
 <0.6mm,0.44mm>*{};<8mm,5mm>*{}**@{-},
   <0mm,0mm>*{};<-8.5mm,5.5mm>*{^1}**@{},
   <0mm,0mm>*{};<-5mm,5.5mm>*{^2}**@{},
   <0mm,0mm>*{};<4.5mm,5.5mm>*{^{m\hspace{-0.5mm}-\hspace{-0.5mm}1}}**@{},
   <0mm,0mm>*{};<9.0mm,5.5mm>*{^m}**@{},
 <-0.6mm,-0.44
 mm>*{};<-8mm,-5mm>*{}**@{-},
 <-0.4mm,-0.7mm>*{};<-4.5mm,-5mm>*{}**@{-},
 <0mm,0mm>*{};<-1mm,-5mm>*{\ldots}**@{},
 <0.4mm,-0.7mm>*{};<4.5mm,-5mm>*{}**@{-},
 <0.6mm,-0.44mm>*{};<8mm,-5mm>*{}**@{-},
   <0mm,0mm>*{};<-8.5mm,-6.9mm>*{^1}**@{},
   <0mm,0mm>*{};<-5mm,-6.9mm>*{^2}**@{},
   <0mm,0mm>*{};<4.5mm,-6.9mm>*{^{n\hspace{-0.5mm}-\hspace{-0.5mm}1}}**@{},
   <0mm,0mm>*{};<9.0mm,-6.9mm>*{^n}**@{},
 \end{xy}}\Ea
\ \ = \ \
 \sum_{[m]=I_1\sqcup I_2\atop
 {|I_1|\geq 0, |I_2|\geq 1}}
 \sum_{[n]=J_1\sqcup J_2\atop
 {|J_1|, |J_2|\geq 0}
}\hspace{0mm}
\pm
\Ba{c}\resizebox{22mm}{!}{ \begin{xy}
 <0mm,0mm>*{\bu};<0mm,0mm>*{}**@{},
 <-0.6mm,0.44mm>*{};<-8mm,5mm>*{}**@{-},
 <-0.4mm,0.7mm>*{};<-4.5mm,5mm>*{}**@{-},
 <0mm,0mm>*{};<0mm,5mm>*{\ldots}**@{},
 <0.4mm,0.7mm>*{};<4.5mm,5mm>*{}**@{-},
 <0.6mm,0.44mm>*{};<12.4mm,4.8mm>*{}**@{-},
     <0mm,0mm>*{};<-2mm,7mm>*{\overbrace{\ \ \ \ \ \ \ \ \ \ \ \ }}**@{},
     <0mm,0mm>*{};<-2mm,9mm>*{^{I_1}}**@{},
 <-0.6mm,-0.44mm>*{};<-8mm,-5mm>*{}**@{-},
 <-0.4mm,-0.7mm>*{};<-4.5mm,-5mm>*{}**@{-},
 <0mm,0mm>*{};<-1mm,-5mm>*{\ldots}**@{},
 <0.4mm,-0.7mm>*{};<4.5mm,-5mm>*{}**@{-},
 <0.6mm,-0.44mm>*{};<8mm,-5mm>*{}**@{-},
      <0mm,0mm>*{};<0mm,-7mm>*{\underbrace{\ \ \ \ \ \ \ \ \ \ \ \ \ \ \
      }}**@{},
      <0mm,0mm>*{};<0mm,-10.6mm>*{_{J_1}}**@{},
 <13mm,5mm>*{};<13mm,5mm>*{\bu}**@{},
 <12.6mm,5.44mm>*{};<5mm,10mm>*{}**@{-},
 <12.6mm,5.7mm>*{};<8.5mm,10mm>*{}**@{-},
 <13mm,5mm>*{};<13mm,10mm>*{\ldots}**@{},
 <13.4mm,5.7mm>*{};<16.5mm,10mm>*{}**@{-},
 <13.6mm,5.44mm>*{};<20mm,10mm>*{}**@{-},
      <13mm,5mm>*{};<13mm,12mm>*{\overbrace{\ \ \ \ \ \ \ \ \ \ \ \ \ \ }}**@{},
      <13mm,5mm>*{};<13mm,14mm>*{^{I_2}}**@{},
 <12.4mm,4.3mm>*{};<8mm,0mm>*{}**@{-},
 <12.6mm,4.3mm>*{};<12mm,0mm>*{\ldots}**@{},
 <13.4mm,4.5mm>*{};<16.5mm,0mm>*{}**@{-},
 <13.6mm,4.8mm>*{};<20mm,0mm>*{}**@{-},
     <13mm,5mm>*{};<14.3mm,-2mm>*{\underbrace{\ \ \ \ \ \ \ \ \ \ \ }}**@{},
     <13mm,5mm>*{};<14.3mm,-4.5mm>*{_{J_2}}**@{},
 \end{xy}}\Ea
\Eeq
where the signs on the r.h.s\ are uniquely fixed for $c+d\in 2\Z$ by the fact that they all equal to $-1$ if $ c$ and $d$ are even integers.
Taking the quotient of $\HoqLBcd$ by the ideal generated by all $(m,0)$-corollas, $m\geq 3$, gives us the minimal model $\HoLBcd$ of the properad $\LBcd$.

\sip

The properads  $\HoLBcd$ and $q\HoLBcd$ with the same parity of $c+d$ are isomorphic to each other up to degree shift,
$$
\HoLBcd=\HoLB_{c+d,0}\{d\}, \ \ \ q\HoLBcd=q\HoLB_{c+d,0}\{d\},
$$
% HoLBcd{N} has generators of degree 1+c+d -cm-dn+ Nm-Nn; set N=-d and get 1+c+d -(c+d)m
%i.e. HoLB_c,d{-d}= HoLB_{c+d,0}
i.e.\ there are essentially two different types of the (quasi-)Lie bialgebra properads, even and odd ones.

\subsection{A short reminder on graph complexes}\label{3: subsec on GC and HGC}
The M.\ Kontsevich graph complexes come in a family $\GC_d$ parameterized by an integer $d\in \Z$. The complex $\GC_d$ for fixed $d$ is generated by arbitrary graphs $\Ga$  with valencies $|v|$ of vertices $v$ of $\Ga$ satisfying $|v|\geq 3$, and with the orientation $or$ defined on each graph $\Ga\in \GC_d$ as an ordering of edges (up to an even permutation) for $d$ even, and an ordering of vertices and half edges (again up to even permutation); each graph $\Ga$ has precisely two different orientations, $or$ and $-or$, and one identifies $(\Ga,or)=-(\Ga,-or)$ and abbreviates the pair $(\Ga,or)$ to $\Ga$.  The cohomological degree of $\Ga\in \GC_d$ is defined by
$$
|\Ga|=d(\# V(\Ga)-1) + (1-d) \#E(\Ga)
$$
The differential $\delta$ on $\GC_d$ is given by an action,
  $
  \delta\Ga=\sum_v \delta_v\Ga
  $,
   on each vertex
$
v= \Ba{c}\resizebox{7mm}{!}{\xy
 (0,0)*{\bullet}="a",
(2.3,4)*{}="1",
(-2.3,4)*{}="2",
(5,1)*{}="3",
(-5,1)*{}="4",
(-3.6,-3.6)*{}="5",
(3.6,-3.6)*{}="6",
(0,-4.5)*{}="7",
\ar @{-} "a";"1" <0pt>
\ar @{-} "a";"2" <0pt>
\ar @{-} "a";"3" <0pt>
\ar @{-} "a";"4" <0pt>
\ar @{-} "a";"5" <0pt>
\ar @{-} "a";"6" <0pt>
\ar @{-} "a";"7" <0pt>
\endxy}
\Ea
 $
 of a graph $\Ga\in \GC_d$
by splitting $v$ into two new vertices connected by an edge, and then re-attaching the edges attached earlier to $v$ to the new vertices in all possible ways,
$$
\delta_v:
 \Ba{c}\resizebox{7mm}{!}{\xy
 (0,0)*{\bullet}="a",
(2.3,4)*{}="1",
(-2.3,4)*{}="2",
(5,1)*{}="3",
(-5,1)*{}="4",
(-3.6,-3.6)*{}="5",
(3.6,-3.6)*{}="6",
(0,-4.5)*{}="7",
\ar @{-} "a";"1" <0pt>
\ar @{-} "a";"2" <0pt>
\ar @{-} "a";"3" <0pt>
\ar @{-} "a";"4" <0pt>
\ar @{-} "a";"5" <0pt>
\ar @{-} "a";"6" <0pt>
\ar @{-} "a";"7" <0pt>
\endxy}
\Ea
\ \lon \ \sum
\Ba{c}\resizebox{7mm}{!}{\xy
%
%(0,-11)*{\underbrace{\ \ \ \ \  \ \ \ \ \ \ \ }_{I''}},
%(0,11)*{\overbrace{ \ \ \  \ \ \ \ \ \ \ }^{I'}},
%
 (0,-2.3)*{\bullet}="a",
 (0,2.3)*{\bullet}="b",
(-7,-5)*{}="1",
(7,-6)*{}="2",
(-3,-8)*{}="3",
(3,-8)*{}="4",
(5,7)*{}="5",
(-5,7)*{}="6",
(0,8)*{}="7",
\ar @{-} "a";"b" <0pt>
\ar @{-} "a";"1" <0pt>
\ar @{-} "a";"2" <0pt>
\ar @{-} "a";"3" <0pt>
\ar @{-} "a";"4" <0pt>
\ar @{-} "b";"5" <0pt>
\ar @{-} "b";"6" <0pt>
\ar @{-} "b";"7" <0pt>
\endxy}
\Ea.
$$
It is very hard to compute the cohomology classes $\GC_d$ explicitly. Here are two examples of degree zero cycles in $\GC_2$
$$
\mathfrak{w}_3 =
 \Ba{c}\resizebox{12mm}{!}{
\xy
 (0,0)*{\bullet}="a",
(0,8)*{\bullet}="b",
(-7.5,-4.5)*{\bullet}="c",
(7.5,-4.5)*{\bullet}="d",
\ar @{-} "a";"b" <0pt>
\ar @{-} "a";"c" <0pt>
\ar @{-} "b";"c" <0pt>
\ar @{-} "d";"c" <0pt>
\ar @{-} "b";"d" <0pt>
\ar @{-} "d";"a" <0pt>
\endxy}
\Ea, \
\mathfrak{w}_5 =
 \Ba{c}\resizebox{13mm}{!}{
\xy
 (0,0)*{\bullet}="0",
(0,8)*{\bullet}="1",
(-8,3)*{\bullet}="5",
(8,3)*{\bullet}="2",
(-5,-7)*{\bullet}="4",
(5,-7)*{\bullet}="3",
\ar @{-} "0";"1" <0pt>
\ar @{-} "0";"2" <0pt>
\ar @{-} "0";"3" <0pt>
\ar @{-} "0";"4" <0pt>
\ar @{-} "0";"5" <0pt>
\ar @{-} "1";"2" <0pt>
\ar @{-} "2";"3" <0pt>
\ar @{-} "3";"4" <0pt>
\ar @{-} "4";"5" <0pt>
\ar @{-} "5";"1" <0pt>
\endxy}\Ea
  \ +\ \frac{5}{2}
 \Ba{c}\resizebox{12mm}{!}{
\xy
 (1,0)*{\bullet}="0",
(0,8)*{\bullet}="1",
(-8,3)*{\bullet}="5",
(8,3)*{\bullet}="2",
(-5,-7)*{\bullet}="4",
(5,-7)*{\bullet}="3",
\ar @{-} "0";"1" <0pt>
\ar @{-} "0";"2" <0pt>
\ar @{-} "0";"4" <0pt>
\ar @{-} "1";"4" <0pt>
\ar @{-} "5";"3" <0pt>
\ar @{-} "1";"2" <0pt>
\ar @{-} "2";"3" <0pt>
\ar @{-} "3";"4" <0pt>
\ar @{-} "4";"5" <0pt>
\ar @{-} "5";"1" <0pt>
\endxy}
\Ea,
$$
which represent non-trivial cohomology classes. It has been proven in \cite{W} that
$H^0(\GC_2)=\grt_1$, the Lie algebra of the Grothendieck-Teichm\"uller group $GRT_1$. Interestingly in the present context, the graph complexes $\GC_{c+d+1}$ control \cite{MW2} the homotopy theory of properads $\HoLBcd$ and $\LBcd$, i.e.\ there is a quasi-isomorphism, up to one rescaling class, of dg Lie algebras
$$
\GC_{c+d+1} \lon \Der(\HoLBcd) \simeq \Def(\LBcd \stackrel{\Id}{\rar} \LBcd)[1]
$$
where $\Der(\HoLBcd)$ is the derivation complex of genus completion of the properad $\HoLBcd$.

\sip

The graph complex $\HGC_d^N$ with $N$ labelled hairs is defined similarly --- the only novelty is that each graph $\Ga$ in $\HGC_d^N$ has precisely $N$ hairs (or legs) attached to its vertex or vertices. Again each vertex must be at least trivalent (with hairs counted), and the differential $\delta$ acts on vertices as before. One can understand hairs as kind of special univalent vertices on which $\delta$ does not act; they are assigned the same cohomological degree $1-d$ as edges.
The haired graph complexes has been  introduced and studied recently in the context of the theory of moduli spaces $\cM_{g,N}$ of algebraic curves of arbitrary genus $g$ with $N$ punctures \cite{CGP}
and the theory of long knots \cite{FTW}. It has been proven in \cite{CGP} that there is an isomorphism of cohomology groups
$$
H^\bu(\HGC_0^N)= \prod_{2g+N\geq 4} W_0 H_c^{\bu - N} \cM_{g,N}
$$
where $W_0 H_c^{\bu} \cM_{g,N}$ stands for the weight zero summand of the compactly supported cohomology of the moduli space $\cM_{g,N}$.

\subsection{Partial twisting of properads under $\LBcd$}\label{3: subsec on reduced twisting of (prop)erads under LBcd} Let $\cP=\{\cP(m,n), \p\}_{m,n\geq 0}$ be a dg properad. We represent its  generic elements pictorially as $(m,n)$-corollas
\Beq\label{3: generic elements of cP as (m,n)-corollas}
\Ba{c}\resizebox{11mm}{!}{ \xy
(0,4.5)*+{...},
(0,-4.5)*+{...},
(0,0)*{\circ}="o",
(-5,5)*{}="1",
(-3,5)*{}="2",
(3,5)*{}="3",
(5,5)*{}="4",
(-3,-5)*{}="5",
(3,-5)*{}="6",
(5,-5)*{}="7",
(-5,-5)*{}="8",
(-5.5,7)*{_1},
(-3,7)*{_2},
(3,6)*{},
(5.9,7)*{m},
(-3,-7)*{_2},
(3,-7)*+{},
(5.9,-7)*{n},
(-5.5,-7)*{_1},
\ar @{-} "o";"1" <0pt>
\ar @{-} "o";"2" <0pt>
\ar @{-} "o";"3" <0pt>
\ar @{-} "o";"4" <0pt>
\ar @{-} "o";"5" <0pt>
\ar @{-} "o";"6" <0pt>
\ar @{-} "o";"7" <0pt>
\ar @{-} "o";"8" <0pt>
\endxy}\Ea
\Eeq
whose white vertex is decorated by an element of $\cP(m,n)$. Properadic compositions in $\cP$ are represented pictorially by gluing out-legs of such decorated corollas to in-legs  of another decorated corollas.

\sip

Assume $\cP$ comes equipped with a non-trivial morphism
\Beq\label{3: i from LBcd to P}
i: \LBcd  \lon  \cP:\ \ \ \
  \Ba{c}\resizebox{8mm}{!}{  \xy
(-5,-6)*+{}="1";
    (-5,-0.2)*{\bu}="L";
  (-8,5)*+{_1}="C";
   (-2,5)*+{_2}="D";
\ar @{-} "D";"L" <0pt>
\ar @{-} "C";"L" <0pt>
\ar @{-} "1";"L" <0pt>
 \endxy}
 \Ea
\stackrel{i}{\rar}
 \Ba{c}\resizebox{8mm}{!}{  \xy
(-5,-6)*+{}="1";
    (-5,-0.2)*{\circledcirc}="L";
  (-8,5)*+{_1}="C";
   (-2,5)*+{_2}="D";
\ar @{-} "D";"L" <0pt>
\ar @{-} "C";"L" <0pt>
\ar @{-} "1";"L" <0pt>
 \endxy}
 \Ea,
 \ \ \ \
\Ba{c}\resizebox{8mm}{!}{  \xy
(-5,6)*{}="1";
    (-5,+1)*{\bu}="L";
  (-8,-5)*+{_1}="C";
   (-2,-5)*+{_2}="D";
\ar @{-} "D";"L" <0pt>
\ar @{-} "C";"L" <0pt>
\ar @{-} "1";"L" <0pt>
 \endxy}
 \Ea
\stackrel{i}{\rar}
 \Ba{c}\resizebox{8mm}{!}{  \xy
(-5,6)*{}="1";
    (-5,+1)*{\circledcirc}="L";
  (-8,-5)*+{_1}="C";
   (-2,-5)*+{_2}="D";
\ar @{-} "D";"L" <0pt>
\ar @{-} "C";"L" <0pt>
\ar @{-} "1";"L" <0pt>
 \endxy}
 \Ea
\Eeq
The images under $i$ of the generators of $\LBcd$ are special elements of $\cP$ and hence we reserve a special notation $\circledcirc$ for the decoration of the associated corollas. In particular, $\cP$ is a properad under $\Lie_d$ and hence can be twisted in the full analogy to the case of operads discussed in the previous section: applying T.\ Willwacher twisting endofunctor  to $(\cP,\p)$ we obtain a dg properad $(\tw \cP, \sd_\centerdot)$ called the {\em partial twisting  of a properad $\cP$ under $\HoLBcd$}. The latter is freely generated by $\cP$ and an extra generator  $\Ba{c}\resizebox{1.2mm}{!}{\begin{xy}
 <0mm,0.5mm>*{};<0mm,4mm>*{}**@{-},
 <0mm,0mm>*{\bu};<0mm,0mm>*{}**@{},
 \end{xy}}\Ea$ of degree $d$. The twisted differential $\pc$  acts on the latter generator by the standard formula (\ref{2: d_c on boxdot}), while its action on elements of $\cP$ is given by the following obvious analogue of (\ref{2: sd_c on generators of TwO}),
 \Beq\label{2: d_centerdot on twP under Lieb}
\sd_\centerdot \Ba{c}\resizebox{14mm}{!}{
 \begin{xy}
 <0mm,0mm>*{\circ};<-8mm,6mm>*{^1}**@{-},
 <0mm,0mm>*{\circ};<-4.5mm,6mm>*{^2}**@{-},
 <0mm,0mm>*{\circ};<0mm,5.5mm>*{\ldots}**@{},
 <0mm,0mm>*{\circ};<3.5mm,5mm>*{}**@{-},
 <0mm,0mm>*{\circ};<8mm,6mm>*{^m}**@{-},
 <0mm,0mm>*{\circ};<-8mm,-6mm>*{_1}**@{-},
 <0mm,0mm>*{\circ};<-4.5mm,-6mm>*{_2}**@{-},
 <0mm,0mm>*{\circ};<0mm,-5.5mm>*{\ldots}**@{},
 <0mm,0mm>*{\circ};<4.5mm,-6mm>*+{}**@{-},
 <0mm,0mm>*{\circ};<8mm,-6mm>*{_n}**@{-},
   \end{xy}}\Ea
=
\sd \Ba{c}\resizebox{14mm}{!}{
 \begin{xy}
 <0mm,0mm>*{\circ};<-8mm,6mm>*{^1}**@{-},
 <0mm,0mm>*{\circ};<-4.5mm,6mm>*{^2}**@{-},
 <0mm,0mm>*{\circ};<0mm,5.5mm>*{\ldots}**@{},
 <0mm,0mm>*{\circ};<3.5mm,5mm>*{}**@{-},
 <0mm,0mm>*{\circ};<8mm,6mm>*{^m}**@{-},
 <0mm,0mm>*{\circ};<-8mm,-6mm>*{_1}**@{-},
 <0mm,0mm>*{\circ};<-4.5mm,-6mm>*{_2}**@{-},
 <0mm,0mm>*{\circ};<0mm,-5.5mm>*{\ldots}**@{},
 <0mm,0mm>*{\circ};<4.5mm,-6mm>*+{}**@{-},
 <0mm,0mm>*{\circ};<8mm,-6mm>*{_n}**@{-},
   \end{xy}}\Ea
+
\overset{m-1}{\underset{i=0}{\sum}}
\Ba{c}\resizebox{15mm}{!}{
\begin{xy}
 %<0mm,0mm>*{\circ};<0mm,0mm>*{}**@{},
 <0mm,0mm>*{\circ};<-8mm,5mm>*{}**@{-},
 <0mm,0mm>*{\circ};<-3.5mm,5mm>*{}**@{-},
 <0mm,0mm>*{\circ};<-6mm,5mm>*{..}**@{},
 <0mm,0mm>*{\circ};<0mm,5mm>*{}**@{-},
  <0mm,13mm>*{\circledcirc};
  <0mm,13mm>*{};<5mm,10mm>*{_\bu}**@{-},
  <0mm,5mm>*{};<0mm,12mm>*{}**@{-},
  <0mm,14mm>*{};<0mm,17mm>*{}**@{-},
  <0mm,5mm>*{};<0mm,19mm>*{^{i\hspace{-0.2mm}+\hspace{-0.5mm}1}}**@{},
<0mm,0mm>*{\circ};<8mm,5mm>*{}**@{-},
<0mm,0mm>*{\circ};<3.5mm,5mm>*{}**@{-},
<6mm,5mm>*{..}**@{},
<-8.5mm,5.5mm>*{^1}**@{},
<-4mm,5.5mm>*{^i}**@{},
<9.0mm,5.5mm>*{^m}**@{},
 <0mm,0mm>*{\circ};<-8mm,-5mm>*{}**@{-},
 <0mm,0mm>*{\circ};<-4.5mm,-5mm>*{}**@{-},
 <-1mm,-5mm>*{\ldots}**@{},
 <0mm,0mm>*{\circ};<4.5mm,-5mm>*{}**@{-},
 <0mm,0mm>*{\circ};<8mm,-5mm>*{}**@{-},
<-8.5mm,-6.9mm>*{^1}**@{},
<-5mm,-6.9mm>*{^2}**@{},
<4.5mm,-6.9mm>*{^{n\hspace{-0.5mm}-\hspace{-0.5mm}1}}**@{},
<9.0mm,-6.9mm>*{^n}**@{},
 \end{xy}}\Ea
 - (-1)^{|a|}
\overset{n-1}{\underset{i=0}{\sum}}
 \Ba{c}\resizebox{14mm}{!}{\begin{xy}
 %<0mm,0mm>*{\circ};
 <0mm,0mm>*{\circ};<-8mm,-5mm>*{}**@{-},
 <0mm,0mm>*{\circ};<-3.5mm,-5mm>*{}**@{-},
 <0mm,0mm>*{\circ};<-6mm,-5mm>*{..}**@{},
 <0mm,0mm>*{\circ};<0mm,-5mm>*{}**@{-},
   <0mm,-11mm>*{\circledcirc};
  <0mm,-12mm>*{};<5mm,-16mm>*{_\bu}**@{-},
  <0mm,-5mm>*{};<0mm,-10mm>*{}**@{-},
  <0mm,-12mm>*{};<0mm,-17mm>*{}**@{-},
  <0mm,-5mm>*{};<0mm,-19mm>*{^{i\hspace{-0.2mm}+\hspace{-0.5mm}1}}**@{},
<0mm,0mm>*{\circ};<8mm,-5mm>*{}**@{-},
<0mm,0mm>*{\circ};<3.5mm,-5mm>*{}**@{-},
 <6mm,-5mm>*{..}**@{},
<-8.5mm,-6.9mm>*{^1}**@{},
<-4mm,-6.9mm>*{^i}**@{},
<9.0mm,-6.9mm>*{^n}**@{},
 <0mm,0mm>*{\circ};<-8mm,5mm>*{}**@{-},
 <0mm,0mm>*{\circ};<-4.5mm,5mm>*{}**@{-},
<-1mm,5mm>*{\ldots}**@{},
 <0mm,0mm>*{\circ};<4.5mm,5mm>*{}**@{-},
 <0mm,0mm>*{\circ};<8mm,5mm>*{}**@{-},
<-8.5mm,5.5mm>*{^1}**@{},
<-5mm,5.5mm>*{^2}**@{},
<4.5mm,5.5mm>*{^{m\hspace{-0.5mm}-\hspace{-0.5mm}1}}**@{},
<9.0mm,5.5mm>*{^m}**@{},
 \end{xy}}\Ea,
\Eeq
 The twisted properad comes equipped with a natural epimorphism of dg properads
 $$
 (\tw\cP, \sd_\centerdot) \lon (\cP, \sd)
 $$
 which sends the MC generator to zero. According to the general twisting machinery, the element  $\Ba{c}\resizebox{7mm}{!}{  \xy
(-5,6)*{}="1";
    (-5,+1)*{\circledcirc}="L";
  (-8,-5)*+{_1}="C";
   (-2,-5)*+{_2}="D";
\ar @{-} "D";"L" <0pt>
\ar @{-} "C";"L" <0pt>
\ar @{-} "1";"L" <0pt>
 \endxy}
 \Ea$ remains a cocycle in $\cP$ even after the twisting of the original differential so that the original morphism $i$ extends to the twisted version by the same formula,
\Beq\label{2: map i from Lie to twP}
\Ba{rccc}
\imath: & (\Lie_d, 0) & \lon & (\tw\cP,\sd_\centerdot)\\
&
\Ba{c}\resizebox{8mm}{!}{  \xy
(-5,6)*{}="1";
    (-5,+1)*{\bu}="L";
  (-8,-5)*+{_1}="C";
   (-2,-5)*+{_2}="D";
\ar @{-} "D";"L" <0pt>
\ar @{-} "C";"L" <0pt>
\ar @{-} "1";"L" <0pt>
 \endxy}
 \Ea
 & \lon & \Ba{c}\resizebox{8mm}{!}{  \xy
(-5,6)*{}="1";
    (-5,+1)*{\circledcirc}="L";
  (-8,-5)*+{_1}="C";
   (-2,-5)*+{_2}="D";
\ar @{-} "D";"L" <0pt>
\ar @{-} "C";"L" <0pt>
\ar @{-} "1";"L" <0pt>
%
 %
 %\ar @{-} "l1";"L" <0pt>
 \endxy}
 \Ea
\Ea.
\Eeq
However the image  of the co-Lie generator in $\cP$ is {\em not}, in general, respected by the  twisted
differential,
$$
\sd_\centerdot\hspace{-3mm}
\Ba{c}\resizebox{8mm}{!}{  \xy
(-5,-6)*+{_1}="1";
    (-5,0)*{\circledcirc}="L";
  (-8,5)*+{_1}="C";
   (-2,5)*+{_2}="D";
\ar @{-} "D";"L" <0pt>
\ar @{-} "C";"L" <0pt>
\ar @{-} "1";"L" <0pt>
 \endxy}
 \Ea
\hspace{-2mm}  = \hspace{-2mm}
 \Ba{c}\resizebox{11mm}{!}{  \xy
(-9,8)*{^2}="1";
    (-9,+2.5)*{\circledcirc}="L";
 (-13,-3)*{\circledcirc}="B";
 (-17,4)*+{^1}="b1";
 (-13,-9)*+{_1}="b2";
  (-5,-3)*{\bu}="C";
\ar @{-} "C";"L" <0pt>
\ar @{-} "B";"L" <0pt>
\ar @{-} "B";"b1" <0pt>
\ar @{-} "B";"b2" <0pt>
\ar @{-} "1";"L" <0pt>
 \endxy}
 \Ea
 +
 (-1)^c\hspace{-2mm}
  \Ba{c}\resizebox{11mm}{!}{  \xy
(-9,8)*{^1}="1";
    (-9,+2.5)*{\circledcirc}="L";
 (-13,-3)*{\circledcirc}="B";
 (-17,4)*+{^2}="b1";
 (-13,-9)*+{_1}="b2";
  (-5,-3)*{\bu}="C";
\ar @{-} "C";"L" <0pt>
\ar @{-} "B";"L" <0pt>
\ar @{-} "B";"b1" <0pt>
\ar @{-} "B";"b2" <0pt>
\ar @{-} "1";"L" <0pt>
 \endxy}
 \Ea
 -(-1)^{c-1}\hspace{-2mm}
  \Ba{c}\resizebox{8mm}{!}{  \xy
(-18,7)*{^1}="1";
    (-14,2.8)*{\circledcirc}="L";
 (-14,-2.8)*{\circledcirc}="B";
 (-18,-8)*{_1}="b1";
 (-10,-8)*{\bu}="b2";
  (-10,7)*{^2}="C";
\ar @{-} "C";"L" <0pt>
\ar @{-} "B";"L" <0pt>
\ar @{-} "B";"b1" <0pt>
\ar @{-} "B";"b2" <0pt>
\ar @{-} "1";"L" <0pt>
 \endxy}
 \Ea\hspace{-2mm}
 =\hspace{-2mm}
  \Ba{c}\resizebox{13mm}{!}{  \xy
(-19,8)*{^1}="1";
    (-19,+3)*{\circledcirc}="L";
 (-14,-2.5)*{\circledcirc}="B";
 (-9,3)*+{^2}="b1";
 (-14,-8)*{\bu}="b2";
  (-23,-3)*+{_1}="C";
\ar @{-} "C";"L" <0pt>
\ar @{-} "B";"L" <0pt>
\ar @{-} "B";"b1" <0pt>
\ar @{-} "B";"b2" <0pt>
\ar @{-} "1";"L" <0pt>
 \endxy}
 \Ea
 +
 (-1)^c\hspace{-2mm}
  \Ba{c}\resizebox{13mm}{!}{  \xy
(-19,8)*{^2}="1";
    (-19,+3)*{\circledcirc}="L";
 (-14,-2.5)*{\circledcirc}="B";
 (-9,3)*+{^1}="b1";
 (-14,-8)*{\bu}="b2";
  (-23,-3)*+{_1}="C";
\ar @{-} "C";"L" <0pt>
\ar @{-} "B";"L" <0pt>
\ar @{-} "B";"b1" <0pt>
\ar @{-} "B";"b2" <0pt>
\ar @{-} "1";"L" <0pt>
 \endxy}
 \Ea.
 %
%
%  \Ba{c}\resizebox{10.5mm}{!}{  \xy
%(-9,8)*{^1}="1";
%    (-9,+3)*{\bu}="L";
% (-14,-3.5)*{\bu}="B";
% (-19,5)*+{^2}="b1";
% (-14,-10)*{\bu}="b2";
%  (-3,-5)*+{_1}="C";
%     %
%\ar @{-} "C";"L" <0pt>
%\ar @{-} "B";"L" <0pt>
%\ar @{-} "B";"b1" <0pt>
%\ar @{-} "B";"b2" <0pt>
%\ar @{-} "1";"L" <0pt>
% \endxy}
% \Ea
 %
 %%%%%%%%%%%%%%%%%%
$$
where we used the image under $i$ of the third relation in (\ref{3: R for LieB}) (and  ordered vertices from bottom to the top).
The first equality equality in the formula just above, the formula (\ref{2: d_c on boxdot}) and the Drinfeld compatibility condition (that is, the bottom relation in (\ref{3: R for qLieB}))  imply
$$
\sd_\centerdot\hspace{-3mm}
\Ba{c}\resizebox{8mm}{!}{  \xy
(-5,-5.5)*{\bu}="1";
    (-5,0)*{\circledcirc}="L";
  (-8,5.5)*+{_1}="C";
   (-2,5.5)*+{_2}="D";
\ar @{-} "D";"L" <0pt>
\ar @{-} "C";"L" <0pt>
\ar @{-} "1";"L" <0pt>
 \endxy}
 \Ea\hspace{-2mm}
  =\hspace{-2mm}
 \Ba{c}\resizebox{11mm}{!}{  \xy
(-9,8)*{^2}="1";
    (-9,+2.5)*{\circledcirc}="L";
 (-13,-3)*{\circledcirc}="B";
 (-17,4)*+{^1}="b1";
 (-13,-9)*{\bu}="b2";
  (-5,-3)*{\bu}="C";
\ar @{-} "C";"L" <0pt>
\ar @{-} "B";"L" <0pt>
\ar @{-} "B";"b1" <0pt>
\ar @{-} "B";"b2" <0pt>
\ar @{-} "1";"L" <0pt>
 \endxy}
 \Ea
 +
 (-1)^c\hspace{-2mm}
  \Ba{c}\resizebox{11mm}{!}{  \xy
(-9,8)*{^1}="1";
    (-9,+2.5)*{\circledcirc}="L";
 (-13,-3)*{\circledcirc}="B";
 (-17,4)*+{^2}="b1";
 (-13,-9)*{\bu}="b2";
  (-5,-3)*{\bu}="C";
\ar @{-} "C";"L" <0pt>
\ar @{-} "B";"L" <0pt>
\ar @{-} "B";"b1" <0pt>
\ar @{-} "B";"b2" <0pt>
\ar @{-} "1";"L" <0pt>
 \endxy}
 \Ea
  -(-1)^{c-1}\hspace{-2mm}
   \Ba{c}\resizebox{8mm}{!}{  \xy
(-18,7)*{^1}="1";
    (-14,2.8)*{\circledcirc}="L";
 (-14,-2.8)*{\circledcirc}="B";
 (-18,-8)*{\bu}="b1";
 (-10,-8)*{\bu}="b2";
  (-10,7)*{^2}="C";
\ar @{-} "C";"L" <0pt>
\ar @{-} "B";"L" <0pt>
\ar @{-} "B";"b1" <0pt>
\ar @{-} "B";"b2" <0pt>
\ar @{-} "1";"L" <0pt>
 \endxy}
 \Ea
 +\frac{(-1)^{c-1}}{2}\hspace{-1mm}
   \Ba{c}\resizebox{8mm}{!}{  \xy
(-18,7)*{^1}="1";
    (-14,2.8)*{\circledcirc}="L";
 (-14,-2.8)*{\circledcirc}="B";
 (-18,-8)*{\bu}="b1";
 (-10,-8)*{\bu}="b2";
  (-10,7)*{^2}="C";
\ar @{-} "C";"L" <0pt>
\ar @{-} "B";"L" <0pt>
\ar @{-} "B";"b1" <0pt>
\ar @{-} "B";"b2" <0pt>
\ar @{-} "1";"L" <0pt>
 \endxy}
 \Ea=0
$$
which in turn implies that the element
$ \hspace{-2mm}\Ba{c}\resizebox{10mm}{!}{  \xy
(-19,8)*{^1}="1";
    (-19,+3)*{\circledcirc}="L";
 (-14,-2.5)*{\circledcirc}="B";
 (-9,3)*+{^2}="b1";
 (-14,-8)*{\bu}="b2";
  (-23,-3)*+{_1}="C";
\ar @{-} "C";"L" <0pt>
\ar @{-} "B";"L" <0pt>
\ar @{-} "B";"b1" <0pt>
\ar @{-} "B";"b2" <0pt>
\ar @{-} "1";"L" <0pt>
 \endxy}
 \Ea\hspace{-2mm}
 \in \tw \cP
$
is a cycle with respect to the twisted differential $\sd_\centerdot$. The linear combination
$$
 \Ba{c}\resizebox{12mm}{!}{  \xy
(-19,8)*{^1}="1";
    (-19,+3)*{\circledcirc}="L";
 (-14,-2.5)*{\circledcirc}="B";
 (-9,3)*+{^2}="b1";
 (-14,-8)*{\bu}="b2";
  (-23,-3)*+{_1}="C";
\ar @{-} "C";"L" <0pt>
\ar @{-} "B";"L" <0pt>
\ar @{-} "B";"b1" <0pt>
\ar @{-} "B";"b2" <0pt>
\ar @{-} "1";"L" <0pt>
 \endxy}
 \Ea
 +\la
 (-1)^c\hspace{-2mm}
  \Ba{c}\resizebox{12mm}{!}{  \xy
(-19,8)*{^2}="1";
    (-19,+3)*{\circledcirc}="L";
 (-14,-2.5)*{\circledcirc}="B";
 (-9,3)*+{^1}="b1";
 (-14,-8)*{\bu}="b2";
  (-23,-3)*+{_1}="C";
\ar @{-} "C";"L" <0pt>
\ar @{-} "B";"L" <0pt>
\ar @{-} "B";"b1" <0pt>
\ar @{-} "B";"b2" <0pt>
\ar @{-} "1";"L" <0pt>
 \endxy}
 \Ea
$$
is a $\sd_\centerdot$-coboundary for $\la=1$, but for other values of the parameter $\la$, say for $\la=-1,$ it represents, in general, a non-trivial cohomology class in $H^\bu(\tw\cP, \sd_\centerdot)$ of cohomological degree $2-c$.

%\sip

%Nevertheless the information about the Lie cobracket generator in $\cP$ is not lost completely as %the following Theorem shows.
%--- it gives rise to a quasi-Lie bialgebra structure in the cohomology $H^\bu(\tw \cP)$ of the %twisted properad.

\subsection{Theorem on partial twisting and quasi-Lie bialgebras}\label{2: qLien to twP}
 {\it Let $\cP$ be a dg properad under $\LBcd$ and $\tw \cP$ the associated twisting of $\cP$ as a properad under $\Lie_d$. Then there is an  explicit morphism of properads
$$
\Ba{rccl}
i^Q: & q\LB_{c-1,d} & \lon &  H^\bu(\tw\cP,\pc) \vspace{1mm}\\
&
 \Ba{c}\resizebox{8mm}{!}{  \xy
(-5,6)*{}="1";
    (-5,+1)*{\bu}="L";
  (-8,-5)*+{_1}="C";
   (-2,-5)*+{_2}="D";
\ar @{-} "D";"L" <0pt>
\ar @{-} "C";"L" <0pt>
\ar @{-} "1";"L" <0pt>
 \endxy}
 \Ea
 &\lon &
 \Ba{c}\resizebox{8mm}{!}{  \xy
(-5,6)*{}="1";
    (-5,+1)*{\circledcirc}="L";
  (-8,-5)*+{_1}="C";
   (-2,-5)*+{_2}="D";
\ar @{-} "D";"L" <0pt>
\ar @{-} "C";"L" <0pt>
\ar @{-} "1";"L" <0pt>
 \endxy}
 \Ea \mod \Img \pc \vspace{1mm} \\
&
\Ba{c}\resizebox{8mm}{!}{  \xy
(-5,-6.5)*+{_1}="1";
    (-5,0)*{\bu}="L";
  (-8.5,6)*+{_1}="C";
   (-1.5,6)*+{_2}="D";
\ar @{-} "D";"L" <0pt>
\ar @{-} "C";"L" <0pt>
\ar @{-} "1";"L" <0pt>
 \endxy}
 \Ea
 & \lon &\Ba{c}\resizebox{9mm}{!}{  \xy
(-5,-6.5)*+{_1}="1";
    (-5,0)*{\boxtimes}="L";
  (-8.5,6)*+{_1}="C";
   (-1.5,6)*+{_2}="D";
\ar @{-} "D";"L" <0pt>
\ar @{-} "C";"L" <0pt>
\ar @{-} "1";"L" <0pt>
 \endxy}
 \Ea:=
 %%%%%%%%%
 \frac{1}{2}
 \left( \hspace{-2mm}
  \Ba{c}\resizebox{13mm}{!}{  \xy
(-19,8)*{^1}="1";
    (-19,+3)*{\circledcirc}="L";
 (-14,-2.5)*{\circledcirc}="B";
 (-9,3)*+{^2}="b1";
 (-14,-8)*{\bu}="b2";
  (-23,-3)*+{_1}="C";
\ar @{-} "C";"L" <0pt>
\ar @{-} "B";"L" <0pt>
\ar @{-} "B";"b1" <0pt>
\ar @{-} "B";"b2" <0pt>
\ar @{-} "1";"L" <0pt>
 \endxy}
 \Ea
 -
 (-1)^c
  \Ba{c}\resizebox{13mm}{!}{  \xy
(-19,8)*{^2}="1";
    (-19,+3)*{\circledcirc}="L";
 (-14,-2.5)*{\circledcirc}="B";
 (-9,3)*+{^1}="b1";
 (-14,-8)*{\bu}="b2";
  (-23,-3)*+{_1}="C";
\ar @{-} "C";"L" <0pt>
\ar @{-} "B";"L" <0pt>
\ar @{-} "B";"b1" <0pt>
\ar @{-} "B";"b2" <0pt>
\ar @{-} "1";"L" <0pt>
 \endxy}
 \Ea
 \hspace{-2mm}
\right)\mod \Img \pc \vspace{2mm} \\
%%%%%%%%%%%%%%%%%%%%%%%%%%
&
\Ba{c}\resizebox{9mm}{!}{ \begin{xy}
 <0mm,-1mm>*{\bu};<-4mm,3mm>*{^{_{{1}}}}**@{-},
 <0mm,-1mm>*{\bu};<0mm,3mm>*{^{_{{2}}}}**@{-},
 <0mm,-1mm>*{\bu};<4mm,3mm>*{^{_{{3}}}}**@{-},
 \end{xy}}\Ea
 & \lon &
 \oint_{123}  \hspace{-2.5mm} \Ba{c}\resizebox{17mm}{!}{  \xy
(0,9)*{^2}="c";
(0,+3)*{\circledcirc}="C";
(-5,-3.5)*{\circledcirc}="d1";
(5,-3.5)*{\circledcirc}="d2";
(-11,4)*+{^1}="u1";
(11,4)*+{^3}="u2";
(-5,-10)*{\bu}="b1";
(5,-10)*{\bu}="b2";
\ar @{-} "d1";"C" <0pt>
\ar @{-} "d2";"C" <0pt>
\ar @{-} "c";"C" <0pt>
\ar @{-} "d1";"u1" <0pt>
\ar @{-} "d2";"u2" <0pt>
\ar @{-} "d1";"b1" <0pt>
\ar @{-} "d2";"b2" <0pt>
 \endxy}
 \Ea
\mod \Img \pc \\
\Ea
$$
}

\sip

The graphs on the r.h.s.\ are cycles in $\tw\cP$ which are not, in general, coboundaries.
What about $\cP=\LBcd$???? Is the image of the bottom graph is a non-trivial cohomology class?

 %{\color{blue} Why is the bottom graph a cycle? Yes.
%{\color{blue} The degree in the r.h.s., bottom line,  $1-d+ 2(1-c) +2d= 3+d-2c$; while on the %l.h.s.
%$1+d -2(c-1)=3+d-2c$. Agreed.}

\begin{proof}
Proof is a straightforward bur rather tedious calculation. Remarkably, the first and and the fourth relations in the list $\cR_Q$ above hold true exactly. However, the remaining third relation holds true only up to $\pc$-exact terms. Let us check it in full details that the map $i^Q$ satisfies
\Beq\label{3: relation (4,0) in qLieb}
\text{Alt}_{\bS_4}^{c-1}\ i^Q\left(\hspace{-1.5mm}
\Ba{c}\resizebox{9.0mm}{!}{
\begin{xy}
(0,0)*{\bu};(-4,5)*{^{1}}**@{-},
(0,0)*{\bu};(0,5)*{^{2}}**@{-},
(0,0)*{\bu};(4,5)*{\bu}**@{-},
(4,5)*{\bu};(1.5,10)*{^{3}}**@{-},
(4,5)*{\bu};(6.5,10)*{^{4}}**@{-},
\end{xy}}
\Ea
\hspace{-1mm}\right)
=0 \bmod \Img\pc .
\Eeq
We have
$$
 i^Q\left(\hspace{-2.8mm}
\Ba{c}\resizebox{9.4mm}{!}{
\begin{xy}
(0,0)*{\bu};(-4,5)*{^{1}}**@{-},
(0,0)*{\bu};(0,5)*{^{2}}**@{-},
(0,0)*{\bu};(4,5)*{\bu}**@{-},
(4,5)*{\bu};(1.5,10)*{^{3}}**@{-},
(4,5)*{\bu};(6.5,10)*{^{4}}**@{-},
\end{xy}}
\Ea
\hspace{-1.8mm}\right)\hspace{-0.5mm}
 =\hspace{-0.5mm}
 \frac{(\Id + (-1)^{c-1}(34))}{2}
 %%%%
\left(\hspace{-3mm}
 \Ba{c}\resizebox{21mm}{!}{  \xy
(0,9)*{^2}="c";
(0,+2)*{\circledcirc}="C";
(-5,-3.5)*{\circledcirc}="d1";
(5,-3.5)*{\circledcirc}="d2";
(15,-3.5)*{\circledcirc}="d3";
(-10,4)*+{^1}="l";
(10,9)*+{^3}="3";
(20,4)*+{^4}="4";
(10,2)*{\circledcirc}="r";
(-5,-10)*{\bu}="b1";
(5,-10)*{\bu}="b2";
(15,-10)*{\bu}="b3";
\ar @{-} "d1";"C" <0pt>
\ar @{-} "d2";"C" <0pt>
\ar @{-} "c";"C" <0pt>
\ar @{-} "d1";"l" <0pt>
\ar @{-} "d2";"r" <0pt>
\ar @{-} "3";"r" <0pt>
\ar @{-} "d1";"b1" <0pt>
\ar @{-} "d2";"b2" <0pt>
\ar @{-} "d3";"r" <0pt>
\ar @{-} "d3";"b3" <0pt>
\ar @{-} "d3";"4" <0pt>
 \endxy}
 \Ea
 \hspace{-1mm}
 +(-1)^{c-1}
 \hspace{-2.7mm}
  %%%%
  \Ba{c}\resizebox{21mm}{!}{  \xy
(0,9)*{^1}="c";
(0,+2)*{\circledcirc}="C";
(-5,-3.5)*{\circledcirc}="d1";
(5,-3.5)*{\circledcirc}="d2";
(15,-3.5)*{\circledcirc}="d3";
(-10,4)*+{^2}="l";
(10,9)*+{^3}="3";
(20,4)*+{^4}="4";
(10,2)*{\circledcirc}="r";
(-5,-10)*{\bu}="b1";
(5,-10)*{\bu}="b2";
(15,-10)*{\bu}="b3";
\ar @{-} "d1";"C" <0pt>
\ar @{-} "d2";"C" <0pt>
\ar @{-} "c";"C" <0pt>
\ar @{-} "d1";"l" <0pt>
\ar @{-} "d2";"r" <0pt>
\ar @{-} "3";"r" <0pt>
\ar @{-} "d1";"b1" <0pt>
\ar @{-} "d2";"b2" <0pt>
\ar @{-} "d3";"r" <0pt>
\ar @{-} "d3";"b3" <0pt>
\ar @{-} "d3";"4" <0pt>
 \endxy}
 \Ea
 + \hspace{-2mm}
  \Ba{c}\resizebox{20mm}{!}{  \xy
(10,4)*{^2}="c";
(0,+3)*{\circledcirc}="C";
(-5,-3.5)*{\circledcirc}="d1";
(5,-3.5)*{\circledcirc}="d2";
(15,9.5)*{\circledcirc}="d3";
(-10,5)*+{^1}="l";
(10,22)*+{^3}="3";
(20,17)*+{^4}="4";
(10,15)*{\circledcirc}="r";
(-5,-10)*{\bu}="b1";
(5,-10)*{\bu}="b2";
(15,2)*{\bu}="b3";
\ar @{-} "d1";"C" <0pt>
\ar @{-} "d2";"C" <0pt>
\ar @{-} "c";"d2" <0pt>
\ar @{-} "d1";"l" <0pt>
\ar @{-} "C";"r" <0pt>
\ar @{-} "3";"r" <0pt>
\ar @{-} "d1";"b1" <0pt>
\ar @{-} "d2";"b2" <0pt>
\ar @{-} "d3";"r" <0pt>
\ar @{-} "d3";"b3" <0pt>
\ar @{-} "d3";"4" <0pt>
 \endxy}
 \Ea
 \hspace{-2mm}\right).
$$
The Jacobi identity for the Lie generator implies the following vanishing
$$
\text{Alt}_{\bS_4}^{c-1} %
 \Ba{c}\resizebox{21mm}{!}{  \xy
(10,4)*{^2}="c";
(0,+3)*{\circledcirc}="C";
(-5,-3.5)*{\circledcirc}="d1";
(5,-3.5)*{\circledcirc}="d2";
(15,9.5)*{\circledcirc}="d3";
(-10,5)*+{^1}="l";
(10,22)*+{^3}="3";
(20,17)*+{^4}="4";
(10,15)*{\circledcirc}="r";
(-5,-10)*{\bu}="b1";
(5,-10)*{\bu}="b2";
(15,2)*{\bu}="b3";
\ar @{-} "d1";"C" <0pt>
\ar @{-} "d2";"C" <0pt>
\ar @{-} "c";"d2" <0pt>
\ar @{-} "d1";"l" <0pt>
\ar @{-} "C";"r" <0pt>
\ar @{-} "3";"r" <0pt>
\ar @{-} "d1";"b1" <0pt>
\ar @{-} "d2";"b2" <0pt>
\ar @{-} "d3";"r" <0pt>
\ar @{-} "d3";"b3" <0pt>
\ar @{-} "d3";"4" <0pt>
 \endxy}
 \Ea
 =0
$$
The symmetry properties of the generators imply, for any $c,d\in \Z$, the equality
$$
 \Ba{c}\resizebox{21mm}{!}{  \xy
(0,9)*{^2}="c";
(0,+2)*{\circledcirc}="C";
(-5,-3.5)*{\circledcirc}="d1";
(5,-3.5)*{\circledcirc}="d2";
(15,-3.5)*{\circledcirc}="d3";
(-10,4)*+{^1}="l";
(10,9)*+{^3}="3";
(20,4)*+{^4}="4";
(10,2)*{\circledcirc}="r";
(-5,-10)*{\bu}="b1";
(5,-10)*{\bu}="b2";
(15,-10)*{\bu}="b3";
\ar @{-} "d1";"C" <0pt>
\ar @{-} "d2";"C" <0pt>
\ar @{-} "c";"C" <0pt>
\ar @{-} "d1";"l" <0pt>
\ar @{-} "d2";"r" <0pt>
\ar @{-} "3";"r" <0pt>
\ar @{-} "d1";"b1" <0pt>
\ar @{-} "d2";"b2" <0pt>
\ar @{-} "d3";"r" <0pt>
\ar @{-} "d3";"b3" <0pt>
\ar @{-} "d3";"4" <0pt>
 \endxy}
 \Ea
 =
  \Ba{c}\resizebox{21mm}{!}{  \xy
(0,9)*{^3}="c";
(0,+2)*{\circledcirc}="C";
(-5,-3.5)*{\circledcirc}="d1";
(5,-3.5)*{\circledcirc}="d2";
(15,-3.5)*{\circledcirc}="d3";
(-10,4)*+{^4}="l";
(10,9)*+{^2}="3";
(20,4)*+{^1}="4";
(10,2)*{\circledcirc}="r";
(-5,-10)*{\bu}="b1";
(5,-10)*{\bu}="b2";
(15,-10)*{\bu}="b3";
\ar @{-} "d1";"C" <0pt>
\ar @{-} "d2";"C" <0pt>
\ar @{-} "c";"C" <0pt>
\ar @{-} "d1";"l" <0pt>
\ar @{-} "d2";"r" <0pt>
\ar @{-} "3";"r" <0pt>
\ar @{-} "d1";"b1" <0pt>
\ar @{-} "d2";"b2" <0pt>
\ar @{-} "d3";"r" <0pt>
\ar @{-} "d3";"b3" <0pt>
\ar @{-} "d3";"4" <0pt>
 \endxy}
 \Ea
$$
Therefore the first two summands in the above formula do not cancel out upon (skew)symmetrization. However one has the equality modulo $\pc$-exact terms,
$$
 \Ba{c}\resizebox{12mm}{!}{  \xy
(-19,8)*{^2}="1";
    (-19,+3)*{\circledcirc}="L";
 (-14,-2.5)*{\circledcirc}="B";
 (-9,3)*+{^x}="b1";
 (-14,-8)*{\bu}="b2";
  (-23,-3)*+{_y}="C";
\ar @{-} "C";"L" <0pt>
\ar @{-} "B";"L" <0pt>
\ar @{-} "B";"b1" <0pt>
\ar @{-} "B";"b2" <0pt>
\ar @{-} "1";"L" <0pt>
 \endxy}
 \Ea
 =(-1)^{c-1}
 \hspace{-2mm}
  \Ba{c}\resizebox{12mm}{!}{  \xy
(-19,8)*{^x }="1";
    (-19,+3)*{\circledcirc}="L";
 (-14,-2.5)*{\circledcirc}="B";
 (-9,3)*+{^2}="b1";
 (-14,-8)*{\bu}="b2";
  (-23,-3)*+{_1}="C";
\ar @{-} "C";"L" <0pt>
\ar @{-} "B";"L" <0pt>
\ar @{-} "B";"b1" <0pt>
\ar @{-} "B";"b2" <0pt>
\ar @{-} "1";"L" <0pt>
 \endxy}
 \Ea
 $$
which can be used to transform  first two terms in the above formula into the third one (up to a permutation) which has been just considered. Hence
$$
\text{Alt}_{\bS_4}^{c-1}  \Ba{c}\resizebox{22mm}{!}{  \xy
(0,10)*{^2}="c";
(0,+3)*{\circledcirc}="C";
(-5,-3.5)*{\circledcirc}="d1";
(5,-3.5)*{\circledcirc}="d2";
(15,-3.5)*{\circledcirc}="d3";
(-13,5)*+{^1}="l";
(10,10)*+{^3}="3";
(22,5)*+{^4}="4";
(10,3)*{\circledcirc}="r";
(-5,-10)*{\bu}="b1";
(5,-10)*{\bu}="b2";
(15,-10)*{\bu}="b3";
\ar @{-} "d1";"C" <0pt>
\ar @{-} "d2";"C" <0pt>
\ar @{-} "c";"C" <0pt>
\ar @{-} "d1";"l" <0pt>
\ar @{-} "d2";"r" <0pt>
\ar @{-} "3";"r" <0pt>
\ar @{-} "d1";"b1" <0pt>
\ar @{-} "d2";"b2" <0pt>
\ar @{-} "d3";"r" <0pt>
\ar @{-} "d3";"b3" <0pt>
\ar @{-} "d3";"4" <0pt>
 \endxy}
 \Ea
 = (-1)^{c-1}
 \text{Alt}_{\bS_4}^{c-1} %
  \Ba{c}\resizebox{22mm}{!}{  \xy
(10,4)*{^2}="c";
(0,+3)*{\circledcirc}="C";
(-5,-3.5)*{\circledcirc}="d1";
(5,-3.5)*{\circledcirc}="d2";
(15,9.5)*{\circledcirc}="d3";
(-13,5)*+{^1}="l";
(10,22)*+{^3}="3";
(22,17)*+{^4}="4";
(10,15)*{\circledcirc}="r";
(-5,-10)*{\bu}="b1";
(5,-10)*{\bu}="b2";
(15,2)*{\bu}="b3";
\ar @{-} "d1";"C" <0pt>
\ar @{-} "d2";"C" <0pt>
\ar @{-} "c";"d2" <0pt>
\ar @{-} "d1";"l" <0pt>
\ar @{-} "C";"r" <0pt>
\ar @{-} "3";"r" <0pt>
\ar @{-} "d1";"b1" <0pt>
\ar @{-} "d2";"b2" <0pt>
\ar @{-} "d3";"r" <0pt>
\ar @{-} "d3";"b3" <0pt>
\ar @{-} "d3";"4" <0pt>
 \endxy}
 \Ea
 \bmod \pc =
  \ \ \ \   0 \bmod \Img\pc
$$
and the claim follows.
\end{proof}

\subsection{Haired graph complexes and $\tw\LBcd$}\label{3: subsec on reduced twisting of (prop)erads under HoLBcd}
%Do we have a map of dg properads
%$$
%q\HoLB_{c-1,d}\lon \tw \HoLBcd ?????
%$$
%$(m,n)$-generators of the l.h.s. have degree $1+(c-1) + d -m(c-1) - dn=(1+c+d -cm-dn) +m-1$. Hence it is hard %to imagine what could shift degree by $m-1$..., say $m=m'+m''$, so we use two generators of degrees
%$$
%1+c+d -c(m'+1)-dn', \ \ \ 1+c+d -c(m'')-d(n''+1),\ \ \text{and some number k of black vertices of degree}\ %d.
%$$
%so
%$$
%m-1= kd -d +c
%$$
The epimorphism
$$
\tw \HoLBcd=\{\tw \HoLBcd(N,M)\} \lon \tw\LBcd=\{\tw\LBcd(N,M)\}
$$
is a quasi-isomorphism for any $M,N\geq 1$ as $\tw$ is an exact functor. It is a straightforward inspection to see the complex $\tw \HoLBcd(N,0)$ is identical to th oriented haired graph complex $\mathsf{HHOGC}^N_{c+d+1}[-dN]$ introduced
in \S 3.4.2 of \cite{AWZ}. One of the main results in that paper says that there is an isomorphism of cohomology groups
$$
H^\bu(\HGC^N_d)\simeq H^\bu(\mathsf{HHOGC}^N_{d+1})
$$
%\subsubsection{\bf Proposition}\label{3: Prop on HGC and twHoLB}
Hence we can conclude that {\em for any natural number $N\geq 1$ one has an isomorphisms of cohomology groups
%
%cohomology group  $H^\bu(\tw\HoLBcd)(m,0)$ of the $(m,0)$ part of the dg properad %$\tw\HoLBcd=\{\tw\HoLBcd(m,n)\}_{m\geq 1, n\geq 0}$ can be identified
%with the cohomology group of the haired graph complex $\mathsf{HGC}_?^S$ with $|S|=m$ introduced in [CGP] %and [AWZ]. More precisely, one has
$$
H^\bu\left(\tw \HoLB_{c,d}(N,0)\right)\simeq H^\bu(\tw\LB_{c,d}(N,0))%\simeq H^\bu(\mathsf{HHGC}_n^S)
\simeq H^{\bu-dN}(\mathsf{HGC}_{c+d}^N).
$$
}
%the generator qHoLBcd((|S|,0)=HoLB_{c+d,0}{d}=\HoLB_{c+d,0}[-d|S|]= HGC_{c+d}^S[-d|S|]
%The cohomology graded vector space $H^\bu(\tw\LB_{c,d})$ can be described explicitly in terms of
%generators and relations...........................
In particular, we have an induced morphism
$$
H^\bu(\HGC^N_{c+d}) \lon H^{\bu + dN}(\tw\cP(N,0))
$$
for any dg properad $\cP\in \PROP_{\HoLBcd}$.

\subsection{An example: (chain) gravity properad} A ribbon graph $\Ga$ is a graph with an extra structure: the set of half-edges attached to each vertex comes equipped with a cyclic ordering
(a detailed definition can be found, e.g., in \S 4.1 of the paper \cite{MW1} to which we refer often in this subsection).  Thickening each vertex $v\in V(\Ga)$ of a ribbon graph $\Ga$ into a closed disk, and every edge $e\in \Ga$ attached to $v$ into a 2-dimensional strip glued to that disk, one associates to $\Ga$ a unique topological $2$-dimensional-surface with boundaries; the set of such boundaries is denoted by $B(\Ga)$. Shrinking 2-strips back into 1-dimensional edges, one represents each boundary $b$ as a closed path comprising some vertices and edges of $\Ga$. We work with {\em connected}\, ribbon graphs only, their genus is defined by
 \Beq\label{4: genus of a ribbon graph}
 g= 1+\frac{1}{2}\left(\# E(\Ga) - \# V(\Ga)- \# B(\Ga)\right).
 \Eeq
Let $\RGra_d(m,n)$, $d\in \Z$, stand for the graded vector space generated by ribbon graphs $\Ga$ such that
\Bi
\item[(i)] $\Ga$ has precisely $n$ vertices and $m$ boundaries which are labelled, i.e.\ some isomorphisms $V(\Ga)\rar [n]$ and $B(\Ga)\rar [\bar{m}]:=\{\bar{1},\ldots, \bar{m}\}$ are fixed;
\item[(ii)] $\Ga$ is equipped with an orientation which is for $d$ even is defined as an ordering of edges (up the sign action of $\bS_{\# E(\Ga)}$), while for $d$ odd it is defined as a choice of the direction on each edge (up to the sign action of $\bS_2$).
\item[(iii)] $\Ga$ is assigned the cohomological degree $(1-d)\# E(\Ga)$.
\Ei
For example,
$$
\Ba{c}\resizebox{8.2mm}{!}{  \xy
 (3.5,4)*{^{\bar{1}}};
 (7,0)*+{_2}*\frm{o}="A";
 (0,0)*+{_1}*\frm{o}="B";
 \ar @{-} "A";"B" <0pt>
\endxy} \Ea\hspace{-2mm} \in \RGra_d(1,2),
 \Ba{c}\resizebox{7mm}{!}{ \xy
 (0.5,1)*{^{{^{\bar{1}}}}},
(0.5,5)*{^{{^{\bar{2}}}}},
 (0,-2)*+{_{_1}}*\frm{o}="A";
"A"; "A" **\crv{(7,7) & (-7,7)};
\endxy}\Ea \hspace{-2mm} \in \RGra_d(2,1),
 \Ba{c}\resizebox{9mm}{!}{ \xy
 (-4,0)*{^{\bar{_1}}};
 (-1,0)*{^{\bar{_2}}};
 (1.5,0)*{^{\bar{_3}}};
  (0,5)*+{_1}*\frm{o}="1";
(0,-4)*+{_2}*\frm{o}="3";
"1";"3" **\crv{(4,0) & (4,1)};
"1";"3" **\crv{(-4,0) & (-4,-1)};
\ar @{-} "1";"3" <0pt>
\endxy}\Ea \hspace{-2mm}\in \RGra_d(3,2),
\Ba{c}\resizebox{9mm}{!}{ \xy
 (-3,1)*{^{\bar{_1}}};
 (0,8)*+{_1}*\frm{o}="1";
%(0,5)*{\circ}="1";
(0,-4)*+{_2}*\frm{o}="3";
"1";"3" **\crv{(-5,2) & (5,2)};
"1";"3" **\crv{(5,2) & (-5,2)};
"1";"3" **\crv{(-7,7) & (-7,-7)};
%\ar @{-} "1";"3" <0pt>
\endxy}\Ea \hspace{-2mm}\in \RGra_d(1,2).
$$
The subspace of $\RGra_d(m,n)$ spanned by ribbon graphs of genus $g$ is denoted by $\RGra_d(g;m,n)$.
The permutation group $\bS_.^{op}\times \bS_n$ acts on $\RGra_d(m,n)$ by relabelling vertices and boundaries. The $\bS$-bimodule
$$
\RGra_d=\{\RGra_d(m,n)\}
$$
has the structure of a properad \cite{MW1} given by substituting a boundary $b$  of one ribbon graph into a vertex $v$ of another one, and reattaching half-edges (attached earlier to $v$) among the vertices belonging to $b$ in all possible ways while respecting the cyclic orders of both sets. One of the main motivation behind this definition of $\RGra_d$ is that it comes with a morphism of operads,
\Beq\label{3: Lieb_dd to RGra_d}
\Ba{rccc}
i: & \LB_{d,d} & \lon & \RGra_d\\
   & \Ba{c}\begin{xy}
 <0mm,0.66mm>*{};<0mm,3mm>*{}**@{-},
 <0.39mm,-0.39mm>*{};<2.2mm,-2.2mm>*{}**@{-},
 <-0.35mm,-0.35mm>*{};<-2.2mm,-2.2mm>*{}**@{-},
 <0mm,0mm>*{\bu};<0mm,0mm>*{}**@{},
   <0.39mm,-0.39mm>*{};<2.9mm,-4mm>*{^{_2}}**@{},
   <-0.35mm,-0.35mm>*{};<-2.8mm,-4mm>*{^{_1}}**@{},
 <0mm,4mm>*{^{_{\bar{1}}}}**@{},
\end{xy}\Ea
&\lon&
\Ba{c}\resizebox{8.2mm}{!}{  \xy
 (3.5,4)*{^{\bar{1}}};
 (7,0)*+{_2}*\frm{o}="A";
 (0,0)*+{_1}*\frm{o}="B";
 \ar @{-} "A";"B" <0pt>
\endxy} \Ea
\vspace{1.1mm}\\
& \Ba{c}\begin{xy}
 <0mm,-0.55mm>*{};<0mm,-2.5mm>*{}**@{-},
 <0.5mm,0.5mm>*{};<2.2mm,2.2mm>*{}**@{-},
 <-0.48mm,0.48mm>*{};<-2.2mm,2.2mm>*{}**@{-},
 <0mm,0mm>*{\bu};<0mm,0mm>*{}**@{},
 <0mm,-0.55mm>*{};<0mm,-3.8mm>*{_1}**@{},
 <0.5mm,0.5mm>*{};<2.7mm,2.8mm>*{^{_{\bar{1}}}}**@{},
 <-0.48mm,0.48mm>*{};<-2.7mm,2.8mm>*{^{_{\bar{2}}}}**@{},
 \end{xy}\Ea
 &\lon &
   \Ba{c}\resizebox{7mm}{!}{ \xy
 (0.5,1)*{^{{^{\bar{1}}}}},
(0.5,5)*{^{{^{\bar{2}}}}},
 (0,-2)*+{_{_1}}*\frm{o}="A";
"A"; "A" **\crv{(7,7) & (-7,7)};
\endxy}\Ea
\Ea
\Eeq
In particular, $\RGra_d$ is a properad under $\Lie_d$ and hence can be twisted: $\tw\RGra_d$
is generated by ribbon graphs with two types of vertices, white ones which are labelled and black ones which are unlabelled and assigned the cohomological degree $d$ (cf.\ \S {\ref{2: Example Graphs_d operad}}), e.g.
$$
\Ba{c}\resizebox{9mm}{!}{ \xy
 (-3,1)*{^{\bar{_1}}};
 (0,8)*+{_1}*\frm{o}="1";
%(0,5)*{\circ}="1";
(0,-4)*{\bu}="3";
"1";"3" **\crv{(-5,2) & (5,2)};
"1";"3" **\crv{(5,2) & (-5,2)};
"1";"3" **\crv{(-7,7) & (-7,-7)};
%\ar @{-} "1";"3" <0pt>
\endxy}\Ea \in \tw\cR\cG ra_d(1,1).
$$
The differential $\delta_\centerdot$ in $\tw\RGra_d$ is determined by its action on vertices as in (\ref{2: d in Graphs}). One of the main results in \cite{Me2} is the proof of the following

\subsubsection{\bf Theorem} (i) {\em For any $g\geq 0$, $m\geq 1$ and $n\geq 0$ with $2g+m+n\geq 3$ one has an isomorphism of $\bS_m^{op}\times \bS_n$-modules,}
$$
H^\bu(\tw\RGra_d(g;m,n))= H^{\bu-m +d(2g-2+m+n)}_c(\cM_{g,m+n}\times \R^m)
$$
where $\cM_{g,m+n}$ is the moduli space of genus algebraic curves with $m+n$ marked points, and
$H^\bu_c$ stands for the compactly supported cohomology functor.
\sip

(ii) {\em For any $g\geq 0$, $m\geq 1$ and $n\geq 0$ with $2g+m+n< 3$ one has}
$$
H^k(\tw\RGra_d(g;m,n))=\left\{\Ba{ll} \K & \text{if} \ g=n=0,m=2, k=(1-d)p \ \text{with} \ p\geq 1\ \& \ p\equiv 2d+1 \bmod 4 \\
0 & \text{otherwise}.
\Ea
\right.
$$
{\em where $\K$ is generated by  the unique polytope-like ribbon graph with $p$ edges and $p$ bivalent vertices which are all black}.

\sip

This result says that the most important part of $\tw\RGra_d$ is the dg sub-properad
$\cC h\Grav_d$ spanned by ribbon graphs with black vertices at least trivalent; it is called the {\em chain gravity properad}. Its cohomology
$$
\Grav_d:=\left\{\prod_{g\geq 0\atop 2g\geq 3-m-n}  H^{\bu-m +d(2g-2+m+n)}_c(\cM_{g,m+n}\times \R^m)\right\}_{m\geq1, n\geq 0}
$$
is called the {\em gravity}\, properad.
The general morphism  $i^Q$ from Theorem {\ref{2: qLien to twP}} reads in this concrete situation as follows
$$
\Ba{rccc}
i^Q: & q\LB_{d-1,d} & \lon & \GRav_d \vspace{1mm}\\
&
\Ba{c}\begin{xy}
 <0mm,0.66mm>*{};<0mm,4mm>*{^{_{\bar{1}}}}**@{-},
 <0.39mm,-0.39mm>*{};<2.2mm,-2.2mm>*{}**@{-},
 <-0.35mm,-0.35mm>*{};<-2.2mm,-2.2mm>*{}**@{-},
 <0mm,0mm>*{\bu};<0mm,0mm>*{}**@{},
   <0.39mm,-0.39mm>*{};<2.9mm,-4mm>*{^{_2}}**@{},
   <-0.35mm,-0.35mm>*{};<-2.8mm,-4mm>*{^{_1}}**@{},
\end{xy}\Ea
 &\lon &
  \Ba{c}\resizebox{10mm}{!}{  \xy
 (3.5,4)*{^{\bar{1}}};
 (7,0)*+{_2}*\frm{o}="A";
 (0,0)*+{_1}*\frm{o}="B";
 \ar @{-} "A";"B" <0pt>
\endxy} \Ea \vspace{1mm} \\
&
\Ba{c}\begin{xy}
 <0mm,-0.55mm>*{};<0mm,-2.5mm>*{}**@{-},
 <0.5mm,0.5mm>*{};<2.2mm,2.2mm>*{}**@{-},
 <-0.48mm,0.48mm>*{};<-2.2mm,2.2mm>*{}**@{-},
 <0mm,0mm>*{\bu};<0mm,0mm>*{}**@{},
 <0mm,-0.55mm>*{};<0mm,-3.8mm>*{_1}**@{},
 <0.5mm,0.5mm>*{};<2.7mm,2.8mm>*{^{_{\bar{2}}}}**@{},
 <-0.48mm,0.48mm>*{};<-2.7mm,2.8mm>*{^{_{\bar{1}}}}**@{},
 \end{xy}\Ea
 & \lon &
 \frac{1}{2}
 \left(\hspace{-1.5mm}
\Ba{c}\resizebox{6mm}{!}{
\mbox{$\xy
 (0.5,0.9)*{^{{^{\bar{1}}}}},
(0.5,5)*{^{{^{\bar{2}}}}},
 (0,-8)*+{_{_1}}*\frm{o}="C";
(0,-2)*{\bu}="A";
(0,-2)*{\bu}="B";
"A"; "B" **\crv{(6,6) & (-6,6)};
 \ar @{-} "A";"C" <0pt>
\endxy$}}
\Ea
- (-1)^d
\Ba{c}\resizebox{6mm}{!}{
\mbox{$\xy
 (0.5,0.9)*{^{{^{\bar{2}}}}},
(0.5,5)*{^{{^{\bar{1}}}}},
 (0,-8)*+{_{_1}}*\frm{o}="C";
(0,-2)*{\bu}="A";
(0,-2)*{\bu}="B";
"A"; "B" **\crv{(6,6) & (-6,6)};
 \ar @{-} "A";"C" <0pt>
\endxy$}}
\Ea
\hspace{-1.5mm}\right) \vspace{2mm} \\
%%%%%%%%%%%%%%%%%%%%%%%%%%
&
\Ba{c}\begin{xy}
 <0mm,-1mm>*{\bu};<-4mm,3mm>*{^{_{\bar{1}}}}**@{-},
 <0mm,-1mm>*{\bu};<0mm,3mm>*{^{_{\bar{2}}}}**@{-},
 <0mm,-1mm>*{\bu};<4mm,3mm>*{^{_{\bar{3}}}}**@{-},
 \end{xy}\Ea
 & \lon &
 \frac{1}{2}\left(\hspace{-1.5mm}
 \Ba{c}\resizebox{13.5mm}{!}{
\mbox{$\xy
 (0.5,0.9)*{^{{^{\bar{1}}}}},
(5,0)*{^{{^{\bar{2}}}}},
 (10.5,0.9)*{^{{^{\bar{3}}}}},
 (10,-2)*{\bu}="B";
(0,-2)*{\bu}="A";
"A"; "A" **\crv{(6,6) & (-6,6)};
"B"; "B" **\crv{(16,6) & (4,6)};
 \ar @{-} "A";"B" <0pt>
\endxy$}}
\Ea
-(-1)^d
 \Ba{c}\resizebox{13.5mm}{!}{
\mbox{$\xy
 (0.5,0.9)*{^{{^{\bar{2}}}}},
(5,0)*{^{{^{\bar{1}}}}},
 (10.5,0.9)*{^{{^{\bar{3}}}}},
 (10,-2)*{\bu}="B";
(0,-2)*{\bu}="A";
"A"; "A" **\crv{(6,6) & (-6,6)};
"B"; "B" **\crv{(16,6) & (4,6)};
 \ar @{-} "A";"B" <0pt>
\endxy$}}
\Ea
\hspace{-1.5mm}\right)
\\
\Ea
$$

We have in $\tw\RGra_d$
$$
\delta_\centerdot \Ba{c}
\mbox{$\xy
 (-3,-0.2)*{^{{^{\bar{1}}}}},
(3,-0.2)*{^{{^{\bar{3}}}}},
(0,3)*{^{{^{\bar{2}}}}},
(0,1)*{\bu}="A1";
(0,1)*{\bu}="A2";
"A1"; "A2" **\crv{(-6,7) & (-6,-5)};
"A1"; "A2" **\crv{(6,7) & (6,-5)};
\endxy$}
\Ea
=
 \Ba{c}\resizebox{13.5mm}{!}{
\mbox{$\xy
 (0.5,0.9)*{^{{^{\bar{1}}}}},
(5,0)*{^{{^{\bar{2}}}}},
 (10.5,0.9)*{^{{^{\bar{3}}}}},
 (10,-2)*{\bu}="B";
(0,-2)*{\bu}="A";
"A"; "A" **\crv{(6,6) & (-6,6)};
"B"; "B" **\crv{(16,6) & (4,6)};
 \ar @{-} "A";"B" <0pt>
\endxy$}}
\Ea
-
 \Ba{c}\resizebox{9mm}{!}{
\xy
(2.0,-3.5)*{^{{^{\bar{2}}}}},
 (-2,-3.5)*{^{{^{\bar{1}}}}},
(0.5,4)*{^{{^{\bar{3}}}}},
 (0,-8)*{\bu}="C";
(0,3)*{\bu}="A1";
(0,3)*{\bu}="A2";
"C"; "A1" **\crv{(-5,-9) & (-5,4)};
"C"; "A2" **\crv{(5,-9) & (5,4)};
 \ar @{-} "A1";"C" <0pt>
\endxy}
\Ea
$$
so that the above map can be re-written exactly in form first found in \cite{Me2} via a completely independent calculation using ribbon graphs only,
$$
i^Q: \Ba{c}\begin{xy}
 <0mm,-1mm>*{\bu};<-4mm,3mm>*{^{_1}}**@{-},
 <0mm,-1mm>*{\bu};<0mm,3mm>*{^{_2}}**@{-},
 <0mm,-1mm>*{\bu};<4mm,3mm>*{^{_3}}**@{-},
 \end{xy}\Ea
 \lon
 \frac{1}{2}\left(-\hspace{-1mm}
 \Ba{c}\resizebox{9mm}{!}{
\xy
(2.0,-3.5)*{^{{^{\bar{2}}}}},
 (-2,-3.5)*{^{{^{\bar{1}}}}},
(0.5,4)*{^{{^{\bar{3}}}}},
 (0,-8)*{\bu}="C";
(0,3)*{\bu}="A1";
(0,3)*{\bu}="A2";
"C"; "A1" **\crv{(-5,-9) & (-5,4)};
"C"; "A2" **\crv{(5,-9) & (5,4)};
 \ar @{-} "A1";"C" <0pt>
\endxy}
\Ea
+\hspace{-1mm}
\Ba{c}\resizebox{9mm}{!}{
\xy
(2.0,-3.5)*{^{{^{\bar{1}}}}},
 (-2,-3.5)*{^{{^{\bar{2}}}}},
(0.5,4)*{^{{^{\bar{3}}}}},
 (0,-8)*{\bu}="C";
(0,3)*{\bu}="A1";
(0,3)*{\bu}="A2";
"C"; "A1" **\crv{(-5,-9) & (-5,4)};
"C"; "A2" **\crv{(5,-9) & (5,4)};
 \ar @{-} "A1";"C" <0pt>
\endxy}
\Ea\right)
$$
This version of the map $i^Q$ was used in \cite{Me2} to show that this map is injective
on infinitely elements of $q\LB_{d-1,d}$ constructing thereby infinitely  many higher genus  cohomology classes in $H^\bu(\cM_{g,m+n})$ from the unique cohomology class in $H^\bu(\cM_{0,3})$ via properadic compositions in $\Grav_d\subset H^\bu(\tw\RGra_d)$.

\subsection{Quasi-Lie bialgebra structures on Hochschild cohomologies of cyclic $A_\infty$-algebras}\label{3: subsec on qLieb on Cyc(A)}  Let $A$ be a graded vector space equipped with a degree $-n$
non-degenerate scalar product
\Beq\label{3: scalar product in A}
\Ba{rccc}
\langle\ ,\ \rangle: & A\odot A & \lon & \K[-n]\\
                     &a\odot b     & \lon & \langle a,b \rangle=(-1)^{|a||b|} \langle a,b \rangle.
\Ea
\Eeq
One has  an associated isomorphism of graded vector spaces,
$$
 A\simeq  A^*[-n]:=\Hom(A,\K)[-n],
$$
and an induced non-degenerate pairing
$$
\Ba{rccc}
\Theta: & \ot^2\left(A[n-1]\right) & \lon & \K[n-2]\equiv \K[1-(3-n)]\\
           & (a'=\fs^{n-1}a, b'= \fs^{n-1} b) & \lon & \fs^{2n-2} \langle a,b \rangle
\Ea
$$
which satisfies the following equation (cf.\ \S 2.3 in \cite{MW1} with $d=3-n$ in the notation of that paper)
\Beqrn
\Theta(b', a') &=& \fs^{2n-2} \langle b,a \rangle \\
&=& (-1)^{|a||b|} \fs^{2n-2} \langle a,b \rangle \\
&=& (-1)^{(|a'|+n-1)(|b'|+n-1)} \Theta(a', b') \\
&=& (-1)^{|a'| + |b'| + (3-n)}\Theta(a', b')
\Eeqrn
where we used the fact that $\Theta(a,b)=0$ unless $|a|+|b|=n$. By Theorem 4.2.2 in \cite{MW1} this symmetry equation
 implies  that the (reduced) space of cyclic word
$$
Cyc(A):=\bigoplus_{p\geq 2} \left(\ot^p(A[n-1])\right)^{\Z_p}\simeq
\bigoplus_{p\geq 2}
\left(\ot^p(A[n-1])\right)^{\Z_p}
$$
carries canonically a representation of the properad $\RGra_{3-n}$ discussed in the previous subsection. In particular this space is a $\LB_{3-n,3-n}$-algebra (see (\ref{3: Lieb_dd to RGra_d})) with the Lie bracket given by a simple formula,
$$
\{(a'_1\ot...\ot b'_k)^{\Z_k}, (b'_1\ot ...\ot b'_l)^{\Z_n}\}:=\hspace{100mm}
$$
$$
\hspace{9mm} \sum_{i=1}^k\sum_{j=1}^l
 \pm
\Theta(a'_i,b'_j) (a'_1\ot ...\ot  a'_{i-1}\ot b'_{j+1}\ot ... \ot b'_l\ot b'_1\ot ... \ot b'_{j-1}\ot a'_{i+1}\ot\ldots\ot a'_k)^{\Z_{k+l-2}}
$$
Maurer Cartan elements elements of this $\Lie_d$-algebra are degree $d=3-n$ elements
 $
 \ga\in Cyc(A)$ such that $\{\ga, \ga\}=0$. There is a one-to-one correspondence
 between such Maurer-Cartan elements and\footnote{One can use this statement as a definition of a degree $n$ cyclic $A_\infty$-algebra structure on $A$.} {\em degree $n$ cyclic strongly homotopy algebra structures in $A$}. The dg Lie algebra
 $$
CH(A):= \left( Cyc(A), d_\ga:=\{\ga,\ \} \right)
 $$
is precisely the (reduced) cyclic Hochschild complex of the cyclic $A_\infty$-algebra $(A,\ga)$.

\sip

By the very definition of the twisting endofunctor $\tw$, the chain gravity properad $\cC h\Grav_{3-n}$
admits a canonical representation in $CH(A)$ for any degree $n$ cyclic $A_\infty$-algebra $A$.
In particular, the gravity properad $\Grav_{3-n}$ acts on its cohomology $H^\bu(CH(A))$ implying, by  Theorem {\ref{2: qLien to twP}}, the following observation.

\subsubsection{\bf Corollary} {\em The Hochschild cohomology $H^\bu(CH(A))$ of any degree $n$ cyclic $A_\infty$ algebra is a quasi-Lie bialgebra; more precisely it carries a representation of the properad $q\LB_{2-n,3-n}$.}

\sip

If $A$ is a dg  Poincare model of some compact $n$-dimensional manifold, then
there is a linear map
$$
\bar{H}_\bu^{S^1}(LM) \lon  H^\bu(CH(A))
$$
from the reduced equivariant homology $\bar{H}_\bu^{S^1}(LM)$ of the free loop space $LM$ of $M$. If $M$ is simply connected, this map is an isomorphism so that the  gravity properad $\Grav_d$ acts
 on $\bar{H}_\bu^{S^1}(LM)$. However the Maurer-Cartan associated to any Poincare model is a relatively simple linear combination of cyclic words in the letters, and the just mentioned action is much trivialized. That ``trivialization" is studied in detail in
 \cite{Me4} where it is shown that the action of $\cC h\Grav_{3-n}$ on  $CH(A)$
factors through a quotient properad $\cS\cT_{3-n}$ which contains $\LB_{d,d}$ \cite{CS}, the gravity operad \cite{Ge, We} and the four $\HoLoB_{d-1}$-operations found in \cite{Me3}.

\sip

One has to consider a less trivial class (comparing to the
class of  Poincare models) of cyclic $A_\infty$-algebras $A$  in order to get a chance to see a less trivial action of the gravity properad  on the associated cyclic Hochschild cohomologies.

\bip

%%%%%%%%%%%%%%%%%%%%%%%%%%%%%%%%%%%%%%%%%%%%%%%%%%%%%%%%%%%%%%

{\Large
\section{\bf A full twisting of properads under $\HoLBcd$}
}

\mip

\subsection{Reminder on the polydifferential functor $\f$}
There is an exact polydifferential functor \cite{MW1}
$$
\Ba{rccc}
\f: & \text{\sf Category of dg props} & \lon & \text{\sf Category of dg operads}\\
&     \cP   &\lon & \f\cP
\Ea
$$
which has the following property:
given any dg prop $\cP$ and its arbitrary representation $\rho: \cP\rar \cE nd_V$
 in a dg vector space $V$, the associated dg operad $\f\cP$ has an associated representation, $\f\rho: \f(\cP)\rar \cE nd_{\odot^\bu V}$,
 in the  graded commutative tensor algebra ${\odot^\bu} V$ given in terms of polydifferential (with respect to the standard multiplication in ${\odot^\bu} V$) operators. Roughly speaking the functor $\f$ symmetrizes all outputs of elements of $\cP$, and splits all inputs into symmetrized blocks; pictorially, if we identify elements of $\cP$ with decorated corollas as in (\ref{3: generic elements of cP as (m,n)-corollas}), then  every element of $\f(\cP)$ can be identified with a decorated corolla
    which is allowed to have the {\em same}\, numerical labels assigned to its different in-legs, and also with the same label 1 assigned to its all outgoing legs\footnote{Treating out- and inputs legs in this procedure on equal footing, one gets a polydifferential functor $\caD$ in the category of dg props such that $\f(\cP)$ a sub-operad of $\caD(\cP)$. It was introduced and studied in \cite{MW3}.}
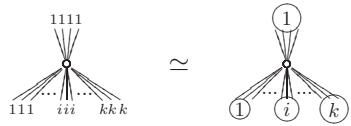
$$
\Ba{c}\resizebox{16mm}{!}{\xy
(-9,-6)*{};
(0,0)*{\circ }
**\dir{-};
(-7.5,-6)*{};
(0,0)*{\circ }
**\dir{-};
(-6,-6)*{};
(0,0)*{\circ }
**\dir{-};
(-1,-6)*{};
(0,0)*{\circ }
**\dir{-};
(0,-6)*{};
(0,0)*{\circ }
**\dir{-};
(1,-6)*{};
(0,0)*{\circ }
**\dir{-};
(9,-6)*{};
(0,0)*{\circ }
**\dir{-};
(7.5,-6)*{};
(0,0)*{\circ }
**\dir{-};
(6,-6)*{};
(0,0)*{\circ }
**\dir{-};
(-3,-5)*{...};
(3,-5)*{...};
(-9,-7.5)*{_1};
(-7.5,-7.5)*{_1};
(-6,-7.5)*{_1};
(-1.1,-7.5)*{_i};
(1.1,-7.5)*{_i};
(0,-7.5)*{_i};
(7.8,-7.5)*{_k};
(6.3,-7.5)*{_k};
(9.6,-7.5)*{_k};
(0,7)*{{ }^{1111}};
(-2,6)*{};
(0,0)*{\circ }
**\dir{-};
(-0.7,6)*{};
(0,0)*{\circ }
**\dir{-};
(0.7,6)*{};
(0,0)*{\circ }
**\dir{-};
(2,6)*{};
(0,0)*{\circ }
**\dir{-};
\endxy}\Ea  \  \ \ {\simeq} \ \ \
\Ba{c}\resizebox{16mm}{!}{\xy
(-9,-6)*{};
(0,0)*{\circ }
**\dir{-};
(-7.5,-6)*{};
(0,0)*{\circ }
**\dir{-};
(-6,-6)*{};
(0,0)*{\circ }
**\dir{-};
(-1,-6)*{};
(0,0)*{\circ }
**\dir{-};
(0,-6)*{};
(0,0)*{\circ }
**\dir{-};
(1,-6)*{};
(0,0)*{\circ }
**\dir{-};
(9,-6)*{};
(0,0)*{\circ }
**\dir{-};
(7.5,-6)*{};
(0,0)*{\circ }
**\dir{-};
(6,-6)*{};
(0,0)*{\circ }
**\dir{-};
(-3,-5)*{...};
(3,-5)*{...};
%(0,9)*{\overbrace{ }^{[m]}};
%
(-2,6)*{};
(0,0)*{\circ }
**\dir{-};
(-0.7,6)*{};
(0,0)*{\circ }
**\dir{-};
(0.7,6)*{};
(0,0)*{\circ }
**\dir{-};
(2,6)*{};
(0,0)*{\circ }
**\dir{-};
(-8.2,-7.9)*+\hbox{${{1}}$}*\frm{o};
(8.2,-7.9)*+\hbox{${{\, k\, }}$}*\frm{o};
(0,-7.9)*+\hbox{${{\, i\, }}$}*\frm{o};
(0,8)*+\hbox{${{\, 1\, }}$}*\frm{o};
\endxy}\Ea
$$

\sip

  Since we want to apply the above construction to dg props $\cP$ under $\HoLBcd$, we are more interested in this paper in its degree shifted version, $\f_{c,d}$, which was also introduced in \cite{MW1},
 $$
 \f_{c,d}\cP:= \f(\cP\{c\})
 $$
 The notation may be slightly misleading as $\f_{c,d}$ does not depend on $d$ but it suits us well in the context of this paper.
 %We have an associated morphism of operads
%$$
%\f_{c,d}i: \f_{c,d}\HoLBcd \rar \f_{c,d}\cP.
%$$
We refer to \cite{MW2} for more details about the functor $\f_{c,d}$ (and to \cite{MW3} for its extension $\caD$) and discuss next the particular example, the dg operad
 $$
 \f_{c,d}\HoLBcd\simeq \f\HoLB_{0,c+d}\equiv\{\f\HoLB_{0,c+d}(k)\}_{k\geq 1}.
 $$
The $\bS_k$-module $\f\HoLB_{0,c+d}(k)$, $k\geq 1$, is generated by graphs $\ga$ constructed from  arbitrary decorated graphs $\Ga$  from $\HoLB_{0,c+d}(m,n)$, $\forall m,n\geq 1$,
as follows:
\Bi
\item[(i)]
draw new $k$ big white vertices labelled from $1$ to $k$  (these will be inputs of $\ga$) and one extra output big white vertex,
\item[(ii)] symmetrize all $m$ outputs legs of $\Ga$ and attach them to the unique output white vertex;
\item[(iii)] partition the set $[n]$ of input legs of $\Ga$ into $k$ ordered disjoint  (not necessary non-empty) subsets
    $$
    [n]=I_1\sqcup \ldots \sqcup I_k, \ \ \ \ \#I_i\geq 0, i\in [k],
    $$
 then symmetrize the legs in each subset $I_i$ and attach them (if any) to the $i$-labelled input white vertex.
\Ei

For example, the element
$$
\Ga=\Ba{c}\resizebox{11mm}{!}{
\xy
(0,0)*{\bu}="o",
(-2,5)*{}="2",
(4,5)*{\bu}="3",
(4,10)*{}="u",
(4,0)*{}="d1",
(7,0)*{}="d2",
(10,0)*{}="d3",
(-1.5,-5)*{}="5",
(1.5,-5)*{}="6",
(4,-5)*{}="7",
(-4,-5)*{}="8",
(-2,7)*{_1},
(4,12)*{_2},
(-1.5,-7)*{_2},
(1.5,-7)*{_3},
(10.4,-1.6)*{_6},
(-4,-7)*{_1},
(4,-1.6)*{_4},
(7,-1.6)*{_5},
\ar @{-} "o";"2" <0pt>
\ar @{-} "o";"3" <0pt>
\ar @{-} "o";"5" <0pt>
\ar @{-} "o";"6" <0pt>
%\ar @{-} "o";"7" <0pt>
\ar @{-} "o";"8" <0pt>
\ar @{-} "3";"u" <0pt>
\ar @{-} "3";"d1" <0pt>
\ar @{-} "3";"d2" <0pt>
\ar @{-} "3";"d3" <0pt>
\endxy}\Ea
 \in \HoLB_{0,c+d}(2,6)
$$
can produce the following generator
$$
\ga=\Ba{c}\resizebox{15mm}{!}{ \xy
(-1.5,5)*{}="1",
(1.5,5)*{}="2",
(9,5)*{}="3",
 (0,0)*{\bu}="A";
  (9,3)*{\bu}="O";
   (5,12)*+{\hspace{2mm}}*\frm{o}="X";
 (-6,-10)*+{_1}*\frm{o}="B";
  (6,-10)*+{_2}*\frm{o}="C";
   (14,-10)*+{_3}*\frm{o}="D";
    (22,-10)*+{_4}*\frm{o}="E";
 "A"; "B" **\crv{(-5,-0)}; % (5,5) = tangent point
 %"A"; "B" **\crv{(-5,-6)};
  "A"; "D" **\crv{(5,-0.5)};
   %"A"; "B" **\crv{(5,-1)};
  "A"; "C" **\crv{(-5,-7)};
   "A"; "O" **\crv{(5,5)};
\ar @{-} "O";"C" <0pt>
\ar @{-} "O";"D" <0pt>
\ar @{-} "O";"X" <0pt>
\ar @{-} "A";"X" <0pt>
\ar @{-} "O";"B" <0pt>
 \endxy}
 \Ea \in \f\HoLB_{0,c+d}(4)\simeq\f_{c,d}\HoLBcd(4)
$$
in the associated polydifferential operad (note that one and the same element $\Ga\in \HoLB_{0,c+d}$ can give rise to several different generators of $\f\HoLB_{0,c+d}$). The labelled white vertices of elements of $\f(\HoLB_{0,c+d}$ are called {\em external}, while unlabelled black vertices (more, precisely, the vertices of the underlying elements of $\HoLB_{0,c+d}$) are called {\em internal}. The same terminology can applied to $\f\cP$ for any dg prop $\cP$.

\sip

 For any $k,l\geq 1$ and $i\in [k]$ the operadic composition
 $$
 \Ba{rccc}
 \circ_i: &\f\HoLB_{0,c+d}(k)\ot \f\HoLB_{0,c+d}(k) & \lon &
 \f\HoLB_{0,c+d}(k+l-1)\\
 & \Ga_1\ot \Ga_2 &\lon & \Ga_1\circ_i \Ga_2
 \Ea
$$
is defined by
\Bi
\item[(i)]
 substituting the graph $\Ga_2$ (with the output external vertex erased so that all  edges connected to that external vertex are hanging at this step loosely) inside the big circle of the $i$-labelled external vertex of $\Ga_1$,
 \item[(ii)] erasing that big $i$-th labelled external circle (so that all edges of $\Ga_1$ connected to that $i$-th external vertex, if any, are also hanging loosely),   and
 \item[(iii)] finally taking the sum over all possible ways to do the following three operations in any order,
\Bi
\item[(a)] glue some (or all or none) hanging edges of $\Ga_2$ to the same number of hanging edges of $\Ga_1$,
\item[(b)] attach some (or all or none) hanging edges of $\Ga_2$ to the output external vertex of $\Ga_1$,
\item[(c)] attach some (or all or none) hanging edges of $\Ga_1$ to the external input vertices of $\Ga_2$,
\Ei
in such a way that no hanging edges are left.
\Ei
 We refer to \cite{MW1,MW3} for concrete examples of such compositions.

\subsubsection{\bf Proposition \cite{MW1}}\label{4: Prop on map from Lie^+ to OHoLB}
{\em There is a  morphism of dg operads
\Beq\label{3: map Holie^+_{c+d} to f_{c,d}(HoLBcd)}
\Holie^+_{c+d} \to \f_{c,d}\HoLBcd
\Eeq
given explicitly on the $(1,1)$-generator by
$$
\begin{xy}
 <0mm,-0.55mm>*{};<0mm,-3mm>*{}**@{-},
 <0mm,0.5mm>*{};<0mm,3mm>*{}**@{-},
 <0mm,0mm>*{\bu};<0mm,0mm>*{}**@{},
 \end{xy}
 \lon
\sum_{m\geq 2} \Ba{c}\resizebox{12mm}{!}{\begin{xy}
 <0mm,0mm>*{\bu};<0mm,0mm>*{}**@{},
 <-0.6mm,0.44mm>*{};<-8mm,5mm>*{^1}**@{-},
 <-0.4mm,0.7mm>*{};<-4.5mm,5mm>*{^1}**@{-},
 <0mm,5mm>*{\ldots},
 <0.4mm,0.7mm>*{};<4.5mm,5mm>*{^1}**@{-},
 <0.6mm,0.44mm>*{};<8mm,5mm>*{^1}**@{-},
     <0mm,9mm>*{\overbrace{\ \ \ \ \ \ \ \ \ \ \ \ \ \ \ \ }},
     <0mm,11mm>*{^m},
 <0.0mm,-0.44mm>*{};<0mm,-5mm>*{}**@{-},
 \end{xy}}\Ea
$$
and on the remaining $(1,n)$-generators with $n\geq 2$ by}
$$
\Ba{c}\resizebox{20mm}{!}{ \xy
(1,-5)*{\ldots},
(-13,-7)*{_1},
(-8,-7)*{_2},
(-3,-7)*{_3},
%(7,-7)*{_{n-1}},
(13,-7)*{_n},
 (0,0)*{\bu}="a",
(0,5)*{}="0",
(-12,-5)*{}="b_1",
(-8,-5)*{}="b_2",
(-3,-5)*{}="b_3",
(8,-5)*{}="b_4",
(12,-5)*{}="b_5",
\ar @{-} "a";"0" <0pt>
\ar @{-} "a";"b_2" <0pt>
\ar @{-} "a";"b_3" <0pt>
\ar @{-} "a";"b_1" <0pt>
\ar @{-} "a";"b_4" <0pt>
\ar @{-} "a";"b_5" <0pt>
\endxy}\Ea
\lon
\sum_{m\geq 1}
\Ba{c}\resizebox{12mm}{!}{  \xy
(0,8)*{^1}="1";
(-2,8)*{^1}="2";
(2,8)*{^1}="3";
(2,-4)*{\ldots};
    (0,+3)*{\bu}="L";
 (-8,-5)*+{_1}*\frm{o}="B";
  (-3,-5)*+{_2}*\frm{o}="C";
   (8,-5)*+{_n}*\frm{o}="D";
    <0mm,12mm>*{\overbrace{ \ \ \ \ \ \ \ \ }},
     <0mm,14.6mm>*{_m},
\ar @{-} "D";"L" <0pt>
\ar @{-} "C";"L" <0pt>
\ar @{-} "B";"L" <0pt>
\ar @{-} "1";"L" <0pt>
\ar @{-} "2";"L" <0pt>
\ar @{-} "3";"L" <0pt>
 \endxy}
 \Ea
$$

Proof is a straightforward calculation (cf.\ \S 5.5 and \S 5.7 in \cite{MW1}).

\sip

%{\tiny
%There is a similar map of dg operads for $c+d\in 2\Z$,
%$$
%\Holie_{c+d}^{\diamond +} \lon \f_{c,d} (\HoLoBcd)
%$$
%}

Using T.\ Willwacher twisting endofuctor discussed in \S 2 one obtains via the morphism
(\ref{3: map Holie^+_{c+d} to f_{c,d}(HoLBcd)}) a dg operad
$
\tw\, \f_{c,d} \HoLBcd %\ \ \ \ Tw\left(\f_{c,d} (\HoLoBcd) \right),
$.
%This construction can be applied straightforwardly to the following class of dg props.

\subsubsection{\bf Properads under $\HoLBcd$} Assume $\cP$ is a dg properad {\em under}\, $\HoLBcd$, i.e.\ the one which comes equipped with a non-trivial morphism
\Beq\label{4: i from HoLB to P}
\Ba{rccc}
i: & \HoLBcd & \lon & \cP \vspace{0mm} \\
&
\Ba{c}\resizebox{10mm}{!}{ \xy
(0,4.5)*+{...},
(0,-4.5)*+{...},
(0,0)*{\bu}="o",
(-5,5)*{}="1",
(-3,5)*{}="2",
(3,5)*{}="3",
(5,5)*{}="4",
(-3,-5)*{}="5",
(3,-5)*{}="6",
(5,-5)*{}="7",
(-5,-5)*{}="8",
(-5.5,7)*{_1},
(-3,7)*{_2},
(3,6)*{},
(5.9,7)*{m},
(-3,-7)*{_2},
(3,-7)*+{},
(5.9,-7)*{n},
(-5.5,-7)*{_1},
\ar @{-} "o";"1" <0pt>
\ar @{-} "o";"2" <0pt>
\ar @{-} "o";"3" <0pt>
\ar @{-} "o";"4" <0pt>
\ar @{-} "o";"5" <0pt>
\ar @{-} "o";"6" <0pt>
\ar @{-} "o";"7" <0pt>
\ar @{-} "o";"8" <0pt>
\endxy}\Ea
 &\lon &
 \Ba{c}\resizebox{10mm}{!}{ \xy
(0,4.5)*+{...},
(0,-4.5)*+{...},
(0,0)*{\circledcirc}="o",
(-5,5)*{}="1",
(-3,5)*{}="2",
(3,5)*{}="3",
(5,5)*{}="4",
(-3,-5)*{}="5",
(3,-5)*{}="6",
(5,-5)*{}="7",
(-5,-5)*{}="8",
(-5.5,7)*{_1},
(-3,7)*{_2},
(3,6)*{},
(5.9,7)*{m},
(-3,-7)*{_2},
(3,-7)*+{},
(5.9,-7)*{n},
(-5.5,-7)*{_1},
\ar @{-} "o";"1" <0pt>
\ar @{-} "o";"2" <0pt>
\ar @{-} "o";"3" <0pt>
\ar @{-} "o";"4" <0pt>
\ar @{-} "o";"5" <0pt>
\ar @{-} "o";"6" <0pt>
\ar @{-} "o";"7" <0pt>
\ar @{-} "o";"8" <0pt>
\endxy}\Ea
\Ea
\Eeq

 Note that corollas on the right hand side (the ones with $\circledcirc$ as the vertex) stand from now on for images of the generators of $\HoLBcd$ under the map $i$ so that {\it some (or all) of them can in fact be equal to zero}.

 \sip

Applying the functor $\f_{c,d}$ and using the above proposition we obtain an associated chain of morphism of dg operads,
$$
\iota: \Holie_{c+d}^+ \lon \f_{c,d} \HoLBcd \lon \f_{c,d}\cP
$$
and hence a morphism of the associated twisted dg operads,
$$
\tw\f(i): \tw\f_{c,d} \HoLBcd\simeq \tw\f(\HoLB_{0,c+d}) \lon \tw\f_{c,d}\cP
$$
%Note that
%$$
%\f_{c,d} (\HoLBcd)\equiv\f(\HoLBcd\{c\})\simeq  \f(\HoLB_{0,c+d}), \ \ \ %\f_{c,d}\cP=\f(\cP\{c\})
%$$
The degree of the generating $(m,n)$-corolla $\Ba{c}\resizebox{7mm}{!}{ \xy
(0,4.5)*+{...},
(0,-4.5)*+{...},
(0,0)*{\bu}="o",
(-5,5)*{}="1",
(-3,5)*{}="2",
(3,5)*{}="3",
(5,5)*{}="4",
(-3,-5)*{}="5",
(3,-5)*{}="6",
(5,-5)*{}="7",
(-5,-5)*{}="8",
\ar @{-} "o";"1" <0pt>
\ar @{-} "o";"2" <0pt>
\ar @{-} "o";"3" <0pt>
\ar @{-} "o";"4" <0pt>
\ar @{-} "o";"5" <0pt>
\ar @{-} "o";"6" <0pt>
\ar @{-} "o";"7" <0pt>
\ar @{-} "o";"8" <0pt>
\endxy}\Ea$ of $\HoLB_{0,c+d}$ is equal to $1+c+d- (c+d)n$, so its $m$ out-legs are carry trivial representation of $\bS_m$ (and are assigned degree $0$), while its in-legs are (skew)symmetrized according to the parity of $c+d\in \Z$ (and assigned degree $(c+d)$; the vertex is assigned the degree $1+c+d$.
{\em Hence it is only the sum $c+d$ of our integer parameters which plays a role in this story}.
Therefore we can assume without loss of generality that
$$
c=0,\ \ \ d\ \text{is an arbitrary integer},
$$
from now on, i.e. work solely with dg props under $\HoLB_{0,d}$.

\subsection{Maurer-Cartan elements of strongly homotopy Lie bialgebras}\label{4: subsec on  MC elements of HoLBcd} Given a $\HoLB_{0,d}$-algebra structure in a dg vector space $(V,\delta)$, i.e.\ a morphism of properads
$$
\rho: \HoLB_{0,d} \lon \cE nd_V.
$$
 Its {\it Maurer-Cartan element}\, is, by definition, a Maurer-Cartan element $\ga\in \odot^{\geq 1}V$ of the associated $\Holie_d^+$ structure induced on $\odot^{\geq 1}V$
via the canonical monomorphism
\[
\Holie^+_{d} \to \f\HoLB_{0,d}
\]
described explicitly in Proposition {\ref{4: Prop on map from Lie^+ to OHoLB}}. Let us describe it in more detail. %(assuming for simplicity of signs that $d$ is even, say $d=0$).
The $\HoLB_{0,d}$-structure on $V$ is given by a collection of  linear maps of cohomological degree $1+d-dn$,
$$
\rho\left(\Ba{c}\resizebox{10mm}{!}{ \xy
(0,4.5)*+{...},
(0,-4.5)*+{...},
(0,0)*{\bu}="o",
(-5,5)*{}="1",
(-3,5)*{}="2",
(3,5)*{}="3",
(5,5)*{}="4",
(-3,-5)*{}="5",
(3,-5)*{}="6",
(5,-5)*{}="7",
(-5,-5)*{}="8",
(-5.5,7)*{_1},
(-3,7)*{_2},
(3,6)*{},
(5.9,7)*{m},
(-3,-7)*{_2},
(3,-7)*+{},
(5.9,-7)*{n},
(-5.5,-7)*{_1},
\ar @{-} "o";"1" <0pt>
\ar @{-} "o";"2" <0pt>
\ar @{-} "o";"3" <0pt>
\ar @{-} "o";"4" <0pt>
\ar @{-} "o";"5" <0pt>
\ar @{-} "o";"6" <0pt>
\ar @{-} "o";"7" <0pt>
\ar @{-} "o";"8" <0pt>
\endxy}\Ea
\right)=:\mu_{m,n}: \otimes^n V \lon \odot^m V
$$
satisfying compatibility conditions. Each such linear map gives rise to a map
$$
\hat{\mu}_{m,n}: \ot^n (\odot^{\geq 1}V) \lon \odot^{\geq 1}V
$$
given, in
%which extends uniquely each action on the tensor factor $V$ to $\odot^{\geq 1}V$ as a %derivation with respect to the natural graded symmetric product on $\odot^{\geq 1}V$. Assume %
arbitrary basis $\{p_\al\}$  of $V$ as follows
%$$
%\mu_{m,n}(p_{\al_1}\ot  p_{\al_2} \ot \ldots\ot  p_{\al_n})=: \sum \frac{1}{m!} p_{\be_1} %p_{\be_2}\ldots p_{\be_m}\, \mu_{\al_1,\al_2,\ldots, \al_n}^{\be_1\be_2\ldots \be_m} ,
%$$
% then its extension $\hat{\mu}$ is given as a polydifferential operator,
 $$
 %\hat{\mu}_{m,n}(f_1,\ldots, f_n):= \sum \pm \frac{1}{m!}p_{\be_1} p_{\be_2}\ldots %p_{\be_m}\, \mu_{\al_1,\al_2,\ldots, \al_n}^{\be_1\be_2\ldots \be_m}
 % \frac{\p f_1}{\p p_{\al_1}} \frac{\p f_2}{\p p_{\al_2}}\ldots \frac{\p f_n}{\p p_{\al_n}},\ %
 \hat{\mu}_{m,n}(f_1,\ldots, f_n):= \sum \pm \mu_{m,n}(p_{\al_1}\ot p_{\al_2}\ot \ldots \ot p_{\al_n})\cdot
  \frac{\p f_1}{\p p_{\al_1}} \frac{\p f_2}{\p p_{\al_2}}\cdots \frac{\p f_n}{\p p_{\al_n}},\ %
  \ \ \forall f_1,\ldots, f_n\in \odot^{\geq 1}V.
 $$
  Then a degree $d$ element $\ga \in \odot^{\geq 1}V$ is a Maurer-Cartan element of the (appropriately filtered or nilpotent)  $\HoLB_{0,d}$-algebra on structure $V$  if and only if the following equation holds,
\Beq\label{4: MC eqn for HoLB0d}
\delta\ga+ \sum_{m,n\geq 1}\frac{1}{n!}\hat{\mu}_{m,n}(\underbrace{\ga,\ldots,\ga}_n)=0.
\Eeq
The operator
\Beq\label{4: twisted by MC differentiail in V}
\delta_\ga v :=\delta v + \sum_{n\geq 1}\frac{1}{n!}{\mu}_{1,n+1}(\underbrace{\ga_1,\ldots,\ga_1}_n, v ),\ \ \ \ \forall\, v\in V,
\Eeq
with $\ga_1$ being the image of $\ga$ under the projection $\odot^{\geq 1}V \rar V$,
is a twisted differential on $V$.

\subsubsection{\bf Combinatorial incarnation}\label{4: subsec on comb incarn of HoLB MC elements}
Maurer-Cartan elements of $\HoLB_{0,d}$-algebras admit a simple combinatorial description as (a representation in $V$ of) an infinite linear combination
$$
\ga \simeq \sum_{m\geq 1}\frac{1}{m!} \Ba{c}\resizebox{13mm}{!}{\begin{xy}
 <0mm,0mm>*{\bu};<0mm,0mm>*{}**@{},
 <-0.6mm,0.44mm>*{};<-8mm,5mm>*{}**@{-},
 <-0.4mm,0.7mm>*{};<-4.5mm,5mm>*{}**@{-},
 <0mm,5mm>*{\ldots},
 <0.4mm,0.7mm>*{};<4.5mm,5mm>*{}**@{-},
 <0.6mm,0.44mm>*{};<8mm,5mm>*{}**@{-},
     <0mm,7mm>*{\overbrace{\ \ \ \ \ \ \ \ \ \ \ \ \ \ \ \ }},
     <0mm,9mm>*{^m},
 %<0.0mm,-0.44mm>*{};<0mm,-5mm>*{}**@{-},
 \end{xy}}\Ea
$$
of degree $d$ $(m,0)$-corollas with symmetrized outgoing legs. One extends the standard differential $\delta$ in $\HoLB_{0,d}$ to such new generating corollas as follows
 %(cf.\ (\ref{2: delta on (1,0) MC corolla})),
\Beq\label{4: d on HoLB MC}
\delta\ \Ba{c}\resizebox{13mm}{!}{\begin{xy}
 <0mm,0mm>*{\bu};<0mm,0mm>*{}**@{},
 <-0.6mm,0.44mm>*{};<-8mm,5mm>*{}**@{-},
 <-0.4mm,0.7mm>*{};<-4.5mm,5mm>*{}**@{-},
 <0mm,5mm>*{\ldots},
 <0.4mm,0.7mm>*{};<4.5mm,5mm>*{}**@{-},
 <0.6mm,0.44mm>*{};<8mm,5mm>*{}**@{-},
     <0mm,7mm>*{\overbrace{\ \ \ \ \ \ \ \ \ \ \ \ \ \ \ \ }},
     <0mm,9mm>*{^m},
 %<0.0mm,-0.44mm>*{};<0mm,-5mm>*{}**@{-},
 \end{xy}}\Ea
 =
- \sum_{k\geq 1, [m]=\sqcup  [m_\bu]\atop {  m_{0}\geq 1, k+m_0\geq 3\atop m_1,...,m_k\geq 0}}\frac{1}{k!}
 \Ba{c}\resizebox{22mm}{!}{  \xy
(-27,8)*{}="1";
(-25,8)*{}="2";
(-23,8)*{}="3";
%(2,-4)*{\ldots};
%
(-18,8)*{}="n11";
(-15,7)*{...};
(-13,8)*{}="n12";
(-11,8)*{}="n21";
(-8.4,7)*{...};
(-6,8)*{}="n22";
(1,8)*{}="nn1";
(3,7)*{...};
(5,8)*{}="nn2";
    (-25,+3)*{\bu}="L";
 (-15,-5)*{\bu}="B";
  (-8,-5)*{\bu}="C";
   (3,-5)*{\bu}="D";
   <-5mm,-10mm>*{\underbrace{ \ \ \ \ \ \ \ \ \ \ \ \ \ \ \ \ \ }_{k}},
   (-3,-5)*{...};
    <-25mm,10.6mm>*{\overbrace{ \ \ \ \ \ \ }^{m_0}},
    <-16mm,10.6mm>*{\overbrace{ \ \ \ \ }^{m_1}},
    <-9mm,10.6mm>*{\overbrace{ \ \ \ \ }^{m_2}},
    <3mm,10.6mm>*{\overbrace{ \ \ \ \ }^{m_n}},
     %<0mm,12mm>*{_m},
\ar @{-} "D";"L" <0pt>
\ar @{-} "C";"L" <0pt>
\ar @{-} "B";"L" <0pt>
\ar @{-} "1";"L" <0pt>
\ar @{-} "2";"L" <0pt>
\ar @{-} "3";"L" <0pt>
\ar @{-} "1";"L" <0pt>
\ar @{-} "2";"L" <0pt>
\ar @{-} "n11";"B" <0pt>
\ar @{-} "n12";"B" <0pt>
\ar @{-} "n21";"C" <0pt>
\ar @{-} "n22";"C" <0pt>
\ar @{-} "nn1";"D" <0pt>
\ar @{-} "nn2";"D" <0pt>
 \endxy}
 \Ea \ \ \ \  \forall\ m\geq 1.
\Eeq
where we take the sum over all partitions of the ordered set $[m]$ into $k+1$ ordered subsets. For $k=1$ we recover the standard formula (cf.\ (\ref{2: delta on (1,0) MC corolla})).

%\sip

% The black square corolla with $m$ legs on the l.h.s.\  has  total degree $c+d$. Then the %r.h.s. has the total degree
%$$
%1+c+d - (c+d)n + n(c+d)= 1+c+d
%$$
%as expected.

\sip

Let us check first that the above definition makes sense.

\subsubsection{\bf Lemma}
$\delta^2\ \Ba{c}\resizebox{13mm}{!}{\begin{xy}
 <0mm,0mm>*{\bu};<0mm,0mm>*{}**@{},
 <-0.6mm,0.44mm>*{};<-8mm,5mm>*{}**@{-},
 <-0.4mm,0.7mm>*{};<-4.5mm,5mm>*{}**@{-},
 <0mm,5mm>*{\ldots},
 <0.4mm,0.7mm>*{};<4.5mm,5mm>*{}**@{-},
 <0.6mm,0.44mm>*{};<8mm,5mm>*{}**@{-},
     <0mm,7mm>*{\overbrace{\ \ \ \ \ \ \ \ \ \ \ \ \ \ \ \ }},
     <0mm,9mm>*{^m},
 %<0.0mm,-0.44mm>*{};<0mm,-5mm>*{}**@{-},
 \end{xy}}\Ea\equiv 0$ \ {\it for any $m\geq 1$}.
\begin{proof} The proof is based on a straightforward calculation which is similar to the one made  in the proof of Lemma {\ref{2: Lemma on d^2 for MC Lie}} above. The only really new phenomenon is the appearance in $\delta^2$ of summands of the form
\Beq\label{4: terms of the for lastochka}
\sum_{[m]= [m'_\bu] \sqcup [m''_\bu]\sqcup [m_0] }\ \ \frac{1}{k'!} \frac{1}{k''!}
 \Ba{c}\resizebox{55mm}{!}{  \xy
(8,8)*{}="1";
(10,8)*{}="2";
(12,8)*{}="3";
%(2,-4)*{\ldots};
%
(17,8)*{}="n11";
(20,7)*{...};
(22,8)*{}="n12";
(24,8)*{}="n21";
(26.4,7)*{...};
(29,8)*{}="n22";
(36,8)*{}="nn1";
(38,7)*{...};
(40,8)*{}="nn2";
(-8,8)*{}="1'";
(-10,8)*{}="2'";
(-12,8)*{}="3'";
(-17,8)*{}="n11'";
(-20,7)*{...};
(-22,8)*{}="n12'";
(-24,8)*{}="n21'";
(-26.4,7)*{...};
(-29,8)*{}="n22'";
(-36,8)*{}="nn1'";
(-38,7)*{...};
(-40,8)*{}="nn2'";
    (10,+3)*{\bu}="R";
 (20,-5)*{\bu}="B1";
  (27,-5)*{\bu}="C1";
   (38,-5)*{\bu}="D1";
   <30mm,-10mm>*{\underbrace{ \ \ \ \ \ \ \ \ \ \ \ \ \ \ \ \ \ }_{k''}},
   (32,-5)*{...};
    <10mm,10.6mm>*{\overbrace{ \ \ \ \ \ \ }^{m''_0}},
    <19mm,10.6mm>*{\overbrace{ \ \ \ \ }^{m''_1}},
    <26mm,10.6mm>*{\overbrace{ \ \ \ \ }^{m''_2}},
    <38mm,10.6mm>*{\overbrace{ \ \ \ \ }^{m''_n}},
 (0,-5)*{\bu}="0";
 (-3,8)*{}="s1";
(-0,7)*{...};
(3,8)*{}="s2";
<0mm,10.6mm>*{\overbrace{ \ \ \ \ \ \ }^{m_0}},
(-10,+3)*{\bu}="L";
 (-20,-5)*{\bu}="B1'";
  (-27,-5)*{\bu}="C1'";
   (-38,-5)*{\bu}="D1'";
   <-30mm,-10mm>*{\underbrace{ \ \ \ \ \ \ \ \ \ \ \ \ \ \ \ \ \ }_{k'}},
   (-32,-5)*{...};
    <-10mm,10.6mm>*{\overbrace{ \ \ \ \ \ \ }^{m'_0}},
    <-19mm,10.6mm>*{\overbrace{ \ \ \ \ }^{m'_1}},
    <-26mm,10.6mm>*{\overbrace{ \ \ \ \ }^{m'_2}},
    <-38mm,10.6mm>*{\overbrace{ \ \ \ \ }^{m'_n}},
%
     %<0mm,12mm>*{_m},
\ar @{-} "D1";"R" <0pt>
\ar @{-} "C1";"R" <0pt>
\ar @{-} "B1";"R" <0pt>
\ar @{-} "1";"R" <0pt>
\ar @{-} "2";"R" <0pt>
\ar @{-} "3";"R" <0pt>
\ar @{-} "1";"R" <0pt>
\ar @{-} "2";"R" <0pt>
\ar @{-} "n11";"B1" <0pt>
\ar @{-} "n12";"B1" <0pt>
\ar @{-} "n21";"C1" <0pt>
\ar @{-} "n22";"C1" <0pt>
\ar @{-} "nn1";"D1" <0pt>
\ar @{-} "nn2";"D1" <0pt>
 \ar @{-} "0";"L" <0pt>
 \ar @{-} "0";"R" <0pt>
  \ar @{-} "0";"s1" <0pt>
   \ar @{-} "0";"s2" <0pt>
 \ar @{-} "D1'";"L" <0pt>
\ar @{-} "C1'";"L" <0pt>
\ar @{-} "B1'";"L" <0pt>
\ar @{-} "1'";"L" <0pt>
\ar @{-} "2'";"L" <0pt>
\ar @{-} "3'";"L" <0pt>
\ar @{-} "1'";"L" <0pt>
\ar @{-} "2'";"L" <0pt>
\ar @{-} "n11'";"B1'" <0pt>
\ar @{-} "n12'";"B1'" <0pt>
\ar @{-} "n21'";"C1'" <0pt>
\ar @{-} "n22'";"C1'" <0pt>
\ar @{-} "nn1'";"D1'" <0pt>
\ar @{-} "nn2'";"D1'" <0pt>
 %
% \ar @{-} "0";"L" <0pt>
% \ar @{-} "0";"R" <0pt>
 \endxy}
 \Ea
\Eeq
which cancel each other for symmetry reasons.
\end{proof}

\subsection{Full twisting of properads under $\HoLB_{c,d}$}\label{4: Subsec on Def of Tw of Tw(HoLB) and Tw(P)}
Let $\cP$ be a properad under $\HoLB_{0,d}$ (as in (\ref{4: i from HoLB to P})). We construct
 the associated {\em fully twisted}\, properad $(\Tw \cP,\pc)$ in  several steps.

\sip

First we define  $\wTW\cP$ %=\{\wTW\cP)(m,n), \sd\}_{m\geq 1, n\geq 0}$
to be
 be the properad generated freely by $\cP$ and a family of new $(m,0)$-corollas, $\Ba{c}\resizebox{11mm}{!}{\begin{xy}
 <0mm,0mm>*{\bu};<0mm,0mm>*{}**@{},
 <-0.6mm,0.44mm>*{};<-8mm,5mm>*{}**@{-},
 <-0.4mm,0.7mm>*{};<-4.5mm,5mm>*{}**@{-},
 <0mm,5mm>*{\ldots},
 <0.4mm,0.7mm>*{};<4.5mm,5mm>*{}**@{-},
 <0.6mm,0.44mm>*{};<8mm,5mm>*{}**@{-},
     <0mm,7mm>*{\overbrace{\ \ \ \ \ \ \ \ \ \ \ \ \ \ \ \ }},
     <0mm,9mm>*{^m} \end{xy}}\Ea$, $m\geq 1$, of cohomological degree $d$ which are called {\it MC generators}.
     %; moreover, we assume that {\em $\wTW\cP$ is completed with the respect to the %filtration by the number of these new MC vertices}.
     We make this properad differential by using
  the original differential $\sd$ on elements of $\cP$ and extending its action on the new generators  by (cf.\ (\ref{4: d on HoLB MC}))
\Beq\label{3: d on MC in TW(P)}
\sd\ \Ba{c}\resizebox{13mm}{!}{\begin{xy}
 <0mm,0mm>*{\bu};<0mm,0mm>*{}**@{},
 <-0.6mm,0.44mm>*{};<-8mm,5mm>*{}**@{-},
 <-0.4mm,0.7mm>*{};<-4.5mm,5mm>*{}**@{-},
 <0mm,5mm>*{\ldots},
 <0.4mm,0.7mm>*{};<4.5mm,5mm>*{}**@{-},
 <0.6mm,0.44mm>*{};<8mm,5mm>*{}**@{-},
     <0mm,7mm>*{\overbrace{\ \ \ \ \ \ \ \ \ \ \ \ \ \ \ \ }},
     <0mm,9mm>*{^m},
 %<0.0mm,-0.44mm>*{};<0mm,-5mm>*{}**@{-},
 \end{xy}}\Ea
 =
- \sum_{k\geq 1, [m]=\sqcup  [m_\bu]\atop {  m_{0}\geq 1, k+m_0\geq 3\atop m_1,...,m_k\geq 0}}\frac{1}{k!}
 \Ba{c}\resizebox{22mm}{!}{  \xy
(-27,8)*{}="1";
(-25,8)*{}="2";
(-23,8)*{}="3";
%(2,-4)*{\ldots};
%
(-18,8)*{}="n11";
(-15,7)*{...};
(-13,8)*{}="n12";
(-11,8)*{}="n21";
(-8.4,7)*{...};
(-6,8)*{}="n22";
(1,8)*{}="nn1";
(3,7)*{...};
(5,8)*{}="nn2";
    (-25,+3)*{\circledcirc}="L";
 (-14,-5)*{\bu}="B";
  (-8,-5)*{\bu}="C";
   (3,-5)*{\bu}="D";
   <-5mm,-9mm>*{\underbrace{ \ \ \ \ \ \ \ \ \ \ \ \ \ \ \ \ \ }_{k}},
   (-3,-5)*{...};
    <-25mm,10.6mm>*{\overbrace{ \ \ \ \ \ \ }^{m_0}},
    <-16mm,10.6mm>*{\overbrace{ \ \ \ \ }^{m_1}},
    <-9mm,10.6mm>*{\overbrace{ \ \ \ \ }^{m_2}},
    <3mm,10.6mm>*{\overbrace{ \ \ \ \ }^{m_n}},
     %<0mm,12mm>*{_m},
\ar @{-} "D";"L" <0pt>
\ar @{-} "C";"L" <0pt>
\ar @{-} "B";"L" <0pt>
\ar @{-} "1";"L" <0pt>
\ar @{-} "2";"L" <0pt>
\ar @{-} "3";"L" <0pt>
\ar @{-} "1";"L" <0pt>
\ar @{-} "2";"L" <0pt>
\ar @{-} "n11";"B" <0pt>
\ar @{-} "n12";"B" <0pt>
\ar @{-} "n21";"C" <0pt>
\ar @{-} "n22";"C" <0pt>
\ar @{-} "nn1";"D" <0pt>
\ar @{-} "nn2";"D" <0pt>
 \endxy}
 \Ea \ \ \ \  \forall\ m\geq 1.
\Eeq
Note that corollas with $\circledcirc$-vertices are images of the geneartors of $\HoLB_{0,d}$
in $\cP$  under the morphism $i$ (and hence some of them can, in principle, be zero). The map (\ref{4: i from HoLB to P}) extends to a morphism
$$
\wTW(i): \wTW\HoLB_{0,d} \lon \wTW\cP
$$
which restricts on the MC generators as the identoty map.

\sip

Next we notice the following surprising result which tells us essentially the MC elements
originating in $\f(\HoLBcd)$ can be used to twist not only the operad $\f(\HoLBcd)$ (which is obvious) but also the
properad $\HoLBcd$ itself!

\subsubsection{\bf Theorem}\label{3: Th on HoLB^+ to Tw(HoLB)}  {\em
There is a canonical  monomorphism of dg props
$$
\HoLB_{0,d}^+ \lon \wTW\HoLB_{0,d}
$$
given on the generating $(m,n)$-corollas with $m,n\geq 1$ as follows,
$$
\Ba{c}\resizebox{12mm}{!}{\begin{xy}
 <0mm,0mm>*{\bu};<0mm,0mm>*{}**@{},
 <-0.6mm,0.44mm>*{};<-8mm,5mm>*{}**@{-},
 <-0.4mm,0.7mm>*{};<-4.5mm,5mm>*{}**@{-},
 <0mm,0mm>*{};<-1mm,5mm>*{\ldots}**@{},
 <0.4mm,0.7mm>*{};<4.5mm,5mm>*{}**@{-},
 <0.6mm,0.44mm>*{};<8mm,5mm>*{}**@{-},
   <0mm,0mm>*{};<-8.5mm,5.5mm>*{^1}**@{},
   <0mm,0mm>*{};<-5mm,5.5mm>*{^2}**@{},
   <0mm,0mm>*{};<4.5mm,5.5mm>*{^{m\hspace{-0.5mm}-\hspace{-0.5mm}1}}**@{},
   <0mm,0mm>*{};<9.0mm,5.5mm>*{^m}**@{},
 <-0.6mm,-0.44
 mm>*{};<-8mm,-5mm>*{}**@{-},
 <-0.4mm,-0.7mm>*{};<-4.5mm,-5mm>*{}**@{-},
 <0mm,0mm>*{};<-1mm,-5mm>*{\ldots}**@{},
 <0.4mm,-0.7mm>*{};<4.5mm,-5mm>*{}**@{-},
 <0.6mm,-0.44mm>*{};<8mm,-5mm>*{}**@{-},
   <0mm,0mm>*{};<-8.5mm,-6.9mm>*{^1}**@{},
   <0mm,0mm>*{};<-5mm,-6.9mm>*{^2}**@{},
   <0mm,0mm>*{};<4.5mm,-6.9mm>*{^{n\hspace{-0.5mm}-\hspace{-0.5mm}1}}**@{},
   <0mm,0mm>*{};<9.0mm,-6.9mm>*{^n}**@{},
 \end{xy}}\Ea
\lon
\sum_{k\geq 0, [m]=\sqcup  [m_\bu],\atop { m_{0}\geq 1, k+m_0+n\geq 3\atop m_1,...,m_k\geq 0}}
\frac{1}{k!}
% \sum_{m=p+m_1+...+m_n, n,p\geq 1, m_{\bu}\geq 0 \atop n+p\geq 3}
\Ba{c}\resizebox{23mm}{!}{  \xy
(-27,8)*{}="1";
(-25,8)*{}="2";
(-23,8)*{}="3";
%(2,-4)*{\ldots};
%
(-18,8)*{}="n11";
(-15,7)*{...};
(-13,8)*{}="n12";
(-11,8)*{}="n21";
(-8.4,7)*{...};
(-6,8)*{}="n22";
(1,8)*{}="nn1";
(3,7)*{...};
(5,8)*{}="nn2";
(-29,-8)*+{_1}="l2";
    (-27,-8)*+{_2}="l1";
 (-24,-6)*{...};
   (-20,-8)*+{_n}="ln";
(-3,-5)*{...};
    (-25,+3)*{\bu}="L";
 (-14,-5)*{\bu}="B";
  (-8,-5)*{\bu}="C";
   (3,-5)*{\bu}="D";
   %<0mm,-10mm>*{\underbrace{ \ \ \ \ \ \ \ \ \ \ \ \ \ \ \ \ \ }_{n}},
 %
    <-25mm,10.6mm>*{\overbrace{ \ \ \ \ \ \ }^{m_0}},
    <-16mm,10.6mm>*{\overbrace{ \ \ \ \ }^{m_1}},
    <-9mm,10.6mm>*{\overbrace{ \ \ \ \ }^{m_2}},
    <3mm,10.6mm>*{\overbrace{ \ \ \ \ }^{m_k}},
     <-5mm,-10mm>*{\underbrace{ \ \ \ \ \ \ \ \ \ \ \ \ \ \ \ \ \ }_{k}},
\ar @{-} "D";"L" <0pt>
\ar @{-} "C";"L" <0pt>
\ar @{-} "B";"L" <0pt>
\ar @{-} "1";"L" <0pt>
\ar @{-} "2";"L" <0pt>
\ar @{-} "3";"L" <0pt>
\ar @{-} "1";"L" <0pt>
\ar @{-} "2";"L" <0pt>
\ar @{-} "n11";"B" <0pt>
\ar @{-} "n12";"B" <0pt>
\ar @{-} "n21";"C" <0pt>
\ar @{-} "n22";"C" <0pt>
\ar @{-} "nn1";"D" <0pt>
\ar @{-} "nn2";"D" <0pt>
 \ar @{-} "l1";"L" <0pt>
\ar @{-} "l2";"L" <0pt>
\ar @{-} "ln";"L" <0pt>
 \endxy}
 \Ea.
$$
}
\vspace{-2mm}
\begin{proof} One has to check that the above explicit map of properads respects their differentials, $\delta^+$ on the l.h.s.\ and $\p$ on the r.h.s. This is a straightforward calculation  which is analogous to (but much more tedious than)  the one used in the proof of Theorem {\ref{2: map from HoLB+ to tilda twA}}. The only really new aspect is again the appearance of terms (\ref{4: terms of the for lastochka}) which cancel out for symmetry reasons. We omit full details.
\end{proof}

Hence for any properad $\cP$ under $\HoLB_{0,d}$ there is an associated morphism
of dg properads

$$
i^+: \HoLB_{0,d}^+ \lon \wTW\cP
$$
which factors through the morphism described in the Theorem just above.

%\subsubsection{\bf Proposition} {\em There is a isomorphism of dg props,}
%$$
%\f(TW(\HoLB_{0,d}))= TW(\f(\HoLBcd\{c\}).
%$$

%\subsubsection{\bf Open problem} What is the cohomology of $Tw(\HoLBcd\{c\})$? There is a %similar question for $Tw(\Holie_d)$... which is probably not interesting.

\subsubsection{\bf Twisting of  the differential} The argument in \S {\ref{2: subsec on twisting d in operads}} about twisting of differentials in operads  extends straightforwardly to properads. Indeed,  given any dg prop $(\cP=\{\cP(m,n)\},\sd)$,
and any  $h\in \cA(1,1)$, there is an associated derivation $D_h$ of the (non-differential) prop $\cP$ by the formula analogous to (\ref{2: delta^+}),
\Beq\label{3: formula for D_h}
D_h(a)= \sum_{i=1}^m h _1\circ_i a - (-1)^{|h||a|} \sum_{j=1}^n a _j\circ_1 h, \ \ \forall \ a\in \cP(m,n).
\Eeq
Moreover, if $|h|=1$ and
$h\circ_1 h=-\sd h$, then the sum
$
\p_h:= \p+ D_h
$
is a differential in $\cP$.

\sip

 Assume we have a morphism of dg props
\Beq\label{3: g^+ map from HoLB}
g^+: (\HoLB_{0,d}^+, \delta^+) \lon (\cP,\sd)
\Eeq
Then the element
\Beq\label{3: blacklozenge (1,1) element}
\begin{xy}
 <0mm,-3mm>*{};<0mm,3mm>*{}**@{-},
 %<0mm,0.5mm>*{};<0mm,3mm>*{}**@{-},
 <0mm,0mm>*{_\blacklozenge};<0mm,0mm>*{}**@{},
 \end{xy}
 :=g^+(\begin{xy}
 <0mm,-0.55mm>*{};<0mm,-3mm>*{}**@{-},
 <0mm,0.5mm>*{};<0mm,3mm>*{}**@{-},
 <0mm,0mm>*{\bullet};<0mm,0mm>*{}**@{},
 \end{xy})=
 \sum_{k=1}^\infty \frac{1}{k!}\ \resizebox{17mm}{!}{  \xy
(-25,8)*{}="1";
(-9,-4)*{...};
 <-11mm,-8mm>*{\underbrace{ \ \ \ \ \ \ \ \ \ \ \ \ \ \ \ \ \ }_{k}},
    (-25,+3)*{\circledast}="L";
    (-25,-3)*{}="N";
 (-19,-4)*{\bu}="B";
  (-13,-4)*{\bu}="C";
   (-2,-4)*{\bu}="D";
\ar @{-} "D";"L" <0pt>
\ar @{-} "C";"L" <0pt>
\ar @{-} "B";"L" <0pt>
\ar @{-} "1";"L" <0pt>
\ar @{-} "N";"L" <0pt>
%
 %
% \ar @{-} "l1";"N" <0pt>
%  \ar @{-} "l2";"N" <0pt>
%   \ar @{-} "ln";"N" <0pt>
 \endxy} \in \wt{\Tw}\cP
\Eeq
 satisfies all the conditions specified above for $h$ so that the sum
 $$
 \sd_{\centerdot}:= \p + D_{\hspace{-1.8mm}\Ba{c}\resizebox{1.8mm}{!}{\begin{xy}
 <0mm,-0.55mm>*{};<0mm,-3mm>*{}**@{-},
 <0mm,0.5mm>*{};<0mm,3mm>*{}**@{-},
 <0mm,0mm>*{\blacklozenge};<0mm,0mm>*{}**@{},
 \end{xy}}\Ea}
 $$
is a differential in $\wt{\Tw}\cP$.

\subsubsection{\bf Main definition} {\em Let $\cP$ be a dg properad under $\HoLB_{0,d}$. The full twisting, $\Tw\cP$, of $\cP$ is a dg properad  defined as the properad $\wt{\Tw}\cP$ equipped with the twisted differential $\pc$.}

\sip

%{\bf (i)} Given any dg prop(erad) $\cP=\{\cP(m,n)\}_{m,n\geq 0}$ under $\HoLB_{0,d}$. The %associated prop(erad) $(\wTW\cP,\sd)$ under $\HoLB_{0,d}^+$. Twisting the differential in the %latter from $\sd$ into $\sd_\centerdot$
% we obtain a new dg prop(erad)
%which is denoted from now on by
%by $\Tw\cP=\{\Tw\cP(m,n)\}_{m\geq 1,n\geq 0}$ and is called the {\em (full) twisting of $%\cP$}.

\sip

Thus
$\Tw\cP$ is identical to $\wt{\Tw}\cP$ as a non-differential properad, i.e.\ it is generated freely  by $\cP$ and the family of extra generators $\Ba{c}\resizebox{11mm}{!}{\begin{xy}
 <0mm,0mm>*{\bu};
 <-0.6mm,0.44mm>*{};<-8mm,5mm>*{^1}**@{-},
 <-0.4mm,0.7mm>*{};<-4.5mm,5mm>*{^2}**@{-},
 <0mm,5mm>*{\ldots},
 <0.4mm,0.7mm>*{};<4.5mm,5mm>*{^{}}**@{-},
 <0.6mm,0.44mm>*{};<8mm,5mm>*{^m}**@{-},
 \end{xy}}\Ea$, $m\geq 1$, of degree $d$.
If we represent elements of $\cP$ as decorated corollas (\ref{3: generic elements of cP as (m,n)-corollas}), then the twisted differential $\sd_\centerdot$ acts on elements of $\cP$  as follows,
\Beq\label{4: d_centerdot on Tw(P)}
\sd_\centerdot \Ba{c}\resizebox{13mm}{!}{
 \begin{xy}
 <0mm,0mm>*{\circ};<-8mm,6mm>*{^1}**@{-},
 <0mm,0mm>*{\circ};<-4.5mm,6mm>*{^2}**@{-},
 <0mm,0mm>*{\circ};<0mm,5.5mm>*{\ldots}**@{},
 <0mm,0mm>*{\circ};<3.5mm,5mm>*{}**@{-},
 <0mm,0mm>*{\circ};<8mm,6mm>*{^m}**@{-},
 <0mm,0mm>*{\circ};<-8mm,-6mm>*{_1}**@{-},
 <0mm,0mm>*{\circ};<-4.5mm,-6mm>*{_2}**@{-},
 <0mm,0mm>*{\circ};<0mm,-5.5mm>*{\ldots}**@{},
 <0mm,0mm>*{\circ};<4.5mm,-6mm>*+{}**@{-},
 <0mm,0mm>*{\circ};<8mm,-6mm>*{_n}**@{-},
   \end{xy}}\Ea
=
\sd \Ba{c}\resizebox{13mm}{!}{
 \begin{xy}
 <0mm,0mm>*{\circ};<-8mm,6mm>*{^1}**@{-},
 <0mm,0mm>*{\circ};<-4.5mm,6mm>*{^2}**@{-},
 <0mm,0mm>*{\circ};<0mm,5.5mm>*{\ldots}**@{},
 <0mm,0mm>*{\circ};<3.5mm,5mm>*{}**@{-},
 <0mm,0mm>*{\circ};<8mm,6mm>*{^m}**@{-},
 <0mm,0mm>*{\circ};<-8mm,-6mm>*{_1}**@{-},
 <0mm,0mm>*{\circ};<-4.5mm,-6mm>*{_2}**@{-},
 <0mm,0mm>*{\circ};<0mm,-5.5mm>*{\ldots}**@{},
 <0mm,0mm>*{\circ};<4.5mm,-6mm>*+{}**@{-},
 <0mm,0mm>*{\circ};<8mm,-6mm>*{_n}**@{-},
   \end{xy}}\Ea
+
\overset{m-1}{\underset{i=0}{\sum}}
\Ba{c}\resizebox{14mm}{!}{
\begin{xy}
 %<0mm,0mm>*{\circ};<0mm,0mm>*{}**@{},
 <0mm,0mm>*{\circ};<-8mm,5mm>*{}**@{-},
 <0mm,0mm>*{\circ};<-3.5mm,5mm>*{}**@{-},
 <0mm,0mm>*{\circ};<-6mm,5mm>*{..}**@{},
 <0mm,0mm>*{\circ};<0mm,5mm>*{}**@{-},
  <0mm,5mm>*{\blacklozenge};
  <0mm,5mm>*{};<0mm,8mm>*{}**@{-},
  <0mm,5mm>*{};<0mm,9mm>*{^{i\hspace{-0.2mm}+\hspace{-0.5mm}1}}**@{},
<0mm,0mm>*{\circ};<8mm,5mm>*{}**@{-},
<0mm,0mm>*{\circ};<3.5mm,5mm>*{}**@{-},
<6mm,5mm>*{..}**@{},
<-8.5mm,5.5mm>*{^1}**@{},
<-4mm,5.5mm>*{^i}**@{},
<9.0mm,5.5mm>*{^m}**@{},
 <0mm,0mm>*{\circ};<-8mm,-5mm>*{}**@{-},
 <0mm,0mm>*{\circ};<-4.5mm,-5mm>*{}**@{-},
 <-1mm,-5mm>*{\ldots}**@{},
 <0mm,0mm>*{\circ};<4.5mm,-5mm>*{}**@{-},
 <0mm,0mm>*{\circ};<8mm,-5mm>*{}**@{-},
<-8.5mm,-6.9mm>*{^1}**@{},
<-5mm,-6.9mm>*{^2}**@{},
<4.5mm,-6.9mm>*{^{n\hspace{-0.5mm}-\hspace{-0.5mm}1}}**@{},
<9.0mm,-6.9mm>*{^n}**@{},
 \end{xy}}\Ea
 - (-1)^{|a|}
\overset{n-1}{\underset{i=0}{\sum}}
 \Ba{c}\resizebox{14mm}{!}{\begin{xy}
 %<0mm,0mm>*{\circ};
 <0mm,0mm>*{\circ};<-8mm,-5mm>*{}**@{-},
 <0mm,0mm>*{\circ};<-3.5mm,-5mm>*{}**@{-},
 <0mm,0mm>*{\circ};<-6mm,-5mm>*{..}**@{},
 <0mm,0mm>*{\circ};<0mm,-5mm>*{}**@{-},
  <0mm,-5mm>*{\blacklozenge};
  <0mm,-5mm>*{};<0mm,-8mm>*{}**@{-},
  <0mm,-5mm>*{};<0mm,-10mm>*{^{i\hspace{-0.2mm}+\hspace{-0.5mm}1}}**@{},
<0mm,0mm>*{\circ};<8mm,-5mm>*{}**@{-},
<0mm,0mm>*{\circ};<3.5mm,-5mm>*{}**@{-},
 <6mm,-5mm>*{..}**@{},
<-8.5mm,-6.9mm>*{^1}**@{},
<-4mm,-6.9mm>*{^i}**@{},
<9.0mm,-6.9mm>*{^n}**@{},
 <0mm,0mm>*{\circ};<-8mm,5mm>*{}**@{-},
 <0mm,0mm>*{\circ};<-4.5mm,5mm>*{}**@{-},
<-1mm,5mm>*{\ldots}**@{},
 <0mm,0mm>*{\circ};<4.5mm,5mm>*{}**@{-},
 <0mm,0mm>*{\circ};<8mm,5mm>*{}**@{-},
<-8.5mm,5.5mm>*{^1}**@{},
<-5mm,5.5mm>*{^2}**@{},
<4.5mm,5.5mm>*{^{m\hspace{-0.5mm}-\hspace{-0.5mm}1}}**@{},
<9.0mm,5.5mm>*{^m}**@{},
 \end{xy}}\Ea
\Eeq
where $\Ba{c}\resizebox{1.8mm}{!}{\begin{xy}
 <0mm,3mm>*{};<0mm,-3mm>*{}**@{-},
 <0mm,0.5mm>*{};<0mm,3mm>*{}**@{-},
 <0mm,0mm>*{\blacklozenge};<0mm,0mm>*{}**@{},
 \end{xy}}\Ea
$ is given by (\ref{3: blacklozenge (1,1) element}).
On the other hand, the action of $\sd_\centerdot$ on the MC generators is given by,
\Beq\label{4: d_centerdot on MC generators of Tw(P)}
 \sd_\centerdot
 \Ba{c}\resizebox{14mm}{!}{
 \begin{xy}
 <0mm,-1mm>*{\bu};
 <0mm,0mm>*{};<-8mm,5mm>*{}**@{-},
 <0mm,0mm>*{};<-4.5mm,5mm>*{}**@{-},
 <0mm,0mm>*{};<-1mm,5mm>*{\ldots}**@{},
 <0mm,0mm>*{};<4.5mm,5mm>*{}**@{-},
 <0mm,0mm>*{};<8mm,5mm>*{}**@{-},
   <0mm,0mm>*{};<-8.5mm,5.5mm>*{^1}**@{},
   <0mm,0mm>*{};<-5mm,5.5mm>*{^2}**@{},
   %<0mm,0mm>*{};<4.5mm,5.5mm>*{^{m\hspace{-0.5mm}-\hspace{-0.5mm}1}}**@{},
   <0mm,0mm>*{};<9.0mm,5.5mm>*{^m}**@{},
 \end{xy}}\Ea:=
\overset{m-1}{\underset{i=0}{\sum}}
\Ba{c}\resizebox{14mm}{!}{
\begin{xy}
 <0mm,-1mm>*{\bu};
 <0mm,0mm>*{};<-8mm,5mm>*{}**@{-},
 <0mm,0mm>*{};<-3.5mm,5mm>*{}**@{-},
 <0mm,0mm>*{};<-6mm,5mm>*{..}**@{},
 <0mm,0mm>*{};<0mm,5mm>*{}**@{-},
  <0mm,5mm>*{\blacklozenge};
  <0mm,5mm>*{};<0mm,8mm>*{}**@{-},
  <0mm,5mm>*{};<0mm,9mm>*{^{i\hspace{-0.2mm}+\hspace{-0.5mm}1}}**@{},
<0mm,0mm>*{};<8mm,5mm>*{}**@{-},
<0mm,0mm>*{};<3.5mm,5mm>*{}**@{-},
 <0mm,0mm>*{};<6mm,5mm>*{..}**@{},
   <0mm,0mm>*{};<-8.5mm,5.5mm>*{^1}**@{},
   <0mm,0mm>*{};<-4mm,5.5mm>*{^i}**@{},
   <0mm,0mm>*{};<9.0mm,5.5mm>*{^m}**@{},
 \end{xy}}\Ea
- \sum_{k\geq 1, [m]=\sqcup  [m_\bu],\atop {  m_{0}\geq 1, k+m_0\geq 3\atop m_1,...,m_k\geq 0}}\frac{1}{k!}
 \Ba{c}\resizebox{21mm}{!}{  \xy
(-27,8)*{}="1";
(-25,8)*{}="2";
(-23,8)*{}="3";
%(2,-4)*{\ldots};
%
(-18,8)*{}="n11";
(-15,7)*{...};
(-13,8)*{}="n12";
(-11,8)*{}="n21";
(-8.4,7)*{...};
(-6,8)*{}="n22";
(1,8)*{}="nn1";
(3,7)*{...};
(5,8)*{}="nn2";
    (-25,+2)*{\circledcirc}="L";
 (-14,-5)*{\bu}="B";
  (-8,-5)*{\bu}="C";
   (3,-5)*{\bu}="D";
   <-5mm,-10mm>*{\underbrace{ \ \ \ \ \ \ \ \ \ \ \ \ \ \ \ \ \ }_{k}},
   (-3,-5)*{...};
    <-25mm,10.6mm>*{\overbrace{ \ \ \ \ \ \ }^{m_0}},
    <-16mm,10.6mm>*{\overbrace{ \ \ \ \ }^{m_1}},
    <-9mm,10.6mm>*{\overbrace{ \ \ \ \ }^{m_2}},
    <3mm,10.6mm>*{\overbrace{ \ \ \ \ }^{m_n}},
     %<0mm,12mm>*{_m},
\ar @{-} "D";"L" <0pt>
\ar @{-} "C";"L" <0pt>
\ar @{-} "B";"L" <0pt>
\ar @{-} "1";"L" <0pt>
\ar @{-} "2";"L" <0pt>
\ar @{-} "3";"L" <0pt>
\ar @{-} "1";"L" <0pt>
\ar @{-} "2";"L" <0pt>
\ar @{-} "n11";"B" <0pt>
\ar @{-} "n12";"B" <0pt>
\ar @{-} "n21";"C" <0pt>
\ar @{-} "n22";"C" <0pt>
\ar @{-} "nn1";"D" <0pt>
\ar @{-} "nn2";"D" <0pt>
 \endxy}
 \Ea \ \ \ \  \forall\ m\geq 1.
\Eeq
Note that for $m\geq 2$ the first sum on the r.h.s.\ of (\ref{4: d_centerdot on MC generators of Tw(P)}) cancels out with all the summands corresponding to $k\geq 2 ,m_0=1,m_i=m-1, i\in [k]$, in the second sum.

\sip

%It is sometimes useful to consider separately the part $\{\Tw\cP(m,0)\}_{m\geq 1}$ of the %full twisted prop $\Tw\cP=\{TW(\cP)(m,n)\}_{m,n\geq 0}$  which is spanned by elements of %$\cP$ with {\em all}\, inputs filled with the MC  generators.

%\sip

%{\bf (ii)} We define $Tw(\cP):=\{TW(\cP)(m,n)\}_{m,n\geq 1}$ as the subprop of $TW(\cP)$ spanned %by elements with at least one input not filled by an $MC$-generator, and equipped with the %twisted differential $\sd_\centerdot$.
% \sip

 By its very construction, this twisted prop $\Tw\cP$ has the following properties:

\Bi
\item[(a)] There is a canonical chain of morphisms of dg prop(erad)s
 $$
 (\HoLB_{0,d},\delta)  \lon (\Tw\HoLB_{0,d}, \delta_\centerdot) \lon (\Tw\cP,\sd_\centerdot).
 $$
 given explicitly by
 \Beq\label{4: map frpm HoLBcd to TwP}
\Ba{c}\resizebox{13mm}{!}{\begin{xy}
 <0mm,0mm>*{\bu};<0mm,0mm>*{}**@{},
 <-0.6mm,0.44mm>*{};<-8mm,5mm>*{}**@{-},
 <-0.4mm,0.7mm>*{};<-4.5mm,5mm>*{}**@{-},
 <0mm,0mm>*{};<-1mm,5mm>*{\ldots}**@{},
 <0.4mm,0.7mm>*{};<4.5mm,5mm>*{}**@{-},
 <0.6mm,0.44mm>*{};<8mm,5mm>*{}**@{-},
   <0mm,0mm>*{};<-8.5mm,5.5mm>*{^1}**@{},
   <0mm,0mm>*{};<-5mm,5.5mm>*{^2}**@{},
   <0mm,0mm>*{};<4.5mm,5.5mm>*{^{m\hspace{-0.5mm}-\hspace{-0.5mm}1}}**@{},
   <0mm,0mm>*{};<9.0mm,5.5mm>*{^m}**@{},
 <-0.6mm,-0.44
 mm>*{};<-8mm,-5mm>*{}**@{-},
 <-0.4mm,-0.7mm>*{};<-4.5mm,-5mm>*{}**@{-},
 <0mm,0mm>*{};<-1mm,-5mm>*{\ldots}**@{},
 <0.4mm,-0.7mm>*{};<4.5mm,-5mm>*{}**@{-},
 <0.6mm,-0.44mm>*{};<8mm,-5mm>*{}**@{-},
   <0mm,0mm>*{};<-8.5mm,-6.9mm>*{^1}**@{},
   <0mm,0mm>*{};<-5mm,-6.9mm>*{^2}**@{},
   <0mm,0mm>*{};<4.5mm,-6.9mm>*{^{n\hspace{-0.5mm}-\hspace{-0.5mm}1}}**@{},
   <0mm,0mm>*{};<9.0mm,-6.9mm>*{^n}**@{},
 \end{xy}}\Ea
\lon
\sum_{k\geq 0,  [m]=\sqcup [m_\bu]\atop { m_{0}\geq 1, k+m_0+n\geq 3\atop m_1,...,m_k\geq 0}}
\frac{1}{k!}
% \sum_{m=p+m_1+...+m_n, n,p\geq 1, m_{\bu}\geq 0 \atop n+p\geq 3}
\Ba{c}\resizebox{21mm}{!}{  \xy
(-27,8)*{}="1";
(-25,8)*{}="2";
(-23,8)*{}="3";
%(2,-4)*{\ldots};
%
(-18,8)*{}="n11";
(-15,7)*{...};
(-13,8)*{}="n12";
(-11,8)*{}="n21";
(-8.4,7)*{...};
(-6,8)*{}="n22";
(1,8)*{}="nn1";
(3,7)*{...};
(5,8)*{}="nn2";
(-29,-8)*+{_1}="l2";
    (-27,-8)*+{_2}="l1";
 (-24,-6)*{...};
   (-20,-8)*+{_n}="ln";
(-3,-5)*{...};
    (-25,+3)*{\circledcirc}="L";
 (-14,-5)*{\bu}="B";
  (-8,-5)*{\bu}="C";
   (3,-5)*{\bu}="D";
   %<0mm,-10mm>*{\underbrace{ \ \ \ \ \ \ \ \ \ \ \ \ \ \ \ \ \ }_{n}},
 %
    <-25mm,10.6mm>*{\overbrace{ \ \ \ \ \ \ }^{m_0}},
    <-16mm,10.6mm>*{\overbrace{ \ \ \ \ }^{m_1}},
    <-9mm,10.6mm>*{\overbrace{ \ \ \ \ }^{m_2}},
    <3mm,10.6mm>*{\overbrace{ \ \ \ \ }^{m_k}},
     <-5mm,-10mm>*{\underbrace{ \ \ \ \ \ \ \ \ \ \ \ \ \ \ \ \ \ }_{k}},
\ar @{-} "D";"L" <0pt>
\ar @{-} "C";"L" <0pt>
\ar @{-} "B";"L" <0pt>
\ar @{-} "1";"L" <0pt>
\ar @{-} "2";"L" <0pt>
\ar @{-} "3";"L" <0pt>
\ar @{-} "1";"L" <0pt>
\ar @{-} "2";"L" <0pt>
\ar @{-} "n11";"B" <0pt>
\ar @{-} "n12";"B" <0pt>
\ar @{-} "n21";"C" <0pt>
\ar @{-} "n22";"C" <0pt>
\ar @{-} "nn1";"D" <0pt>
\ar @{-} "nn2";"D" <0pt>
 \ar @{-} "l1";"L" <0pt>
\ar @{-} "l2";"L" <0pt>
\ar @{-} "ln";"L" <0pt>
 \endxy}
 \Ea,\ \ \ \ \ \ m,n\geq 1, m+n\geq 3
\Eeq

%\item[(b)] There is an isomorphism of dg operads (commutativity of the endofunctor $Tw$ %with the functor $\f$), {\color{blue} THIS IS NOT TRUE: the differential on the l.h.s. %is $\sd_\centerdot \simeq (1,1)$ while on the r.h.s. $\sum_m (m,1)+\bu$ }
%$$
%\f (Tw(\cP))= Tw(\f(\cP))
%$$

\item[(b)] There is a canonical epimorphism of of dg prop(erad)s,
$$
(\Tw\cP, \sd_\centerdot) \lon (\cP, \sd)
$$
which sends all the MC generators $\Ba{c}\resizebox{11mm}{!}{\begin{xy}
 <0mm,0mm>*{\bu};<0mm,0mm>*{}**@{},
 <-0.6mm,0.44mm>*{};<-8mm,5mm>*{}**@{-},
 <-0.4mm,0.7mm>*{};<-4.5mm,5mm>*{}**@{-},
 <0mm,5mm>*{\ldots},
 <0.4mm,0.7mm>*{};<4.5mm,5mm>*{}**@{-},
 <0.6mm,0.44mm>*{};<8mm,5mm>*{}**@{-},
     <0mm,7mm>*{\overbrace{\ \ \ \ \ \ \ \ \ \ \ \ \ \ \ \ }},
     <0mm,9mm>*{^m} \end{xy}}\Ea$, $m\geq 1$, to zero. %It is well-defined as the ideal generators all such $(m,0)$-corollas is {\em differential}.
     Note that the natural inclusion of $\bS$-bimodules $\cP\rar \Tw\cP$ is not, in general, a morphism of {\em dg}\, properads.
\Ei

\subsubsection{\bf Remark} Assume $\cP$ is a properad under $\LB_{0,d}$, i.e.\
all corollas on the r.h.s.\ of the map (\ref{4: i from HoLB to P}) vanish except the following two, $\Ba{c}\resizebox{7mm}{!}{  \xy
(-5,-6)*+{}="1";
    (-5,-0.2)*{\circledcirc}="L";
  (-8,5)*+{_1}="C";
   (-2,5)*+{_2}="D";
\ar @{-} "D";"L" <0pt>
\ar @{-} "C";"L" <0pt>
\ar @{-} "1";"L" <0pt>
 \endxy}
 \Ea$ and
$ \Ba{c}\resizebox{7mm}{!}{  \xy
(-5,6)*{}="1";
    (-5,+1)*{\circledcirc}="L";
  (-8,-5)*+{_1}="C";
   (-2,-5)*+{_2}="D";
\ar @{-} "D";"L" <0pt>
\ar @{-} "C";"L" <0pt>
\ar @{-} "1";"L" <0pt>
 \endxy}
 \Ea$. Then the associated twisted properad $\Tw\cP$ is, in general, a properad under the minimal resolution $\HoLB_{0,d}$ of $\LB_{0,d}$, not just under $\LB_{0,d}$! Put another way, {\em the full twisting of $\cP$ produces, in general,  higher homotopy Lie bialgebras operations},
% $$
% \Tw: \PROP_{\LBcd} \lon \PROP_{\HoLBcd},
% $$
  a new phenomenon comparing to what we get in $\tw\cP$ under the partial twisting of $\cP$.

\subsubsection{\bf Full twisting for general values of the integer parameters $c$ and $d$}
 {\em The full twisting $\Tw\cP$ of a dg properad $\cP$ under $\HoLBcd$}\, is defined as $(\Tw\cP\{c\})\{-c\}$; note that $\cP\{c\}$ is a dg properad under $\HoLB_{0,c+d}$ so that the above twisting functor $\Tw$ applies.
Thus the full twisting $\Tw\cP$ of $\cP\in \PROP_{\HoLBcd}$ is generated freely by $\cP$ and extra MC generators
$$
\Ba{c}\resizebox{14mm}{!}{\begin{xy}
 <0mm,0mm>*{\bu};
 <-0.6mm,0.44mm>*{};<-8mm,5mm>*{^1}**@{-},
 <-0.4mm,0.7mm>*{};<-4.5mm,5mm>*{^2}**@{-},
 <0mm,5mm>*{\ldots},
 <0.4mm,0.7mm>*{};<4.5mm,5mm>*{^{}}**@{-},
 <0.6mm,0.44mm>*{};<8mm,5mm>*{^m}**@{-},
 \end{xy}}\Ea =
 (-1)^{c|\sigma|}
 \Ba{c}\resizebox{17mm}{!}{\begin{xy}
 <0mm,0mm>*{\bu};
 <-0.6mm,0.44mm>*{};<-11mm,7mm>*{^{\sigma(1)}}**@{-},
 <-0.4mm,0.7mm>*{};<-4.5mm,7mm>*{^{\sigma(2)}}**@{-},
 <0mm,5mm>*{\ldots},
 <0.4mm,0.7mm>*{};<4.5mm,6mm>*{}**@{-},
 <0.6mm,0.44mm>*{};<10mm,7mm>*{^{\sigma(m)}}**@{-},
 \end{xy}}\Ea \ \ \ \forall \sigma\in \bS_m, \ m\geq 1,
$$
of cohomological degree  $(1-m)c+d$. It comes equipped with a canonical morphism
$
\HoLBcd \rar \Tw\cP
$
given by (\ref{4: map frpm HoLBcd to TwP}).

% In this case one has a chain of maps
%$$
%\Tw(i): (\HoLBcd,\delta)  \lon (\Tw\HoLBcd, \delta_\centerdot) \lon %(\Tw\cP,\sd_\centerdot).
% $$
%given explicitly on the generating corollas of $\HoLBcd$ by

%and the isomorphism
%$$
%\f_{c,d} (Tw(\cP))= Tw(\f_{c,d}(\cP)).
%$$
%QUESTION: Does $Tw(\cP)$ remember the differential in $\cP$? WE have only a quotient map %to $\cP$...

%\sip
%{\bf (iii)} The part $(\{TW(\cP)(m,0)\}_{m\geq 0}, \sd_\centerdot$ with no inputs is a %collection of complexes (which is also a module over $Tw(\cP)$ which we denote by %$Tw_\bu(\cP)$. In applications it is often useful to study both parts, $Tw(\cP)$ and %$Tw_\bu(\cP)$, of $(TW(\cP), \sd_\centerdot$ separately (see \S 4 for concrete examples).

\sip

Next we show that quasi-isomorphisms (\ref{2: quasi-iso twHoLB to HoLB, twLB to LB}) extend to their full twisting analogues.
\subsubsection{\bf Theorem}\label{4: theorem on H(TwHoLBcd)}  {\em The canonical projection $\pi: \Tw\HoLBcd \rar \HoLBcd$ is a quasi-isomorphism, i.e.\
$H^\bu(\Tw\HoLBcd)=\LBcd$.}
\begin{proof}  For any $m\geq 2$ the r.h.s.\ of formula (\ref{4: d_centerdot on MC generators of Tw(P)}) applied to $\cP=\HoLBcd$ contains a unique summand of the form
$$
 \sd_\centerdot
 \Ba{c}\resizebox{14mm}{!}{
 \begin{xy}
 <0mm,-1mm>*{\bu};
 <0mm,0mm>*{};<-8mm,5mm>*{}**@{-},
 <0mm,0mm>*{};<-4.5mm,5mm>*{}**@{-},
 <0mm,0mm>*{};<-1mm,5mm>*{\ldots}**@{},
 <0mm,0mm>*{};<4.5mm,5mm>*{}**@{-},
 <0mm,0mm>*{};<8mm,5mm>*{}**@{-},
   <0mm,0mm>*{};<-8.5mm,5.5mm>*{^1}**@{},
   <0mm,0mm>*{};<-5mm,5.5mm>*{^2}**@{},
   %<0mm,0mm>*{};<4.5mm,5.5mm>*{^{m\hspace{-0.5mm}-\hspace{-0.5mm}1}}**@{},
   <0mm,0mm>*{};<9.0mm,5.5mm>*{^m}**@{},
 \end{xy}}\Ea:=
-
\Ba{c}\resizebox{14mm}{!}{
 \begin{xy}
 <0mm,-5mm>*{\bu};<0mm,0mm>*{\bu}**@{-},
 <0mm,-0mm>*{\bu};
 <0mm,0mm>*{};<-8mm,5mm>*{}**@{-},
 <0mm,0mm>*{};<-4.5mm,5mm>*{}**@{-},
 <0mm,0mm>*{};<-1mm,5mm>*{\ldots}**@{},
 <0mm,0mm>*{};<4.5mm,5mm>*{}**@{-},
 <0mm,0mm>*{};<8mm,5mm>*{}**@{-},
   <0mm,0mm>*{};<-8.5mm,5.5mm>*{^1}**@{},
   <0mm,0mm>*{};<-5mm,5.5mm>*{^2}**@{},
   %<0mm,0mm>*{};<4.5mm,5.5mm>*{^{m\hspace{-0.5mm}-\hspace{-0.5mm}1}}**@{},
   <0mm,0mm>*{};<9.0mm,5.5mm>*{^m}**@{},
 \end{xy}}\Ea
 + \ldots.
$$
Let us call the unique edge of such a summand {\em special}, and consider a filtration of $\Tw\HoLBcd$ by the number of non-special edges plus the total number of MC generators. On the initial page $E_0$ of the associated spectral sequence  the induced differential acts only on MC generators with $m\geq 2$ by the formula given just above (with no additional terms). Hence the next page $E_1$ of the spectral sequence is equal to the quotient subcomplex $\tw \HoLBcd'$ of the partially twisted properad $\tw\HoLBcd$  by the differential ideal generated by graphs with at least one special edge; the induced differential acts only on the generators
of $\HoLBcd$ by the standard formula (\ref{3: d in qLBcd_infty}). We consider next a filtration of $E_1$ by the total number of paths connecting in-legs and univalent MC generators $\Ba{c}\resizebox{2.5mm}{!}{ \xy
 (0,0)*{\bu}="a",
(0,3.5)*{}="0",
\ar @{-} "a";"0" <0pt>
\endxy}\Ea$  to the out-legs of elements of $E_1$. The induced differential $d$ on the associated graded complex $gr E_1$ is precisely the $\frac{1}{2}$-prop differential\footnote{The notion of $\frac{1}{2}$-prop (as well the closely related notion of path filtration) was introduced by Maxim Kontsevich \cite{Ko3}. A nice exposition of this theory can be found in \cite{MaVo}.}  in $\HoLBcd$ given explicitly by those summands in  (\ref{3: d in qLBcd_infty}) whose lower (or upper) corolla has  type $(1,p\geq 2)$ (or, resp., $(p\geq 2,1))$ only. The point is that such summands never create new special edges which have to be set to zero by hand, i.e.\ the fact that we have to take the quotient by graphs with at least one special edge does not complicate the action of the induced differential any more. Hence the next page $E_2  \sim  H^\bu(gr E_1)$ of our spectral sequence is spanned by graphs generated by the following three corollas
$$
\Ba{c}\resizebox{8mm}{!}{  \xy
(-5,4)*{}="1";
    (-5,0)*{\bu}="L";
  (-8,-5)*+{_1}="C";
   (-2,-5)*+{_2}="D";
\ar @{-} "D";"L" <0pt>
\ar @{-} "C";"L" <0pt>
\ar @{-} "1";"L" <0pt>
 \endxy}
 \Ea
 = (-1)^d
 \Ba{c}\resizebox{8mm}{!}{  \xy
(-5,4)*{}="1";
    (-5,0)*{\bu}="L";
  (-8,-5)*+{_2}="C";
   (-2,-5)*+{_1}="D";
\ar @{-} "D";"L" <0pt>
\ar @{-} "C";"L" <0pt>
\ar @{-} "1";"L" <0pt>
 \endxy}
 \Ea
, \ \ \ \ \
\Ba{c}\resizebox{8mm}{!}{  \xy
(-5,-5)*+{_1}="1";
    (-5,0)*{\bu}="L";
  (-8,4.5)*+{_1}="C";
   (-2,4.5)*+{_2}="D";
\ar @{-} "D";"L" <0pt>
\ar @{-} "C";"L" <0pt>
\ar @{-} "1";"L" <0pt>
 \endxy}
 \Ea
 =(-1)^c
 \Ba{c}\resizebox{8mm}{!}{  \xy
(-5,-5)*+{_1}="1";
    (-5,0)*{\bu}="L";
  (-8,4.5)*+{_2}="C";
   (-2,4.5)*+{_1}="D";
\ar @{-} "D";"L" <0pt>
\ar @{-} "C";"L" <0pt>
\ar @{-} "1";"L" <0pt>
 \endxy}
 \Ea,
 \ \ \ \Ba{c}\resizebox{2.5mm}{!}{ \xy
 (0,0)*{\bu}="a",
(0,5)*{}="0",
\ar @{-} "a";"0" <0pt>
\endxy}\Ea,
$$
subject to the relations
$$
\oint_{123}\hspace{-1mm} \Ba{c}\resizebox{7mm}{!}{
\begin{xy}
 <0mm,0mm>*{\bu};<0mm,0mm>*{}**@{},
 <0mm,-0.49mm>*{};<0mm,-3.0mm>*{}**@{-},
 <0.49mm,0.49mm>*{};<1.9mm,1.9mm>*{}**@{-},
 <-0.5mm,0.5mm>*{};<-1.9mm,1.9mm>*{}**@{-},
 <-2.3mm,2.3mm>*{\bu};<-2.3mm,2.3mm>*{}**@{},
 <-1.8mm,2.8mm>*{};<0mm,4.9mm>*{}**@{-},
 <-2.8mm,2.9mm>*{};<-4.6mm,4.9mm>*{}**@{-},
   <0.49mm,0.49mm>*{};<2.7mm,2.3mm>*{^3}**@{},
   <-1.8mm,2.8mm>*{};<0.4mm,5.3mm>*{^2}**@{},
   <-2.8mm,2.9mm>*{};<-5.1mm,5.3mm>*{^1}**@{},
 \end{xy}}\Ea
=0, \ \  \ \
%%%%%%%%%%%%%% Lie %%%%%%%%%%%%%%%%%%%%%%%%
\oint_{123}\hspace{-1mm} \Ba{c}\resizebox{8.4mm}{!}{ \begin{xy}
 <0mm,0mm>*{\bu};<0mm,0mm>*{}**@{},
 <0mm,0.69mm>*{};<0mm,3.0mm>*{}**@{-},
 <0.39mm,-0.39mm>*{};<2.4mm,-2.4mm>*{}**@{-},
 <-0.35mm,-0.35mm>*{};<-1.9mm,-1.9mm>*{}**@{-},
 <-2.4mm,-2.4mm>*{\bu};<-2.4mm,-2.4mm>*{}**@{},
 <-2.0mm,-2.8mm>*{};<0mm,-4.9mm>*{}**@{-},
 <-2.8mm,-2.9mm>*{};<-4.7mm,-4.9mm>*{}**@{-},
    <0.39mm,-0.39mm>*{};<3.3mm,-4.0mm>*{^3}**@{},
    <-2.0mm,-2.8mm>*{};<0.5mm,-6.7mm>*{^2}**@{},
    <-2.8mm,-2.9mm>*{};<-5.2mm,-6.7mm>*{^1}**@{},
 \end{xy}}\Ea
=0, \ \ \ \Ba{c}\resizebox{5mm}{!}{\begin{xy}
 <0mm,2.47mm>*{};<0mm,0.12mm>*{}**@{-},
 <0.5mm,3.5mm>*{};<2.2mm,5.2mm>*{}**@{-},
 <-0.48mm,3.48mm>*{};<-2.2mm,5.2mm>*{}**@{-},
 <0mm,3mm>*{\bu};<0mm,3mm>*{}**@{},
  <0mm,-0.8mm>*{\bu};<0mm,-0.8mm>*{}**@{},
<-0.39mm,-1.2mm>*{};<-2.2mm,-3.5mm>*{}**@{-},
 <0.39mm,-1.2mm>*{};<2.2mm,-3.5mm>*{}**@{-},
     <0.5mm,3.5mm>*{};<2.8mm,5.7mm>*{^2}**@{},
     <-0.48mm,3.48mm>*{};<-2.8mm,5.7mm>*{^1}**@{},
   <0mm,-0.8mm>*{};<-2.7mm,-5.2mm>*{^1}**@{},
   <0mm,-0.8mm>*{};<2.7mm,-5.2mm>*{^2}**@{},
\end{xy}}\Ea=0, \ \ \ \Ba{c}\resizebox{8mm}{!}{  \xy
(-5,-5)*{\bu}="1";
    (-5,-1)*{\bu}="L";
  (-8,4.5)*+{_1}="C";
   (-2,4.5)*+{_2}="D";
\ar @{-} "D";"L" <0pt>
\ar @{-} "C";"L" <0pt>
\ar @{-} "1";"L" <0pt>
 \endxy}
 \Ea=0.
$$
The induced differential acts only on the MC generator by the standard formula
$$
 \Ba{c}\resizebox{1.4mm}{!}{\begin{xy}
 <0mm,0.5mm>*{};<0mm,5mm>*{}**@{-},
 <0mm,0mm>*{\bu};
 \end{xy}}\Ea
 \lon
\frac{1}{2}
\Ba{c}\resizebox{6.1mm}{!}{  \xy
(-5,6)*{}="1";
    (-5,+1)*{\bu}="L";
  (-8,-4)*{\bu}="C";
   (-2,-4)*{\bu}="D";
\ar @{-} "D";"L" <0pt>
\ar @{-} "C";"L" <0pt>
\ar @{-} "1";"L" <0pt>
%
 %
 %\ar @{-} "l1";"L" <0pt>
 \endxy}
 \Ea.
 $$
 As $H^\bu(\tw\Lie_d)=\Lie_d$, we conclude that the cohomology is spanned by the standard two generators of $\LBcd$ modulo the above three relations. Hence $H^\bu(\Tw\HoLBcd)=\LBcd$ and the Theorem is proven.
\end{proof}

%%%%%%%%%%%%%%%%%%%%%%%%%%%%%%%%%%%%%%%%%%%%%%%%%%%%%%%%%%%%%%%%%%%%%%%%%%%%%%%%%%%

\subsection{An action of the deformation complex  $\Def(\HoLBcd\rar \cP)$  on $\Tw\cP$}
  Given a dg properad $\cP$ under $\HoLBcd$  (see (\ref{4: i from HoLB to P})), one can consider $\cP$ as a dg properad under $\HoLBcd^+$ using the composition
 $$
 i^+: \HoLBcd^+ \lon \HoLBcd \stackrel{i}{\lon} \cP
 $$
 where the first arrow is the unique morphism which sends the $(1,1)$-generator to zero and is the identity on all other generators. Following \cite{MW1} we define the deformation complex of the morphism $i$ in (\ref{4: i from HoLB to P})  as the deformation complex of the morphism $i^+$,
 $$
 \Def\left(\HoLBcd^+ \stackrel{i^+}{\lon} \cP\right)=\prod_{m,n\geq 1} \cP(m,n)\ot_{\bS_m^{op}\times \bS_n}\left(\sgn_m^{|c|}\ot \sgn_n^{|d|}\right)[c(1-m)+d(1-n)].
$$
Note that even in the case when the properad $\cP$ is generated by $(m,n)$-operations with $m,n\geq 1$ and $m+n\geq 3$ (as, e.g., in the case $\cP=\HoLBcd$) the term $\cP(1,1)$ is often non-zero  and should not be ignored  in the deformation theory (cf.\ \cite{MW3}); hence the $+$-extension. As in \cite{MW2}, we abuse notation and re-denote $\Def(\HoLBcd^+ \stackrel{i^+}{\rar} \cP)$ by
$\Def(\HoLBcd \stackrel{i}{\rar} \cP)$.

\sip

%Its Maurer-Cartan elements are in 1-1 correspondence with morphisms $\HoLBcd\rar \cP$ which %are deformations of $i$.
Following \cite{MV} one can describe the differential (and Lie brackets) in
  $\Def\left(\HoLBcd \stackrel{i}{\rar} \cP\right)$ very explicitly. Let us represent a generic element of this complex  pictorially as a collection of $(m,n)$-corollas,
$$
\ga=\left\{
\Ba{c}\resizebox{12mm}{!}{  \xy
(0,9)*{\overbrace{\ \ \ \ \ \ \ \ \ \ \ }^m};
(0,-9)*{\underbrace{\ \ \ \ \ \ \ \ \ \ \ }_n};
(-7,7)*+{}="U1";
(-3,7)*+{}="U2";
(2,5)*{...};
(7,7)*+{}="U3";
(0,0)*{\circledast}="C";
(-7,-7)*+{}="L1";
(-3,-7)*+{}="L2";
(2,-5)*{...};
(7,-7)*+{}="L3";
\ar @{-} "C";"L1" <0pt>
\ar @{-} "C";"L2" <0pt>
\ar @{-} "C";"L3" <0pt>
\ar @{-} "C";"U1" <0pt>
\ar @{-} "C";"U2" <0pt>
\ar @{-} "C";"U3" <0pt>
 \endxy}
 \Ea
 \right\}_{m,n\geq 1},
 \ \ \text{or as a formal sum}\
 \ga=\sum_{m,n\geq 1}
 \Ba{c}\resizebox{12mm}{!}{  \xy
(-7,7)*+{_1}="U1";
(-3,7)*+{_2}="U2";
(2,5)*{...};
(7,7)*+{_m}="U3";
(0,0)*{\circledast}="C";
(-7,-7)*+{_1}="L1";
(-3,-7)*+{_2}="L2";
(2,-5)*{...};
(7,-7)*+{_n}="L3";
\ar @{-} "C";"L1" <0pt>
\ar @{-} "C";"L2" <0pt>
\ar @{-} "C";"L3" <0pt>
\ar @{-} "C";"U1" <0pt>
\ar @{-} "C";"U2" <0pt>
\ar @{-} "C";"U3" <0pt>
 \endxy}
 \Ea,
$$
%they have degree |c| - c(1-m) - d(1-n), where |c$ is the degree as c an element of P(m,n)
%A formal sum of such corollas is homogeneous of degree $k$ if and only f the degree |c| of %each (m,n)-corolla in P in  is equal to k+c(1-m)+d(1-n)
of corollas whose vertices are decorated with  elements\footnote{A formal sum $\ga$ of such $(m,n)$-corollas is homogeneous of degree $k$ as an element of the deformation complex if and only if their degrees $|\circledast|$
as elements of $\cP(m,n)\ot_{\bS_m^{op}\times \bS_n} (\sgn_m^{|c|}\ot \sgn_n^{|d|})$ are equal to  $k+c(1-m)+d(1-n)$. This explains the grading consistency of the explicit formulae
shown
in Theorem {\ref{3: Theorem on Def action on TwP}} below.} of $\cP(m,n)\ot_{\bS_m^{op}\times \bS_n} (\sgn_m^{|c|}\ot \sgn_n^{|d|})$ and denoted by $\circledast$ in order to distinguish them from the generic elements of $\cP$ which are represented as corollas (\ref{3: generic elements of cP as (m,n)-corollas})
 and also from  the images of $\Holie_d$-generators under $i$ which are represented pictorially as $\circledcirc$-vertex corollas. Since labels of in- and out legs are (skew)symmetrized, one cam omit them in pictures. The differential in the deformation complex is given by
\Beq\label{3: differential in Def(Lieb to P)}
\delta \hspace{-3mm}
\Ba{c}\resizebox{12mm}{!}{  \xy
(0,9)*{\overbrace{\ \ \ \ \ \ \ \ \ \ \ }^m};
(0,-9)*{\underbrace{\ \ \ \ \ \ \ \ \ \ \ }_n};
(-7,7)*+{}="U1";
(-3,7)*+{}="U2";
(2,5)*{...};
(7,7)*+{}="U3";
(0,0)*{\circledast}="C";
(-7,-7)*+{}="L1";
(-3,-7)*+{}="L2";
(2,-5)*{...};
(7,-7)*+{}="L3";
\ar @{-} "C";"L1" <0pt>
\ar @{-} "C";"L2" <0pt>
\ar @{-} "C";"L3" <0pt>
\ar @{-} "C";"U1" <0pt>
\ar @{-} "C";"U2" <0pt>
\ar @{-} "C";"U3" <0pt>
 \endxy}
 \Ea
 %%%%%%%%%%%%%%%%%%%%%%%%%%%%%%%%%%%%
 =
 %%%%%%%%%%%%%%%%%%%%%%%%%%%%%%%%%%%
\sd \hspace{-3mm}
\Ba{c}\resizebox{12mm}{!}{  \xy
(0,9)*{\overbrace{\ \ \ \ \ \ \ \ \ \ \ }^m};
(0,-9)*{\underbrace{\ \ \ \ \ \ \ \ \ \ \ }_n};
(-7,7)*+{}="U1";
(-3,7)*+{}="U2";
(2,5)*{...};
(7,7)*+{}="U3";
(0,0)*{\circledast}="C";
(-7,-7)*+{}="L1";
(-3,-7)*+{}="L2";
(2,-5)*{...};
(7,-7)*+{}="L3";
\ar @{-} "C";"L1" <0pt>
\ar @{-} "C";"L2" <0pt>
\ar @{-} "C";"L3" <0pt>
\ar @{-} "C";"U1" <0pt>
\ar @{-} "C";"U2" <0pt>
\ar @{-} "C";"U3" <0pt>
 \endxy}
 \Ea
+
 %%%%%%%%%%%%%%%%%%%%%%%%%%%%%%%%%%%%%%%%%%%%%%%
\sum_{[n]=[n']\sqcup [n'']\atop n'\geq 2,n''\geq 1}\left(
\pm
\Ba{c}\resizebox{14mm}{!}{  \xy
(-10,8)*{}="11";
(-7,8)*{}="12";
(0,8)*{}="13";
(-3,-5)*{...};
(-13,-13)*{...};
 <-15mm,-17mm>*{\underbrace{ \ \ \ \ \  \ \ \ \ \ \ }_{n'}},
  <-16mm,5mm>*{\overbrace{ \ \ \ \ \ \ \ \ \ \  }^{m'}},
 <-3mm,-9mm>*{\underbrace{ \ \ \ \  \ \ \ \ \  }_{n''}},
 <-5mm,12mm>*{\overbrace{ \ \ \ \ \ \  \ \ \ \ \  }^{m''}},
    (-5,+2)*{\circledcirc}="L";
%%%
 (-14,-6)*{\circledast}="B";
 (-20,-13)*{}="b1";
 (-17,-13)*{}="b2";
 (-10,-13)*{}="b3";
 (-20,1)*{}="a1";
 (-17,1)*{}="a2";
 (-11,1)*{}="a3";
  (-8,-5)*{}="C";
   (1,-5)*{}="D";
     %<0mm,12mm>*{_m},
\ar @{-} "D";"L" <0pt>
\ar @{-} "C";"L" <0pt>
\ar @{-} "B";"L" <0pt>
\ar @{-} "B";"b1" <0pt>
\ar @{-} "B";"b2" <0pt>
\ar @{-} "B";"b3" <0pt>
\ar @{-} "B";"a1" <0pt>
\ar @{-} "B";"a2" <0pt>
\ar @{-} "B";"a3" <0pt>
\ar @{-} "11";"L" <0pt>
\ar @{-} "12";"L" <0pt>
\ar @{-} "13";"L" <0pt>
%
 %
 %\ar @{-} "l1";"L" <0pt>
 \endxy}
 \Ea
 %%%%%%%%%%%%%%%%%%%%%%%%%%%%%%%%%%%%%%%%%%%%%
\mp
\Ba{c}\resizebox{14mm}{!}{  \xy
(-10,8)*{}="11";
(-7,8)*{}="12";
(0,8)*{}="13";
(-3,-5)*{...};
(-13,-13)*{...};
 <-15mm,-17mm>*{\underbrace{ \ \ \ \ \  \ \ \ \ \ \ }_{n'}},
  <-16mm,5mm>*{\overbrace{ \ \ \ \ \ \ \ \ \ \  }^{m'}},
 <-3mm,-9mm>*{\underbrace{ \ \ \ \  \ \ \ \ \  }_{n''}},
 <-5mm,12mm>*{\overbrace{ \ \ \ \ \ \  \ \ \ \ \  }^{m''}},
    (-5,+2)*{\circledast}="L";
%%%
 (-14,-6)*{\circledcirc}="B";
 (-20,-13)*{}="b1";
 (-17,-13)*{}="b2";
 (-10,-13)*{}="b3";
 (-20,1)*{}="a1";
 (-17,1)*{}="a2";
 (-11,1)*{}="a3";
  (-8,-5)*{}="C";
   (1,-5)*{}="D";
     %<0mm,12mm>*{_m},
\ar @{-} "D";"L" <0pt>
\ar @{-} "C";"L" <0pt>
\ar @{-} "B";"L" <0pt>
\ar @{-} "B";"b1" <0pt>
\ar @{-} "B";"b2" <0pt>
\ar @{-} "B";"b3" <0pt>
\ar @{-} "B";"a1" <0pt>
\ar @{-} "B";"a2" <0pt>
\ar @{-} "B";"a3" <0pt>
\ar @{-} "11";"L" <0pt>
\ar @{-} "12";"L" <0pt>
\ar @{-} "13";"L" <0pt>
%
 %
 %\ar @{-} "l1";"L" <0pt>
 \endxy}
 \Ea
 %%%%%%%%%%%%%%%%%%%%%%%%%%%%%%%%%%%%%%%%%%%%%%%
 \right)
\Eeq
where the rule of signs depends on $d$ and is read from (\ref{2: d in Lie_infty}); for $d$ even the first ambiguous sign symbol on the r.h.s.\ is $+1$, while the second one is $-(-1)^{|\circledast|}$.
%the differential $\sd$ in the first term above acts on the decorations of $\circledast$-corollas.
%In fact the above equation should be understood as an infinite collection of equations, one for for all %summands with precisely  $N\geq 2$  incoming legs.

\sip

Let $(\Der(\Tw\cP),\ [\ ,\ ])$ be the  Lie algebra of derivations of the properad $\Tw\cP$. The differential $\sd_\centerdot$ is its MC element making $\Der(\Tw\cP)$ into a dg Lie algebra with the differential $[\sd_\centerdot,\ ]$.
% Recall that any element $h\in \Tw\cP(1,1)$ defines a derivation $D_h$ of $\Tw\cP$
% (viewed as a non-differential prop(erad)) given by the formula (\ref{3: formula for D_h}).

\subsubsection{\bf Theorem}\label{3: Theorem on Def action on TwP}{\it
There is a morphism of dg Lie algebras
\Beq\label{3: Def to Der(P)}
\Ba{rccc}
\Phi: & \Def\left(\HoLBcd \stackrel{i}{\rar} \cP\right) &\lon & \Der(\Tw\cP)\\
      &    \ga & \lon &    \Phi_\ga
\Ea
\Eeq
where the derivation $\Phi_\ga$ is given on the generators by
$$
\Ba{rccc}
\Phi_\ga: &\Tw\cP & \lon &  \Tw\cP \ \ \ \ \ \ \ \ \ \ \  \ \ \ \ \ \\
&
 \Ba{c}\resizebox{14mm}{!}{
 \begin{xy}
 <0mm,-1mm>*{\bu};
 <0mm,0mm>*{};<-8mm,5mm>*{}**@{-},
 <0mm,0mm>*{};<-4.5mm,5mm>*{}**@{-},
 <0mm,0mm>*{};<-1mm,5mm>*{\ldots}**@{},
 <0mm,0mm>*{};<4.5mm,5mm>*{}**@{-},
 <0mm,0mm>*{};<8mm,5mm>*{}**@{-},
   <0mm,0mm>*{};<-8.5mm,5.5mm>*{^1}**@{},
   <0mm,0mm>*{};<-5mm,5.5mm>*{^2}**@{},
   %<0mm,0mm>*{};<4.5mm,5.5mm>*{^{m\hspace{-0.5mm}-\hspace{-0.5mm}1}}**@{},
   <0mm,0mm>*{};<9.0mm,5.5mm>*{^m}**@{},
 \end{xy}}\Ea
& \lon &
\displaystyle
 \overset{m-1}{\underset{i=0}{\sum}}
\Ba{c}\resizebox{15mm}{!}{
\begin{xy}
 <0mm,0mm>*{\bu};<-8mm,5mm>*{}**@{-},
 <0mm,0mm>*{\bu};<-3.5mm,5mm>*{}**@{-},
 <0mm,0mm>*{\bu};<-6mm,5mm>*{..}**@{},
 <0mm,0mm>*{\bu};<0mm,5mm>*{}**@{-},
  <0mm,5mm>*{\lozenge};
  <0mm,5mm>*{};<0mm,8mm>*{}**@{-},
  <0mm,5mm>*{};<0mm,9mm>*{^{i\hspace{-0.2mm}+\hspace{-0.5mm}1}}**@{},
<0mm,0mm>*{\bu};<8mm,5mm>*{}**@{-},
<0mm,0mm>*{\bu};<3.5mm,5mm>*{}**@{-},
<6mm,5mm>*{..}**@{},
<-8.5mm,5.5mm>*{^1}**@{},
<-4mm,5.5mm>*{^i}**@{},
<9.0mm,5.5mm>*{^m}**@{},
 \end{xy}}\Ea
 +
 \sum_{k\geq 1, m=\sum m_\bu,\atop {  m_{0}\geq 1, k+m_0\geq 3\atop m_1,...,m_k\geq 0}}\frac{1}{k!}
 \Ba{c}\resizebox{22mm}{!}{  \xy
(-28,12)*{}="1";
(-25,12)*{}="2";
(-22,12)*{}="3";
%(2,-4)*{\ldots};
%
(-18,8)*{}="n11";
(-15,7)*{...};
(-13,8)*{}="n12";
(-11,8)*{}="n21";
(-8.4,7)*{...};
(-6,8)*{}="n22";
(1,8)*{}="nn1";
(3,7)*{...};
(5,8)*{}="nn2";
    (-25,+4)*{\circledast}="L";
 (-15,-5)*{\bu}="B";
  (-8,-5)*{\bu}="C";
   (3,-5)*{\bu}="D";
   <-5mm,-10mm>*{\underbrace{ \ \ \ \ \ \ \ \ \ \ \ \ \ \ \ \ \ }_{k}},
   (-3,-5)*{...};
    <-25mm,14.6mm>*{\overbrace{ \ \ \ \ \ \ }^{m_0}},
    <-16mm,10.6mm>*{\overbrace{ \ \ \ \ }^{m_1}},
    <-9mm,10.6mm>*{\overbrace{ \ \ \ \ }^{m_2}},
    <3mm,10.6mm>*{\overbrace{ \ \ \ \ }^{m_n}},
     %<0mm,12mm>*{_m},
\ar @{-} "D";"L" <0pt>
\ar @{-} "C";"L" <0pt>
\ar @{-} "B";"L" <0pt>
\ar @{-} "1";"L" <0pt>
\ar @{-} "2";"L" <0pt>
\ar @{-} "3";"L" <0pt>
\ar @{-} "1";"L" <0pt>
\ar @{-} "2";"L" <0pt>
\ar @{-} "n11";"B" <0pt>
\ar @{-} "n12";"B" <0pt>
\ar @{-} "n21";"C" <0pt>
\ar @{-} "n22";"C" <0pt>
\ar @{-} "nn1";"D" <0pt>
\ar @{-} "nn2";"D" <0pt>
 \endxy}
 \Ea\vspace{3mm}
 \\
 &
 \Ba{c}\resizebox{14mm}{!}{
 \begin{xy}
 <0mm,0mm>*{\circ};<-8mm,6mm>*{^1}**@{-},
 <0mm,0mm>*{\circ};<-4.5mm,6mm>*{^2}**@{-},
 <0mm,0mm>*{\circ};<0mm,5.5mm>*{\ldots}**@{},
 <0mm,0mm>*{\circ};<3.5mm,5mm>*{}**@{-},
 <0mm,0mm>*{\circ};<8mm,6mm>*{^m}**@{-},
 <0mm,0mm>*{\circ};<-8mm,-6mm>*{_1}**@{-},
 <0mm,0mm>*{\circ};<-4.5mm,-6mm>*{_2}**@{-},
 <0mm,0mm>*{\circ};<0mm,-5.5mm>*{\ldots}**@{},
 <0mm,0mm>*{\circ};<4.5mm,-6mm>*+{}**@{-},
 <0mm,0mm>*{\circ};<8mm,-6mm>*{_n}**@{-},
   \end{xy}}\Ea
& \lon &
\overset{m-1}{\underset{i=0}{\sum}}
\Ba{c}\resizebox{14mm}{!}{
\begin{xy}
 %<0mm,0mm>*{\circ};<0mm,0mm>*{}**@{},
 <0mm,0mm>*{\circ};<-8mm,5mm>*{}**@{-},
 <0mm,0mm>*{\circ};<-3.5mm,5mm>*{}**@{-},
 <0mm,0mm>*{\circ};<-6mm,5mm>*{..}**@{},
 <0mm,0mm>*{\circ};<0mm,5mm>*{}**@{-},
  <0mm,5mm>*{\lozenge};
  <0mm,5mm>*{};<0mm,8mm>*{}**@{-},
  <0mm,5mm>*{};<0mm,9mm>*{^{i\hspace{-0.2mm}+\hspace{-0.5mm}1}}**@{},
<0mm,0mm>*{\circ};<8mm,5mm>*{}**@{-},
<0mm,0mm>*{\circ};<3.5mm,5mm>*{}**@{-},
<6mm,5mm>*{..}**@{},
<-8.5mm,5.5mm>*{^1}**@{},
<-4mm,5.5mm>*{^i}**@{},
<9.0mm,5.5mm>*{^m}**@{},
 <0mm,0mm>*{\circ};<-8mm,-5mm>*{}**@{-},
 <0mm,0mm>*{\circ};<-4.5mm,-5mm>*{}**@{-},
 <-1mm,-5mm>*{\ldots}**@{},
 <0mm,0mm>*{\circ};<4.5mm,-5mm>*{}**@{-},
 <0mm,0mm>*{\circ};<8mm,-5mm>*{}**@{-},
<-8.5mm,-6.9mm>*{^1}**@{},
<-5mm,-6.9mm>*{^2}**@{},
%<4.5mm,-6.9mm>*{^{n\hspace{-0.5mm}-\hspace{-0.5mm}1}}**@{},
<9.0mm,-6.9mm>*{^n}**@{},
 \end{xy}}\Ea
 - (-1)^{|a|}
\overset{n-1}{\underset{i=0}{\sum}}
 \Ba{c}\resizebox{14mm}{!}{\begin{xy}
 %<0mm,0mm>*{\circ};
 <0mm,0mm>*{\circ};<-8mm,-5mm>*{}**@{-},
 <0mm,0mm>*{\circ};<-3.5mm,-5mm>*{}**@{-},
 <0mm,0mm>*{\circ};<-6mm,-5mm>*{..}**@{},
 <0mm,0mm>*{\circ};<0mm,-5mm>*{}**@{-},
  <0mm,-5mm>*{\lozenge};
  <0mm,-5mm>*{};<0mm,-8mm>*{}**@{-},
  <0mm,-5mm>*{};<0mm,-10mm>*{^{i\hspace{-0.2mm}+\hspace{-0.5mm}1}}**@{},
<0mm,0mm>*{\circ};<8mm,-5mm>*{}**@{-},
<0mm,0mm>*{\circ};<3.5mm,-5mm>*{}**@{-},
 <6mm,-5mm>*{..}**@{},
<-8.5mm,-6.9mm>*{^1}**@{},
<-4mm,-6.9mm>*{^i}**@{},
<9.0mm,-6.9mm>*{^n}**@{},
 <0mm,0mm>*{\circ};<-8mm,5mm>*{}**@{-},
 <0mm,0mm>*{\circ};<-4.5mm,5mm>*{}**@{-},
<-1mm,5mm>*{\ldots}**@{},
 <0mm,0mm>*{\circ};<4.5mm,5mm>*{}**@{-},
 <0mm,0mm>*{\circ};<8mm,5mm>*{}**@{-},
<-8.5mm,5.5mm>*{^1}**@{},
<-5mm,5.5mm>*{^2}**@{},
%<4.5mm,5.5mm>*{^{m\hspace{-0.5mm}-\hspace{-0.5mm}1}}**@{},
<9.0mm,5.5mm>*{^m}**@{},
 \end{xy}}\Ea
\Ea
$$
where (cf.\ (40))
\Beq\label{3: (1,1) part Def action on TwP}
\begin{xy}
 <0mm,-3mm>*{};<0mm,3mm>*{}**@{-},
 %<0mm,0.5mm>*{};<0mm,3mm>*{}**@{-},
 <0mm,0mm>*{_\lozenge};<0mm,0mm>*{}**@{},
 \end{xy}:= -
 \sum_{k=1}^\infty \frac{1}{k!}\ \resizebox{18mm}{!}{  \xy
(-25,8)*{}="1";
(-9,-4)*{...};
 <-11mm,-8mm>*{\underbrace{ \ \ \ \ \ \ \ \ \ \ \ \ \ \ \ \ \ }_{k}},
    (-25,+3)*{\circledast}="L";
    (-25,-3)*{}="N";
 (-19,-4)*{\bu}="B";
  (-13,-4)*{\bu}="C";
   (-2,-4)*{\bu}="D";
\ar @{-} "D";"L" <0pt>
\ar @{-} "C";"L" <0pt>
\ar @{-} "B";"L" <0pt>
\ar @{-} "1";"L" <0pt>
\ar @{-} "N";"L" <0pt>
%
 %
% \ar @{-} "l1";"N" <0pt>
%  \ar @{-} "l2";"N" <0pt>
%   \ar @{-} "ln";"N" <0pt>
 \endxy}
\Eeq
}
% If ga is homogeneous of degree $k$, then each circldust-corolla, viewed as an element of %P(1,k+1)  has degree |c|= k + d(1-(k+1))=k-dk, so that every summand in this sum is %homogeneous of degree k as expected.
% Similarly, for MC generators with m outputs: the k-th summand has degree k+c(1-m_0)+ %d(1-k) + (1-m_1-1)c+d+...+(1-m_k-1)c + d= k + c(1-m) as expected!
\hspace{-2mm}
\begin{proof} (A sketch). Any derivation of $\Tw \cP$ (viewed as a non-differential properad) is uniquely determined by its values on the MC generators  $\Ba{c}\resizebox{10mm}{!}{
 \begin{xy}
 <0mm,-1mm>*{\bu};
 <0mm,0mm>*{};<-8mm,5mm>*{}**@{-},
 <0mm,0mm>*{};<-4.5mm,5mm>*{}**@{-},
 <0mm,0mm>*{};<-1mm,5mm>*{\ldots}**@{},
 <0mm,0mm>*{};<4.5mm,5mm>*{}**@{-},
 <0mm,0mm>*{};<8mm,5mm>*{}**@{-},
   <0mm,0mm>*{};<-8.5mm,5.5mm>*{^1}**@{},
   <0mm,0mm>*{};<-5mm,5.5mm>*{^2}**@{},
   <0mm,0mm>*{};<4.5mm,5.5mm>*{^{m\hspace{-0.5mm}-\hspace{-0.5mm}1}}**@{},
   <0mm,0mm>*{};<9.0mm,5.5mm>*{^m}**@{},
 \end{xy}}\Ea$ and on arbitrary elements of $\cP$.
  The first values can be chosen arbitrary, while the second ones must be compatible with the properad compositions; as the second values are, by the definition, of the form (\ref{3: formula for D_h}), we conclude that the above formulae do define a derivation of $\TW\cP$ as a non-differential prop(erad). Hence the main point is to show that $\Phi_\ga$ respects differentials in both dg Lie algebras, i.e. satisfies the equation
 \Beq\label{3: compatibility of Def to Der with d}
 \Phi_{\delta\ga}=[\p_\centerdot, \Phi_\ga].
\Eeq
Consider first a simpler morphism of non-differential graded Lie algebras,
$$
\Ba{rccc}
\wt{\Phi}: & \Def\left(\HoLBcd \stackrel{i}{\rar} \cP\right) &\lon & \Der(\TW\cP)\\
      &    \ga & \lon &    \wt{\Phi}_\ga
\Ea
$$
where the derivation $\wt{\Phi}_\ga\in \Der(\TW \cP)$ is given on the generators by
$$
\wt{\Phi}_\ga\left(\hspace{-2mm}\Ba{c}\resizebox{14mm}{!}{
 \begin{xy}
 <0mm,-1mm>*{\bu};
 <0mm,0mm>*{};<-8mm,5mm>*{}**@{-},
 <0mm,0mm>*{};<-4.5mm,5mm>*{}**@{-},
 <0mm,0mm>*{};<-1mm,5mm>*{\ldots}**@{},
 <0mm,0mm>*{};<4.5mm,5mm>*{}**@{-},
 <0mm,0mm>*{};<8mm,5mm>*{}**@{-},
   <0mm,0mm>*{};<-8.5mm,5.5mm>*{^1}**@{},
   <0mm,0mm>*{};<-5mm,5.5mm>*{^2}**@{},
   <0mm,0mm>*{};<4.5mm,5.5mm>*{^{m\hspace{-0.5mm}-\hspace{-0.5mm}1}}**@{},
   <0mm,0mm>*{};<9.0mm,5.5mm>*{^m}**@{},
 \end{xy}}\Ea\hspace{-1mm}\right)
 =
%\displaystyle
%
 \sum_{k\geq 1, m=\sum m_\bu,\atop {  m_{0}\geq 1, k+m_0\geq 3\atop m_1,...,m_k\geq 0}}\frac{1}{k!}
 \Ba{c}\resizebox{20mm}{!}{  \xy
(-28,12)*{}="1";
(-25,12)*{}="2";
(-22,12)*{}="3";
%(2,-4)*{\ldots};
%
(-18,8)*{}="n11";
(-15,7)*{...};
(-13,8)*{}="n12";
(-11,8)*{}="n21";
(-8.4,7)*{...};
(-6,8)*{}="n22";
(1,8)*{}="nn1";
(3,7)*{...};
(5,8)*{}="nn2";
    (-25,+4)*{\circledast}="L";
 (-15,-5)*{\bu}="B";
  (-8,-5)*{\bu}="C";
   (3,-5)*{\bu}="D";
   <-5mm,-10mm>*{\underbrace{ \ \ \ \ \ \ \ \ \ \ \ \ \ \ \ \ \ }_{k}},
   (-3,-5)*{...};
    <-25mm,14.6mm>*{\overbrace{ \ \ \ \ \ \ }^{m_0}},
    <-16mm,10.6mm>*{\overbrace{ \ \ \ \ }^{m_1}},
    <-9mm,10.6mm>*{\overbrace{ \ \ \ \ }^{m_2}},
    <3mm,10.6mm>*{\overbrace{ \ \ \ \ }^{m_n}},
     %<0mm,12mm>*{_m},
\ar @{-} "D";"L" <0pt>
\ar @{-} "C";"L" <0pt>
\ar @{-} "B";"L" <0pt>
\ar @{-} "1";"L" <0pt>
\ar @{-} "2";"L" <0pt>
\ar @{-} "3";"L" <0pt>
\ar @{-} "1";"L" <0pt>
\ar @{-} "2";"L" <0pt>
\ar @{-} "n11";"B" <0pt>
\ar @{-} "n12";"B" <0pt>
\ar @{-} "n21";"C" <0pt>
\ar @{-} "n22";"C" <0pt>
\ar @{-} "nn1";"D" <0pt>
\ar @{-} "nn2";"D" <0pt>
 \endxy}
 \Ea, \ \ \ \ \ \
\wt{\Phi}_\ga\left(\hspace{-1.5mm} \Ba{c}\resizebox{14mm}{!}{
 \begin{xy}
 <0mm,0mm>*{\circ};<-8mm,6mm>*{^1}**@{-},
 <0mm,0mm>*{\circ};<-4.5mm,6mm>*{^2}**@{-},
 <0mm,0mm>*{\circ};<0mm,5.5mm>*{\ldots}**@{},
 <0mm,0mm>*{\circ};<3.5mm,5mm>*{}**@{-},
 <0mm,0mm>*{\circ};<8mm,6mm>*{^m}**@{-},
 <0mm,0mm>*{\circ};<-8mm,-6mm>*{_1}**@{-},
 <0mm,0mm>*{\circ};<-4.5mm,-6mm>*{_2}**@{-},
 <0mm,0mm>*{\circ};<0mm,-5.5mm>*{\ldots}**@{},
 <0mm,0mm>*{\circ};<4.5mm,-6mm>*+{}**@{-},
 <0mm,0mm>*{\circ};<8mm,-6mm>*{_n}**@{-},
   \end{xy}}\Ea\hspace{-1.5mm}\right) =  0.
$$
The map $\wt{\Phi}$ respects the Lie bracket while the obstruction for this map to respect
 the differentials is given by the derivation of type  (\ref{3: formula for D_h}),
$$
[\sd_\centerdot, \wt{\Phi}_\ga] - \wt{\Phi}_\ga = D_{\ga_1}, \ \ \
\ga_1=\wt{\Phi}_\ga\left(\hspace{-1mm} \Ba{c}\resizebox{2mm}{!}{\begin{xy}
 <0mm,-0.55mm>*{};<0mm,-3mm>*{}**@{-},
 <0mm,0.5mm>*{};<0mm,3mm>*{}**@{-},
 <0mm,0mm>*{\blacklozenge};<0mm,0mm>*{}**@{},
 \end{xy}}\Ea
 \hspace{-1mm}
 \right) \in \TW\cP(1,1),
$$
with $\Ba{c}\resizebox{1.8mm}{!}{\begin{xy}
 <0mm,-0.55mm>*{};<0mm,-3mm>*{}**@{-},
 <0mm,0.5mm>*{};<0mm,3mm>*{}**@{-},
 <0mm,0mm>*{\blacklozenge};<0mm,0mm>*{}**@{},
 \end{xy}}\Ea
$ given by (\ref{3: blacklozenge (1,1) element}). It is a straightforward calculation to check
that the adjustment of the derivation $\wt{\Phi}_\ga$  with an extra term of the type (\ref{3: formula for D_h}),
$$
\wt{\Phi}_\ga \lon \Phi_\ga = \wt{\Phi}_\ga + D_{\begin{xy}
 <0mm,-2.5mm>*{};<0mm,2.5mm>*{}**@{-},
 %<0mm,0.5mm>*{};<0mm,3mm>*{}**@{-},
 <0mm,0mm>*{_\lozenge};<0mm,0mm>*{}**@{},
 \end{xy}}
$$
solves the problem of the compatibility with the differentials.
\end{proof}

\subsection{Grothendieck-Teichm\"uller group and twisted properads} Let $\wHoLBcd$ be the genus completion of the properad $\HoLBcd$. It was proven in \cite{MW2} that for any $c,d\in \Z$ there is a morphism of dg Lie algebras
$$
F \colon \GC_{c+d+1}^{or}\to \Der(\wHoLBcd)
$$
which is a quasi-isomorphism up to one rescaling class (which controls the
automorphism of $\HoLBcd$ given by rescaling each $(m,n)$ generator by $\la^{m+n-2}$ for any $\la \in
K^*$). Here $\GC_{c+d+1}^{or}$ stands for the oriented version of the Kontsevich graph complex from \S {\ref{3: subsec on GC and HGC}} which was studied in \cite{W2} and where it was proven that
$$
H^\bu(\GC_{3}^{or})=H^\bu(\GC_{2})= \fg\fr\ft_1,
$$
 This result implies that for any $c,d\in \Z$ with $c+d=2$, one has an isomorphism of Lie algebras,
$$
H^0(\Der(\wHoLBcd))=\grt
$$
where $\fg\fr\ft$ is the Lie algebra of the ``full" Grothendieck-Teichm\"uller group $GRT_1$
\cite{Dr2}.

\sip

Let $\wh{\cP}$ be a dg properad under $\wHoLBcd$ and let $\TW\wh{\cP}$ be the associated twisted properad. One has morphisms
 %constructed as explained in \S {\ref{3: subsec on twsting of P under Lieb}}.
% It comes equipped with a morphism of dg properads
$$
\TW(i): (\wHoLBcd,\delta) \lon (\TW\wh{\cP},\sd_\centerdot), \ \ \ \
\Phi:\Def\left(\wHoLBcd \stackrel{i}{\rar} \wh{\cP}\right) \lon  \Der(\TW\wh{\cP})
 $$
given explicitly by the same formulae as in (\ref{4: map frpm HoLBcd to TwP}) and in Theorem {\ref{3: Theorem on Def action on TwP}}.

\subsubsection{\bf Proposition}  {\em For any dg properad  $\wh{\cP}$ under $\wHoLBcd$ there is an associated morphism of complexes,
$$
\cF: \GC_{c+d+1}^{or} \lon \Der(\TW\wh{\cP})[1],
$$
where $\GC_{c+d+1}^{or}$ is the oriented version of the Kontsevich graph complex. If $c+d=2$, there is an associated linear map
$
\grt\lon H^1\left(\Der(\TW\wh{\cP})\right).
$
}

\begin{proof} The morphism $\TW(i)$ induces a morphism of dg Lie algebras,
$$
\Def(\wHoLBcd \stackrel{\Id}{\rar} \wHoLBcd) \lon \Def(\wHoLBcd \stackrel{\TW(i)}{\lon} \TW\wh{\cP})
$$
The l.h.s.\ can be identified as a complex (but not as a Lie algebra) with the degree shifted
derivation complex $\Der(\wHoLBcd)[-1]$ while the r.h.s.\ can be mapped, according to
Theorem  {\ref{3: Theorem on Def action on TwP}}, into the complex $\Der(\TW(\wh{\cP})$. Thus we obtain a chain of morphisms of complexes
$$
\cF: \GC_{c+d+1}^{or} \stackrel{F}{\lon} \Der(\wHoLBcd) \lon  \Der(\TW\wh{\cP})[1]
$$
which proves the claim.
\end{proof}

 Thus  fully twisted completed properads under $\HoLBcd$ can have potentially a highly non-trivial homotopy theory depending on the properties of the above map $\cF$ at the cohomology level.

\subsection{From representations of $\cP$ to representations of $\Tw\cP$}
Assume a dg properad $\cP$ under $\HoLBcd$ admits a representation in a dg space $(V,d)$. Then the graded vector space $V[-c]$ is a $\HoLB_{0,c+d}$-algebra. For any Maurer-Cartan element $\ga\in \odot^{\geq 1}(V[-c])$ , that is, for any a solution of the equation (\ref{4: MC eqn for HoLB0d}), we obtain  a presentation of $\Tw\cP$ in $V$ equipped with the twisted differential  (\ref{4: twisted by MC differentiail in V}).  Let us consider such twisted representations in the case $\cP=\HoLBcd$ in  detail.

\subsubsection{\bf Twisted $\HoLBcd$-algebra structures: an explicit description} Let $(V,d)$ be a graded vector space equipped with a basis $\{e_\al\}$ and $V^*$ its dual equipped with the dual basis $\{e^\al\}$. Consider a graded
commutative tensor algebra
$$
\odot^{\geq 1}(V[-c])\otimes\odot^{\geq 1}(V^*[-d]) \subset \odot^{\bu}\left(V[-c]\oplus V^*[-d]\right)\simeq \K[x^\al, p_\al]
$$
where $x^\al:=\fs^{-d} e^\al$, $p_\al:=\fs^{-c} e_\al$.
The paring $V[-c]\ot V^*[-d]\rar \K[-c-d]$ makes this space into a Lie algebra with respect to the Poisson type brackets $\{\ ,\ \}$ (of degree $-c-d$).

\sip

There is a 1-1 correspondence representations of $\HoLBcd$ in $V$,
$$
\rho: \HoLBcd \lon \cE nd_V   \Rightarrow \left\{\rho\left(\hspace{-2mm} \Ba{c}\resizebox{8mm}{!}{ \xy
(0,4.5)*+{...},
(0,-4.5)*+{...},
(0,0)*{\bu}="o",
(-5,5)*{}="1",
(-3,5)*{}="2",
(3,5)*{}="3",
(5,5)*{}="4",
(-3,-5)*{}="5",
(3,-5)*{}="6",
(5,-5)*{}="7",
(-5,-5)*{}="8",
(-5.5,7)*{_1},
(-3,7)*{_2},
(3,6)*{},
(5.9,7)*{m},
(-3,-7)*{_2},
(3,-7)*+{},
(5.9,-7)*{n},
(-5.5,-7)*{_1},
\ar @{-} "o";"1" <0pt>
\ar @{-} "o";"2" <0pt>
\ar @{-} "o";"3" <0pt>
\ar @{-} "o";"4" <0pt>
\ar @{-} "o";"5" <0pt>
\ar @{-} "o";"6" <0pt>
\ar @{-} "o";"7" <0pt>
\ar @{-} "o";"8" <0pt>
\endxy}\Ea \hspace{-2mm} \right)=: \pi_n^m(x,p)\in \odot^n(V^*[-d]) \otimes \odot^m(V[-c])\right\}_{m,n\geq 1, m+n\geq 3},
$$
and Maurer-Cartan elements of the Lie algebra $(\odot^{\geq 1}(V[c])\bigotimes\odot^{\geq 1}(V^*[d]), \{\ , \ \})$, that is, with degree $1+c+d$ elements
$$
\pi(x,p)=\sum_{n,m\geq 1} \pi_n^m(x,p)=\sum_{m\geq 1}\sum_{\al_\bu, \be_\bu} \frac{1}{m!n!} \pi^{\al_1\ldots \al_m}_{\be_1\ldots\be_n} p_{\al_1}\cdots p_{\al_m} x^{\be_1}\cdots x^{\be_n}
$$
such that
$$
\{\pi,\pi\}=2\sum_{\al} \pm \frac{\p \pi}{\p x^\al}\frac{\p \pi}{\p p_\al}=0
$$
and the degree $1$ summand $\pi_1^1\in V\ot V^*$ is precisely the given differential $d$ in $V$.
The representation $\rho$ induces a representation of the associated polydifferential operad
$$
\f_{c,d}(\rho):\ \ \ \f_{c,d}(\HoLBcd) \lon \cE nd_{\odot^{\bu} (V[-c])}
$$
and hence a $\Holie_{c+d}$-algebra structure on $\cE nd_{\odot^\bu (V\{-c\})}$ via the composition of  $\f_{c,d}(\rho)$ with the map (\ref{3: map Holie^+_{c+d} to f_{c,d}(HoLBcd)}). Let
$$
\ga(p)=\sum_{m\geq 1} \ga_m(p),\ \ \ \ga_m\in \odot^{m} (V[-c])
\subset \K[p_\al]
$$
be a Maurer-Cartan element of that $\Holie_{c+d}$-algebra structure, that is, a degree $c+d$ solution of the following explicit coordinate incarnation of the Maurer-Cartan equation
(\ref{4: MC eqn for HoLB0d}),
$$
\sum_{n\geq 1} \pm \frac{1}{n!}\frac{\p^n \pi}{\p x^{\al_1}\ldots \p x^{\al_n}}|_{x=0}
\frac{\p \ga(p)}{\p p_{\al_1}}\ldots\frac{\p \ga(p)}{\p p_{\al_n}}=0.
$$
 Then the data
$\pi(x,p)$ and $\ga(p)$   give rise to a representation,
$$
\Ba{rccc}
 \rho^{Tw}: & \TW\HoLBcd & \lon & \cE nd_V\vspace{1mm}\\
            &\Ba{c}\resizebox{11mm}{!}{ \xy
(0,4.5)*+{...},
(0,-4.5)*+{...},
(0,0)*{\bu}="o",
(-5,5)*{}="1",
(-3,5)*{}="2",
(3,5)*{}="3",
(5,5)*{}="4",
(-3,-5)*{}="5",
(3,-5)*{}="6",
(5,-5)*{}="7",
(-5,-5)*{}="8",
(-5.5,7)*{_1},
(-3,7)*{_2},
(3,6)*{},
(5.9,7)*{m},
(-3,-7)*{_2},
(3,-7)*+{},
(5.9,-7)*{n},
(-5.5,-7)*{_1},
\ar @{-} "o";"1" <0pt>
\ar @{-} "o";"2" <0pt>
\ar @{-} "o";"3" <0pt>
\ar @{-} "o";"4" <0pt>
\ar @{-} "o";"5" <0pt>
\ar @{-} "o";"6" <0pt>
\ar @{-} "o";"7" <0pt>
\ar @{-} "o";"8" <0pt>
\endxy}\Ea &\lon & \pi_n^m  \vspace{1mm}\\
&
\Ba{c}\resizebox{12mm}{!}{\begin{xy}
 <0mm,0mm>*{\bu};<0mm,0mm>*{}**@{},
 <-0.6mm,0.44mm>*{};<-8mm,5mm>*{}**@{-},
 <-0.4mm,0.7mm>*{};<-4.5mm,5mm>*{}**@{-},
 <0mm,5mm>*{\ldots},
 <0.4mm,0.7mm>*{};<4.5mm,5mm>*{}**@{-},
 <0.6mm,0.44mm>*{};<8mm,5mm>*{}**@{-},
     <0mm,7mm>*{\overbrace{\ \ \ \ \ \ \ \ \ \ \ \ \ \ \ \ }},
     <0mm,9mm>*{^m},
 \end{xy}}\Ea
&\lon&
\ga_m
\Ea
$$
of the twisted prop $\TW\HoLBcd$
in $V=\text{span}\langle e_\al\rangle$ equipped with the deformed differential
\Beq\label{4: differential twisted in V}
d_\centerdot=d + \sum_{k\geq 1} \pm e_\be \frac{1}{k!}\frac{\p^{k+2} \pi}{\p p_\be\p x^{\al_0}x^{\al_1}\ldots \p x^{\al_k}}_{x=p=0}
\frac{\p \ga(p)}{\p p_{\al_1}}|_{p=0}\ldots \frac{\p \ga(p)}{\p p_{\al_k}}|_{p=0} \frac{\p}{\p e_{\al_0}}
\Eeq
The associated twisted  $\HoLBcd$ structure on $V$ is given explicitly by (cf.\ (\ref{4: map frpm HoLBcd to TwP}))
$$
\pi^{\Tw}%=\pi^\ga(x,p) - \pi^\ga(0,p) \ \ \text{with}\ \pi^\ga
:=\hspace{-1mm}\sum_{m,n\geq 1\atop \al_\bu, \be_\bu} \frac{1}{m!n!} \pi_{\al_1\ldots \al_n}^{\be_1\ldots\be_m} p_{\be_1}\ldots p_{\be_m} (x^{\al_1}+ \frac{\p \ga}{\p p_{\al_1}}) \ldots (x^{\al_n}+ \frac{\p \ga}{\p p_{\al_n}})
%
%\sum_{n=1}^N \pi^{a_1\ldots a_n}(x)(p_{a_1}+ \frac{\p \ga}{\p x^{a_1}}) \cdots (p_{a_n} %+\frac{\p \ga}{\p x^{a_1}})
$$
The MC equation for $\ga$ ensures that  $\pi^\Tw|_{x=0}=0$. As $\pi^{\Tw}$ is produced from $\pi$ by the change of variables $x^{\al} \rar x^{\al} + \frac{\p \ga(p)}{\p p_{\al}}$, it is easy to check --- using the vanishing of the sum
$$
\sum \pm  \frac{\p^2 \ga(p)}{\p p_{\al}\p p_{\be}})\frac{\p \pi(x,p)}{\p x^{\al}}\frac{\p \pi(x,p)}{\p x^{\be}}\equiv 0
$$
solely for degree+symmetry reasons --- that
 that the equation
$
\{\pi^{\Tw}, \pi^{\Tw}\}=0$ holds true indeed. Finally, one notices that the $(1,1)$ summand in $\pi^\ga$ (which is responsible for the differential on $V$) is precisely the twisted differential (\ref{4: differential twisted in V}) or, equivalently,  $d + \sum_{k\geq 2}\hat{\mu}_{k,1}$ in the notation of {\S \ref{4: subsec on  MC elements of HoLBcd})}. This gives a short and independent ``local coordinate" check of many ``properadic" claims made above.
%Thus we get a $\ga$-twisted $\HoLBcd$-algebra on the complex $(V, d_\centerdot)$.

% {\tiny A generic representation of $Tw(\HoLBcd)$ in a dg vector space $V$ is given by a %collection of linear maps
%$$
%\left\{\pi_n^m: \odot^{\geq 1}(V[d]) \lon  \odot^{m}(V[-c])[1+c+d],\ \ \ga_m \in %\odot^m(V[-c]))\ \text{with}\ |\ga_m|=c+d-mc
%\right\}_{m,n\geq 1}
%$$
%such that....}

\subsection{Homotopy triangular Lie bialgebras and Lie trialgebras}\label{4: subsec on trinagular} Assume $V$ is a vector space concentrated in degree zero, say, $V=\K^N$ for some $N\in \N$, and let $\HoLB_{1,1}$
be the prop of ordinary Lie bialgebras. Any representation $\rho$ of  $\TW\HoLB_{1,1}$ in $V$ is uniquely determined by its values on generators of cohomological degree zero only, i.e.\ only on the following three generators of $\TW\HoLB_{1,1}$,
$$
\rho\, (\hspace{-3mm}
\Ba{c}\resizebox{7mm}{!}{  \xy
(-5,-5)*+{}="1";
    (-5,0)*{\bu}="L";
  (-8,4.5)*+{}="C";
   (-2,4.5)*+{}="D";
\ar @{-} "D";"L" <0pt>
\ar @{-} "C";"L" <0pt>
\ar @{-} "1";"L" <0pt>
 \endxy}
 \Ea\hspace{-2mm}
 ): V\rar \wedge^2 V,
 \ \ \ \
\rho\, (\hspace{-3mm}
\Ba{c}\resizebox{7mm}{!}{  \xy
(-5,5)*+{}="1";
    (-5,0)*{\bu}="L";
  (-8,-4.5)*+{}="C";
   (-2,-4.5)*+{}="D";
\ar @{-} "D";"L" <0pt>
\ar @{-} "C";"L" <0pt>
\ar @{-} "1";"L" <0pt>
 \endxy}
 \Ea\hspace{-2mm}
 ): \wedge^2 V \rar V,
 \ \ \ \
 \rho\, (\hspace{-1mm}
 \Ba{c}\resizebox{3.8mm}{!}{\begin{xy}
 <0mm,0.5mm>*{};<-3mm,5mm>*{}**@{-},
 <0mm,0.5mm>*{};<3mm,5mm>*{}**@{-},
 <0mm,0mm>*{\bu};
 \end{xy}}\Ea \hspace{-1mm} ) \in \wedge^2V
$$
which satisfy the standard relations (\ref{3: R for LieB}) as well as the following one,
\Beq\label{3: Tw(LB) relation}
 \Ba{c}\resizebox{8mm}{!}{  \xy
(-5,-6)*{\bu}="1";
    (-5,0)*{\bu}="L";
      (1,-0)*+{_3}="R";
  (-8,5)*{^1}="C";
   (-2,5)*{^2}="D";
\ar @{-} "D";"L" <0pt>
\ar @{-} "C";"L" <0pt>
\ar @{-} "1";"L" <0pt>
\ar @{-} "1";"R" <0pt>
 \endxy}
 \Ea
 +
  \Ba{c}\resizebox{8mm}{!}{  \xy
(-5,-6)*{\bu}="1";
    (-5,0)*{\bu}="L";
      (1,-0)*+{_2}="R";
  (-8,5)*{^3}="C";
   (-2,5)*{^1}="D";
\ar @{-} "D";"L" <0pt>
\ar @{-} "C";"L" <0pt>
\ar @{-} "1";"L" <0pt>
\ar @{-} "1";"R" <0pt>
 \endxy}
 \Ea
 +
  \Ba{c}\resizebox{8mm}{!}{  \xy
(-5,-6)*{\bu}="1";
    (-5,0)*{\bu}="L";
      (1,-0)*+{_1}="R";
  (-8,5)*{^2}="C";
   (-2,5)*{^3}="D";
\ar @{-} "D";"L" <0pt>
\ar @{-} "C";"L" <0pt>
\ar @{-} "1";"L" <0pt>
\ar @{-} "1";"R" <0pt>
 \endxy}
 \Ea
 +
 (-1)^c
 \left(
 \Ba{c}\resizebox{10mm}{!}{  \xy
(-5,6)*{^1}="1";
    (-5,+1)*{\bu}="L";
  (-8,-5)*{\bu}="C";
   (-2,-5)*{\bu}="D";
   (-12,+1)*{^2}="l";
   (+2,+1)*{^3}="r";
\ar @{-} "D";"L" <0pt>
\ar @{-} "D";"r" <0pt>
\ar @{-} "C";"L" <0pt>
\ar @{-} "C";"l" <0pt>
\ar @{-} "1";"L" <0pt>
%
 %
 %\ar @{-} "l1";"L" <0pt>
 \endxy}
 \Ea
  +
 \Ba{c}\resizebox{10mm}{!}{  \xy
(-5,6)*{^1}="1";
    (-5,+1)*{\bu}="L";
  (-8,-5)*{\bu}="C";
   (-2,-5)*{\bu}="D";
   (-12,+1)*{^2}="l";
   (+2,+1)*{^3}="r";
\ar @{-} "D";"L" <0pt>
\ar @{-} "D";"r" <0pt>
\ar @{-} "C";"L" <0pt>
\ar @{-} "C";"l" <0pt>
\ar @{-} "1";"L" <0pt>
%
 %
 %\ar @{-} "l1";"L" <0pt>
 \endxy}
 \Ea
  +
 \Ba{c}\resizebox{10mm}{!}{  \xy
(-5,6)*{^3}="1";
    (-5,+1)*{\bu}="L";
  (-8,-5)*{\bu}="C";
   (-2,-5)*{\bu}="D";
   (-12,+1)*{^1}="l";
   (+2,+1)*{^2}="r";
\ar @{-} "D";"L" <0pt>
\ar @{-} "D";"r" <0pt>
\ar @{-} "C";"L" <0pt>
\ar @{-} "C";"l" <0pt>
\ar @{-} "1";"L" <0pt>
%
 %
 %\ar @{-} "l1";"L" <0pt>
 \endxy}
 \Ea
 \right)=0.
\Eeq

If $\rho\, (\hspace{-3mm}
\Ba{c}\resizebox{7mm}{!}{  \xy
(-5,5)*+{}="1";
    (-5,0)*{\bu}="L";
  (-8,-4.5)*+{}="C";
   (-2,-4.5)*+{}="D";
\ar @{-} "D";"L" <0pt>
\ar @{-} "C";"L" <0pt>
\ar @{-} "1";"L" <0pt>
 \endxy}
 \Ea\hspace{-2mm}
 )$ happens to be zero, then the new relation (\ref{3: Tw(LB) relation}) reduced to the classical Yang-Baxter equation so that  associated $\TW\HoLB_{1,1}$-algebra structure in $V$ becomes precisely a so called {\it  triangular Lie bialgebra structure}\, on $V$ \cite{Dr1}. Thus a generic $\TW\HoLB_{1,1}$-algebra structure on $\K^N$  is a version of that notion in which  $V$ has two Lie bialgebra structures, one is given by  the pair
 $\rho\, (\hspace{-3mm}
\Ba{c}\resizebox{7mm}{!}{  \xy
(-5,5)*+{}="1";
    (-5,0)*{\bu}="L";
  (-8,-4.5)*+{}="C";
   (-2,-4.5)*+{}="D";
\ar @{-} "D";"L" <0pt>
\ar @{-} "C";"L" <0pt>
\ar @{-} "1";"L" <0pt>
 \endxy}
 \Ea\hspace{-2mm}
 )$ and $\rho\, (\hspace{-3mm}
\Ba{c}\resizebox{7mm}{!}{  \xy
(-5,-5)*+{}="1";
    (-5,0)*{\bu}="L";
  (-8,4.5)*+{}="C";
   (-2,4.5)*+{}="D";
\ar @{-} "D";"L" <0pt>
\ar @{-} "C";"L" <0pt>
\ar @{-} "1";"L" <0pt>
 \endxy}
 \Ea\hspace{-2mm}
 )$
 and one is given by a pair
 $$
\rho \left(\hspace{-2mm}
\Ba{c}\resizebox{8mm}{!}{  \xy
(-5,4)*{}="1";
    (-5,0)*{\bu}="L";
  (-8,-5)*+{_1}="C";
   (-2,-5)*+{_2}="D";
\ar @{-} "D";"L" <0pt>
\ar @{-} "C";"L" <0pt>
\ar @{-} "1";"L" <0pt>
 \endxy}
 \Ea\hspace{-2mm}
\right) \ \ \ \text{and}\ \ \
\rho\left(\hspace{-2mm}
\Ba{c}\resizebox{8mm}{!}{  \xy
(-5,-5)*+{_1}="1";
    (-5,0)*{\bu}="L";
  (-8,4.5)*+{_1}="C";
   (-2,4.5)*+{_2}="D";
\ar @{-} "D";"L" <0pt>
\ar @{-} "C";"L" <0pt>
\ar @{-} "1";"L" <0pt>
 \endxy}
 \Ea
 +
   \Ba{c}\resizebox{12.5mm}{!}{  \xy
(-18,8)*{^1}="1";
    (-18,+2.5)*{\bu}="L";
 (-14,-2.5)*{\bu}="B";
 (-9,4)*+{^2}="b1";
  (-23,-4)*+{_1}="C";
\ar @{-} "C";"L" <0pt>
\ar @{-} "B";"L" <0pt>
\ar @{-} "B";"b1" <0pt>
\ar @{-} "1";"L" <0pt>
 \endxy}
 \Ea
 -
    \Ba{c}\resizebox{12.5mm}{!}{  \xy
(-18,8)*{^2}="1";
    (-18,+2.5)*{\bu}="L";
 (-14,-2.5)*{\bu}="B";
 (-9,4)*+{^1}="b1";
  (-23,-4)*+{_1}="C";
\ar @{-} "C";"L" <0pt>
\ar @{-} "B";"L" <0pt>
\ar @{-} "B";"b1" <0pt>
\ar @{-} "1";"L" <0pt>
 \endxy}
 \Ea
 \hspace{-1.5mm}\right)
 $$
in which the Lie cobracket is twisted by the coboundary term.

 \sip

 Motivated by the above observation, we introduce a properad of {\em  Lie trialgebras}\,
$\LBcd^\vee$ which is generated by the $\bS$-bimodule $T=\{T(m,n)\}_{m,n\geq 0}$ with
 all $T(m,n)=0$ except
  $$
T(2,1):=\id_1\ot \sgn_2^{|c|}[c-1]=\mbox{span}\left\langle\hspace{-2mm}
\Ba{c}\resizebox{8mm}{!}{  \xy
(-5,-5)*+{_1}="1";
    (-5,0)*{\bu}="L";
  (-8,4.5)*+{_1}="C";
   (-2,4.5)*+{_2}="D";
\ar @{-} "D";"L" <0pt>
\ar @{-} "C";"L" <0pt>
\ar @{-} "1";"L" <0pt>
 \endxy}
 \Ea
 =(-1)^c
 \Ba{c}\resizebox{8mm}{!}{  \xy
(-5,-5)*+{_1}="1";
    (-5,0)*{\bu}="L";
  (-8,4.5)*+{_2}="C";
   (-2,4.5)*+{_1}="D";
\ar @{-} "D";"L" <0pt>
\ar @{-} "C";"L" <0pt>
\ar @{-} "1";"L" <0pt>
 \endxy}
 \Ea
 \hspace{-1mm}  \right\rangle
$$
$$
T(1,2):= \sgn_2^{|d|}\ot \id_1[d-1]=\mbox{span}\left\langle\hspace{-2mm}
 \Ba{c}\resizebox{8mm}{!}{  \xy
(-5,4)*{}="1";
    (-5,0)*{\bu}="L";
  (-8,-5)*+{_1}="C";
   (-2,-5)*+{_2}="D";
\ar @{-} "D";"L" <0pt>
\ar @{-} "C";"L" <0pt>
\ar @{-} "1";"L" <0pt>
 \endxy}
 \Ea
 = (-1)^d
 \Ba{c}\resizebox{8mm}{!}{  \xy
(-5,4)*{}="1";
    (-5,0)*{\bu}="L";
  (-8,-5)*+{_2}="C";
   (-2,-5)*+{_1}="D";
\ar @{-} "D";"L" <0pt>
\ar @{-} "C";"L" <0pt>
\ar @{-} "1";"L" <0pt>
 \endxy}
 \Ea\hspace{-1mm}
\right\rangle
 $$
   $$
T(2,0):=\sgn_2^{|c|}[c-d]=\mbox{span}\left\langle
  \Ba{c}\resizebox{5.0mm}{!}{\begin{xy}
 <0mm,0.5mm>*{};<-3mm,6mm>*{^1}**@{-},
 <0mm,0.5mm>*{};<3mm,6mm>*{^2}**@{-},
 <0mm,0mm>*{\bu};
 \end{xy}}\Ea
 =(-1)^c
  \Ba{c}\resizebox{5.8mm}{!}{\begin{xy}
 <0mm,0.5mm>*{};<-3mm,6mm>*{^2}**@{-},
 <0mm,0.5mm>*{};<3mm,6mm>*{^1}**@{-},
 <0mm,0mm>*{\bu};
 \end{xy}}\Ea
   \right\rangle
$$
modulo relations (\ref{3: R for LieB}) and
(\ref{3: Tw(LB) relation}). This properad comes equipped with two morphisms from $\LBcd$, the one which is identity on the generators of $\LBcd$ and the twisted one  given by
\Beq\label{4: map LBcd to TwLBcd}
\Ba{c}\resizebox{8mm}{!}{  \xy
(-5,4)*{}="1";
    (-5,0)*{\bu}="L";
  (-8,-5)*+{_1}="C";
   (-2,-5)*+{_2}="D";
\ar @{-} "D";"L" <0pt>
\ar @{-} "C";"L" <0pt>
\ar @{-} "1";"L" <0pt>
 \endxy}
 \Ea
\rar
 \Ba{c}\resizebox{8mm}{!}{  \xy
(-5,4)*{}="1";
    (-5,0)*{\bu}="L";
  (-8,-5)*+{_1}="C";
   (-2,-5)*+{_2}="D";
\ar @{-} "D";"L" <0pt>
\ar @{-} "C";"L" <0pt>
\ar @{-} "1";"L" <0pt>
 \endxy}
 \Ea
, \ \ \ \
 \Ba{c}\resizebox{8mm}{!}{  \xy
(-5,-5)*+{_1}="1";
    (-5,0)*{\bu}="L";
  (-8,4.7)*+{_1}="C";
   (-2,4.7)*+{_2}="D";
\ar @{-} "D";"L" <0pt>
\ar @{-} "C";"L" <0pt>
\ar @{-} "1";"L" <0pt>
 \endxy}
 \Ea
 \rar
 \Ba{c}\resizebox{8mm}{!}{  \xy
(-5,-5)*+{_1}="1";
    (-5,0)*{\bu}="L";
  (-8,4.7)*+{_1}="C";
   (-2,4.7)*+{_2}="D";
\ar @{-} "D";"L" <0pt>
\ar @{-} "C";"L" <0pt>
\ar @{-} "1";"L" <0pt>
 \endxy}
 \Ea
 +
   \Ba{c}\resizebox{12.5mm}{!}{  \xy
(-18,8)*{^1}="1";
    (-18,+2.5)*{\bu}="L";
 (-14,-2.5)*{\bu}="B";
 (-9,4)*+{^2}="b1";
  (-23,-4)*+{_1}="C";
\ar @{-} "C";"L" <0pt>
\ar @{-} "B";"L" <0pt>
\ar @{-} "B";"b1" <0pt>
\ar @{-} "1";"L" <0pt>
 \endxy}
 \Ea
 + (-1)^c
    \Ba{c}\resizebox{12.5mm}{!}{  \xy
(-18,8)*{^2}="1";
    (-18,+2.5)*{\bu}="L";
 (-14,-2.5)*{\bu}="B";
 (-9,4)*+{^1}="b1";
  (-23,-4)*+{_1}="C";
\ar @{-} "C";"L" <0pt>
\ar @{-} "B";"L" <0pt>
\ar @{-} "B";"b1" <0pt>
\ar @{-} "1";"L" <0pt>
 \endxy}
 \Ea
\Eeq
The full twisting construction gives us a minimal resolution of $\LBcd^\vee$ as follows. Consider a quotient dg properad
$$
\HoLBcd^\vee :=\Tw\HoLBcd/ \langle  \Ba{c}\resizebox{1.6mm}{!}{\begin{xy}
 <0mm,0.5mm>*{};<0mm,4mm>*{}**@{-},
 <0mm,0mm>*{_\bu};
 \end{xy}}\Ea\rangle
$$
by the ideal generated by the univalent MC generator.

\subsubsection{\bf Theorem} {\em The canonical projection $\HoLBcd^\vee \rar \LBcd^\vee$ is a quasi-isomorphism.}
\begin{proof} Consider a filtration of $\HoLBcd^\vee$ by the number of MC generators. The differential $d$ in the associated graded complex $gr \HoLBcd^\vee$ acts on the generators coming from $
\HoLBcd$ by the standard formula (\ref{3: d in qLBcd_infty}) while on the MC generators by
$$
d  \Ba{c}\resizebox{5.0mm}{!}{\begin{xy}
 <0mm,0.5mm>*{};<-3mm,6mm>*{^1}**@{-},
 <0mm,0.5mm>*{};<3mm,6mm>*{^2}**@{-},
 <0mm,0mm>*{\bu};
 \end{xy}}\Ea=0, \ \ \  d
 \Ba{c}\resizebox{13mm}{!}{
 \begin{xy}
 <0mm,-1mm>*{\bu};
 <0mm,0mm>*{};<-8mm,5mm>*{}**@{-},
 <0mm,0mm>*{};<-4.5mm,5mm>*{}**@{-},
 <0mm,0mm>*{};<-1mm,5mm>*{\ldots}**@{},
 <0mm,0mm>*{};<4.5mm,5mm>*{}**@{-},
 <0mm,0mm>*{};<8mm,5mm>*{}**@{-},
   <0mm,0mm>*{};<-8.5mm,5.5mm>*{^1}**@{},
   <0mm,0mm>*{};<-5mm,5.5mm>*{^2}**@{},
   <0mm,0mm>*{};<4.5mm,5.5mm>*{^{m\hspace{-0.5mm}-\hspace{-0.5mm}1}}**@{},
   <0mm,0mm>*{};<9.0mm,5.5mm>*{^m}**@{},
 \end{xy}}\Ea:=
 - \sum_{[m]=[m_0]\sqcup [m_1]\atop {\# m_0= 2, \# m_1\geq 1}}
 \Ba{c}\resizebox{11mm}{!}{  \xy
(-27,8)*{}="1";
(-22,8)*{}="3";
(-18,8)*{}="n11";
(-13,8)*{}="n12";
(-15.6,7.1)*{...};
(-24.9,7.2)*{...};
    (-25,+2)*{\bu}="L";
 (-15.5,-5)*{_\bu}="B";
    <-25mm,10.6mm>*{\overbrace{ \ \ \ \ \ \ }^{m_0}},
    <-16mm,10.6mm>*{\overbrace{ \ \ \ \ }^{m_1}},
\ar @{-} "B";"L" <0pt>
\ar @{-} "1";"L" <0pt>
\ar @{-} "3";"L" <0pt>
\ar @{-} "1";"L" <0pt>
\ar @{-} "n11";"B" <0pt>
\ar @{-} "n12";"B" <0pt>
 \endxy}
 \Ea \ \ \ \forall m\geq 3.
$$
Since the number of the MC generators is preserved, we can assume that they are distinguished, say, labelled by integers. Then the direct summand of $gr \HoLBcd^\vee$
with, say, $k$ MC generators (labelled by integers from $[k]$) can by identified with a direct summand in $\HoLBcd$ whose first in-legs (labelled by integers from $[k]$) are attached to ``operadic type" $(m_i,1)$-corollas with $m_i\geq 2$, $i\in [k]$. The cohomology of this summand is spanned by trivalent corollas only; trivalent $(2,1)$ corollas whose unique in-legs are labelled by integers from $[k]$  correspond in this approach precisely to $k$ copies of the MC generator $\Ba{c}\resizebox{5.0mm}{!}{\begin{xy}
 <0mm,0.5mm>*{};<-3mm,6mm>*{^1}**@{-},
 <0mm,0.5mm>*{};<3mm,6mm>*{^2}**@{-},
 <0mm,0mm>*{\bu};
 \end{xy}}\Ea$. This result proves the claim.
\end{proof}

\subsubsection{\bf Properad of triangular Lie bialgebras and its minimal resolution} Triangular Lie bialgebras appear naturally in the representation theory of the twisted properad $\TW \HoLBcd$. Consider a quotient properad
%$\LB^\triangle_{c,d}$ be a prop(erad) defined, for any pair of integer  $c,d\in \Z$, is %defined as the quotient
$$
\LB^\triangle_{c,d}:= \LB^\vee_{c,d}/I
$$
of the defined above properad $\LB^\vee_{c,d}$ by the ideal $I$ generated by
the coLie corolla $\hspace{-3mm}
\Ba{c}\resizebox{7mm}{!}{  \xy
(-5,-4)*+{}="1";
    (-5,0)*{\bu}="L";
  (-8,3.5)*+{}="C";
   (-2,3.5)*+{}="D";
\ar @{-} "D";"L" <0pt>
\ar @{-} "C";"L" <0pt>
\ar @{-} "1";"L" <0pt>
 \endxy}
 \Ea\hspace{-2mm}$. Thus $\LB^\triangle_{c,d}$ governs two operations of degrees $1-d$ and $d-c$ respectively,
 $$
 \Ba{c}\resizebox{8mm}{!}{  \xy
(-5,4)*{}="1";
    (-5,0)*{\bu}="L";
  (-8,-5)*+{_1}="C";
   (-2,-5)*+{_2}="D";
\ar @{-} "D";"L" <0pt>
\ar @{-} "C";"L" <0pt>
\ar @{-} "1";"L" <0pt>
 \endxy}
 \Ea
 = (-1)^d
 \Ba{c}\resizebox{8mm}{!}{  \xy
(-5,4)*{}="1";
    (-5,0)*{\bu}="L";
  (-8,-5)*+{_2}="C";
   (-2,-5)*+{_1}="D";
\ar @{-} "D";"L" <0pt>
\ar @{-} "C";"L" <0pt>
\ar @{-} "1";"L" <0pt>
 \endxy}
 \Ea
, \ \ \
  \Ba{c}\resizebox{5.0mm}{!}{\begin{xy}
 <0mm,0.5mm>*{};<-3mm,6mm>*{^1}**@{-},
 <0mm,0.5mm>*{};<3mm,6mm>*{^2}**@{-},
 <0mm,0mm>*{\bu};
 \end{xy}}\Ea
 =(-1)^c
  \Ba{c}\resizebox{5.8mm}{!}{\begin{xy}
 <0mm,0.5mm>*{};<-3mm,6mm>*{^2}**@{-},
 <0mm,0.5mm>*{};<3mm,6mm>*{^1}**@{-},
 <0mm,0mm>*{\bu};
 \end{xy}}\Ea
 $$
%$\cF ree\langle E^\triangle\rangle/\langle\cR^\triangle\rangle$,
%of the free prop generated by an  $\bS$-bimodule $E^\triangle=\{E^\triangle(m,n)\}_{m,n\geq %0}$ with
% all $E^\triangle(m,n)=0$ except
%  $$
%E^\triangle(2,0):=\id_1\ot \sgn_2^{c}[c-1]=\mbox{span}\left\langle
%\right\rangle
%E^\triangle(1,2):= \sgn_2^{d}\ot \id_1[d-1]=\mbox{span}\left\langle
%\right\rangle
%
which are subject to the following relations
\Beq\label{3: R for triangular LieB}
\cR^\triangle:\left\{
\Ba{c}
%%%%%%%%%%%%%% Lie %%%%%%%%%%%%%%%%%%%%%%%%
\Ba{c}\resizebox{8.4mm}{!}{ \begin{xy}
 <0mm,0mm>*{\bu};<0mm,0mm>*{}**@{},
 <0mm,0.69mm>*{};<0mm,3.0mm>*{}**@{-},
 <0.39mm,-0.39mm>*{};<2.4mm,-2.4mm>*{}**@{-},
 <-0.35mm,-0.35mm>*{};<-1.9mm,-1.9mm>*{}**@{-},
 <-2.4mm,-2.4mm>*{\bu};<-2.4mm,-2.4mm>*{}**@{},
 <-2.0mm,-2.8mm>*{};<0mm,-4.9mm>*{}**@{-},
 <-2.8mm,-2.9mm>*{};<-4.7mm,-4.9mm>*{}**@{-},
    <0.39mm,-0.39mm>*{};<3.3mm,-4.0mm>*{^3}**@{},
    <-2.0mm,-2.8mm>*{};<0.5mm,-6.7mm>*{^2}**@{},
    <-2.8mm,-2.9mm>*{};<-5.2mm,-6.7mm>*{^1}**@{},
 \end{xy}}\Ea
 +
\Ba{c}\resizebox{8.4mm}{!}{ \begin{xy}
 <0mm,0mm>*{\bu};<0mm,0mm>*{}**@{},
 <0mm,0.69mm>*{};<0mm,3.0mm>*{}**@{-},
 <0.39mm,-0.39mm>*{};<2.4mm,-2.4mm>*{}**@{-},
 <-0.35mm,-0.35mm>*{};<-1.9mm,-1.9mm>*{}**@{-},
 <-2.4mm,-2.4mm>*{\bu};<-2.4mm,-2.4mm>*{}**@{},
 <-2.0mm,-2.8mm>*{};<0mm,-4.9mm>*{}**@{-},
 <-2.8mm,-2.9mm>*{};<-4.7mm,-4.9mm>*{}**@{-},
    <0.39mm,-0.39mm>*{};<3.3mm,-4.0mm>*{^2}**@{},
    <-2.0mm,-2.8mm>*{};<0.5mm,-6.7mm>*{^1}**@{},
    <-2.8mm,-2.9mm>*{};<-5.2mm,-6.7mm>*{^3}**@{},
 \end{xy}}\Ea
 +
\Ba{c}\resizebox{8.4mm}{!}{ \begin{xy}
 <0mm,0mm>*{\bu};<0mm,0mm>*{}**@{},
 <0mm,0.69mm>*{};<0mm,3.0mm>*{}**@{-},
 <0.39mm,-0.39mm>*{};<2.4mm,-2.4mm>*{}**@{-},
 <-0.35mm,-0.35mm>*{};<-1.9mm,-1.9mm>*{}**@{-},
 <-2.4mm,-2.4mm>*{\bu};<-2.4mm,-2.4mm>*{}**@{},
 <-2.0mm,-2.8mm>*{};<0mm,-4.9mm>*{}**@{-},
 <-2.8mm,-2.9mm>*{};<-4.7mm,-4.9mm>*{}**@{-},
    <0.39mm,-0.39mm>*{};<3.3mm,-4.0mm>*{^1}**@{},
    <-2.0mm,-2.8mm>*{};<0.5mm,-6.7mm>*{^3}**@{},
    <-2.8mm,-2.9mm>*{};<-5.2mm,-6.7mm>*{^2}**@{},
 \end{xy}}\Ea=0,
 \ \ \
  \Ba{c}\resizebox{10mm}{!}{  \xy
(-5,6)*{^1}="1";
    (-5,+1)*{\bu}="L";
  (-8,-5)*{\bu}="C";
   (-2,-5)*{\bu}="D";
   (-12,+1)*{^2}="l";
   (+2,+1)*{^3}="r";
\ar @{-} "D";"L" <0pt>
\ar @{-} "D";"r" <0pt>
\ar @{-} "C";"L" <0pt>
\ar @{-} "C";"l" <0pt>
\ar @{-} "1";"L" <0pt>
%
 %
 %\ar @{-} "l1";"L" <0pt>
 \endxy}
 \Ea
  +
 \Ba{c}\resizebox{10mm}{!}{  \xy
(-5,6)*{^1}="1";
    (-5,+1)*{\bu}="L";
  (-8,-5)*{\bu}="C";
   (-2,-5)*{\bu}="D";
   (-12,+1)*{^2}="l";
   (+2,+1)*{^3}="r";
\ar @{-} "D";"L" <0pt>
\ar @{-} "D";"r" <0pt>
\ar @{-} "C";"L" <0pt>
\ar @{-} "C";"l" <0pt>
\ar @{-} "1";"L" <0pt>
%
 %
 %\ar @{-} "l1";"L" <0pt>
 \endxy}
 \Ea
  +
 \Ba{c}\resizebox{10mm}{!}{  \xy
(-5,6)*{^3}="1";
    (-5,+1)*{\bu}="L";
  (-8,-5)*{\bu}="C";
   (-2,-5)*{\bu}="D";
   (-12,+1)*{^1}="l";
   (+2,+1)*{^2}="r";
\ar @{-} "D";"L" <0pt>
\ar @{-} "D";"r" <0pt>
\ar @{-} "C";"L" <0pt>
\ar @{-} "C";"l" <0pt>
\ar @{-} "1";"L" <0pt>
%
 %
 %\ar @{-} "l1";"L" <0pt>
 \endxy}
 \Ea=0.
\Ea
\right.
\Eeq
Its representations in a dg vector space $V$ are precisely degree shifted triangular Lie bialgebra structures in $V$, the case $c=d=1$ corresponding to he ordinary triangular Lie bialgebras \cite{Dr1}. There is a morphism of properads
$$
f: \LBcd \lon \LBcd^\triangle
$$
given on the generators by
$$
f\, ( \hspace{-2mm} \Ba{c}\resizebox{8mm}{!}{  \xy
(-5,4)*{}="1";
    (-5,0)*{\bu}="L";
  (-8,-5)*+{_1}="C";
   (-2,-5)*+{_2}="D";
\ar @{-} "D";"L" <0pt>
\ar @{-} "C";"L" <0pt>
\ar @{-} "1";"L" <0pt>
 \endxy}
 \Ea  \hspace{-2mm} )
 = \hspace{-2mm}
 \Ba{c}\resizebox{8mm}{!}{  \xy
(-5,4)*{}="1";
    (-5,0)*{\bu}="L";
  (-8,-5)*+{_1}="C";
   (-2,-5)*+{_2}="D";
\ar @{-} "D";"L" <0pt>
\ar @{-} "C";"L" <0pt>
\ar @{-} "1";"L" <0pt>
 \endxy}
 \Ea
, \ \ \ \ \
f(\hspace{-3mm} \Ba{c}\resizebox{7mm}{!}{  \xy
(-5,-5)*+{_1}="1";
    (-5,0)*{\bu}="L";
  (-8,4.7)*+{_1}="C";
   (-2,4.7)*+{_2}="D";
\ar @{-} "D";"L" <0pt>
\ar @{-} "C";"L" <0pt>
\ar @{-} "1";"L" <0pt>
 \endxy}
 \Ea
\hspace{-2mm})
=\hspace{-2mm}
   \Ba{c}\resizebox{12.5mm}{!}{  \xy
(-18,8)*{^1}="1";
    (-18,+2.5)*{\bu}="L";
 (-14,-2.5)*{\bu}="B";
 (-9,4)*+{^2}="b1";
  (-23,-4)*+{_1}="C";
\ar @{-} "C";"L" <0pt>
\ar @{-} "B";"L" <0pt>
\ar @{-} "B";"b1" <0pt>
\ar @{-} "1";"L" <0pt>
 \endxy}
 \Ea
 + (-1)^c
    \Ba{c}\resizebox{12.5mm}{!}{  \xy
(-18,8)*{^2}="1";
    (-18,+2.5)*{\bu}="L";
 (-14,-2.5)*{\bu}="B";
 (-9,4)*+{^1}="b1";
  (-23,-4)*+{_1}="C";
\ar @{-} "C";"L" <0pt>
\ar @{-} "B";"L" <0pt>
\ar @{-} "B";"b1" <0pt>
\ar @{-} "1";"L" <0pt>
 \endxy}
 \Ea
$$
Consider an ideal $I^{\triangle}$ in $\TW \HoLBcd$ generated by all $(m,n)$-corollas
$\Ba{c}\resizebox{10mm}{!}{  \xy
(-7,7)*+{_1}="U1";
(-3,7)*+{_2}="U2";
(2,5)*{...};
(7,7)*+{_m}="U3";
(0,0)*{\circledast}="C";
(-7,-7)*+{_1}="L1";
(-3,-7)*+{_2}="L2";
(2,-5)*{...};
(7,-7)*+{_n}="L3";
\ar @{-} "C";"L1" <0pt>
\ar @{-} "C";"L2" <0pt>
\ar @{-} "C";"L3" <0pt>
\ar @{-} "C";"U1" <0pt>
\ar @{-} "C";"U2" <0pt>
\ar @{-} "C";"U3" <0pt>
 \endxy}
 \Ea$ with $m\geq 2$, $n\geq 1$. This ideal is differential, and, moreover, the quotient
properad $\TW \HoLBcd/I^{\triangle}$
is a dg {\em free}\, properad with generators
$$
\Ba{c}\resizebox{14mm}{!}{ \xy
(1,-5)*{\ldots},
%(-13,-7)*{_1},
(-8,-7)*{_1},
(-3,-7)*{_3},
(8,-7)*{_{n}},
%(13,-7)*{_n},
%
 (0,0)*{\bu}="a",
(0,5)*{}="0",
(-8,-5)*{}="b_2",
(-3,-5)*{}="b_3",
(8,-5)*{}="b_4",
\ar @{-} "a";"0" <0pt>
\ar @{-} "a";"b_2" <0pt>
\ar @{-} "a";"b_3" <0pt>
\ar @{-} "a";"b_4" <0pt>
\endxy}\Ea
\hspace{-3mm}
=(-1)^d \hspace{-2.6mm}
\Ba{c}\resizebox{17mm}{!}{\xy
(1,-6)*{\ldots},
(-7.5,-7)*{_{\sigma(1)}},
(11,-7)*{_{\sigma(n)}},
 (0,0)*{\bu}="a",
(0,5)*{}="0",
(-8,-5)*{}="b_2",
(-3,-5)*{}="b_3",
(8,-5)*{}="b_4",
\ar @{-} "a";"0" <0pt>
\ar @{-} "a";"b_2" <0pt>
\ar @{-} "a";"b_3" <0pt>
\ar @{-} "a";"b_4" <0pt>
\endxy}\Ea \hspace{-3mm}\forall \sigma\in \bS_{n\geq 2}\ \ , \ \
%%%%%%%%%%
\Ba{c}\resizebox{13mm}{!}{\begin{xy}
 <0mm,0mm>*{\bu};
 <-0.6mm,0.44mm>*{};<-8mm,5mm>*{^1}**@{-},
 <-0.4mm,0.7mm>*{};<-4.5mm,5mm>*{^2}**@{-},
 <0mm,5mm>*{\ldots},
 <0.4mm,0.7mm>*{};<4.5mm,5mm>*{^{}}**@{-},
 <0.6mm,0.44mm>*{};<8mm,5mm>*{^m}**@{-},
 \end{xy}}\Ea\hspace{-2mm} =
 (-1)^{c|\tau|}\hspace{-1mm}
 \Ba{c}\resizebox{16mm}{!}{\begin{xy}
 <0mm,0mm>*{\bu};
 <-0.6mm,0.44mm>*{};<-11mm,7mm>*{^{\tau(1)}}**@{-},
 <-0.4mm,0.7mm>*{};<-4.5mm,7mm>*{^{\tau(2)}}**@{-},
 <0mm,5mm>*{\ldots},
 <0.4mm,0.7mm>*{};<4.5mm,6mm>*{}**@{-},
 <0.6mm,0.44mm>*{};<10mm,7mm>*{^{\tau(m)}}**@{-},
 \end{xy}}\Ea  \forall \tau\in \bS_{m\geq 1},
$$
of degrees $1+d-nd$ and $d+c-mc$ respectively. This is an extension of $\Tw\Holie_d$ by the MC $(m,0)$-generators with $m\geq 2$. The induced differential acts on the unique $(1,2)$- and $(2,0)$-generators trivially, i.e.\ they are cohomology classes. In fact every cohomology class in $\TW \HoLBcd/I^{\triangle}$ is generated by this pair via properadic compositions. Indeed, consider a further quotient of $\TW \HoLBcd/I^{\triangle}$ by the  ideal generated by  $\Ba{c}\resizebox{2.5mm}{!}{ \xy
 (0,0)*{\bu}="a",
(0,3.5)*{}="0",
\ar @{-} "a";"0" <0pt>
\endxy}\Ea$ and denote that quotient by $\HoLBcd^\triangle$; notice that the induced differential in $\HoLBcd^\triangle$ is much simplified: this properad is freely generated by the operad $\Holie_d$ equipped with standard differential (\ref{2: d in Lie_infty}) and the MC elements $\Ba{c}\resizebox{12.5mm}{!}{\begin{xy}
 <0mm,0mm>*{\bu};
 <-0.6mm,0.44mm>*{};<-8mm,5mm>*{^1}**@{-},
 <-0.4mm,0.7mm>*{};<-4.5mm,5mm>*{^2}**@{-},
 <0mm,5mm>*{\ldots},
 <0.4mm,0.7mm>*{};<4.5mm,5mm>*{^{}}**@{-},
 <0.6mm,0.44mm>*{};<8mm,5mm>*{^m}**@{-},
 \end{xy}}\Ea$, $m\geq 2$ with the differential given by
 $$
  \pc  \Ba{c}\resizebox{5.8mm}{!}{\begin{xy}
 <0mm,0.5mm>*{};<-3mm,6mm>*{^1}**@{-},
 <0mm,0.5mm>*{};<3mm,6mm>*{^2}**@{-},
 <0mm,0mm>*{\bu};
 \end{xy}}\Ea
 =0
 , \ \ \
 \sd_\centerdot\hspace{-2mm}
 \Ba{c}\resizebox{14mm}{!}{
 \begin{xy}
 <0mm,-1mm>*{\bu};
 <0mm,0mm>*{};<-8mm,5mm>*{}**@{-},
 <0mm,0mm>*{};<-4.5mm,5mm>*{}**@{-},
 <0mm,0mm>*{};<-1mm,5mm>*{\ldots}**@{},
 <0mm,0mm>*{};<4.5mm,5mm>*{}**@{-},
 <0mm,0mm>*{};<8mm,5mm>*{}**@{-},
   <0mm,0mm>*{};<-8.5mm,5.5mm>*{^1}**@{},
   <0mm,0mm>*{};<-5mm,5.5mm>*{^2}**@{},
   %<0mm,0mm>*{};<4.5mm,5.5mm>*{^{m\hspace{-0.5mm}-\hspace{-0.5mm}1}}**@{},
   <0mm,0mm>*{};<9.0mm,5.5mm>*{^m}**@{},
 \end{xy}}\Ea=
- \hspace{-2mm} \sum_{k\geq 2, [m]=\sqcup  [m_\bu],\atop {  m_{0}=1, m_1,...,m_k\geq 1}}\frac{1}{k!}
 \Ba{c}\resizebox{21mm}{!}{  \xy
(-25,8)*{}="1";
(-18,8)*{}="n11";
(-15,7)*{...};
(-13,8)*{}="n12";
(-11,8)*{}="n21";
(-8.4,7)*{...};
(-6,8)*{}="n22";
(1,8)*{}="nn1";
(3,7)*{...};
(5,8)*{}="nn2";
    (-25,+2)*{\bu}="L";
 (-14,-5)*{\bu}="B";
  (-8,-5)*{\bu}="C";
   (3,-5)*{\bu}="D";
   <-5mm,-10mm>*{\underbrace{ \ \ \ \ \ \ \ \ \ \ \ \ \ \ \ \ \ }_{k}},
   (-3,-5)*{...};
    <-25mm,10.6mm>*{\overbrace{ \ \ \ \ \ \ }^{m_0}},
    <-16mm,10.6mm>*{\overbrace{ \ \ \ \ }^{m_1}},
    <-9mm,10.6mm>*{\overbrace{ \ \ \ \ }^{m_2}},
    <3mm,10.6mm>*{\overbrace{ \ \ \ \ }^{m_n}},
     %<0mm,12mm>*{_m},
\ar @{-} "D";"L" <0pt>
\ar @{-} "C";"L" <0pt>
\ar @{-} "B";"L" <0pt>
\ar @{-} "1";"L" <0pt>
\ar @{-} "n11";"B" <0pt>
\ar @{-} "n12";"B" <0pt>
\ar @{-} "n21";"C" <0pt>
\ar @{-} "n22";"C" <0pt>
\ar @{-} "nn1";"D" <0pt>
\ar @{-} "nn2";"D" <0pt>
 \endxy}
 \Ea \ \ \forall\ m\geq 3.
$$
% i.e.
%the canonical projection
%$
%\HoLBcd^\triangle \rar \LBcd^\triangle
%$
%is a quasi-isomorphism.
Let show that the projection
$$
\TW \HoLBcd/I^{\triangle} \lon \HoLBcd^\triangle
$$
is a quasi-isomorphism. Consider a filtration of both sides by the number of MC $(m,0)$ generators with $m\geq 2$ (this number can not decreese). The induced differential on the asscoaied graded of the r.h.s.\ acts only on $\Holie_d$ generators so that its cohomology is a properad, say $P$, generated freely by $\Lie_d$ and the MC generators with $m\geq 2$. On the other hand, the associated graded of the l.h.s.\ is isomorphic to tensor products (modulo the action of finite permutation roups) of complexes $\tw\Holie_d$ and the trivial complex spanned by $(m,0)$ generators with $m\geq 2$. According to \cite{DW}, the $H^\bu(\tw\Holie_d)=\Lie_d$ so that on the l.h.s.\ we get the same properad $P$. Thus at the second page of the spectral sequence the above map becomes the identity map implying, by the Comparison Theorem, the resquired quasi-isomorphism.

\sip

Using Gr\"obner basis techniques A.\ Khoroshkin has proven \cite{Kh} that $\HoLBcd^\triangle$ is a minimal resolution of $\LBcd^\triangle$ at the dioperadic level. perhaps this result holds true at the properadic level as well.
%
%We can summarize our discussion in the following
%\subsubsection{\bf Proposition} {\em The canonical epimorphism
%$
%\TW \HoLBcd/I^{\triangle} \lon \LBcd^\triangle.
%$
%is a quasi-isomorphism.}
%\sip
%
Homotopy triangular Lie bialgebras (in a different, non-properadic, context) have been studied in \cite{LST} where their relation to homotopy  Rota-Baxter Lie algebras has been established.

%%%%%%%%%%%%%%%%%%%%%%%%%%%%%%%%%%%%%%%%%%%%%%%%%%%%%%%%%%%%%%%%%%%%%%%%%%%%%%%%%%%%%%%
\subsection{Full twisting of properads under $\LB_{c,d}$} Assume $\cP=\LB_{c,d}$ and the map $i: \HoLB_{c,d}\rar \LBcd$ is the canonical projection. The associated twisted dg prop(erad) $(\TW\LBcd, \delta_\centerdot)$ is generated by
the standard corollas
$$
\Ba{c}\resizebox{8mm}{!}{  \xy
(-5,4)*{}="1";
    (-5,0)*{\bu}="L";
  (-8,-5)*+{_1}="C";
   (-2,-5)*+{_2}="D";
\ar @{-} "D";"L" <0pt>
\ar @{-} "C";"L" <0pt>
\ar @{-} "1";"L" <0pt>
 \endxy}
 \Ea
 = (-1)^d
 \Ba{c}\resizebox{8mm}{!}{  \xy
(-5,4)*{}="1";
    (-5,0)*{\bu}="L";
  (-8,-5)*+{_2}="C";
   (-2,-5)*+{_1}="D";
\ar @{-} "D";"L" <0pt>
\ar @{-} "C";"L" <0pt>
\ar @{-} "1";"L" <0pt>
 \endxy}
 \Ea
, \ \ \ \ \
\Ba{c}\resizebox{8mm}{!}{  \xy
(-5,-5)*+{_1}="1";
    (-5,0)*{\bu}="L";
  (-8,4.5)*+{_1}="C";
   (-2,4.5)*+{_2}="D";
\ar @{-} "D";"L" <0pt>
\ar @{-} "C";"L" <0pt>
\ar @{-} "1";"L" <0pt>
 \endxy}
 \Ea
 =(-1)^c
 \Ba{c}\resizebox{8mm}{!}{  \xy
(-5,-5)*+{_1}="1";
    (-5,0)*{\bu}="L";
  (-8,4.5)*+{_2}="C";
   (-2,4.5)*+{_1}="D";
\ar @{-} "D";"L" <0pt>
\ar @{-} "C";"L" <0pt>
\ar @{-} "1";"L" <0pt>
 \endxy}
 \Ea
$$
of degrees $1-c$ and $1-d$ respectively modulo relations (\ref{3: R for LieB}), as well by the
family of extra generators,
$$
\Ba{c}\resizebox{12.5mm}{!}{\begin{xy}
 <0mm,0mm>*{\bu};
 <-0.6mm,0.44mm>*{};<-8mm,5mm>*{^1}**@{-},
 <-0.4mm,0.7mm>*{};<-4.5mm,5mm>*{^2}**@{-},
 <0mm,5mm>*{\ldots},
 <0.4mm,0.7mm>*{};<4.5mm,5mm>*{^{}}**@{-},
 <0.6mm,0.44mm>*{};<8mm,5mm>*{^m}**@{-},
 \end{xy}}\Ea =
 (-1)^{c|\sigma|}
 \Ba{c}\resizebox{15mm}{!}{\begin{xy}
 <0mm,0mm>*{\bu};
 <-0.6mm,0.44mm>*{};<-11mm,7mm>*{^{\sigma(1)}}**@{-},
 <-0.4mm,0.7mm>*{};<-4.5mm,7mm>*{^{\sigma(2)}}**@{-},
 <0mm,5mm>*{\ldots},
 <0.4mm,0.7mm>*{};<4.5mm,6mm>*{}**@{-},
 <0.6mm,0.44mm>*{};<10mm,7mm>*{^{\sigma(m)}}**@{-},
 \end{xy}}\Ea \ \ \ \forall \sigma\in \bS_m, \ m\geq 1,
$$
of cohomological degree $(1-m)c+d$. The twisted differential $\delta_\centerdot$ is given explicitly on the first pair of generators by
$$
\delta_\centerdot\hspace{-3mm}
\Ba{c}\resizebox{8mm}{!}{  \xy
(-5,6)*{}="1";
    (-5,+1)*{\bu}="L";
  (-8,-5)*+{_1}="C";
   (-2,-5)*+{_2}="D";
\ar @{-} "D";"L" <0pt>
\ar @{-} "C";"L" <0pt>
\ar @{-} "1";"L" <0pt>
 \endxy}
 \Ea =0,
 \ \ \
\delta_\centerdot\hspace{-3mm}
\Ba{c}\resizebox{8mm}{!}{  \xy
(-5,-6)*+{_1}="1";
    (-5,0)*{\bu}="L";
  (-8,5)*+{_1}="C";
   (-2,5)*+{_2}="D";
\ar @{-} "D";"L" <0pt>
\ar @{-} "C";"L" <0pt>
\ar @{-} "1";"L" <0pt>
 \endxy}
 \Ea
\hspace{-2mm}
 =\hspace{-2mm}
  \Ba{c}\resizebox{11mm}{!}{  \xy
(-19,8)*{^1}="1";
    (-19,+3)*{\bu}="L";
 (-14,-2.5)*{\bu}="B";
 (-9,3)*+{^2}="b1";
 (-14,-8)*{\bu}="b2";
  (-23,-3)*+{_1}="C";
\ar @{-} "C";"L" <0pt>
\ar @{-} "B";"L" <0pt>
\ar @{-} "B";"b1" <0pt>
\ar @{-} "B";"b2" <0pt>
\ar @{-} "1";"L" <0pt>
 \endxy}
 \Ea
 +
 (-1)^c\hspace{-2mm}
  \Ba{c}\resizebox{11mm}{!}{  \xy
(-19,8)*{^2}="1";
    (-19,+3)*{\bu}="L";
 (-14,-2.5)*{\bu}="B";
 (-9,3)*+{^1}="b1";
 (-14,-8)*{\bu}="b2";
  (-23,-3)*+{_1}="C";
\ar @{-} "C";"L" <0pt>
\ar @{-} "B";"L" <0pt>
\ar @{-} "B";"b1" <0pt>
\ar @{-} "B";"b2" <0pt>
\ar @{-} "1";"L" <0pt>
 \endxy}
 \Ea.
 %
%
%  \Ba{c}\resizebox{10.5mm}{!}{  \xy
%(-9,8)*{^1}="1";
%    (-9,+3)*{\bu}="L";
% (-14,-3.5)*{\bu}="B";
% (-19,5)*+{^2}="b1";
% (-14,-10)*{\bu}="b2";
%  (-3,-5)*+{_1}="C";
%     %
%\ar @{-} "C";"L" <0pt>
%\ar @{-} "B";"L" <0pt>
%\ar @{-} "B";"b1" <0pt>
%\ar @{-} "B";"b2" <0pt>
%\ar @{-} "1";"L" <0pt>
% \endxy}
% \Ea
 %
 %%%%%%%%%%%%%%%%%%
$$
where we used  relations (\ref{3: R for LieB}) (and an ordering of vertices in the second line just above goes from the bottom to the top).
 The action of $\sd_\centerdot$  on the remaining MC generators is given by
$$
 \p_\centerdot \Ba{c}\resizebox{1.6mm}{!}{\begin{xy}
 <0mm,0.5mm>*{};<0mm,4mm>*{}**@{-},
 <0mm,0mm>*{_\bu};
 \end{xy}}\Ea= +
\frac{1}{2}
\Ba{c}\resizebox{6.1mm}{!}{  \xy
(-5,6)*{}="1";
    (-5,+1)*{\bu}="L";
  (-8,-4)*{\bu}="C";
   (-2,-4)*{\bu}="D";
\ar @{-} "D";"L" <0pt>
\ar @{-} "C";"L" <0pt>
\ar @{-} "1";"L" <0pt>
%
 %
 %\ar @{-} "l1";"L" <0pt>
 \endxy}
 \Ea
 \ \ \ \text{and}\ \ \ \
  \p_\centerdot
  \Ba{c}\resizebox{6.0mm}{!}{\begin{xy}
 <0mm,0.5mm>*{};<-3mm,6mm>*{^1}**@{-},
 <0mm,0.5mm>*{};<3mm,6mm>*{^2}**@{-},
 <0mm,0mm>*{_\bu};
 \end{xy}}\Ea=
  -
 \Ba{c}\resizebox{7.5mm}{!}{  \xy
(-5,-6)*{_\bu}="1";
    (-5,-0.5)*{\bu}="L";
  (-8,5)*+{_1}="C";
   (-2,5)*+{_2}="D";
\ar @{-} "D";"L" <0pt>
\ar @{-} "C";"L" <0pt>
\ar @{-} "1";"L" <0pt>
 \endxy}
 \Ea
 $$
and, for $m\geq 3$, by
% already corrected  BY $\sd_\centerdot$, i.e. some terms already cancelled out,  cf (\ref{3: %d_centerdot on Tw(P)})

$$
 \sd_\centerdot
 \Ba{c}\resizebox{13.5mm}{!}{
 \begin{xy}
 <0mm,-1mm>*{\bu};
 <0mm,0mm>*{};<-8mm,5mm>*{}**@{-},
 <0mm,0mm>*{};<-4.5mm,5mm>*{}**@{-},
 <0mm,0mm>*{};<-1mm,5mm>*{\ldots}**@{},
 <0mm,0mm>*{};<4.5mm,5mm>*{}**@{-},
 <0mm,0mm>*{};<8mm,5mm>*{}**@{-},
   <0mm,0mm>*{};<-8.5mm,5.5mm>*{^1}**@{},
   <0mm,0mm>*{};<-5mm,5.5mm>*{^2}**@{},
   <0mm,0mm>*{};<4.5mm,5.5mm>*{^{m\hspace{-0.5mm}-\hspace{-0.5mm}1}}**@{},
   <0mm,0mm>*{};<9.0mm,5.5mm>*{^m}**@{},
 \end{xy}}\Ea
 =
\sum_{[m]=[m_0]\sqcup [m_1]\atop {\# m_0= 2, \# m_1\geq 1}}(-1)^{1+c\sigma'}
 \Ba{c}\resizebox{11mm}{!}{  \xy
(-27,8)*{}="1";
(-22,8)*{}="3";
(-18,8)*{}="n11";
(-13,8)*{}="n12";
(-15,7)*{...};
    (-25,+2)*{\bu}="L";
 (-15.5,-5)*{\bu}="B";
    <-25mm,10.6mm>*{\overbrace{ \ \ \ \ \ \ }^{m_0}},
    <-16mm,10.6mm>*{\overbrace{ \ \ \ \ }^{m_1}},
\ar @{-} "B";"L" <0pt>
\ar @{-} "1";"L" <0pt>
\ar @{-} "3";"L" <0pt>
\ar @{-} "1";"L" <0pt>
\ar @{-} "n11";"B" <0pt>
\ar @{-} "n12";"B" <0pt>
 \endxy}
 \Ea
-
\sum_{[m]=[m_0]\sqcup [m_1]\sqcup [m_2] \atop {\# m_0=1, \# m_1,\# m_2\geq 1}}\frac{(-1)^{c\sigma''}}{2}
 \Ba{c}\resizebox{17mm}{!}{  \xy
(-27,8)*{}="1";
(-25,8)*{}="2";
(-23,8)*{}="3";
(-19,8)*{}="n11";
(-16.8,7)*{...};
(-14,8)*{}="n12";
(-36,8)*{}="n21";
(-33.4,7)*{...};
(-31,8)*{}="n22";
    (-25,+2)*{\bu}="L";
 (-17,-5)*{\bu}="B";
  (-33,-5)*{\bu}="C";
   %(-3,-5)*{...};
 %
    <-25mm,10.6mm>*{\overbrace{ \ \ \ \ \ \ }^{m_0}},
    <-17mm,10.6mm>*{\overbrace{ \ \ \ \ }^{m_2}},
    <-34mm,10.6mm>*{\overbrace{ \ \ \ \ }^{m_1}},
\ar @{-} "C";"L" <0pt>
\ar @{-} "B";"L" <0pt>
\ar @{-} "2";"L" <0pt>
\ar @{-} "n11";"B" <0pt>
\ar @{-} "n12";"B" <0pt>
\ar @{-} "n21";"C" <0pt>
\ar @{-} "n22";"C" <0pt>
 \endxy}
 \Ea,
$$
where $\sigma'$ (resp. $\sigma''$) is the parity of the permutation $[m]\rar[m_0]\sqcup [m_1]$ (resp.\ $[m]\rar[m_0]\sqcup [m_1]\sqcup [m_2]$) associated with the partition of the ordered set $[m]$ into two (resp.\ three) disjoint ordered subsets.

% (1,1)-terms cancel out.
These relations imply
$$
\sd_\centerdot\left(\hspace{-3mm}
\Ba{c}\resizebox{8mm}{!}{  \xy
(-5,-6)*+{_1}="1";
    (-5,-0.2)*{\bu}="L";
  (-8,5)*+{_1}="C";
   (-2,5)*+{_2}="D";
\ar @{-} "D";"L" <0pt>
\ar @{-} "C";"L" <0pt>
\ar @{-} "1";"L" <0pt>
 \endxy}
 \Ea
 +
   \Ba{c}\resizebox{13.5mm}{!}{  \xy
(-19,9)*{^1}="1";
    (-19,+3)*{\circledcirc}="L";
 (-14,-3.5)*{\bu}="B";
 (-7,5)*+{^2}="b1";
  (-25,-5)*+{_1}="C";
\ar @{-} "C";"L" <0pt>
\ar @{-} "B";"L" <0pt>
\ar @{-} "B";"b1" <0pt>
\ar @{-} "1";"L" <0pt>
 \endxy}
 \Ea
 +
 (-1)^c
   \Ba{c}\resizebox{13.5mm}{!}{  \xy
(-19,9)*{^2}="1";
    (-19,+3)*{\bu}="L";
 (-14,-3.5)*{\bu}="B";
 (-7,5)*+{^1}="b1";
  (-25,-5)*+{_1}="C";
\ar @{-} "C";"L" <0pt>
\ar @{-} "B";"L" <0pt>
\ar @{-} "B";"b1" <0pt>
\ar @{-} "1";"L" <0pt>
 \endxy}
 \Ea
 \right)=0
$$
which is in agreement with the general result saying that $\Tw\LBcd$ is a properad under $\HoLBcd$; the morphism (\ref{4: map frpm HoLBcd to TwP}) takes in this case the following form
 \Beq\label{4: map frpm HoLBcd to TwP for P under LBcd}
\Ba{c}\resizebox{13mm}{!}{\begin{xy}
 <0mm,0mm>*{\bu};<0mm,0mm>*{}**@{},
 <-0.6mm,0.44mm>*{};<-8mm,5mm>*{}**@{-},
 <-0.4mm,0.7mm>*{};<-4.5mm,5mm>*{}**@{-},
 <0mm,0mm>*{};<-1mm,5mm>*{\ldots}**@{},
 <0.4mm,0.7mm>*{};<4.5mm,5mm>*{}**@{-},
 <0.6mm,0.44mm>*{};<8mm,5mm>*{}**@{-},
   <0mm,0mm>*{};<-8.5mm,5.5mm>*{^1}**@{},
   <0mm,0mm>*{};<-5mm,5.5mm>*{^2}**@{},
   <0mm,0mm>*{};<4.5mm,5.5mm>*{^{m\hspace{-0.5mm}-\hspace{-0.5mm}1}}**@{},
   <0mm,0mm>*{};<9.0mm,5.5mm>*{^m}**@{},
 <-0.6mm,-0.44
 mm>*{};<-8mm,-5mm>*{}**@{-},
 <-0.4mm,-0.7mm>*{};<-4.5mm,-5mm>*{}**@{-},
 <0mm,0mm>*{};<-1mm,-5mm>*{\ldots}**@{},
 <0.4mm,-0.7mm>*{};<4.5mm,-5mm>*{}**@{-},
 <0.6mm,-0.44mm>*{};<8mm,-5mm>*{}**@{-},
   <0mm,0mm>*{};<-8.5mm,-6.9mm>*{^1}**@{},
   <0mm,0mm>*{};<-5mm,-6.9mm>*{^2}**@{},
   <0mm,0mm>*{};<4.5mm,-6.9mm>*{^{n\hspace{-0.5mm}-\hspace{-0.5mm}1}}**@{},
   <0mm,0mm>*{};<9.0mm,-6.9mm>*{^n}**@{},
 \end{xy}}\Ea
\lon
\left\{\Ba{cl}
0 & \text{if}\, n\geq 2\ \text{and}\ m+n> 3,\\
\Ba{c}\resizebox{8mm}{!}{  \xy
(-5,6)*{}="1";
    (-5,+1)*{\bu}="L";
  (-8,-5)*+{_1}="C";
   (-2,-5)*+{_2}="D";
\ar @{-} "D";"L" <0pt>
\ar @{-} "C";"L" <0pt>
\ar @{-} "1";"L" <0pt>
 \endxy}
 \Ea  & \text{if}\ n=2,m=1,
 \\
\Ba{c}\resizebox{8mm}{!}{  \xy
(-5,-6)*+{_1}="1";
    (-5,-0.2)*{\bu}="L";
  (-8,5)*+{_1}="C";
   (-2,5)*+{_2}="D";
\ar @{-} "D";"L" <0pt>
\ar @{-} "C";"L" <0pt>
\ar @{-} "1";"L" <0pt>
 \endxy}
 \Ea
 +
   \Ba{c}\resizebox{13.5mm}{!}{  \xy
(-19,9)*{^1}="1";
    (-19,+3)*{\bu}="L";
 (-14,-3.5)*{\bu}="B";
 (-7,5)*+{^2}="b1";
  (-25,-5)*+{_1}="C";
\ar @{-} "C";"L" <0pt>
\ar @{-} "B";"L" <0pt>
\ar @{-} "B";"b1" <0pt>
\ar @{-} "1";"L" <0pt>
 \endxy}
 \Ea
 +
 (-1)^c
   \Ba{c}\resizebox{13.5mm}{!}{  \xy
(-19,9)*{^2}="1";
    (-19,+3)*{\bu}="L";
 (-14,-3.5)*{\bu}="B";
 (-7,5)*+{^1}="b1";
  (-25,-5)*+{_1}="C";
\ar @{-} "C";"L" <0pt>
\ar @{-} "B";"L" <0pt>
\ar @{-} "B";"b1" <0pt>
\ar @{-} "1";"L" <0pt>
 \endxy}
 \Ea    & \text{if}\ n=1,m=2,
\\
  \sum_{[m]=[m_0]\sqcup [m_1]\atop {\# m_0= 1, \# m_1\geq 1}}(-1)^{c\sigma_{m_0,m_1}}
 \Ba{c}\resizebox{11mm}{!}{  \xy
(-25,8)*{}="1";
%(-22,8)*{}="3";
%
(-18,8)*{}="n11";
(-13,8)*{}="n12";
(-15,7)*{...};
    (-25,+2)*{\bu}="L";
 (-15.5,-5)*{\bu}="B";
 (-27,-5)*{}="C";
    <-25mm,10.6mm>*{\overbrace{  \ \ }^{m_0}},
    <-16mm,10.6mm>*{\overbrace{ \ \ \ \ }^{m_1}},
\ar @{-} "B";"L" <0pt>
\ar @{-} "C";"L" <0pt>
\ar @{-} "1";"L" <0pt>
%\ar @{-} "3";"L" <0pt>
\ar @{-} "1";"L" <0pt>
\ar @{-} "n11";"B" <0pt>
\ar @{-} "n12";"B" <0pt>
 \endxy}
 \Ea
 &
\text{if}\ n=1, m\geq 3.
\Ea
\right.
\Eeq

 Note that the quotient of the dg operad $\TW\LBcd$ by the (differential) ideal generated by corollas $\Ba{c}\resizebox{11mm}{!}{
 \begin{xy}
 <0mm,-1mm>*{\bu};
 <0mm,0mm>*{};<-8mm,5mm>*{}**@{-},
 <0mm,0mm>*{};<-4.5mm,5mm>*{}**@{-},
 <0mm,0mm>*{};<-1mm,5mm>*{\ldots}**@{},
 <0mm,0mm>*{};<4.5mm,5mm>*{}**@{-},
 <0mm,0mm>*{};<8mm,5mm>*{}**@{-},
   <0mm,0mm>*{};<-8.5mm,5.5mm>*{^1}**@{},
   <0mm,0mm>*{};<-5mm,5.5mm>*{^2}**@{},
   <0mm,0mm>*{};<4.5mm,5.5mm>*{^{m\hspace{-0.5mm}-\hspace{-0.5mm}1}}**@{},
   <0mm,0mm>*{};<9.0mm,5.5mm>*{^m}**@{},
 \end{xy}}\Ea$ with $m\neq 2$
 gives us precisely the properad of Lie trialgebras $\LBcd^\vee$.
 Theorem {\ref{4: theorem on H(TwHoLBcd)}} impli that the canonical
 projection
 $
 \Tw\LBcd \rar \LBcd
 $
is a quasi-isomorphism.

% \sip

Similarly one can describe explicitly the twisted properad $\Tw\cP$ associated to any properad $\cP\in \PROP_{\LBcd}$.
Note that $\Tw\cP$ is {\em not}\, in general a properad under $\LBcd$ as higher homotopy operations of type $(m\geq 3,1)$ can be non-trivial.

\bip

%%%%%%%%%%%%%%%%%%%%%%%%%%%%%%%%%%%%%%%%%%%%%%%%%%%%%%%%%%%%%%

{\Large
\section{\bf Full twisting endofunctor in the case of involutive Lie bialgebras}
}
\mip

\subsection{Introduction}
This section adopts the full twisting endofunctor $\Tw$ in the category of properads under $\HoLBcd$ to the case when (strongly homotopy) Lie bialgebras satisfy the {\em involutivity or diamond}\, condition (which is often satisfied in applications). The corresponding twisting endofunctor
$$
\Tw^\diamond: \PROP_{\HoLoBcd} \lon \PROP_{\HoLoBcd}
$$
 admits a much shorter and nicer formulation than $\Tw$ due to the equivalence of $\HoLoBcd$-algebra structures and the so called homotopy $\cB\cV^{com}$-structures which are used heavily in the Batalin-Vilkovisky formalism of the mathematical physics and QFT.

\sip

We show many formulae explicitly but omit all the calculations proving them because proofs are much analogous to the ones given in the previous sections, and because this {\em diamond}\, version $\Tw^\diamond$ of $\Tw$ gives us essentially nothing new --- upon representations of the twisted properad $\Tw^\diamond\HoLoBcd$ defined below one recovers the well-known constructions from
\cite{CFL,NW}. Thus the only small novelty is that we arrive to these beautiful results using a different language, the properadic one.

\mip

\subsection{Reminder on involutive Lie bialgebras} Given any pair of integers $c,d$ of the same parity, $c=d\bmod 2\Z$, the properad $\LoBcd$ of {\em involutive}\, Lie bialgebras is defined as the quotient of $\LBcd$ by the ideal generated by the involutivity, or ``diamond", relation
$$
\Ba{c}\resizebox{5.5mm}{!}{
\xy
 (0,0)*{\bu}="a",
(0,6)*{\bu}="b",
(3,3)*{}="c",
(-3,3)*{}="d",
 (0,9)*{}="b'",
(0,-3)*{}="a'",
\ar@{-} "a";"c" <0pt>
\ar @{-} "a";"d" <0pt>
\ar @{-} "a";"a'" <0pt>
\ar @{-} "b";"c" <0pt>
\ar @{-} "b";"d" <0pt>
\ar @{-} "b";"b'" <0pt>
\endxy}
\Ea=0.
$$
Note that this relation is void in $\LBcd$ for $c$ and $d$ of opposite parities.

\sip

It was proven in \cite{CMW} that the minimal resolution $\HoLoBcd$ of the properad $\LoBcd$  is a free properad generated
by the following (skew)symmetric corollas of degree $1+c(1-m-a) + d(1-n-a)$
\Beq\label{5: generators of HoLoBcd}
\Ba{c}\resizebox{13mm}{!}{\xy
(-9,-6)*{};
(0,0)*+{a}*\cir{}
**\dir{-};
(-5,-6)*{};
(0,0)*+{a}*\cir{}
**\dir{-};
(9,-6)*{};
(0,0)*+{a}*\cir{}
**\dir{-};
(5,-6)*{};
(0,0)*+{a}*\cir{}
**\dir{-};
(0,-6)*{\ldots};
(-10,-8)*{_1};
(-6,-8)*{_2};
(10,-8)*{_n};
(-9,6)*{};
(0,0)*+{a}*\cir{}
**\dir{-};
(-5,6)*{};
(0,0)*+{a}*\cir{}
**\dir{-};
(9,6)*{};
(0,0)*+{a}*\cir{}
**\dir{-};
(5,6)*{};
(0,0)*+{a}*\cir{}
**\dir{-};
(0,6)*{\ldots};
(-10,8)*{_1};
(-6,8)*{_2};
(10,8)*{_m};
\endxy}\Ea
=
(-1)^{(d+1)(\sigma+\tau)}
\Ba{c}\resizebox{16mm}{!}{\xy
(-9,-6)*{};
(0,0)*+{a}*\cir{}
**\dir{-};
(-5,-6)*{};
(0,0)*+{a}*\cir{}
**\dir{-};
(9,-6)*{};
(0,0)*+{a}*\cir{}
**\dir{-};
(5,-6)*{};
(0,0)*+{a}*\cir{}
**\dir{-};
(0,-6)*{\ldots};
(-12,-8)*{_{\tau(1)}};
(-6,-8)*{_{\tau(2)}};
(12,-8)*{_{\tau(n)}};
(-9,6)*{};
(0,0)*+{a}*\cir{}
**\dir{-};
(-5,6)*{};
(0,0)*+{a}*\cir{}
**\dir{-};
(9,6)*{};
(0,0)*+{a}*\cir{}
**\dir{-};
(5,6)*{};
(0,0)*+{a}*\cir{}
**\dir{-};
(0,6)*{\ldots};
(-12,8)*{_{\sigma(1)}};
(-6,8)*{_{\sigma(2)}};
(12,8)*{_{\sigma(m)}};
\endxy}\Ea\ \ \ \forall \sigma\in \bS_m, \forall \tau\in \bS_n,
\Eeq
where $m+n+ a\geq 3$, $m\geq 1$, $n\geq 1$, $a\geq 0$. The differential in
$\HoLoBd$ is given on the generators by
\Beq\label{5: d on HoLoBcd}
\delta
\Ba{c}\resizebox{13mm}{!}{\xy
(-9,-6)*{};
(0,0)*+{a}*\cir{}
**\dir{-};
(-5,-6)*{};
(0,0)*+{a}*\cir{}
**\dir{-};
(9,-6)*{};
(0,0)*+{a}*\cir{}
**\dir{-};
(5,-6)*{};
(0,0)*+{a}*\cir{}
**\dir{-};
(0,-6)*{\ldots};
(-10,-8)*{_1};
(-6,-8)*{_2};
(10,-8)*{_n};
(-9,6)*{};
(0,0)*+{a}*\cir{}
**\dir{-};
(-5,6)*{};
(0,0)*+{a}*\cir{}
**\dir{-};
(9,6)*{};
(0,0)*+{a}*\cir{}
**\dir{-};
(5,6)*{};
(0,0)*+{a}*\cir{}
**\dir{-};
(0,6)*{\ldots};
(-10,8)*{_1};
(-6,8)*{_2};
(10,8)*{_m};
\endxy}\Ea
=
\sum_{l\geq 1}\sum_{a=b+c+l-1}\sum_{[m]=I_1\sqcup I_2\atop
[n]=J_1\sqcup J_2} \pm
\Ba{c}
%
%
%%%%%%%%%%%%%%%% two vertex graph with l internal edges %%%%%%%%%%
\Ba{c}\resizebox{18mm}{!}{\xy
(0,0)*+{b}*\cir{}="b",
(10,10)*+{c}*\cir{}="c",
%
%%%%%%%%%% edges to b %%%%%%%%%%%%
(-9,6)*{}="1",
(-7,6)*{}="2",
(-2,6)*{}="3",
(-3.5,5)*{...},
(-4,-6)*{}="-1",
(-2,-6)*{}="-2",
(4,-6)*{}="-3",
(1,-5)*{...},
(0,-8)*{\underbrace{\ \ \ \ \ \ \ \ }},
(0,-11)*{_{J_1}},
(-6,8)*{\overbrace{ \ \ \ \ \ \ }},
(-6,11)*{_{I_1}},
%%%%%%%%%% edges to c %%%%%%%%%%%%
(6,16)*{}="1'",
(8,16)*{}="2'",
(14,16)*{}="3'",
(11,15)*{...},
(11,6)*{}="-1'",
(16,6)*{}="-2'",
(18,6)*{}="-3'",
(13.5,6)*{...},
(15,4)*{\underbrace{\ \ \ \ \ \ \ }},
(15,1)*{_{J_2}},
(10,18)*{\overbrace{ \ \ \ \ \ \ \ \ }},
(10,21)*{_{I_2}},
%
%%%%%%%%%%% internal curved edges %%%%%%%%%%%
(0,2)*-{};(8.0,10.0)*-{}
**\crv{(0,10)};
(0.5,1.8)*-{};(8.5,9.0)*-{}
**\crv{(0.4,7)};
(1.5,0.5)*-{};(9.1,8.5)*-{}
**\crv{(5,1)};
(1.7,0.0)*-{};(9.5,8.6)*-{}
**\crv{(6,-1)};
(5,5)*+{...};
\ar @{-} "b";"1" <0pt>
\ar @{-} "b";"2" <0pt>
\ar @{-} "b";"3" <0pt>
\ar @{-} "b";"-1" <0pt>
\ar @{-} "b";"-2" <0pt>
\ar @{-} "b";"-3" <0pt>
\ar @{-} "c";"1'" <0pt>
\ar @{-} "c";"2'" <0pt>
\ar @{-} "c";"3'" <0pt>
\ar @{-} "c";"-1'" <0pt>
\ar @{-} "c";"-2'" <0pt>
\ar @{-} "c";"-3'" <0pt>
\endxy}\Ea
%%%%%%%%%%%%%%%%%%%%%%%%%%%%%%%%%%%%%%%%%%%
\Ea
\Eeq
where the summation parameter $l$ counts the number of internal edges connecting the two vertices
on the r.h.s., and the signs are  fixed by the fact that they all equal to $-1$ for $c$ and $d$
odd integers.

\sip

The ``plus" extension (see \S {\ref{2: subsection on plus functor}}), $\HoLoBpcd$, of this properad looks especially natural --- one adds just one extra $(1,1)$-generator (which we denote from nowon by $\xy
(0,-4)*{};
(0,0)*+{_0}*\cir{}
**\dir{-};
(0,4)*{};
(0,0)*+{_0}*\cir{}
**\dir{-};
\endxy$) to the list while keeping the differential (\ref{5: d on HoLoBcd}) formally the same.

\sip

Let $\hbar$ be a formal variable of degree $c+d$ and, for a vector space $V$. let $V[[\hbar]]$ stand for the topological vector space of formal power series with coefficients in $V$; it is a module over the topological ring $\K[[\hbar]]$ of formal power series in $\hbar$. Consider a dg properad $\HoLBcd^{\hbar+}$ which is identical to
$\HoLBcd^+[[\hbar]]$ as a topological $\K[[\hbar]]$-module  but is equipped with a different $\hbar$-dependent  differential
\Beq\label{5: d_hbar differential in HoLoB[[h]]}
\delta
\Ba{c}\resizebox{13mm}{!}{\xy
(-9,-6)*{};
(0,0)*+{\bu}%*\cir{}
**\dir{-};
(-5,-6)*{};
(0,0)*+{}%*\cir{}
**\dir{-};
(9,-6)*{};
(0,0)*+{}%*\cir{}
**\dir{-};
(5,-6)*{};
(0,0)*+{}%*\cir{}
**\dir{-};
(0,-6)*{\ldots};
(-10,-8)*{_1};
(-6,-8)*{_2};
(10,-8)*{_m};
(-9,6)*{};
(0,0)*+{}%*\cir{}
**\dir{-};
(-5,6)*{};
(0,0)*+{}%*\cir{}
**\dir{-};
(9,6)*{};
(0,0)*+{}%*\cir{}
**\dir{-};
(5,6)*{};
(0,0)*+{}%*\cir{}
**\dir{-};
(0,6)*{\ldots};
(-10,8)*{_1};
(-6,8)*{_2};
(10,8)*{_n};
\endxy}\Ea
=
\sum_{l\geq 1}\sum_{[m]=I_1\sqcup I_2\atop
[n]=J_1\sqcup J_2}\pm \hbar^{l-1}
\Ba{c}
%
%
%%%%%%%%%%%%%%%% two vertex graph with l internal edges %%%%%%%%%%
\Ba{c}\resizebox{18mm}{!}{\xy
(0,0)*+{\bu}="b",
(10,10)*+{\bu}="c",
%
%%%%%%%%%% edges to b %%%%%%%%%%%%
(-9,6)*{}="1",
(-7,6)*{}="2",
(-2,6)*{}="3",
(-3.5,5)*{...},
(-4,-6)*{}="-1",
(-2,-6)*{}="-2",
(4,-6)*{}="-3",
(1,-5)*{...},
(0,-8)*{\underbrace{\ \ \ \ \ \ \ \ }},
(0,-11)*{_{J_1}},
(-6,8)*{\overbrace{ \ \ \ \ \ \ }},
(-6,11)*{_{I_1}},
%%%%%%%%%% edges to c %%%%%%%%%%%%
(6,16)*{}="1'",
(8,16)*{}="2'",
(14,16)*{}="3'",
(11,15)*{...},
(11,6)*{}="-1'",
(16,6)*{}="-2'",
(18,6)*{}="-3'",
(13.5,6)*{...},
(15,4)*{\underbrace{\ \ \ \ \ \ \ }},
(15,1)*{_{J_2}},
(10,18)*{\overbrace{ \ \ \ \ \ \ \ \ }},
(10,21)*{_{I_2}},
%
%%%%%%%%%%% internal curved edges %%%%%%%%%%%
(0,2)*-{};(8.0,10.0)*-{}
**\crv{(0,10)};
(0.5,1.8)*-{};(8.5,9.0)*-{}
**\crv{(0.4,7)};
(1.5,0.5)*-{};(9.1,8.5)*-{}
**\crv{(5,1)};
(1.7,0.0)*-{};(9.5,8.6)*-{}
**\crv{(6,-1)};
(5,5)*+{...};
\ar @{-} "b";"1" <0pt>
\ar @{-} "b";"2" <0pt>
\ar @{-} "b";"3" <0pt>
\ar @{-} "b";"-1" <0pt>
\ar @{-} "b";"-2" <0pt>
\ar @{-} "b";"-3" <0pt>
\ar @{-} "c";"1'" <0pt>
\ar @{-} "c";"2'" <0pt>
\ar @{-} "c";"3'" <0pt>
\ar @{-} "c";"-1'" <0pt>
\ar @{-} "c";"-2'" <0pt>
\ar @{-} "c";"-3'" <0pt>
\endxy}\Ea
%%%%%%%%%%%%%%%%%%%%%%%%%%%%%%%%%%%%%%%%%%%
\Ea
\Eeq
where $l$ counts the number of internal edges connecting the two vertices on the r.h.s. The symbol $\pm$ stands for $-1$ in the case $c,d\in 2\Z$. There is a morphism of dg properads (cf.\ \cite{CMW})
$$
\cF^+:  \HoLBcd^{\hbar+} \lon \HoLoBpcd[[\hbar]]
$$
given on the generators as follows (cf.\ \cite{CMW})
\Beq\label{5: formal power series of generators}
\cF^+:\
\Ba{c}\resizebox{13mm}{!}{\begin{xy}
 <0mm,0mm>*{\bu};<0mm,0mm>*{}**@{},
 <-0.6mm,0.44mm>*{};<-8mm,5mm>*{}**@{-},
 <-0.4mm,0.7mm>*{};<-4.5mm,5mm>*{}**@{-},
 <0mm,0mm>*{};<-1mm,5mm>*{\ldots}**@{},
 <0.4mm,0.7mm>*{};<4.5mm,5mm>*{}**@{-},
 <0.6mm,0.44mm>*{};<8mm,5mm>*{}**@{-},
   <0mm,0mm>*{};<-8.5mm,5.5mm>*{^1}**@{},
   <0mm,0mm>*{};<-5mm,5.5mm>*{^2}**@{},
   %<0mm,0mm>*{};<4.5mm,5.5mm>*{^{m\hspace{-0.5mm}-\hspace{-0.5mm}1}}**@{},
   <0mm,0mm>*{};<9.0mm,5.5mm>*{^m}**@{},
 <-0.6mm,-0.44mm>*{};<-8mm,-5mm>*{}**@{-},
 <-0.4mm,-0.7mm>*{};<-4.5mm,-5mm>*{}**@{-},
 <0mm,0mm>*{};<-1mm,-5mm>*{\ldots}**@{},
 <0.4mm,-0.7mm>*{};<4.5mm,-5mm>*{}**@{-},
 <0.6mm,-0.44mm>*{};<8mm,-5mm>*{}**@{-},
   <0mm,0mm>*{};<-8.5mm,-6.9mm>*{^1}**@{},
   <0mm,0mm>*{};<-5mm,-6.9mm>*{^2}**@{},
   %<0mm,0mm>*{};<4.5mm,-6.9mm>*{^{n\hspace{-0.5mm}-\hspace{-0.5mm}1}}**@{},
   <0mm,0mm>*{};<9.0mm,-6.9mm>*{^n}**@{},
 \end{xy}}\Ea
 \lon
 \sum_{a=0}^\infty \hbar^{a}
\Ba{c}\resizebox{14mm}{!}{ \xy
(-9,-6)*{};
(0,0)*+{a}*\cir{}
**\dir{-};
(-5,-6)*{};
(0,0)*+{a}*\cir{}
**\dir{-};
(9,-6)*{};
(0,0)*+{a}*\cir{}
**\dir{-};
(5,-6)*{};
(0,0)*+{a}*\cir{}
**\dir{-};
(0,-6)*{\ldots};
(-10,-8)*{_1};
(-6,-8)*{_2};
(10,-8)*{_m};
(-9,6)*{};
(0,0)*+{a}*\cir{}
**\dir{-};
(-5,6)*{};
(0,0)*+{a}*\cir{}
**\dir{-};
(9,6)*{};
(0,0)*+{a}*\cir{}
**\dir{-};
(5,6)*{};
(0,0)*+{a}*\cir{}
**\dir{-};
(0,6)*{\ldots};
(-10,8)*{_1};
(-6,8)*{_2};
(10,8)*{_n};
\endxy}\Ea
\ \ \ \ \ \ \ \ \ \ \ \ \ \forall \ m,n\geq 1.
\Eeq
There is obviously a 1-1 correspondence between morphisms of dg properads
$\HoLoBcd \lon \cP
$
in the category of graded vector spaces over $\K$, and continuous morphisms
of dg properads
$\HoLBcd^{\hbar+} \lon \cP[[\hbar]]
$
in the category of topological $\K[[\hbar]]$-modules.
%Every morphism of dg properads
%$\HoLoBpcd \rar \cP$ gives us a continuous morphism
%$\HoLBcd^{\hbar+} \rar \cP[[\hbar]]$ and vice versa. Moreover, one has the following useful %observation.

%\subsubsection{\bf Lemma \cite{MW1}}\label{5: lemma on holiebb-h} {\em There is a 1-1 %correspondence between morphisms of dg properads
%$
%f: \HoLoBcd \lon \cP
%$
%in the category of graded vector spaces over $\K$, and continuous morphisms
%of dg properads
%$
%f^\hbar: \HoLBcd^{\hbar+} \lon \cP[[\hbar]]
%$
%in the category of topological $\K[[\hbar]]$-modules satisfying the condition
%\Beq\label{2: admissible maps condition}
%$
%\rho_\hbar\left(\begin{xy}
% <0mm,0mm>*{\bu};
% <0mm,-0.5mm>*{};<0mm,-2.5mm>*{}**@{-},
% <0mm,0.5mm>*{};<0mm,2.5mm>*{}**@{-},
% \end{xy}\right)\in \hbar\cP[[\hbar]]
%$
%}

\sip

Let $\Holie_d^\diamond$  be the quotient of $\HoLoBcd$ by the (differential) ideal generated
by all corollas with the number of outgoing legs $\geq 2$.
It is generated by  the following family of  (skew)symmetric corollas with $a\geq 0$, $n\geq 1$ and $a+n\geq 2$,
$$
\Ba{c}\resizebox{13mm}{!}{ \xy
(-7.5,-8.6)*{_{_1}};
(-4.1,-8.6)*{_{_2}};
%(4.5,-8.6)*{_{_{k\hspace{-0.3mm}-\hspace{-0.3mm}1}}};
(9.0,-8.5)*{_{_{n}}};
(0.0,-6)*{...};
(0,6)*{};
(0,0)*+\hbox{$_{{a}}$}*\frm{o}
**\dir{-};
(-4,-7)*{};
(0,0)*+\hbox{$_{{a}}$}*\frm{o}
**\dir{-};
(-7,-7)*{};
(0,0)*+\hbox{$_{{a}}$}*\frm{o}
**\dir{-};
(8,-7)*{};
(0,0)*+\hbox{$_{{a}}$}*\frm{o}
**\dir{-};
(4,-7)*{};
(0,0)*+\hbox{$_{{a}}$}*\frm{o}
**\dir{-};
\endxy}\Ea
=(-1)^{d|\sigma|}
\Ba{c}\resizebox{15.2mm}{!}{ \xy
(-8.5,-8.6)*{_{_{\sigma(1)}}};
(-3.1,-8.6)*{_{_{\sigma(2)}}};
%(4.5,-8.6)*{_{_{k\hspace{-0.3mm}-\hspace{-0.3mm}1}}};
(9.0,-8.5)*{_{_{\sigma(n)}}};
(0.0,-6)*{...};
(0,6)*{};
(0,0)*+\hbox{$_{{a}}$}*\frm{o}
**\dir{-};
(-4,-7)*{};
(0,0)*+\hbox{$_{{a}}$}*\frm{o}
**\dir{-};
(-7,-7)*{};
(0,0)*+\hbox{$_{{a}}$}*\frm{o}
**\dir{-};
(8,-7)*{};
(0,0)*+\hbox{$_{{a}}$}*\frm{o}
**\dir{-};
(4,-7)*{};
(0,0)*+\hbox{$_{{a}}$}*\frm{o}
**\dir{-};
\endxy}\Ea\ \ \ \forall\ \sigma\in \bS_n,
$$
which  are assigned  degree $1 - d(n-1 +a)$;  the induced differential acts as follows
$$
\delta\Ba{c}
\resizebox{12mm}{!}{ \xy
(-7.5,-8.6)*{_{_1}};
(-4.1,-8.6)*{_{_2}};
%(4.5,-8.6)*{_{_{k\hspace{-0.3mm}-\hspace{-0.3mm}1}}};
(9.0,-8.5)*{_{_{n}}};
(0.0,-6)*{...};
(0,6)*{};
(0,0)*+\hbox{$_{{a}}$}*\frm{o}
**\dir{-};
(-4,-7)*{};
(0,0)*+\hbox{$_{{a}}$}*\frm{o}
**\dir{-};
(-7,-7)*{};
(0,0)*+\hbox{$_{{a}}$}*\frm{o}
**\dir{-};
(8,-7)*{};
(0,0)*+\hbox{$_{{a}}$}*\frm{o}
**\dir{-};
(4,-7)*{};
(0,0)*+\hbox{$_{{a}}$}*\frm{o}
**\dir{-};
\endxy}\Ea
=
%%%%%%%%%%%%%%%%%%%%%%%%%%%%%
\sum_{a=p+q\atop [n]=I_1\sqcup I_2}\pm
\Ba{c}
%
%
%%%%%%%%%%%%%%%% two vertex graph with 1 internal edge %%%%%%%%%%
\resizebox{15mm}{!}{ \xy
(0,0)*+{p}*\cir{}="b",
(10,10)*+{q}*\cir{}="c",
%
%%%%%%%%%% edges to b %%%%%%%%%%%%
(-4,-6)*{}="-1",
(-2,-6)*{}="-2",
(4,-6)*{}="-3",
(1,-5)*{...},
(0,-8)*{\underbrace{\ \ \ \ \ \ \ \ }},
(0,-11)*{_{I_1}},
%%%%%%%%%% edges to c %%%%%%%%%%%%
%(6,16)*{}="1'",
(10,16)*{}="2'",
%(14,16)*{}="3'",
%(11,15)*{...},
(11,4)*{}="-1'",
(16,4)*{}="-2'",
(18,4)*{}="-3'",
(13.5,4)*{...},
(15,2)*{\underbrace{\ \ \ \ \ \ \ }},
(15,-1)*{_{I_2}},
%
%%%%%%%%%%% internal curved edges %%%%%%%%%%%
\ar @{-} "b";"c" <0pt>
\ar @{-} "b";"-1" <0pt>
\ar @{-} "b";"-2" <0pt>
\ar @{-} "b";"-3" <0pt>
%
%\ar @{-} "c";"1'" <0pt>
\ar @{-} "c";"2'" <0pt>
%\ar @{-} "c";"3'" <0pt>
\ar @{-} "c";"-1'" <0pt>
\ar @{-} "c";"-2'" <0pt>
\ar @{-} "c";"-3'" <0pt>
\endxy}
%%%%%%%%%%%%%%%%%%%%%%%%%%%%%%%%%%%%%%%%%%%
\Ea
$$
Representations, $\rho: \caH o lie^\diamond_{d}\rar \cE nd_V$, of this operad in a dg vector space $(V,\p)$ can be identified
with continuous representations of  the  topological operad $\caH o lie_{d}[[\hbar]]$
in the topological vector space $V[[\hbar]]$ equipped with the differential
$$
\p + \sum_{p\geq 1}\hbar^p \Delta_p, \ \ \ \Delta_p:=\rho \left(\Ba{c}
\resizebox{4mm}{!}{ \xy
(0,5)*{};
(0,0)*+{_a}*\cir{}
**\dir{-};
(0,-5)*{};
(0,0)*+{_a}*\cir{}
**\dir{-};
\endxy}\Ea\right).
$$
Here the formal parameter $\hbar$ is assumed to have homological degree $d$.

\sip

It is easy to see that the quotient of the dg properad
 $\HoLBcd^{\hbar+}$ by the (differential) ideal generated
by all corollas with the number of outgoing legs $\geq 2$ is identical to $\Holie_d^+[[\hbar]]$ as a dg properad.
Hence we obtain from (\ref{5: formal power series of generators}) a canonical morphism of dg properads
\Beq\label{5: f^+ from Holie^h to Holie[h]}
f^+:  \Holie_d^{+}[[\hbar]] \lon \Holie_d^{\diamond+}[[\hbar]]
\Eeq
It gives us a compact presentation of any morphism $\Holie_d^{\diamond+}\rar \cP$ as an associated continuous morphism of properads
$\Holie_d[[\hbar]] \rar \cP[[\hbar]]$ in the category of topological $\K[[\hbar]]$-modules.

\subsubsection{\bf Proposition}\label{4: Prop on map from diamond Lie^+ to OHoLB}
{\em There is a  morphism of dg operads
$$
F^+: \Holie^{\diamond +}_{c+d} \to \f_{c,d}\HoLoBcd
$$
given explicitly on the $(1,1)$-generators by
$$
\Ba{c}
\resizebox{3.5mm}{!}{ \xy
(0,5)*{};
(0,0)*+{_0}*\cir{}
**\dir{-};
(0,-5)*{};
(0,0)*+{_0}*\cir{}
**\dir{-};
\endxy}\Ea
 \lon
\sum_{m\geq 2} \Ba{c}\resizebox{9mm}{!}{\begin{xy}
(0,-6)*{};
(0,0)*+{_0}*\cir{}
**\dir{-};
(-5,6)*{};
(0,0)*+{_0}*\cir{}
**\dir{-};
(5,6)*{};
(0,0)*+{_0}*\cir{}
**\dir{-};
(-2,6)*{};
(0,0)*+{_0}*\cir{}
**\dir{-};
(2,6)*{};
(0,0)*+{_0}*\cir{}
**\dir{-};
     <0mm,8mm>*{\overbrace{ \ \ \ \ \ \ \ \ \ \ }},
     <0mm,10mm>*{^m},
 \end{xy}}\Ea,
$$
and on the remaining $(1,n)$-generators with $a+n\geq 2$ by}
\Beq\label{5: map Holie^diom+c+d to fc,d}
\Ba{c}
\resizebox{12mm}{!}{ \xy
(-7.5,-8.6)*{_{_1}};
(-4.1,-8.6)*{_{_2}};
%(4.5,-8.6)*{_{_{k\hspace{-0.3mm}-\hspace{-0.3mm}1}}};
(9.0,-8.5)*{_{_{n}}};
(0.0,-6)*{...};
(0,6)*{};
(0,0)*+\hbox{$_{{a}}$}*\frm{o}
**\dir{-};
(-4,-7)*{};
(0,0)*+\hbox{$_{{a}}$}*\frm{o}
**\dir{-};
(-7,-7)*{};
(0,0)*+\hbox{$_{{a}}$}*\frm{o}
**\dir{-};
(8,-7)*{};
(0,0)*+\hbox{$_{{a}}$}*\frm{o}
**\dir{-};
(4,-7)*{};
(0,0)*+\hbox{$_{{a}}$}*\frm{o}
**\dir{-};
\endxy}\Ea
 \lon
 %%%%%%%%%%%%%
\sum_{m\geq 1, l_i\geq 1 \atop a=c+\sum_{i=1}^n(l_i-1)}
\Ba{c}\resizebox{16mm}{!}{\xy
(0,0)*+{_1}*\cir{}="b",
(10,10)*+{c}*\cir{}="c",
(20,0)*+{_n}*\cir{}="r",
%
%%%%%%%%%% edges to c %%%%%%%%%%%%
(6,16)*{}="1'",
(8,16)*{}="2'",
(14,16)*{}="3'",
(11,15)*{...},
(10,0)*{\ldots},
(10,18)*{\overbrace{ \ \ \ \ \ \ \ \ }},
(10,21)*{_{m}},
%
%%%%%%%%%%% internal left curved edges %%%%%%%%%%%
(0,2)*-{};(8.0,10.0)*-{}
**\crv{(0,10)};
(0.5,1.8)*-{};(8.5,9.0)*-{}
**\crv{(0.4,7)};
(1.8,0.6)*-{};(9.1,8.5)*-{}
**\crv{(5,1)};
(2.0,0.1)*-{};(9.5,8.6)*-{}
**\crv{(6,-1)};
(5,5)*+{_{l_1}};
%%%%%%%%%%%%% internal right curved edges %%%%%
(20,2)*-{};(12.0,10.0)*-{}
**\crv{(20,10)};
(19.5,1.8)*-{};(11.5,9.0)*-{}
**\crv{(20.4,7)};
(17.9,0.6)*-{};(10.9,8.5)*-{}
**\crv{(15,1)};
(1.7,0.0)*-{};(9.5,8.6)*-{}
**\crv{(6,-1)};
(16,5)*+{_{l_n}};
\ar @{-} "c";"1'" <0pt>
\ar @{-} "c";"2'" <0pt>
\ar @{-} "c";"3'" <0pt>
\endxy}\Ea.
\Eeq
Proof is a straightforward direct calculation (cf.\  \S{\ref{4: Prop on map from Lie^+ to OHoLB}}).
The existence of such a map  follows also from  Proposition 5.4.1 and Lemma B.4.1 proven in \cite{CMW}.
% where the above map is obtained (essentially) as the composition
%$$
%\Holie_{c+d}^+ \lon \cB\cV_\infty^{com} \lon \f_{c,d}\HoLBcd
%$$
%where $\cB\cV_\infty^{com}$ is a minimal resolution of the operad of so called (degree %shifted)  $\cB\cV^{com}$-algebras (see \S 5.3 in \cite{CMW} for full details).

\sip

Given a $\HoLoBcd$-algebra structure,
$$
\rho: \HoLoBcd \lon \cE nd_V,
$$
$$
\mu^a_{m,n}:=\rho\left(\Ba{c}\resizebox{13mm}{!}{ \xy
(-9,-6)*{};
(0,0)*+{a}*\cir{}
**\dir{-};
(-5,-6)*{};
(0,0)*+{a}*\cir{}
**\dir{-};
(9,-6)*{};
(0,0)*+{a}*\cir{}
**\dir{-};
(5,-6)*{};
(0,0)*+{a}*\cir{}
**\dir{-};
(0,-6)*{\ldots};
(-10,-8)*{_1};
(-6,-8)*{_2};
(10,-8)*{_m};
(-9,6)*{};
(0,0)*+{a}*\cir{}
**\dir{-};
(-5,6)*{};
(0,0)*+{a}*\cir{}
**\dir{-};
(9,6)*{};
(0,0)*+{a}*\cir{}
**\dir{-};
(5,6)*{};
(0,0)*+{a}*\cir{}
**\dir{-};
(0,6)*{\ldots};
(-10,8)*{_1};
(-6,8)*{_2};
(10,8)*{_n};
\endxy}\Ea
\right): \odot^n (V[-c]) \lon (\odot^m (V[-c])) [1+(c+d)(1-n-a)],
$$
in a graded vector space $V$, there is an associated  $\f_{c,d}\HoLBcd$-algebra structure
in $\odot^\bu(V[-c])$ given in terms of polydifferential operators, and hence a
continuous $\Holie_d^{+}[[\hbar]]$-algebra structure on  $\odot^\bu(V[-c])[[\hbar]]$ given by the composition
$$
\Holie_d^{+}[[\hbar]] \stackrel{f^+}{\rar} \Holie_d^{\diamond+}[[\hbar]] \stackrel{F^+}{\rar} \f_{c,d}\HoLBcd[[\hbar]] \rar \cE nd_{\odot^\bu(V[-c])}[[\hbar]]
$$
% Hence there is an induced
%$\Holie_{c+d}^+$-algebra structure in  $\odot^\bu(V[-c])$.
Assuming that the latter is nilpotent (or appropriately filtered which is often the case in applications), one defines a {\em Maurer-Cartan element $\ga$ of the given
$\HoLBcd$-algebra structure in $V$}\, as a Maurer-Cartan of the induced continuous $\Holie_{c+d}^{\hbar+}$-algebra structure in  $\odot^\bu(V[-c])[[\hbar]]$. Using
 (\ref{5: map Holie^diom+c+d to fc,d}) one can describe such an MC element as a homogeneous
 (of degree $c+d$)  formal power series
 \Beq\label{5: ga-MC as h-series}
 \ga=\sum_{a\geq 0, m\geq 0}\hbar^a\ga_{a,m}\in \odot^{\bu\geq 1}(V[-c])[[\hbar]], \ \
 \ga_{a,m}\in \odot^m (V[-c]),
 \Eeq
 satisfying the equation
 \Beq\label{5: MC eqn Delta(e^g)}
\Delta_\rho\left(e^{\frac{\ga}{\hbar}}\right)=0
\Eeq
where $\Delta_\rho$ is a degree $+1$ polydifferential operator on $\odot^{\bu}(V[-c])[[\hbar]]$ given, in an arbitrary basis $\{p_\al\}$ of $V[-c]$ as a sum (cf. \S {\ref{4: subsec on  MC elements of HoLBcd}})
\Beq\label{5: Delta_rho}
\Delta_\rho:= \sum_{a\geq 0\atop
m,n\geq 1}\pm \hbar^{a+n-1} \mu^a_{m,n}(p_{\al_1}\ot ...\ot p_{\al_n})\frac{\p^n}{\p p_{\al_1}\ldots \p p_{\al_n}}
\Eeq
Here the differential in $V$ is encoded as $\mu^0_{1,1}$. The operator $\Delta_\rho$ encodes fully the given $\HoLoBcd$-algebra structure $\rho$ in $V$: there is a {\em one-to-one correspondence}\,  \cite{CMW,Me3} between $\HoLoBcd$-algebra structures in $V$ and degree
$1$ operators on $\odot^\bu(V[-c])[[\hbar]]$ of the form
$$
\Delta=\sum_{a\geq 0} \hbar^a \Delta_a
$$
such that $\Delta_a$ is a derivation of the graded commutative algebra $\odot^\bu(V[-c])$ of order $\leq a+1$ (such structures are often called $\cB\cV_\infty^{com}$-algebra structures in the literatures).

\sip

The sum
\Beq\label{5: twisted by MC differentiail in V}
\delta_\ga v :=\sum_{n\geq 0}\frac{1}{n!}{\mu}^0_{1,n+1}(\underbrace{\ga_1,\ldots,\ga_1}_n, v ),
\Eeq
 is a twisted  differential on $V$ which is used in the following definition-proposition.

\subsection{Diamond twisting endofunctor} The {\em diamond twisting}, $(\Tw^\diamond \cP,\pc)$, of a dg properad $(\cP,\p)$ under $\HoLoBcd$
 is, by definition, a properad freely generated by $\cP$ and the $\bS$-bimodule $M=\{M(m,n)\}$ such that $M(m,n)=0$ for $n\geq 1$ and $M(m,0)=\oplus_{a\geq 0} M_a(0,m)$  with $M_a(0,m)$ being the following 1-dimensional representations of $\bS_m$,
$$
M_a(0,m):=\sgn_m^{|c|}[(c+d)(1-a)-cm]=\text{span}\left\langle
\Ba{c}\resizebox{9mm}{!}{\begin{xy}
(0,6)*{...};
(-5.5,7)*{^1};
(0,0)*+{_a}*\cir{}
**\dir{-};
(5.5,7)*{^m};
(0,0)*+{_a}*\cir{}
**\dir{-};
(-2.5,7)*{^2};
(0,0)*+{_a}*\cir{}
**\dir{-};
(2.5,6)*{};
(0,0)*+{_a}*\cir{}
**\dir{-};
%
 %   <0mm,8mm>*{\overbrace{ \ \ \ \ \ \ \ \ \ \ }},
 %   <0mm,10mm>*{^m},
 \end{xy}}\Ea
\right\rangle
$$
The differential in $\Tw^\diamond \cP$ is defined on the generators as follows
\Beq\label{5: d_centerdot on Tw^diamond P}
\sd_\centerdot \Ba{c}\resizebox{13mm}{!}{
 \begin{xy}
 <0mm,0mm>*{\circ};<-8mm,6mm>*{^1}**@{-},
 <0mm,0mm>*{\circ};<-4.5mm,6mm>*{^2}**@{-},
 <0mm,0mm>*{\circ};<0mm,5.5mm>*{\ldots}**@{},
 <0mm,0mm>*{\circ};<3.5mm,5mm>*{}**@{-},
 <0mm,0mm>*{\circ};<8mm,6mm>*{^m}**@{-},
 <0mm,0mm>*{\circ};<-8mm,-6mm>*{_1}**@{-},
 <0mm,0mm>*{\circ};<-4.5mm,-6mm>*{_2}**@{-},
 <0mm,0mm>*{\circ};<0mm,-5.5mm>*{\ldots}**@{},
 <0mm,0mm>*{\circ};<4.5mm,-6mm>*+{}**@{-},
 <0mm,0mm>*{\circ};<8mm,-6mm>*{_n}**@{-},
   \end{xy}}\Ea
=
\sd \Ba{c}\resizebox{13mm}{!}{
 \begin{xy}
 <0mm,0mm>*{\circ};<-8mm,6mm>*{^1}**@{-},
 <0mm,0mm>*{\circ};<-4.5mm,6mm>*{^2}**@{-},
 <0mm,0mm>*{\circ};<0mm,5.5mm>*{\ldots}**@{},
 <0mm,0mm>*{\circ};<3.5mm,5mm>*{}**@{-},
 <0mm,0mm>*{\circ};<8mm,6mm>*{^m}**@{-},
 <0mm,0mm>*{\circ};<-8mm,-6mm>*{_1}**@{-},
 <0mm,0mm>*{\circ};<-4.5mm,-6mm>*{_2}**@{-},
 <0mm,0mm>*{\circ};<0mm,-5.5mm>*{\ldots}**@{},
 <0mm,0mm>*{\circ};<4.5mm,-6mm>*+{}**@{-},
 <0mm,0mm>*{\circ};<8mm,-6mm>*{_n}**@{-},
   \end{xy}}\Ea
+
\overset{m-1}{\underset{i=0}{\sum}}
\Ba{c}\resizebox{14mm}{!}{
\begin{xy}
 %<0mm,0mm>*{\circ};<0mm,0mm>*{}**@{},
 <0mm,0mm>*{\circ};<-8mm,5mm>*{}**@{-},
 <0mm,0mm>*{\circ};<-3.5mm,5mm>*{}**@{-},
 <0mm,0mm>*{\circ};<-6mm,5mm>*{..}**@{},
 <0mm,0mm>*{\circ};<0mm,5mm>*{}**@{-},
  <0mm,5mm>*{\blacklozenge};
  <0mm,5mm>*{};<0mm,8mm>*{}**@{-},
  <0mm,5mm>*{};<0mm,9mm>*{^{i\hspace{-0.2mm}+\hspace{-0.5mm}1}}**@{},
<0mm,0mm>*{\circ};<8mm,5mm>*{}**@{-},
<0mm,0mm>*{\circ};<3.5mm,5mm>*{}**@{-},
<6mm,5mm>*{..}**@{},
<-8.5mm,5.5mm>*{^1}**@{},
<-4mm,5.5mm>*{^i}**@{},
<9.0mm,5.5mm>*{^m}**@{},
 <0mm,0mm>*{\circ};<-8mm,-5mm>*{}**@{-},
 <0mm,0mm>*{\circ};<-4.5mm,-5mm>*{}**@{-},
 <-1mm,-5mm>*{\ldots}**@{},
 <0mm,0mm>*{\circ};<4.5mm,-5mm>*{}**@{-},
 <0mm,0mm>*{\circ};<8mm,-5mm>*{}**@{-},
<-8.5mm,-6.9mm>*{^1}**@{},
<-5mm,-6.9mm>*{^2}**@{},
<4.5mm,-6.9mm>*{^{n\hspace{-0.5mm}-\hspace{-0.5mm}1}}**@{},
<9.0mm,-6.9mm>*{^n}**@{},
 \end{xy}}\Ea
 - (-1)^{|a|}
\overset{n-1}{\underset{i=0}{\sum}}
 \Ba{c}\resizebox{14mm}{!}{\begin{xy}
 %<0mm,0mm>*{\circ};
 <0mm,0mm>*{\circ};<-8mm,-5mm>*{}**@{-},
 <0mm,0mm>*{\circ};<-3.5mm,-5mm>*{}**@{-},
 <0mm,0mm>*{\circ};<-6mm,-5mm>*{..}**@{},
 <0mm,0mm>*{\circ};<0mm,-5mm>*{}**@{-},
  <0mm,-5mm>*{\blacklozenge};
  <0mm,-5mm>*{};<0mm,-8mm>*{}**@{-},
  <0mm,-5mm>*{};<0mm,-10mm>*{^{i\hspace{-0.2mm}+\hspace{-0.5mm}1}}**@{},
<0mm,0mm>*{\circ};<8mm,-5mm>*{}**@{-},
<0mm,0mm>*{\circ};<3.5mm,-5mm>*{}**@{-},
 <6mm,-5mm>*{..}**@{},
<-8.5mm,-6.9mm>*{^1}**@{},
<-4mm,-6.9mm>*{^i}**@{},
<9.0mm,-6.9mm>*{^n}**@{},
 <0mm,0mm>*{\circ};<-8mm,5mm>*{}**@{-},
 <0mm,0mm>*{\circ};<-4.5mm,5mm>*{}**@{-},
<-1mm,5mm>*{\ldots}**@{},
 <0mm,0mm>*{\circ};<4.5mm,5mm>*{}**@{-},
 <0mm,0mm>*{\circ};<8mm,5mm>*{}**@{-},
<-8.5mm,5.5mm>*{^1}**@{},
<-5mm,5.5mm>*{^2}**@{},
<4.5mm,5.5mm>*{^{m\hspace{-0.5mm}-\hspace{-0.5mm}1}}**@{},
<9.0mm,5.5mm>*{^m}**@{},
 \end{xy}}\Ea
\Eeq

\Beq\label{5: d_centerdot on involutive MC generators of Tw^diamondP}
 \sd_\centerdot
 \Ba{c}\resizebox{9mm}{!}{\begin{xy}
(0,6)*{...};
(-5.5,7)*{^1};
(0,0)*+{_a}*\cir{}
**\dir{-};
(5.5,7)*{^m};
(0,0)*+{_a}*\cir{}
**\dir{-};
(-2.5,7)*{^2};
(0,0)*+{_a}*\cir{}
**\dir{-};
(2.5,6)*{};
(0,0)*+{_a}*\cir{}
**\dir{-};
%
 %   <0mm,8mm>*{\overbrace{ \ \ \ \ \ \ \ \ \ \ }},
 %   <0mm,10mm>*{^m},
 \end{xy}}\Ea
:=
\overset{m-1}{\underset{i=0}{\sum}}
\Ba{c}\resizebox{15mm}{!}{
\begin{xy}
 <0mm,-1mm>*+{_a}*\cir{};
 <-0.9mm,0.2mm>*{};<-8mm,5mm>*{}**@{-},
 <-0.3mm,0.4mm>*{};<-3.5mm,5mm>*{}**@{-},
 <-6mm,5mm>*{..}**@{},
 <0mm,0.5mm>*{};<0mm,5mm>*{}**@{-},
  <0mm,5mm>*{\blacklozenge};
  <0mm,5mm>*{};<0mm,8mm>*{}**@{-},
  <0mm,5mm>*{};<0mm,9mm>*{^{i\hspace{-0.2mm}+\hspace{-0.5mm}1}}**@{},
<0.9mm,0.2mm>*{};<8mm,5mm>*{}**@{-},
<0.3mm,0.4mm>*{};<3.5mm,5mm>*{}**@{-},
 <0mm,0mm>*{};<6mm,5mm>*{..}**@{},
   <0mm,0mm>*{};<-8.5mm,5.7mm>*{^1}**@{},
   <0mm,0mm>*{};<-4mm,5.7mm>*{^i}**@{},
   <0mm,0mm>*{};<9.0mm,5.7mm>*{^m}**@{},
 \end{xy}}\Ea
-
\sum_{k\geq 1, [m]=\sqcup  [m_\bu]\atop
%{m_{0}\geq 1, k+m_0\geq 3 \atop m_1,...,m_k\geq 0}} \atop
a=b+\sum_{i=1}^k(c_i+l_i-1)}
\frac{1}{k!}
\Ba{c}\resizebox{21mm}{!}{\xy
(0,0)*+{_{c_1}}*\cir{}="b",
(10,10)*+{b}*\cir{}="c",
(20,0)*+{_{c_k}}*\cir{}="r",
%
%%%%%%%%%% edges to c %%%%%%%%%%%%
(6,16)*{}="1'",
(8,16)*{}="2'",
(14,16)*{}="3'",
(11,15)*{...},
(10,0)*{\ldots},
(10,18)*{\overbrace{ \ \ \ \ \ \ \ \ }},
(10,21)*{_{m_0}},
%
%%%%%%%%%%% internal left curved edges %%%%%%%%%%%
(0,2)*-{};(8.0,10.0)*-{}
**\crv{(0,10)};
(0.5,1.8)*-{};(8.5,9.0)*-{}
**\crv{(0.4,7)};
(1.8,0.6)*-{};(9.1,8.5)*-{}
**\crv{(5,1)};
(2.0,0.1)*-{};(9.5,8.6)*-{}
**\crv{(6,-1)};
(5,5)*+{_{l_1}};
%%%%%%%%%%%%% internal right curved edges %%%%%
(20,2)*-{};(12.0,10.0)*-{}
**\crv{(20,10)};
(19.5,1.8)*-{};(11.5,9.0)*-{}
**\crv{(20.4,7)};
(17.9,0.6)*-{};(10.9,8.5)*-{}
**\crv{(15,1)};
(1.7,0.0)*-{};(9.5,8.6)*-{}
**\crv{(6,-1)};
(16,5)*+{_{l_k}};
(-7,10)*{}="1l";
(-5,10)*{}="2l";
(-3,10)*{}="3l";
(-5,12)*{\overbrace{  \ \ \ \ \ \ }},
(-5,14)*{_{m_1}},
(27,10)*{}="1r";
(25,10)*{}="2r";
(23,10)*{}="3r";
(25,12)*{\overbrace{  \ \ \ \ \ \ }},
(25,14)*{_{m_k}},
\ar @{-} "c";"1'" <0pt>
\ar @{-} "c";"2'" <0pt>
\ar @{-} "c";"3'" <0pt>
\ar @{-} "b";"1l" <0pt>
\ar @{-} "b";"2l" <0pt>
\ar @{-} "b";"3l" <0pt>
\ar @{-} "r";"1r" <0pt>
\ar @{-} "r";"2r" <0pt>
\ar @{-} "r";"3r" <0pt>
\endxy}\Ea\ \ \ \  \forall\ a,m\geq 1.
\Eeq
where $\Ba{c}\resizebox{1.8mm}{!}{\begin{xy}
 <0mm,-0.55mm>*{};<0mm,-3mm>*{}**@{-},
 <0mm,0.5mm>*{};<0mm,3mm>*{}**@{-},
 <0mm,0mm>*{\blacklozenge};<0mm,0mm>*{}**@{},
 \end{xy}}\Ea
$ is given by

\Beq\label{5: involutive blacklozenge (1,1) element}
\begin{xy}
 <0mm,-3mm>*{};<0mm,3mm>*{}**@{-},
 %<0mm,0.5mm>*{};<0mm,3mm>*{}**@{-},
 <0mm,0mm>*{_\blacklozenge};<0mm,0mm>*{}**@{},
 \end{xy}
 :=
 \sum_{k=1}^\infty \frac{1}{k!}\ \resizebox{17mm}{!}{  \xy
(-25,9)*{}="1";
(-9,-4)*{...};
 <-11mm,-8mm>*{\underbrace{ \ \ \ \ \ \ \ \ \ \ \ \ \ \ \ \ \ }_{k}},
(-25,3)*+{_0}*\cir{}="L";
%(-25,+3)*{\circledast}="L";
    (-25,-3)*{}="N";
 (-19,-4)*+{_0}*\cir{}="B";
  (-13,-4)*+{_0}*\cir{}="C";
   (-2,-4)*+{_0}*\cir{}="D";
\ar @{-} "D";"L" <0pt>
\ar @{-} "C";"L" <0pt>
\ar @{-} "B";"L" <0pt>
\ar @{-} "1";"L" <0pt>
\ar @{-} "N";"L" <0pt>
%
 %
% \ar @{-} "l1";"N" <0pt>
%  \ar @{-} "l2";"N" <0pt>
%   \ar @{-} "ln";"N" <0pt>
 \endxy}
\Eeq
Note that for $m+a\geq 1$ the first sum on the r.h.s.\ of (\ref{4: d_centerdot on MC generators of Tw(P)}) cancels out with all the summands corresponding to $k\geq 2 ,m_0=1,m_i=m-1, c_i=a, i\in [k]$, in the second sum.

\subsubsection{\bf Theorem}\label{5: Th on HoLB^h to Tw(HoLB)[[h]]}\label{5: Theorem on map from HoLoBcd to TwP}  {\em
For any dg properad $\cP$ equipped with a map
 $$
 f: \HoLoB_{c,d} \lon \cP
$$
 %$\HoLoBcd$
 there is an associated map of dg properads
$$
\Tw f: \HoLoB_{c,d} \lon \Tw^\diamond\cP
$$
given explicitly by
$$
\Ba{c}\resizebox{13mm}{!}{\xy
(-9,-6)*{};
(0,0)*+{a}*\cir{}
**\dir{-};
(-5,-6)*{};
(0,0)*+{a}*\cir{}
**\dir{-};
(9,-6)*{};
(0,0)*+{a}*\cir{}
**\dir{-};
(5,-6)*{};
(0,0)*+{a}*\cir{}
**\dir{-};
(0,-6)*{\ldots};
(-10,-8)*{_1};
(-6,-8)*{_2};
(10,-8)*{_n};
(-9,6)*{};
(0,0)*+{a}*\cir{}
**\dir{-};
(-5,6)*{};
(0,0)*+{a}*\cir{}
**\dir{-};
(9,6)*{};
(0,0)*+{a}*\cir{}
**\dir{-};
(5,6)*{};
(0,0)*+{a}*\cir{}
**\dir{-};
(0,6)*{\ldots};
(-10,8)*{_1};
(-6,8)*{_2};
(10,8)*{_m};
\endxy}\Ea
\lon
\sum_{k\geq 1, [m]=\sqcup  [m_\bu]
\atop
%{m_{0}\geq 1, k+m_0\geq 3 \atop m_1,...,m_k\geq 0}} \atop
a=b+\sum_{i=1}^k(c_i+l_i-1)}
\pm
\frac{1}{k!}
\Ba{c}\resizebox{21mm}{!}{\xy
(-3,0)*+{_{c_1}}*\cir{}="b",
(10,10)*+{b}*\cir{}="c",
(23,0)*+{_{c_k}}*\cir{}="r",
%
%%%%%%%%%% edges to c %%%%%%%%%%%%
(6,16)*{}="1'",
(8,16)*{}="2'",
(14,16)*{}="3'",
(11,15)*{...},
(8,2)*{}="1''",
(10,2)*{}="2''",
(12,2)*{}="3''",
(3,-1)*{...},
(17,-1)*{...},
(10,18)*{\overbrace{ \ \ \ \ \ \ \ \ }},
(10,21)*{_{m_0}},
(10,1)*{\underbrace{ \  }},
(10,-1)*{_{n}},
%
%%%%%%%%%%% internal left curved edges %%%%%%%%%%%
(-3,2)*-{};(8.0,10.0)*-{}
**\crv{(0,10)};
(-2.5,1.8)*-{};(8.5,9.0)*-{}
**\crv{(0.4,7)};
(-1.2,0.6)*-{};(9.1,8.5)*-{}
**\crv{(5,1)};
%(-1.0,0.1)*-{};(9.5,8.6)*-{}
%**\crv{(6,-1)};
(5,5)*+{_{l_1}};
%%%%%%%%%%%%% internal right curved edges %%%%%
(23,2)*-{};(12.0,10.0)*-{}
**\crv{(20,10)};
(22.5,1.8)*-{};(11.5,9.0)*-{}
**\crv{(20.4,7)};
(20.9,0.6)*-{};(10.9,8.5)*-{}
**\crv{(15,1)};
%(1.7,0.0)*-{};(9.5,8.6)*-{}
%**\crv{(6,-1)};
(18,5)*+{_{l_k}};
(-7,10)*{}="1l";
(-5,10)*{}="2l";
(-3,10)*{}="3l";
(-5,12)*{\overbrace{  \ \ \ \ \ \ }},
(-5,14)*{_{m_1}},
(27,10)*{}="1r";
(25,10)*{}="2r";
(23,10)*{}="3r";
(25,12)*{\overbrace{  \ \ \ \ \ \ }},
(25,14)*{_{m_k}},
\ar @{-} "c";"1'" <0pt>
\ar @{-} "c";"2'" <0pt>
\ar @{-} "c";"3'" <0pt>
\ar @{-} "c";"1''" <0pt>
\ar @{-} "c";"2''" <0pt>
\ar @{-} "c";"3''" <0pt>
\ar @{-} "b";"1l" <0pt>
\ar @{-} "b";"2l" <0pt>
\ar @{-} "b";"3l" <0pt>
\ar @{-} "r";"1r" <0pt>
\ar @{-} "r";"2r" <0pt>
\ar @{-} "r";"3r" <0pt>
\endxy}\Ea
$$
}

(In the case $c,d\in 2\Z$ the symbol $\pm$ above can be replaced by $+1$.) %Note that the image of the $(1,1)$ generator under this map has no $\hbar^0$ %summand????????????????

To prove this statement one has to check the compatibility of $\Tw f$ with the differentials on both sides.
This can be done either by a direct (but tedious) computation or by studying  generic representations of both properads involved in the above statement as is done briefly  in 
\S {\ref{5: subsec on repr of Tw^diamond HoLoB}}    below  in the most important and illustrative case $\cP=\HoLoBcd$.

\sip

In a full analogy to \S {\ref{3: Theorem on Def action on TwP}}   the deformation complex of any morphism $f$ as above acts on $\Tw^\diamond \cP$ by derivations, that is,
there is a morphism of dg Lie algebras
$$
 \Def\left(\HoLoBcd \stackrel{f}{\rar} \cP\right) \lon  \Der(\Tw\cP)
$$

\subsection{Representations of $\Tw^\diamond \HoLoBcd$}\label{5: subsec on repr of Tw^diamond HoLoB} Let $\rho: \HoLoBcd\rar \cE nd_V$ be a homotopy involutive Lie bialgebra structure in a graded vector space $V$ and let $\Delta_\rho$ be its
equivalent incarnation as a differential operator (\ref{5: Delta_rho}) on $\odot^\bu(V[-c][[\hbar]]$. Assume a Maurer-Cartan element $\ga$ of this $\HoLoBcd$-structure is fixed, that is, a formal power series (\ref{5: ga-MC as h-series}) satisfying the equation
(\ref{5: MC eqn Delta(e^g)}). These data $(\rho, \ga)$ gives us

\Bi
\item[(i)]
 a representation of
$\Tw\HoLoBcd$ in $V$ which sends the MC generators,
$$
\Ba{c}\resizebox{9mm}{!}{\begin{xy}
(0,6)*{...};
(-5.5,7)*{^1};
(0,0)*+{_a}*\cir{}
**\dir{-};
(5.5,7)*{^m};
(0,0)*+{_a}*\cir{}
**\dir{-};
(-2.5,7)*{^2};
(0,0)*+{_a}*\cir{}
**\dir{-};
(2.5,6)*{};
(0,0)*+{_a}*\cir{}
**\dir{-};
%
 %   <0mm,8mm>*{\overbrace{ \ \ \ \ \ \ \ \ \ \ }},
 %   <0mm,10mm>*{^m},
 \end{xy}}\Ea \lon \ga_{a,m} \in \odot^m(V[-c])
$$
to the corresponding summands of the MC series (\ref{5: ga-MC as h-series}).
\item[(ii)] a {\em twisted}\, $\HoLoBcd$-algebra structure on $V$ which can be encoded as the following
$\ga$-twisted differential operator
$$
\Delta_\ga:= e^{-\frac{\ga}{\hbar}}\circ \Delta_\rho \circ e^{\frac{\ga}{\hbar}}=:\sum_{a\geq 0} \hbar^a\Delta_{(a)\ga}
$$
\Ei
The MC equation (\ref{5: MC eqn Delta(e^g)}) guarantees that summands  $\Delta_{(a)\ga}$  are differential operators of order $\leq a+1$ so that $\Delta_\ga$ induces some $\HoLoBcd$-algebra structure on $V$ indeed. A straightforward combinatorial inspection  of $\Delta_\ga$ recovers the universal properadic formula shown in Theorem  \S {\ref{5: Theorem on map from HoLoBcd to TwP}}.

\sip

Thus, contrary to the twisting endofunctor $\Tw$ introduced in the previous section, its diamond version $\Tw^\diamond$ gives us essentially nothing new --- it reproduces in a different language the well-known twisting construction introduced in \S 9 of \cite{CFL} in terms of generic representations of the properads $\HoLoBcd$ and $\Tw^\diamond\HoLoBcd$. In the special class of representations of $\LoBcd$ on the spaces of cyclic words, the MC equation (\ref{5: MC eqn Delta(e^g)}) has been introduced and studied by S.\ Barannikov \cite{Ba1,Ba2} in the context of the deformation theory of modular operads and its applications in the theory of Kontsevich moduli spaces.

\sip

A beautiful concrete solution $\Ga$ of the MC equation (\ref{5: MC eqn Delta(e^g)}) has been constructed by F.\  N\"aef and T.\ Willwacher in \cite{NW} when studying string topology of not necessarily simply connected manifolds $M$; that MC element $\ga$ has been obtained in \cite{NW} from the so called {\em partition function}\, $Z_M$ on $M$ which has been constructed earlier by  R.\ Campos and T.\ Willwacher in \cite{CW} when studying new graph models of configuration spaces of manifolds. Thus the $\HoLoB_{3-n}$-algebra structure constructed in \cite{NW} on the space of cyclic words $Cyc(H^\bu(M)[1])$ of the de Rham cohomology $H^\bu(M)$ of an $n$-dimensional closed manifold $M$ gives us an example of the action of the twisting endofunctor $\Tw^\diamond$ on standard $\LoB_{3-n}$-algebra structure on $Cyc(H^\bu(M)[1])$.

%%%%%%%%%%%%%%%%%%%%%%%%%%%%%%%%%%%%%%%%%%%%%%%%%%%%%%%%%%%%%%%%%%%

\bip

\def\cprime{$'$}

\end{document}